\documentclass{amsbook}

\usepackage{setspace}
\usepackage{amssymb}
\usepackage[all]{xy}

\usepackage[bookmarksopen=true]{hyperref}

\newcommand{\JTO}{$J^T_{\mc{O},\ve}(\;)$}
\newcommand{\JTchi}{$J^T_{\chi,\ve}(\;)$}

\newcommand{\indexs}{\index{symbols}}
\newcommand{\indexi}{\index{index}}
\newcommand{\RR}{\mathcal{R}}
\newcommand{\sig}{{}^\sigma\!}

\newcommand{\cp}{\tilde{\pi}}
\newcommand{\OO}{\mathcal{O}}
\newcommand{\reg}{{\rm reg}}
\newcommand{\rep}{representation } 
\newcommand{\intw}{intertwining operator }

\newcommand{\rR}{{\rm R}}
\newcommand{\ol}{\overline}

\newcommand{\Af}{\mathbb{A}_F}
\newcommand{\Ae}{\mathbb{A}_E}

\newcommand{\adl}{\mathbb{A}}
\newcommand{\Gal}{{\rm Gal}}
\newcommand{\diag}{{\rm diag}}

\newcommand{\blockdiag}{{\rm blockdiag}}
\newcommand{\Vr}{V^{\rm un}(E/F)}

\newcommand{\R}{\mathcal{R}}

\newcommand{\E}{\mathcal{E}}
\newcommand{\Z}{\mathcal{Z}}
\newcommand{\mct}{\mathcal{T}}

\newcommand{\X}{\mathcal{X}}
\newcommand{\U}{\mathcal{U}}
\newcommand{\K}{\mathcal{K}}
\newcommand{\mco}{\mathcal{O}}
\newcommand{\CC}{\mathbb{C}}
\newcommand{\Gm}{\mathbb{G}_m}

\newcommand{\Hom}{{\rm Hom}}
\newcommand{\fhat}{f^{\vee}}

\newcommand{\lp}{\left(}
\newcommand{\rp}{\right)}
\newcommand{\la}{\left\langle}
\newcommand{\ra}{\right\rangle}

\newcommand{\lmx}{\begin{matrix}}
\newcommand{\rmx}{\end{matrix}}
\newcommand{\lsm}{\begin{smallmatrix}}
\newcommand{\rsm}{\end{smallmatrix}}

\newcommand{\tr}{{\rm tr}\;}
\newcommand{\str}{{\rm tr}^{\rm stable}\;}
\newcommand{\ve}{\varepsilon}
\newcommand{\ep}{\epsilon}
\newcommand{\vp}{\varpi}
\newcommand{\fee}{\varphi}
\newcommand{\GL}{{\rm GL}}

\newcommand{\PGL}{{\rm PGL}}

\newcommand{\PGSp}{{\rm PGSp}}
\newcommand{\GSp}{{\rm GSp}}
\newcommand{\SL}{{\rm SL}}

\newcommand{\rH}{H}
\newcommand{\HH}{\mathbf{H}(\Af)}
\newcommand{\GG}{\mathbf{G}(\Af)}
\newcommand{\ZZ}{\mathbb{Z}}
\newcommand{\Zo}{\mathbf{Z}_0(\Af)}
\newcommand{\N}{{\rm N}}
\newcommand{\ord}{{\rm ord}\,}
\newcommand{\bs}{\backslash}
\newcommand{\Int}{{\rm Int}}

\newcommand{\EA}{F_A}
\newcommand{\ED}{F_D}
\newcommand{\EAD}{F_{AD}}

\newcommand{\spl}{{\rm \bf spl}}
\newcommand{\idele}{id\`{e}le }

\newcommand{\adeles}{ad\`{e}les }
\newcommand{\FH}{Frobenius-Hecke }

\newcommand{\A}{\sqrt{A}}

\newcommand{\D}{\sqrt{D}}

\newcommand{\vpmm}{\lp\begin{smallmatrix}1&0\\0&\vp^{-m}\end{smallmatrix}\rp}

\providecommand{\bg}[1]{{\boldsymbol #1}}
\providecommand{\mb}[1]{\mathbf{#1}}
\providecommand{\mc}[1]{\mathcal{#1}}
\providecommand{\mf}[1]{\mathfrak{#1}}
\providecommand{\mbb}[1]{\mathbb{#1}}

\providecommand{\abs}[1]{\left\vert{#1}\right\vert}

\providecommand{\floor}[1]{\lfloor{#1}\rfloor}
\providecommand{\iH}[1]{\iota(G, H_{#1})}

\providecommand{\Acon}[1]{\mathcal{A}_{#1/\mathbf{G}}}
\providecommand{\CL}[1]{{\rm CL}_{{\rm ss}}(#1)}

\providecommand{\idc}[1]{C_{#1}}
\providecommand{\idcc}[1]{\widehat{C_{#1}}}
\providecommand{\dmtwo}[2]{\lp\lsm {#1}&\\&{#2}\rsm\rp}
\providecommand{\dmfour}[4]{\lp\lsm {#1}&&&\\&{#2}&&\\&&{#3}&\\&&&{#4}\rsm\rp}
\providecommand{\ddmfour}[4]{\diag\lp{#1}, {#2}, {#3},{#4}\rp}
\providecommand{\C}{C}
\providecommand{\EL}[1]{E(u_1, u_2, {#1})}
\providecommand{\ELw}[1]{E(w_1, w_2, {#1})}
\providecommand{\St}[1]{{\rm St}_{2, {#1}}}
\providecommand{\one}[1]{1_{2, {#1}}}

\newcommand{\dt}{}
\newcommand{\dts}{}

\newcommand{\dps}{\displaystyle}

\makeatletter
\def\cleardoublepage{\clearpage\if@twoside \ifodd\c@page\else
\hbox{}
\vspace*{\fill}
\vspace{\fill}
\thispagestyle{empty}
\newpage
\if@twocolumn\hbox{}\newpage\fi\fi\fi}
\makeatother

 \let\printindex\relax
 \usepackage{makeidx}
 \usepackage{multind}
\usepackage{times}

\swapnumbers
\newtheorem{generic}{}[chapter]
\newtheorem{thm}[generic]{Theorem}

\newtheorem*{conjecture}{Conjecture}
\newtheorem{lemma}[generic]{Lemma}
\newtheorem{prop}[generic]{Proposition}
\newtheorem{claim}[generic]{Claim}
\newtheorem{define}[generic]{Definition}

\newtheorem{corollary}[generic]{Corollary}

\numberwithin{equation}{section}

\newtheorem{results}{}[chapter]
\newtheorem{rprop}[results]{Proposition}
\newtheorem{rprop'}[results]{Proposition*}
\newtheorem{rprop*}{Proposition}
\newtheorem{rprop*'}{Proposition}

\theoremstyle{remark}
\newtheorem*{remark}{\sc Remark}

\numberwithin{section}{chapter}
\numberwithin{equation}{chapter}

\begin{document}
\frontmatter
\title{Invariant Representations of $\GSp(2)$ under Tensor Product
with a Quadratic Character}

\author{Ping Shun Chan}
\address{Department of Mathematics, The Ohio State University, Columbus,
Ohio 43210}
\email{pschan@math.ohio-state.edu}


\subjclass[2000]{Primary 11F70, 11F72, 11F85}
\keywords{$\GSp(2)$, automorphic representations, $p$-adic harmonic analysis}

\begin{abstract}
Let $F$ be a number field or a $p$-adic field.  We introduce in Chapter \ref{chap:endoscopy} of this work two 
reductive rank one $F$-groups, $\mb{H_1}$, $\mb{H_2}$, which are twisted endoscopic groups of $\GSp(2)$
with respect to a fixed quadratic character $\ve$ of the id\`ele class group of $F$ if $F$ is global, $F^\times$
if $F$ is local.
When $F$ is global, Langlands functoriality predicts that there exists a canonical lifting of the
automorphic representations of $\mb{H_1}$, $\mb{H_2}$ to those of $\GSp(2)$.  In Chapter \ref{chap:global}, 
we establish this lifting in terms of the Satake parameters which parametrize the automorphic
representations.  By means of this lifting we provide a classification
of the discrete spectrum automorphic representations of $\GSp(2)$ which are invariant under tensor product with 
$\ve$.  

The techniques  through which we arrive at our results are inspired by those of Kazhdan's in \cite{K},
which introduced the trace formula twisted by a character.  
In particular, our techniques involve comparing the spectral sides of the trace formulas for the groups under consideration.  
We make use of the twisted extension of Arthur's trace formula, and Kottwitz-Shelstad's stabilization of the 
elliptic part of the geometric side of the twisted trace formula.

When $F$ is local, in Chapter \ref{chap:local} we provide a classification of the irreducible 
admissible representations of $\GSp(2, F)$ which are invariant under tensor product with the quadratic character
$\ve$ of $F^\times$.
More precisely, we use the global results from Chapter \ref{chap:global} to express the twisted characters of these
invariant representations in terms of the characters of the admissible representations of $\mb{H}_i(F)$ ($i = 1, 2$).
These (twisted) character identities provide candidates for the liftings predicted by the conjectural 
local Langlands functoriality. The proofs are inspired by \cite{K} and
rely on Sally-Tadi\'c's classification of the irreducible admissible representations of $\GSp(2, F)$,
and on Flicker's results on the lifting from $\PGSp(2)$ to $\PGL(4)$.

In both the global and local cases, our results depend on the conjectural statement that there exist smooth functions
on the groups which have matching (twisted) orbital integrals.
In particular, these results depend on the validity of the (twisted) Fundamental Lemma in the context of $\GSp(2)$ and
$\mb{H_1}$, $\mb{H_2}$.  We prove in Appendix \ref{chap:fundlemma} the Fundamental Lemma in our case, 
making use of Flicker's work on the twisted Fundamental Lemma for $\GL(4)$ and $\GSp(2)$, 
and Hales' reduction of the Fundamental Lemma to unit elements.

\end{abstract}

\maketitle
\setcounter{page}{4}

\tableofcontents

\chapter*{Technical Preface}
\section*{$\ve$-Endoscopic Groups}
Let $F$ be a number field.  Let $\ve$ be a quadratic character of the \idele class group $\idc{F}$ of $F$.
Let $E$ be the quadratic extension of $F$ which corresponds to $\ve$ via global class field theory.
The quasi-split $\ve$-endoscopic groups of $\mb{G}$ over $F$ are as follows:
\begin{enumerate}
\item
$
\displaystyle
\mb{H_1} = \lp \GL(2) \times \rR_{E/F}\mbb{G}_m\rp',
$
where $\rR_{E/F}$ is the restriction of scalars (from $E$ to $F$) functor.  The prime 
denotes the condition that the determinant of the $\GL(2)$-factor
is equal to the norm of the $\rR_{E/F}\mbb{G}_m$-factor in $\mbb{G}_m$.  In particular, the group of $F$-points
of $\mb{H_1}$ is
\[\dps
\mb{H_1}(F)
= \left\{(g, x) \in \GL(2, F)\times E^\times : \det g = \N_{E/F}x\right\}.
\]
\item
$\mb{H_2}$ is the unique quasi-split reductive $F$-group whose group of $F$-points is 
\[\dps
\mb{H_2}(F) = \lp \GL(2, E) \times F^\times\rp/ \{\lp\diag(z, z),\N_{E/F}z^{-1}\rp:z\in E^\times\}.
\]
\end{enumerate}
\section*{Global Lifting}



We say that an automorphic representation $\pi$ of $\GL(2, \Af)$ is $E$-{\bf monomial} if it is the monomial
representation $\pi(\theta)$ associated with a character $\theta$ of $\idc{E}$.
For a representation $\pi$, let $\omega_\pi$ denote its central character.
Let $B_{E/F}\pi$ denote the $\GL(2, \Ae)$-module which is the base change of a $\GL(2, \Af)$-module $\pi$.
\begin{rprop*}
Let $\pi$ be an irreducible automorphic $\GL(2, \Af)$-module which is either
cuspidal non-$E$-monomial or one dimensional.
The packets $\pi\otimes_1 1$ of $\;\mb{H_1}(\Af)$ and $B_{E/F}\pi\otimes_2 \omega_{\pi}$ of $\;\mb{H_2}(\Af)$ 
lift to an unstable $($quasi-$)$packet $\{\Pi\}$ of $\GSp(2, \Af)$, which in turn lifts to the induced representation $I_{(2,2)}(\pi, \ve\pi)$ of $\GL(4, \Af)$.
\begin{remark} If $\pi$ is one dimensional, $\{\Pi\}$ belongs to a class of quasi-packets constructed by Howe, Piatetski-Shapiro.
\end{remark}
\end{rprop*}

For an automorphic representation $\tau$ of $\GL(2, \Ae)$, 
let $\pi(\tau)$ denote the $\GL(4, \Af)$-module which is obtained from
$\tau$ via the twisted endoscopic lifting from $\rR_{E/F}\GL(2)$ to $\GL(4)$. 
let ${}^\sigma \tau$ denote the representation
\[
{}^\sigma \tau : (g_{ij}) \mapsto \tau(\sigma g_{ij}),\quad \forall (g_{ij}) \in \GL(2, \Ae).
\]
\begin{rprop*'}
Let $\pi$ be an irreducible, cuspidal, automorphic representation of $\GL(2, \Ae)$,
$\mu$ a character of $\idc{F}$, such that $\pi \neq {}^\sigma \pi$ and $\omega_{\pi} = \mu\circ\N_{E/F}$.
The representation $\pi\otimes_2 \mu$ of $\mb{H_2}(\Af)$ lifts to a stable packet $\{\Pi\}$, which in turn lifts
to the cuspidal automorphic representation $\pi(\tau)$ of $\GL(4, \Af)$.
No discrete spectrum packet of $\mb{H_1}(\Af)$ lifts to the  packet $\{\Pi\}$.
\end{rprop*'}

\begin{rprop*'}
Suppose $\{\Pi\}$ is a global packet of $\GSp(2, \Af)$ which is the lift of some automorphic 
representation/packet of one of the $\ve$-endoscopic groups.  Then, $\{\Pi\}$ contains at least
one representation which is $\ve$-invariant.
If $\{\Pi\}$ is a discrete spectrum global packet of $\GSp(2, \Af)$ which is not the lift of 
any discrete spectrum representation/packet of the $\ve$-endoscopic groups, then $\{\Pi\}$
does not contain any discrete spectrum automorphic representation which is $\ve$-invariant.
\end{rprop*'}
\section*{Local Lifting}
Let $F$, $E$ now be local nonarchimedean fields, with $[E:F] = 2$.
Let $\ve$ be the nontrivial quadratic character of $F^\times$ corresponding to $\K/F$.
Let $G = \GSp(2, F)$, $H_1 = \mb{H_1}(F)$, $H_2 = \mb{H_2}(F)$.


Let $Z$ be the center of $G$.
For an irreducible $\ve$-invariant admissible representation $\pi$ of $G$, and a smooth function $f$ on $G$, put
\[
\la \pi, f\ra_A := \tr \int_{Z\bs G} \pi(g)f(g)A\;dg,
\]
where $A$ is an intertwining operator in $\Hom_G(\pi, \ve\pi)$ such that $A^2 = 1$.

For $i = 1, 2$, let $Z_0(H_i) $ denote the maximal $F$-split component of the center of $H_i$.
For an admissible representation $\pi_i$ of $H_i$, and a smooth function $f_i$ on $H_i$, put
\[
\la \pi_i, f_i\ra := \tr \int_{Z_0(H_i)\bs H_i} \pi_i(h_i)f_i(h_i)\;dh_i.
\]


\begin{rprop*} 
 Let $\tau$ be an irreducible, cuspidal, non-$\K$-monomial, admissible $\GL(2, F)$-module whose central character
$\omega_\tau$ is trivial.
Let $\{\pi^+, \pi^-\}$ be the local packet of $G$ which lifts to the induced representation
$I_{(2, 2)}(\tau, \ve\tau)$ of $\PGL(4, F)$ according to \cite{F1}.  
There exist intertwining operators $A^+$ in $\Hom_G(\pi^+, \ve\pi^+)$, and $A^-$ in $\Hom_G(\pi^-, \ve\pi^-)$,
such that the following holds for matching functions:
\[
\begin{split}
\la \pi^+, f\ra_{A^+} + \la \pi^-, f\ra_{A^-} &= \la \tau\otimes_1 1,\, f_1\ra,\\
\la \pi^+, f\ra_{A^+} - \la \pi^-, f\ra_{A^-} &= \la B_{\K/F}\tau\otimes_2 \omega_\tau,\, f_2\ra.
\end{split}
\]
\end{rprop*}

\begin{rprop*}
Let $\tau$ be an irreducible, cuspidal, non-$\K$-monomial, admissible representation of \dt $\GL(2, F)$,
$\chi$ be a character of $\K^\times$, such that: There does not exist a character $\ve'$ of $F^\times$ such that
${}^\sigma\chi/\chi = \ve'\circ\N_{\K/F}$ and $\tau \cong \ve'\tau$.
Then, there exist a collection $[\pi]$ of distinct $\ve$-invariant, irreducible, cuspidal, admissible 
$G$-modules, and an operator $A(\pi') \in \Hom_G(\pi',\,\ve\pi')$ for each $\pi' \in [\pi]$,
such that the following holds for matching functions:
\[
\sum_{\pi' \in [\pi]}\la \pi', f\ra_{A(\pi')} = \la \tau\otimes_1\chi,\,f_1\ra.
\]
\end{rprop*}

\begin{rprop*}
Let $\tau$ be an irreducible, cuspidal, admissible representation of $\GL(2, \K)$, not fixed by the action of $\Gal(\K/F)$,
such that its central character $\omega_\tau$ is equal to $\mu\circ\N_{\K/F}$ for some character $\mu$ of $F^\times$.
Then, there exist a collection $[\pi]$ of distinct $\ve$-invariant, irreducible, cuspidal, admissible 
$G$-modules, and an operator $A(\pi')$ in $\Hom_G(\pi',\,\ve\pi')$ for each $\pi' \in [\pi]$, 
such that the following holds for matching functions:
\[
\sum_{\pi' \in [\pi]}\la \pi', f\ra_{A(\pi')} = \la \tau\otimes_2\mu,\,f_2\ra.
\]
\end{rprop*}

\aufm{Ping Shun Chan\\\today\\Columbus, Ohio, USA}

\mainmatter

\chapter{Introduction}\label{chap:intro}
\section{An Overview}
\subsection{Langlands Functoriality}
In the 1960's, Robert Langlands conjectured that there exists a canonical (not necessarily one-to-one) correspondence
between Galois representations and automorphic representations of reductive algebraic groups, and that this correspondence
coincides with the Artin reciprocity law when the underlying group is $\GL(1)$.  
More precisely, for every reductive group $\mb{H}$ over a number field $F$ with Langlands dual ${}^L\! H$, 
it is believed that each homomorphism from the Langlands group $L_F$---a conjectural group generalizing the Galois and 
Weil groups---to ${}^L\! H$ parametrizes an $L$-packet 
of automorphic representations of $\mb{H}$.

Suppose there are two reductive groups $\mb{H}$, $\mb{G}$ with $L$-groups ${}^L\!H$, ${}^LG$, respectively, such that
there is a homomorphism from ${}^L\! H$ to ${}^L G$.  Then, each homomorphism from the Langlands group to ${}^L\! H$ induces 
one from the Langlands group to ${}^L G$.  
Consequently, one should expect a canonical lifting of the automorphic representations of
$\mb{H}$ to those of $\mb{G}$.  A more precise formulation of this principle is as follows:

\begin{conjecture}[Langlands Functoriality]
Let $\mb{H}$, $\mb{G}$ be reductive groups over a number field such that there is an $L$-homomorphism
$\xi$ of the $L$-group of ${}^L\!H$ of $\mb{H}$ into that of $\mb{G}$.  Then, for every automorphic representation 
$\pi$ of $\mb{H}$, there is a packet $\{\Pi\}$ of automorphic representations of 
$\mb{G}$, parametrized by the same \FH classes at all but finitely many places, 
such that the \FH classes in ${}^L\! H$ parametrizing $\pi$ are mapped by $\xi$ to those parametrizing $\{\Pi\}$ in ${}^L G$. 
\end{conjecture}


Over the years, various examples of Langlands functorial lifting have been proved via the comparison of trace formulas,
a technique introduced by Langlands and his collaborators (\cite{JL}, \cite{L1}, \cite{LL}).  
We list here a few notable cases which are of particular interest in relation to our work.  They are all examples of
(twisted) endoscopic liftings (see \cite{KS}).
\begin{itemize}
\item 
The lifting from $\rR_{E/F} \GL(1)$ to $\GL(2)$, where $[E:F] = 2$, 
proved by Jacquet and Langlands \cite{JL}.
\item
Cyclic base change lifting for $\GL(2)$, proved by Langlands \cite{L1}, and for $\GL(n)$, by Arthur and Clozel
\cite{AC}.
\item
The lifting from elliptic algebraic tori to $\SL(2)$, by Labesse and Langlands \cite{LL}.
\item
The lifting from $\rR_{E/F} \GL(1)$ to $\GL(n)$, $E/F$ cyclic of order $n$, by Kazhdan \cite{K}.
\item
The lifting from $\rR_{E/F}\GL(m)$ to $\GL(n)$, where $m [E:F] = n$, by Waldspurger \cite{WAbasechange}.
\item
The lifting from $\PGSp(2)$ to $\PGL(4)$, by Flicker \cite{F1}.
\item
Endoscopic lifting for $\GSp(2)$, by Arthur \cite{A}.
\end{itemize}

The first half of this work uses the trace formula technique to establish the lifting of automorphic representations 
of two rank one twisted endoscopic groups $\mb{H_1}$ and $\mb{H_2}$ to the symplectic group of similitudes $\GSp(2)$ 
($4 \times 4$ matrices). 
At the same time, we show that the lifting provides a classification of the
automorphic representations of $\GSp(2)$ which are invariant under tensor product with a fixed quadratic
character.  

To a great extent, the conception of this project and the methodology employed are inspired by the 
work of Kazhdan's \cite{K}, which uses the trace formula technique both to establish Langlands functorial lifting
from $\rR_{E/F} \mb{G}_m$ to $\GL(n)$, where $E/F$ is cyclic of degree $n$,
and to classify the automorphic representations of $\GL(n, \Af)$ 
which are invariant under tensor product with a fixed 
character of order $n$.  In particular, he showed that each id\`ele class character $\theta$ of $E$
lifts to an automorphic representation $\pi(\theta)$ of $\GL(n, \Af)$.  
Moreover, an automorphic representation of $\GL(n, \Af)$ is invariant
under tensor product with an order $n$ character $\ve$ if and only if it is equal to $\pi(\theta)$ for 
some character $\theta$ of the cyclic field extension $E$ of $F$ defined by $\ve$ via class field theory.
That Kazhdan was able to obtain these spectacular results is a testament to the power of the trace formula technique.

In this work, we have an analogous situation for $\GSp(2)$, except that in this case there are nonsingleton 
(quasi-)$L$-packets, and there are two twisted endoscopic groups instead of one.
 Modulo the conjectural existence of matching functions and certain technical conditions on the local 
components of the automorphic representations, we show that the 
automorphic representations of the two twisted endoscopic groups ($[E:F] = 2$)
\[\mb{H_1} =\lp \GL(2)\times\rR_{E/F}\mb{G}_m\rp'\quad\text{ and }\quad
\mb{H_2} = \lp \rR_{E/F}\GL(2) \times \mb{G}_m\rp/\rR_{E/F}\mb{G}_m
\]
lift to the (quasi-)packets of $\GSp(2)$ containing at 
least one automorphic representation which is invariant under tensoring with the quadratic character
$\ve$ associated with $E/F$.  
Moreover, all such packets are the lifts of the automorphic representations of the twisted endoscopic groups.  
Naturally, an important question is whether every representation in such a packet is invariant under
tensoring with $\ve$, but that remains unanswered in this work.

The presence of nonsingleton packets for $\GSp(2)$ also raises interesting questions regarding the stability of the packets.
Thanks to the classification results of \cite{A} and \cite{F1}, we show that a packet containing an invariant 
representation is stable if and only if it is the lift of exactly one of the twisted endoscopic groups.  

Another interesting phenomenon is that each discrete spectrum automorphic representation $\pi$ of $\GSp(2)$ which is invariant
under tensoring with $\ve$ in turn lifts to an automorphic representation $\Pi$ of $\GL(4)$ which is invariant under
tensor product with the same character.  Moreover, $\Pi$ is induced if and only if $\pi$ belongs to an unstable
(quasi-)packet.

It is also worth mentioning that the lifting from $\mb{H_2}$ to $\GSp(2)$ may be viewed as follows:
Let $E/F$ be the quadratic field extension which corresponds to the character $\ve$.
Each automorphic representation $\pi$ of $\mb{H_2}(\Af)$ may be constructed from an automorphic representation
$\tau$ of $\GL(2, \Ae)$ whose central character is invariant under the action of $\Gal(E/F)$.  Let 
$\varphi : L_E \rightarrow \GL(2, \CC)$ be the representation of the Langlands group of $E$ which corresponds to
$\tau$.  Then, the representation of $L_F$ (the Langlands group of $F$) induced from $\varphi$
factors through $\GSp(2)$:
\[
\xymatrix{
{\rm Ind}_{L_E}^{L_F}\varphi:\mspace{-42mu}&  L_F \ar[rr] \ar[dr] & & \GL(4, \CC)\\
& & \GSp(2) \ar[ur] &
},
\]
and the automorphic representation of $\GL(4, \Af)$ corresponding to ${\rm Ind}_{L_E}^{L_F}\varphi$ is the lift of
$\pi$.  This is beautifully analogous to the construction of the monomial representations $\pi(\theta)$
of $\GL(2)$ in \cite{JL} and \cite{K}.  

It should be noted here that automorphic representations of $\GSp(2)$ have received the attentions of 
mathematicians for quite some time.  In addition to those cited above, 
it suffices to mention Siegel, Shimura, Weissauer, Howe, Piatetski-Shapiro,
Vigneras, D. Prasad, B. Roberts, and R. Schmidt.
Note that many mathematicians prefer to use ${\rm GSp}(4)$ rather than $\GSp(2)$ to denote 
the same group.

\subsection{Local Lifting}
Langlands Functoriality as stated above has a local analogue.  Namely, it is conjectured that, for any reductive group
$\mb{H}$ over a local $p$-adic field $F$, the homomorphisms from the Weil group $W_F$ (or $L_F$) to ${}^L H$ correspond
canonically to packets of admissible representations of $\mb{H}(F)$.  Consequently, analogously to the global case,
a homomorphism between $L$-groups should parallel the lifting of admissible representations from one $p$-adic 
group to another. 

Since each square integrable admissible representation is a local component of an automorphic representation, 
it is a general philosophy that one can derive a large class of local liftings from global liftings.
Works in the past which establish global lifting using trace formulas typically also derive local lifting.
More precisely, they express the local lifts in the form of trace identities among admissible representations.
For example, for each admissible representation $\pi$ of $\GL(n, F)$ which is invariant under
tensor product with a character $\ve$ of $F^\times$ of order $n$,  
Kazhdan establishes in \cite{K} an expression of the $\ve$-twisted Harish-Chandra character of $\pi$ in terms of 
the Galois orbit of a character of $E^\times$, 
where $E$ is a degree $n$ cyclic extension of $F$ determined by $\ve$ via local class field theory.

The second half of this work is occupied with deriving local twisted trace identities among admissible representations of
$\GSp(2)$ and those of its twisted endoscopic groups.  We arrive at these identities using our global lifting results.

\section{$\ve$-Invariant Automorphic Representations}
\index{index}{automorphic!representation}%
Let $F$ be a number field.  Let $\mbb{A}_F$ denote the ring of \adeles of $F$ and
$\idc{F}$ the \idele class group $F^\times\bs\mbb{A}_F^\times$.
\indexs{A@$\mbb{A}_F$}\indexs{C@$C_F$}%
Let $V$ denote the set of places of $F$.
\index{symbols}{V@$V$}%
For each finite place $v$ of $F$, let $\mathcal{O}_v$ denote the
\index{symbols}{O@$\mc{O}_v$}%
ring of integers of $F_v$.
Let $\mb{H}$ be a reductive algebraic group over $F$.
Put $H := \mb{H}(F),\; H_v := \mb{H}(F_v)$ for any place $v$ of $F$.

For any element $\gamma$ of any group, let subscript $\gamma$ denote the centralizer of
$\gamma$ in that group.  For example, $H_\gamma$ denotes the centralizer in $H$ of $\gamma$.
Let $\mb{Z}$ be the center of $\mb{H}$.  
\index{symbols}{Z@$\mb{Z}$}%
Let $\mb{Z}_0$ be the maximal $F$-split component of $\mb{Z}$.
\indexs{Z@$\mb{Z}_0$}%
Put $Z := \mb{Z}(F),\;Z_0 := \mb{Z}_0(F)$.
For the groups $\GSp(2)$ and $\GL(4)$, which are $F$-split, there is no distinction between
$\mb{Z}$ and $\mb{Z}_0$.

Let $\omega$ be a character of $Z_0 \bs \mb{Z}_0(\Af)$.  Let
\index{symbols}{omega@$\omega$}%
$L(\HH, \omega)$ be the space of measurable functions
\indexs{L@$L(\;,\omega)$}%
$\phi$ on $H\bs\HH$ such that $\phi(zg) = \omega(z)\phi(g)$ for all $z \in \Zo$ and
\[
\int_{\Zo H \bs\HH}\abs{\phi(h)}^2 dh < \infty.
\]
Let $\rho$, or $\rho_\omega$ to emphasize the dependence on $\omega$,
\index{symbols}{rho@$\rho$}%
be the right regular representation of $\HH$ on $L(\HH, \omega)$, that is:
\index{index}{representation!right regular}%
\[
\lp\rho_\omega(h)\phi\rp(g) = \phi(gh), \quad \forall h, g \in \HH, \phi \in L(\HH, \omega).
\]
In particular, the restriction of the central character of $\rho_\omega$ to $\Zo$ is equal to $\omega$.
We say that a representation of $\HH$ is {\bf automorphic} \index{index}{automorphic!representation} if it occurs in 
$\rho_\omega$ for some character $\omega$ of $Z_0\bs\Zo$.
Let $\mathcal{A}(\mb{H}, \omega)$ denote the set of 
\index{symbols}{A@$\mathcal{A}(\mb{H}, \omega)$}%
equivalence classes of irreducible automorphic representations 
$\pi$ of $\mb{H}(\Af)$ with the property that their central characters restrict to $\omega$ on $\Zo$.

Write $\pi \in \mc{A}(\mb{H}, \omega)$ as a tensor product $\otimes_{v \in V}\pi_v$ of representations 
$\pi_v$ of $H_v$.
Since $\rho_\omega$ is continuous, $\pi$ too is continuous with respect to the topology on $\mb{H}(\Af)$ (see \cite{PR}).
Consequently, at almost every finite place $v$, the local representation $\pi_v$ is unramified 
\index{index}{representation!unramified}, i.e.
it contains a nonzero vector fixed by the hyperspecial maximal compact subgroup $\mb{H}(\mc{O}_v)$.

As an $\HH$-module, $L(\HH, \omega)$ decomposes into a direct sum
\[
\displaystyle L_d(\HH, \omega)\oplus L_c(\HH, \omega).
\]
Here, $L_d(\HH, \omega)$ is the closed span of the 
\indexs{L@$L(\;,\omega)$!$L_d(\;, \omega)$}%
\indexs{L@$L(\;,\omega)$!$L_c(\;, \omega)$}%
irreducible, closed, invariant subspaces of $L(\HH, \omega)$ and $L_c(\HH, \omega)$ is its orthogonal
complement.  Let $\rho_d, \rho_c$ denote the restrictions of
$\rho = \rho_\omega$ to $L_d(\HH, \omega)$, $L_c(\HH, \omega)$, respectively.
We call the $\mb{H}(\Af)$-module 
$L_d(\HH, \omega)$ (resp. $L_c(\HH, \omega)$) the {\bf discrete} (resp. {\bf continuous}) {\bf spectrum} of $\HH$.
\index{index}{continuous spectrum}%

We say that the algebraic group $\mb{H}$ has the {\bf multiplicity one}  
\index{index}{multiplicity one property}%
property if each automorphic representation of
$\mb{H}(\Af)$ with central character $\omega$ on $\mb{Z}_0(\Af)$ occurs {\it at most once} in the discrete spectrum of $\rho_\omega$.  
\index{index}{discrete spectrum}%
By the multiplicity one theorem for $\GL(2)$ and a result in \cite{LL}, the multiplicity one property holds for the
twisted endoscopic groups studied in this work.
In the case of $\GSp(2)$, the multiplicity one property for automorphic representations with trivial central characters is 
established in \cite{F1}.  For arbitrary central characters, Arthur has announced in \cite{A} that
multiplicity one holds at least for those representations which are associated with semisimple Arthur parameters.
\begin{remark}
An automorphic representation is really a representation which occurs in the space of
slowly increasing functions called {\bf automorphic forms \index{index}{automorphic!forms}}.   This space is not the same as
$L^2(\HH, \omega)$.  However, the discrete spectrum of this space (as an $\HH$-module) consists of
$L^2$ functions; thus, it coincides with the discrete spectrum of $L^2(\HH, \omega)$;
see \cite{BJ}.
\end{remark}



For any place $v \in V$,
let $C(H_v, \omega_v)$ be the space of smooth functions $f_v$ on $H_v$
\index{symbols}{C@$C(H_v, \omega_v)$}%
such that $f_v$ is compactly supported
modulo $Z_{0,v}$ and $f_v(zh) = \omega_v^{-1}(z)f(h)$ for all $z$ in $Z_{0,v}$ and  $h$ in $H_v$.

If $v$ is finite, and $\mb{H}$ is defined over $\mc{O}_v$,
let $K_v$ denote the hyperspecial maximal compact subgroup $\mb{H}(\mathcal{O}_v)$ of
\index{symbols}{K@$K_v$}%
$H_v$.
Let $\mathcal{H}(H_v, \omega_v)$ denote the Hecke algebra \index{index}{Hecke algebra} of $K_v$-biinvariant functions
\index{symbols}{H@$\mathcal{H}(H_v, \omega_v)$}%
in $C(H_v, \omega_v)$.  The algebra
$\mathcal{H}(H_v, \omega_v)$ is nonzero if and only if $\omega_v$ is trivial on
$Z_{0,v} \cap K_v$.
If $v$ is archimedean, we fix a maximal compact subgroup $K_v$ of $H_v$.  Let
$\mathcal{H}(H_v, \omega_v)$  be the set of $K_v$-finite functions in $C(H_v, \omega_v)$.

Let $C(\HH, \omega)$ be the span of the smooth functions $f = \otimes_{v \in V}f_v$ on $\HH$ which are 
\indexs{C@$C(\HH,\omega)$}%
compactly supported modulo $\mb{Z}_0(\Af)$ such that $f_v \in C(H_v, \omega_v)$
for all $v \in V$ and $f_v$ is the unit element in the Hecke algebra $\mc{H}(H_v, \omega_v)$ for almost all finite $v$.

Fix a Tamagawa measure 
 (see \cite{PR}) 
on $\mb{H}(\Af)$.
For any $f \in C(\mb{H}(\Af), \omega)$ and
$\pi \in \mc{A}(\mb{H}, \omega)$, 
let $\pi(f)$ denote the following convolution operator on the space of $\pi$:
\[\begin{split}
\pi(f) &:= \int_{\Zo\bs \HH}\pi(h)f(h)\;dh\\
&= \int_{\mb{Z}(\Af)\bs\mb{H}(\Af)}\pi(h)\lp\int_{\mb{Z}_0(\Af)\bs\mb{Z}(\Af)}\omega(z)f(zh)\;dz\rp dh.
\end{split}\]
\indexi{convolution operator}%
It has finite rank.  Let $\la \pi, f \ra$ denote the trace of $\pi(f)$.
\index{symbols}{b@$\la\;,\;\ra$}%

\subsection{The $\ve$-Twisted Trace Formula} \index{index}{trace formula!epsilon@$\ve$-twisted}%
Recall that $\rho_\omega$ is the right regular representation of $\mb{H}(\Af)$ in $L(\mb{H}(\Af), \omega)$.
Let $\ve$ be a homomorphism in $\Hom\lp\HH,\CC^\times\rp$ which is trivial on $\mb{Z}(\Af)H$.
Let $\pi$ be an automorphic representation in $\mc{A}(\HH, \omega)$.  Let $V_\pi$ be the space
\index{symbols}{V@$V_\pi$}%
of $\pi$.  By definition, $\pi$ is a subquotient of $\rho_\omega$.
Let $\ve\pi$ be the representation of $\HH$ on the space $V_\pi$ of $\pi$ defined as follows:
\[
\ve\pi := \ve\otimes \pi : h \mapsto \ve(h)\pi(h),\quad\forall h \in \HH.
\]
Note that the space $V_{\ve\pi} = V_\pi$ of $\ve\pi$ is a subquotient of $L(\mb{H}(\Af), \omega)$, 
but $\ve\pi$ is not the right regular representation of $\mb{H}(\Af)$ on $V_{\ve\pi}$.

We say that $\pi$ is {\bf $\ve$-invariant} 
\index{index}{representation!$\ve$-invariant}%
if $\pi \cong \ve\pi$; i.e.  there exists
an automorphism $A$ on $V_\pi$ such that $A\pi(h) = \ve\pi(h)A$ for all $h \in \HH$.
By Schur's lemma, $A^2$ is a scalar, and we normalize $A$ (multiplying it by $(A^2)^{-1/2}$ if necessary) 
so that $A^2 = 1$.  By this choice of normalization, the operator $A$ is unique up to a sign.


Define an operator $\rho_\omega(\ve)$, or $\rho(\ve)$ for brevity, on $L(\HH, \omega)$ as follows:
\[
\lp \rho_\omega(\ve)\phi\rp\!(h) := \phi(h)\ve(h),\quad \forall \phi \in L(\HH, \omega),\; h \in \HH.
\index{symbols}{rho@$\rho_\omega(\ve)$}%
\]
The restriction of $\rho(\ve)$ to the space $V_\pi$ is a nontrivial intertwining operator \index{index}{intertwining operator} of $\ve\pi$ into $\rho$.
In particular, $\ve\pi$ is equivalent to the automorphic representation $\rho|_{\rho(\ve) V_\pi}$ of $\HH$.

Assuming that multiplicity one holds for $\mb{H}$,
any two irreducible inequivalent submodules of $\rho_\omega$ in $L(\mb{H}(\Af), \omega)$ 
are either orthogonal or equal to each other,
each being generated by $\rho(\mb{H}(\Af))$ from a single vector.

Suppose $\pi$ is irreducible, and $\pi\cong\ve\pi$. 
In principle, $\rho_\omega(\ve)$ needs not intertwine $\pi$ with $\ve\pi$, namely $\rho_\omega(\ve)V_\pi$
may be equivalent to $V_\pi$ as an $\mb{H}(\Af)$-module but different from it as a subspace of $L(\mb{H}, \omega)$.
However, since we assume that the multiplicity one property holds, we have $\rho_\omega(\ve)V_\pi = V_\pi$.
It then follows that the restriction of $\rho_\omega(\ve)$ to $V_\pi$ is a nontrivial intertwining operator in
$\Hom_{\mb{H}(\Af)}(\pi, \ve\pi)$.  By Schur's lemma, $\rho_\omega(\ve)|_{V_\pi}$ is a 
nontrivial scalar multiple of $A$; hence, it is $A$ or $-A$.

For an irreducible automorphic representation $\pi$ and a function $f$ in $C(\GG, \omega)$, put
\[
\la \pi, f \ra_\ve := \tr \pi(f)\rho_\omega(\ve) = \tr \rho_\omega(\ve)\pi(f).
\]
\index{symbols}{b@$\la\;,\;\ra_\ve$}%
If $\pi \ncong \ve\pi$, then the spaces $V_\pi$ and  $\rho_\omega(\ve)V_\pi$ are disjoint,
and hence $\la \pi, f\ra_\ve = 0$.
Conversely, if $\pi$ is $\ve$-invariant, we shall show in Lemma \ref{lemma:venondegeneracy} that
the distribution $f \mapsto \la \pi, f\ra_\ve$ is nonzero.


By definition, $\rho_\omega$ is the direct sum of all automorphic representations with $\mb{Z}_0(\Af)$-character $\omega$.
If we could express $\la \rho_\omega, f\ra_\ve$ as the sum (discrete plus continuous)
\dt \[\sum_{\pi \in \mc{A}(\mb{H}, \omega)} \la \pi, f\ra_\ve,\] then those automorphic representations
which are not $\ve$-invariant would drop out from the sum.  Thus, at least formally,
the distribution $\la \rho_\omega, \cdot \ra_\ve : f \mapsto \la\rho_\omega, f\ra_\ve$ on the space $C(\HH, \omega)$
detects those automorphic representations which are $\ve$-invariant.
However, the situation is more delicate.

For $f \in C(\HH, \omega)$, $\rho_\omega(f)\rho_\omega(\ve)$ is an integral operator on $L(\HH, \omega)$.
Let $K_{f,\ve} : \HH \times \HH \rightarrow \mbb{C}$  be the kernel of $\rho_\omega(f)\rho_\omega(\ve)$.
\index{symbols}{K@$K_{f, \ve}$}%
The trace $\la \rho_\omega, f \ra_\ve$ is obtained by integrating $K_{f,\ve}(h, h)$ over $h \in \Zo H \bs \HH$.  
The function $h \mapsto K_{f,\ve}(h, h)$
may be expressed as a sum over conjugacy classes in $H$,
and also as a sum over automorphic representations of $\mb{H}(\Af)$.
For applications, we would like to change the order of summation and integration 
in $\int_{\Zo H \bs \HH} K_{f,\ve}(h, h)\;dh$.
However, this change of order is not justified unless $\mb{H}$ is anisotropic over $F$.
Extending Arthur's \index{index}{Arthur, James} results in \cite{A1}, the authors of \cite{CLL} introduce a truncation
\indexi{truncation}%
$K_{f,\ve}^T(h)$ of $K_{f,\ve}(h, h)$.
\index{symbols}{K@$K_{f,\ve}^T(h)$}%
They show that we \emph{may} pull the sum out of the integral in 
$\int_{\Zo H \bs \HH}K_{f,\ve}^T(h)$.  
The {\bf $\ve$-twisted Arthur trace formula} \index{index}{Arthur, James} 
\index{index}{trace formula!epsilon@$\ve$-twisted}
\index{index}{Arthur's trace formula}%
(or $\ve$-trace formula) 
is the following equality of two different ways to express 
the integral $ \int_{\mb{Z}_0(\Af) H \bs \HH} K_{f,\ve}^T(h)\;dh$:
\[
\sum_{\{\mc{O}\}}J^T_{\mc{O},\ve}(f) = \sum_{\{\chi\}}J^T_{\chi,\ve}(f).
\]
\indexs{J@\JTO}%
\indexs{J@\JTchi}%
The left hand side of the equation is a sum of (weighted) orbital integrals
over the set of semisimple conjugacy classes in $H$.
It is called the {\bf $\mc{O}$-expansion} or the {\bf geometric side} \index{index}{trace formula!geometric side} of the trace formula.
The right hand side is a sum of traces of automorphic representations
over a set of spectral data of $\HH$.  These are data which catalogue
the $\ve$-invariant automorphic representations of $\HH$.  We call $\sum_{\{\chi\}}J^T_{\chi,\ve}(f)$ 
the {\bf fine $\chi$-expansion} or the {\bf spectral side}
\index{index}{trace formula!spectral side} of the trace formula.

\subsection{The Case of $\GSp(2)$}
Let $J = \lp\lsm & w \\ -w&\rsm\rp \in \GL(4)$, where $w = \lp\lsm & 1\\1&\rsm\rp \in \GL(2)$.
\indexs{J@$J$}\indexs{w@$w$}%
Let $\mb{G} = \GSp(2)$, the {\bf symplectic group of similitudes of rank 2},
\index{symbols}{G@$\GSp(2)$}%
be the algebraic group over $\mbb{Z}$ defined as follows:
\[
\GSp(2) = \left\{g \in \GL(4): {}^t g J g = \lambda(g) J
\text{ for some } \lambda(g) \in \mathbb{G}_m\right\}.
\]
We call $\lambda$ the {\bf similitude character}.
\indexi{similitude character}%

Let $F$ be a number field, and consider $\mb{G}$ as a reductive $F$-group.  For any $g \in \mb{G}$, 
we call $\lambda(g)$ its {\bf similitude factor \index{index}{similitude factor}}.
The center of $\mb{G}$ is
$\mb{Z} = \{\diag(z, z, z, z):z \in \mbb{G}_m\}$.  It is $F$-split.
Fix a character $\omega$ of $Z\bs\mb{Z}(\Af)$.  Let $\ve$ be a quadratic character of $G\bs\GG$ which
factors through the similitude factor.
Since $\ve$ is quadratic, and $\lambda(\diag(z, z, z, z)) = z^2$, $\ve$ is trivial on $Z\bs\mb{Z}(\Af)$.
Consider the $\ve$-twisted Arthur trace formula for $\mb{G}(\Af)$.
As a distribution on $C(\GG, \omega)$, the $\mc{O}$-expansion $\sum_{\{\mc{O}\}}J_{\mc{O},\ve}^T(f)$ 
of the $\ve$-trace formula is not invariant with respect to stable conjugacy.
In \cite{KS}, Kottwitz and Shelstad express the elliptic regular part
of $\sum_{\{\mc{O}\}}J_{\mc{O},\ve}^T(f)$ 
(i.e. the part of the sum which is indexed by the elliptic regular conjugacy classes $\mc{O}$)
\index{index}{conjugacy class!elliptic}\index{index}{conjugacy class!regular}%
\index{index}{trace formula!geometric side!elliptic part}%
in terms of stably invariant $\mc{O}$-expansions of lower rank groups.  We call these
lower rank groups the $\ve$-endoscopic groups \index{index}{endoscopic!group} of $\mb{G}$.  We compute them in Chapter \ref{chap:endoscopy}.

On the other hand, Arthur's trace formula for the endoscopic groups equate their $\mc{O}$-expansions with their
fine $\chi$-expansions.  
We impose in Section \ref{sec:simplifyO} a condition on the test function $f \in C(\GG, \omega)$ so that
all nonelliptic or singular terms in the geometric expansion of the $\ve$-trace formula vanish.
Using Kottwitz-Shelstad's formula, we obtain a global trace identity relating the fine $\chi$-expansion 
\index{index}{Kottwitz-Shelstad's trace formula}%
of the $\ve$-trace formula of $\mb{G}(\Af)$ with the fine $\chi$-expansions of 
the trace formulas of the $\ve$-endoscopic groups.
We compute in Section \ref{sec:finechiexp} the fine $\chi$-expansions of the groups, and we
obtain in Section \ref{sec:appKS} a global trace identity relating these spectral expansions.
We then use the global trace identity to deduce global lifting results in Sections 
\ref{sec:contributions} and \ref{sec:somegloballifts}

\section{Local Character Identities}
Let $k$ now be a local $p$-adic field. 
Let $G = \mb{G}(k)$ denote the group of $F$-points of $\mb{G}$.
Let $\ve$ be a quadratic character of $k^\times$.  Then $\ve$ defines a character of $G$ as follows:
\indexs{epsilon@$\ve$}%
\[
\ve(g) := \ve(\lambda(g)),\quad\forall g \in G.
\]
\indexs{epsilon@$\ve$}%
For any admissible (\cite{BZ}) representation $\pi$ of $G$, define a representation $\ve\pi$ on the space of $\pi$ by
\index{index}{representation!admissible}%
\[
\ve\pi := \ve\otimes \pi : g \mapsto \ve(g)\pi(g), \quad \forall g \in G.
\]
\index{symbols}{epsilon@$\ve\pi$}%
We are interested in irreducible admissible representations of $G$
which are {\bf $\ve$-invariant}:
\index{index}{representation!$\ve$-invariant}%
A representation $\pi$ of $G$ is said to be $\ve$-invariant if $\pi \cong \ve\pi$,
namely there is an invertible linear operator $A : V_\pi \rightarrow V_\pi$ on the space $V_\pi$ of $\pi$ such that
$A\pi(g) = \ve(g)\pi(g)A$ for all $g \in G$.

We fix a character $\omega$ of the center $Z$ of $G$.  Let $C(G, \omega)$ be the space of smooth
functions $f$ on $G$ which are compactly supported modulo $Z$ and satisfy:
\index{symbols}{omega@$\omega$}%
\index{symbols}{C@$C(\;,\omega)$}%
\index{symbols}{z@$Z$}%
\[
f(zg) = \omega^{-1}(z)f(g),\quad \forall z \in Z, g \in G.
\]
We fix once and for all a right-invariant Haar measure \index{index}{Haar measure}(\cite{BZ}) $dg$ on $G$.
For any admissible representation $\pi$ of $G$ with central character $\omega$, the following
convolution operator on the space of $\pi$ has finite rank:
\[
\pi(f) := \int_{Z\bs G}\pi(g)f(g)\;dg.
\]
\indexi{convolution operator}%
Suppose $\pi$ is an irreducible, admissible  $\ve$-invariant representation.  Let $A$ be a nontrivial
intertwining operator in $\Hom_G(\pi, \ve\pi)$.  By Schur's lemma (\cite{BZ}), $A^2$ is a scalar multiplication on $\pi$.
Replacing $A$ with $\lp A^2\rp^{-1/2}\!A$ if necessary, we assume that $A$ satisfies $A^2 = 1$.  Then,
$A$ is unique up to a sign.

Put $\la \pi, f\ra_A := \tr \int_{Z\bs G}\pi(g)f(g)A\;dg$.  We call the distribution $f \mapsto \la \pi, f\ra_A$
\indexs{b@$\la\;,\;\ra_A$}%
the {\bf $\ve$-character 
\index{index}{epsilon@$\ve$-character}}%
of $\pi$.  The $\ve$-character depends on the choice of the intertwining 
operator $A$.  Since $A^2 = 1$, the $\ve$-character is unique up to
a sign.

Chapter \ref{chap:local} is devoted to the classification of the $\ve$-invariant representations of $G$.  
More precisely, we utilize the global lifting results established in Chapter \ref{chap:global} to 
express the $\ve$-characters of these representations in terms of the (nontwisted) 
characters of representations of 
\index{index}{endoscopic!group}%
the twisted endoscopic groups.

\section{Statement of Main Results}
\subsection{$\ve$-Endoscopic Groups}
Let $F$ be a number field, let $\ve$ be a quadratic character of the \idele class group $\idc{F}$ of $F$.
Let $E$ be the quadratic extension of $F$ which corresponds to $\ve$ via global class field theory.
In Chapter \ref{chap:endoscopy}, we compute the quasi-split $\ve$-endoscopic groups of $\mb{G}$ over $F$.  They are as follows:
\begin{enumerate}
\item
$
\displaystyle
\mb{H_1} = \lp \GL(2) \times \rR_{E/F}\mbb{G}_m\rp',
$
\index{symbols}{H@$\mb{H_1}$}%
where $\rR_{E/F}$ is the restriction of scalars (from $E$ to $F$) functor.  The prime 
\index{index}{restriction of scalars}%
denotes the condition that the determinant of the $\GL(2)$-factor
is equal to the norm of the $\rR_{E/F}\mbb{G}_m$-factor in $\mbb{G}_m$.  In particular, the group of $F$-points
\index{index}{norm}%
of $\mb{H_1}$ is
\[
\mb{H_1}(F)
= \left\{(g, x) \in \GL(2, F)\times E^\times : \det g = \N_{E/F}x\right\}.
\]
\item
$\mb{H_2}$: the unique quasi-split reductive $F$-group whose group of $F$-points is
\[
\mb{H_2}(F) = \lp \GL(2, E) \times F^\times\rp/ \{\lp\diag(z, z),\N_{E/F}z^{-1}\rp:z\in E^\times\}.
\]
\index{symbols}{H@$\mb{H_2}$}%
\end{enumerate}
\subsection{Global Lifting}
The global lifting results are expressed in terms of {\bf global $($quasi-$)$packets}.  
A global (quasi-)packet consists of\index{index}{packet}\index{index}{packet!quasi-}\index{index}{packet!global}\index{index}{packet!local}%
ad\`elic representations which share the same unramified local components 
(with respect to fixed maximal compact subgroups) at almost every finite place.  Equivalently,
a global (quasi-)packet may be expressed as a
restricted tensor product of {\bf local} {\bf $($quasi-$)$packets}.
We say that a global packet is a {\bf discrete spectrum packet} if it contains a representation
\index{index}{packet!discrete spectrum}%
which occurs in the discrete spectrum.
We say that a discrete spectrum global packet is {\bf stable} if all its members occur with \index{index}{packet!stable}
the same nonzero multiplicity in the discrete spectrum.  Otherwise we say that it is {\bf unstable}.\index{index}{packet!unstable}

Associated with any algebraic $F$-group $\mb{H}$ is an {\bf $L$-group}, denoted ${}^L H$ (see \cite{B}).
\index{index}{L@$L$-!group}%
It is a split extension of a complex Lie group by the Weil group $W_F$ of $F$ (see \cite{T}).
\index{index}{Weil group}%
Each unramified local component 
of a global (quasi-)packet of $\mb{H}(\Af)$
is parametrized by a conjugacy class in ${}^L H$, called a {\bf \FH class} (see \cite{B}). 
\index{index}{Frobenius-Hecke class}%
We compute these conjugacy classes for the $\ve$-endoscopic groups in Chapter \ref{FHclass}.

As we shall discuss in Chapter \ref{chap:endoscopy},
associated with each $\ve$-endoscopic group $\mb{H}$ is an $L$-group embedding $\xi_H : {}^L H \rightarrow {}^L G$.
\indexs{xi@$\xi$}%
We say that a global packet $\{\Pi_H\}$ of an $\ve$-endoscopic group $\mb{H}(\Af)$ {\bf lifts} to a packet
\indexi{lift}%
$\{\Pi\}$ of $\mb{G}(\Af)$, if $\xi_H$ maps the \FH classes in ${}^L H$ parametrizing $\{\Pi_H\}$ to the \FH
classes in ${}^L G$ parametrizing $\{\Pi\}$.
\indexi{packet}%

Let $\pi$ be an irreducible automorphic representation of $\GL(2, \Af)$, $\chi$ a character of $\idc{E}$.
Let $V$ be the set of places of $F$.
For any place $v \in V$, let
$\pi_v\otimes_1\chi_v$ denote the representation of $\mb{H_1}(F_v)$ defined as follows.
\indexs{times@$\otimes_1$}%
\[
\pi_v\otimes_1 \chi_v : (g, x) \mapsto \chi_v(x)\pi_v(g),\quad\forall (g, x) \in \mb{H_1}(F_v).
\]
Then, $\pi_v\otimes_1\chi_v$ may be a reducible representation of $\mb{H_1}(F_v)$ with at most two 
irreducible constituents (\cite{LL}).  
The set of irreducible constituents of $\pi_v\otimes_1\chi_v$ form a local packet.
If $\pi_v\otimes_1 \chi_v$ is reducible, we name its two constituents
$\pi_v^+, \pi_v^-$.
\indexs{pi@$\pi^+$, $\pi^-$}%
For each finite place $v$, let $K_{1,v}$ be the hyperspecial maximal compact subgroup
$\mb{H_1}(\mc{O}_v)$, where $\mc{O}_v$ is the ring of integers in $F_v$.
\indexs{K@$K_{1, v}$}%
We say that a representation of $\mb{H_1}(F_v)$ is unramified if it contains a nonzero $K_{1,v}$-fixed vector.
\index{index}{representation!unramified}%
If $\pi_v\otimes_1 \chi_v$ is unramified,
we let $\pi_v^+$ denote its unique irreducible unramified constituent.
If $\pi_v\otimes_1\chi_v$ is irreducible, we put $\pi_v^+ := \pi_v\otimes_1\chi_v$ and
$\pi_v^- := 0$.
From \cite{LL}, a global packet of $\mb{H_1}(\Af)$ has the form
\[
\pi\otimes_1 \chi
:= \left\{ \otimes_v\pi_v : \pi_v \in \{\pi_v^+, \pi_v^-\},\;\forall v;\;
\pi_v = \pi_v^+ \text{ for almost all } v\right\}.
\]
For $f_1 \in C(\mb{H_1}(\Af), \omega)$,
we have $\displaystyle \la \pi\otimes_1\chi, f_1\ra = \sum_{\pi' \in \pi\otimes_1\chi} \la \pi', f_1\ra$.
\indexs{b@$\la\;,\;\ra$}%

The group $\mb{H_2}(\Af)$ is a quotient of $\GL(2, \Ae)\times \Af^\times$.   Since rigidity, or strong multiplicity one
theorem, holds for $\GL(2)$ and $\GL(1)$, it also holds for $\mb{H_2}$.
Consequently, each global packet of $\mb{H_2}(\Af)$ consists of a single irreducible
automorphic representation.  For any automorphic representation $\pi$ of $\GL(2, \Ae)$ and character $\mu$ of $\idc{F}$ 
such that the central character $\omega_\pi$ of $\pi$ is equal to $\mu\circ\N_{E/F}$, let $\pi\otimes_2 \mu$ denote the
representation
\indexs{times@$\otimes_2$}%
\[
\pi\otimes_2 \mu : (g, x) \mapsto \mu(x)\pi(g),\quad\forall (g, x) \in \mb{H_2}(\Af).
\]
\subsubsection{Summary of Global Lifting}
\index{index}{lifting!global}%
We now give a summary of our global lifting results.  For simplicity's sake,  
some of these statements are accurate only up to certain nuanced conditions.  We mark them with an *.
The full, accurate versions of these statements may be found in Sections \ref{sec:contributions}, 
\ref{sec:somegloballifts}, and \ref{sec:finalwords}.

A more concise version of this summary is in Appendix \ref{chap:globaltables}.

We classify packets of $\GSp(2, \Af)$ by
stating what automorphic representations of $\GL(4, \Af)$ they lift to via the natural $L$-group embedding 
$\GSp(2, \mbb{C}) \rightarrow \GL(4, \mbb{C})$.
This classification is not one to one.  There are cases in which two inequivalent
packets of $\GSp(2, \Af)$ lift to the same representation in $\GL(4, \Af)$.


We say that an automorphic representation $\pi$ of $\GL(2, \Af)$ is $E$-{\bf monomial} if it is the monomial 
\index{index}{representation!monomial}%
representation $\pi(\theta)$ associated with a character $\theta$ of $\idc{E}$, or more precisely with 
the induced two-dimensional representation ${\rm Ind}_{W_{E/E}}^{\smash[b]{W_{E/F}}}(\theta)$ of $W_F$ (see \cite{JL}, \cite{K}).
\indexs{I@${\rm Ind}_{W_{E/E}}^{\smash[b]{W_{E/F}}}$}%
The first global lifting result which we examine involves global packets of $\mb{H_1}$ of the form 
$\pi_1 = \pi\otimes_1 \chi$, where $\pi$ is a cuspidal non-$E$-monomial or one dimensional automorphic representation 
of $\GL(2, \Af)$, and $\chi = \mu\circ\N_{E/F}$ for some character $\mu$ of $\idc{F}$.
From the way $\mb{H_1}$ is defined, $\pi\otimes_1 \mu\circ\N_{E/F}$ is equal to $\mu\pi\otimes_11$;
hence, we assume without loss of generality that $\pi_1 = \pi\otimes_1 1$.  Let $\omega_\pi$ be the central character of $\pi$.
Let $B_{E/F}\pi$ be the representation of $\GL(2, \Ae)$ which is the base change of $\pi$  
\index{index}{base change}%
(see \cite{L1}, \cite{F4}).  Then,
$\omega_{B_{E/F}\pi} = \omega_\pi\circ\N_{E/F}$.
\begin{rprop}
The packets $\pi\otimes_1 1$ of $\;\mb{H_1}(\Af)$ and $B_{E/F}\pi\otimes_2 \omega_{\pi}$ of $\;\mb{H_2}(\Af)$ 
lift to a packet $\{\Pi\}$ of $\GSp(2, \Af)$, which in turn lifts to the induced representation $I_{(2,2)}(\pi, \ve\pi)$ of
$\GL(4, \Af)$.   
\begin{remark}(i) 
According to {\rm \cite[Sect. 5]{A}}, the packet $\{\Pi\}$ is unstable.  In the case where the central character
of $\{\Pi\}$ is trivial, the instability of $\{\Pi\}$ is proven in {\rm \cite[V. 10.]{F1}}.
(ii)
In \cite{A}, Arthur denotes $\{\Pi\}$ by the symbol $(\pi\boxtimes 1)\boxplus(\ve\pi\boxtimes 1)$ if
\indexs{times@$(\;\boxtimes\;)\boxplus(\;\boxtimes\;)$}%
$\pi$ is cuspidal non-$E$-monomial and calls it a Yoshida type packet.  If $\pi = \mu 1_{\GL(2, \Af)}$
\index{index}{Yoshida}%
for some character $\mu$ of $\idc{F}$, Arthur calls  $\{\Pi\}$ a Howe, Piatetski-Shapiro type packet
\index{index}{Howe, Piatetski-Shapiro}%
and denotes it by the symbol $(\mu\boxtimes\nu(2))\boxplus (\ve\mu\boxtimes\nu(2))$, 
where $\nu$ is the (normalized) absolute value function on $\Af$.
\end{remark}
\end{rprop}

Let $\sigma$ be the generator of $\Gal(E/F)$.
For any character $\theta$ of $\idc{E}$, put ${}^\sigma \theta := \theta\circ\sigma$.
\begin{rprop}
Let $\theta, \chi$ be characters of $\idc{E}$ such that none of $\theta, \chi, \theta\chi, \theta\,\sig\chi$ 
is fixed by the action of $\Gal(E/F)$.
Let $\pi(\theta)$, $\pi(\chi)$  be the cuspidal monomial representations of $\GL(2, \Af)$ associated 
with $\theta$, $\chi$, respectively.
The global packets 
\begin{itemize}
\item 
$\pi(\theta)\otimes_1\chi$, 
\item
$\pi(\theta)\otimes_1 {}^\sigma \chi$, 
\item
$\pi(\chi)\otimes_1\theta$,
\item 
$\pi(\chi)\otimes_1{}^\sigma \theta$ 
\end{itemize}
all lift to an unstable packet $\{\Pi\}$, which in turn lifts to the induced representation
\dt \[I_{(2, 2)}(\pi(\theta\chi), \pi(\theta\;{}^\sigma\!\chi))\] of $\GL(4, \Af)$.
Moreover, no discrete spectrum representation of $\mb{H_2}(\Af)$ lifts to $\{\Pi\}$.
\begin{remark}
Under the notation of \cite{A}, $\{\Pi\}$ is the unstable Yoshida type packet denoted by
$(\pi(\theta\chi)\boxtimes 1)\boxplus(\pi(\theta\;{}^\sigma\!\chi)\boxtimes 1)$.\end{remark}
\end{rprop}

For an automorphic representation $\tau$ of $\GL(2, \Ae)$, 
let $\pi(\tau)$ denote the automorphic representation of $\GL(4, \Af)$ which is obtained from
$\tau$ via the twisted endoscopic lifting from $\rR_{E/F}\GL(2)$ to $\GL(4)$ (see \cite[Sect. 3.6]{AC}).
\begin{rprop'}
Let $\chi$ be a character of $\idc{E}$ not fixed by $\Gal(E/F)$.  Let $\pi$ be an irreducible,
cuspidal, automorphic representation of $\GL(2, \Af)$ which is not $E$-monomial.
Suppose there does not exist a character $\ve'$ of $\idc{F}$ such that
${}^\sigma\chi/\chi$ is equal to $\ve'\circ\N_{E/F}$ and $\pi \cong \ve'\pi$.
The packet $\pi\otimes_1 \chi$ lifts to a stable packet $\{\Pi\}$ which lifts to 
the cuspidal automorphic representation $\pi(\chi B_{E/F}\tau)$ of $\GL(4, \Af)$.
No discrete spectrum representation of $\mb{H_2}(\Af)$ lifts to $\{\Pi\}$.
\end{rprop'}
\begin{rprop'}
Let $\chi$ be a character of $\idc{E}$ not fixed by $\Gal(E/F)$.
The one dimensional representation $1_{\GL(2, \Af)}\otimes_1 \chi$ lifts to a stable quasi-packet, denoted 
by \[\{L(\nu\ve, \nu^{-1/2}\pi(\chi))\},\] 
which lifts to the Langlands quotient representation $J(\nu^{1/2}\pi(\chi), \nu^{-1/2}\pi(\chi))$ of
$\GL(4, \Af)$.  No discrete spectrum representation of $\mb{H_2}(\Af)$ lifts to $\{L(\nu\ve, \nu^{-1/2}\pi(\chi))\}$.
\end{rprop'}

For a representation $\pi$ of $\GL(2, \Ae)$, let ${}^\sigma \pi$ denote the representation
\[
{}^\sigma \pi : (g_{ij}) \mapsto \pi(\sigma g_{ij}),\quad \forall (g_{ij}) \in \GL(2, \Ae).
\]
\begin{rprop'}
Let $\pi$ be an irreducible discrete spectrum automorphic representation of $\GL(2, \Ae)$,
$\mu$ a character of $\idc{F}$, such that $\pi \neq {}^\sigma \pi$ and $\omega_{\pi} = \mu\circ\N_{E/F}$.
\begin{itemize}
\item
If $\pi$ is cuspidal,
the representation $\pi\otimes_2 \mu$ of $\mb{H_2}(\Af)$ lifts to a stable packet $\{\Pi\}$ which in turn lifts
to the cuspidal automorphic representation $\pi(\tau)$ of $\GL(4, \Af)$.
\item
Suppose $\pi = \chi 1_{\GL(2, \Ae)}$, where $\chi$ is a character of $\idc{E}$ such that $\chi \neq {}^\sigma \chi$ and
${{}^\sigma\chi}/{\chi} = \ve'\circ\N_{E/F}$ for some nontrivial 
quadratic character $\ve'$ of $\idc{F}$.  Let $\zeta$ be either $\ve'$ or $\ve\ve'$.  Then,
$\pi\otimes_2 \chi|_{\Af^\times}\cdot\zeta$ lifts to a stable quasi-packet $\{\Pi\}$
which in turn lifts to the Langlands quotient representation $J(\nu^{1/2}\pi(\chi), \nu^{-1/2}\pi(\chi))$ of $\GL(4, \Af)$.
We denote $\{\Pi\}$ by $\{L(\nu\ve\zeta, \nu^{-1/2}\pi(\chi))\}$.
\end{itemize}
In either case, no discrete spectrum packet of $\mb{H_1}(\Af)$ lifts to the $($quasi-$)$packet $\{\Pi\}$.
\end{rprop'}
Additional global lifting results involving parabolically induced representations of \dt $\mb{G}(\Af)$ are stated in
Section \ref{sec:somegloballifts}.
\begin{rprop'}
Suppose $\{\Pi\}$ is a global packet of $\GSp(2, \Af)$ which is the lift of some automorphic 
representation/packet of one of the $\ve$-endoscopic groups.  Then, $\{\Pi\}$ contains at least
one representation which is $\ve$-invariant.

If $\{\Pi\}$ is a discrete spectrum global packet of $\GSp(2, \Af)$ which is not the lift of 
any discrete spectrum representation/packet of the $\ve$-endoscopic groups, then $\{\Pi\}$
does not contain any discrete spectrum automorphic representation which is $\ve$-invariant.
\end{rprop'}
\subsection{Local Lifting}
Let $k$ be a local nonarchimedean field.  Let $\ve$ be a nontrivial quadratic character of $k^\times$.
Let $\K$ be the quadratic extension of $k$ which corresponds via local class field theory to $\ve$.
Let  $G = \GSp(2, k)$.
\indexs{G@$\GSp(2)$}%
Let $Z = \{\diag(z, z, z, z):z\in k^\times\}$ be the center of $G$.
Let
\begin{enumerate}
\item
$\displaystyle
H_1 = 
 \lp \GL(2, k)\times \K^\times\rp' 
= \left\{(g, x) \in \GL(2,  k)\times \K^\times: \det g = N_{\K/k}x\right\},
$
\indexs{H@$\mb{H_1}$}%
\item
$\displaystyle
H_2 
= \lp \GL(2, \K)\times k^\times\rp/
\left\{\lp\diag(z, z), \N_{\K/k}z^{-1}\rp : z \in \K^\times\right\}.
$
\indexs{H@$\mb{H_2}$}%
\end{enumerate}
Recall that the character $\ve$ defines a character of $G$ via $\ve(g) := \ve(\lambda(g)),\;\forall g \in G$.  
For any representation $\pi$ of $G$, 
we say that $\pi$ is $\ve$-invariant if $\pi \cong \ve\pi$.
\index{index}{representation!epsilon@$\ve$-invariant}%

Fix a character $\omega$ of $Z$.  Recall that $C(G, \omega)$ is the space of smooth functions
on $G$ which are compactly supported modulo $Z$ and transform under $Z$ via $\omega^{-1}$.
We fix once and for all a Haar measure $dg$ on $G$.  
For any irreducible $\ve$-invariant admissible representation $\pi$ of $G$ with central character $\omega$,
\index{index}{representation!admissible}%
and for any function $f \in C(G, \omega)$, recall that
\[
\la \pi, f\ra_A := \tr \int_{Z\bs G} \pi(g)f(g)A\;dg,
\]
\indexs{b@$\la\;,\;\ra_A$}%
where $A$ is an intertwining operator \index{index}{intertwining operator} in $\Hom_G(\pi, \ve\pi)$ such that $A^2 = 1$.

For $i = 1, 2$, let $Z_0(H_i) $ denote the maximal $k$-split component of the center of $H_i$.
We shall see in Section \ref{sec:normmapping} that $Z_0(H_i)$ is isomorphic to $Z$; 
thus, $\omega$ is defines a character on $Z_0(H_i)$.  Let $C(H_i, \omega)$ denote the space of smooth functions 
on $H_i$ which are compactly supported modulo $Z_0(H_i)$ and transform under $Z_0(H_i)$ via $\omega^{-1}$.
Fix once and for all a Haar measure $dh_i$ on $H_i$.  For any admissible representation
$\pi_i$ of $H_i$, and for any function $f_i \in C(H_i, \omega)$, put
\[
\la \pi_i, f_i\ra := \tr \int_{Z_0(H_i)\bs H_i} \pi_i(h_i)f_i(h_i)\;dh_i.
\]
\indexs{b@$\la\;,\;\ra$}%
Using global lifting results,
we deduce in Chapter \ref{chap:local} local character identities
relating quantities of the form $\la \pi, f\ra_A$ and $\la \pi_i, f_i\ra$ for {\bf matching} functions $f, f_i$.
\index{index}{matching functions}%
Functions are matching if their orbital integrals are compatible with the norm correspondence
\index{index}{norm correspondence}%
between $G$ and its $\ve$-endoscopic groups.  We shall discuss norm correspondence in Section \ref{sec:normmapping}.
In general, the existence of matching functions is conjectural.  
It is related to a conjecture called the Fundamental Lemma.
\index{index}{Fundamental Lemma}%
We prove the Fundamental Lemma in the context of $G$ and its $\ve$-endoscopic groups in Appendix \ref{chap:fundlemma}.  

\subsubsection{Notations:}
For a character $\theta$ of $\K^\times$,
let $\pi(\theta)$ denote the cuspidal monomial representation of $\GL(2, k)$ associated with $\theta$, or more
\indexs{pi@$\pi({\rm char.})$}%
precisely with the representation ${\rm Ind}_{W_{\K/\K}}^{\smash[b]{W_{\K/k}}}(\theta)$ of $W_k$.
If $\theta$ is invariant under the action of $\Gal(\K/k)$, $\pi(\theta)$ is induced; otherwise, 
$\pi(\theta)$ is cuspidal (\cite{JL}).

Let $\nu$ be the normalized absolute value function on $k$; that is,
$\nu(x) = q^{-\ord x}$, where $q$ is the cardinality of the residue field of $k$ and $\ord x$ is the
$p$-adic valuation of $x$.
\indexs{nu@$\nu$}%
\indexs{o@$\ord$}\indexs{q@$q$}%

Let $P_0$ be the minimal upper triangular parabolic subgroup of $G$.  It has a decomposition
\indexs{P@$P_0$}%
$P_0 = TN$, where $T$ is the maximal diagonal torus of $G$ and $N$ is the unipotent component of $P_0$.
For characters $\mu_1, \mu_2, \mu$ of $k^\times$, let $\mu_1\times\mu_2\rtimes\mu$ denote the
\indexs{t@$\;\times\;\rtimes\;$}%
representation of $\GL(2, k)$ (normalizedly) parabolically induced from the following representation of $P_0$:
\[
\mu_1\otimes\mu_2\otimes\mu : \dmfour{a}{b}{\lambda/b}{\lambda/a}n\mapsto \mu_1(a)\mu_2(b)\mu(\lambda) ,
\quad \forall a, b, \lambda \in k^\times, n \in N.
\]

Let $P_\beta$  be the standard Heisenberg parabolic subgroup of $G$ containing $P_0$.  Its Levi component is
\indexs{P@$P_\beta$}%
\index{index}{parabolic subgroup!Heisenberg}%
\[
M_\beta = \left\{\lp\lmx a &&\\&g&\\&&\frac{\det g}{a}\rmx\rp:a \in k^\times, g \in \GL(2, k)\right\}.
\indexs{M@$M_\beta$}%
\]
Let $N_\beta$ be the unipotent component of $P_\beta$.  It is a Heisenberg group.  For a character $\mu$ of $k^\times$ and an
admissible representation $\pi$ of $\GL(2, k)$, let $\mu\rtimes\pi$ denote the representation of $G$ (normalizedly)
\indexs{times@$\rtimes$}%
parabolically induced from the following representation of $P_\beta$:
\[
\mu\otimes\pi : \lp\lsm a &&\\&g&\\&&\frac{\det g}{a}\rsm\rp n \mapsto \mu(a)\pi(g),
\quad\forall a \in k^\times, g \in \GL(2, k), n \in N_\beta.
\]
\subsubsection{Summary of Local Character Identities}
Below is a summary of our twisted character identities for the $\ve$-invariant representations of $G = \GSp(2, k)$.\\
\quad\\
1.  Let $\tau$ be an irreducible, cuspidal, non-$\K$-monomial, admissible $\GL(2, k)$-module whose central character
$\omega_\tau$ is trivial.
Let $\{\pi^+, \pi^-\}$ be the local packet of $G$ which lifts to the induced representation
$I_{(2, 2)}(\tau, \ve\tau)$ of $\PGL(4, k)$ according to \cite{F1}.  
There exist intertwining operators $A^+$ in $\Hom_G(\pi^+, \ve\pi^+)$, and $A^-$ in $\Hom_G(\pi^-, \ve\pi^-)$,
such that the following holds for matching functions:
\[
\begin{split}
\la \pi^+, f\ra_{A^+} + \la \pi^-, f\ra_{A^-} &= \la \tau\otimes_1 1,\, f_1\ra,\\
\la \pi^+, f\ra_{A^+} - \la \pi^-, f\ra_{A^-} &= \la B_{\K/k}\tau\otimes_2 \omega_\tau,\, f_2\ra.
\end{split}
\]
2. Let $\theta$, $\chi$ be characters of $\K^\times$ such that none of $\theta$, $\chi$, $\theta\chi$, $\theta\,\sig\chi$
is fixed by the action of $\Gal(\K/k) = \la \sigma \ra$.
\indexs{sigma@$\sigma$}%
Suppose $(\theta\chi)|_{k^\times} = (\theta\,\sig\chi)|_{k^\times} = \ve$, then
the distinct $\K$-monomial representations $\pi(\theta\chi)$, $\pi(\theta\,\sig\chi)$ 
are cuspidal with trivial central characters.

Let $\{\pi^+, \pi^-\}$ be the local packet which lifts to the induced $\PGL(4, k)$-module
\dt \[I_{(2, 2)}(\pi(\theta\chi),\,\pi(\theta\,\sig\chi)).\]  
There exist intertwining operators $A^+$ in \dt $\Hom_G(\pi^+,\,\ve\pi^+)$, $A^-$ in $\Hom_G(\pi^-,\,\ve\pi^-)$,
and a constant $\ep = \pm 1$, such that the following holds for matching functions:
\[
\begin{split}
\la \pi^+, f\ra_{A^+} + \la \pi^-, f\ra_{A^-} &= \;\;\la \pi(\theta)\otimes_1 \chi,\, f_1\ra,\\
\la \pi^+, f\ra_{A^+} - \la \pi^-, f\ra_{A^-} &= \ep\la \pi(\chi)\otimes_1 \theta,\, f_1\ra.
\end{split}
\]
3. Let $\chi$ be a character of $\K^\times$ not fixed by $\Gal(\K/k)$.  Let $\pi^+$, $\pi^-$ be the two 
inequivalent, tempered, irreducible subrepresentations of $1 \rtimes \pi(\chi)$ (\cite{ST}).  
There exist intertwining operators $A^+$ in 
$\Hom_G(\pi^+,\, \ve\pi^+)$ and $A^-$ in $\Hom_G(\pi^-,\, \ve\pi^-)$, such that the following holds for matching functions:
\[
\la \pi^+, f\ra_{A^+} + \la \pi^-, f\ra_{A^-} = \la \pi(\chi)\otimes_1 1,\, f_1\ra.
\]
4.  
Let $1_{2, l}$ ($l = k, \K$) denote the nontempered trivial representation of $\GL(2, l)$.
\indexs{1@$1_2$}%
Let ${\rm St}_{2,l}$ denote the square integrable Steinberg representation of $\GL(2, l)$.  
\indexs{S@${\rm St}_2$}%

Suppose $\ve$ is unramified.  Let $\xi$ be a (possibly trivial) quadratic character of $k^\times$.
Let $\delta = \delta(\nu^{1/2}\ve{\rm St}_{2, k}, \nu^{-1/2}\xi)$ be the unique
\indexs{delta@$\delta({\rm repn., char.})$}%
square integrable constituent of the induced  $G$-module $\nu^{1/2}\ve{\rm St}_{2, k}\rtimes\nu^{-1/2}\xi$.
Let $\delta^- = \delta^-(\nu^{1/2}\ve{\rm St}_{2, k}, \nu^{-1/2}\xi)$ be the cuspidal member of the local packet containing 
\indexs{delta@$\delta^-({\rm repn., char.})$}%
$\delta$.
In particular, the local packet $\{\delta, \delta^-\}$ lifts to $I_{(2, 2)}(\xi\,{\rm St}_{2, k},\, \ve\xi\,{\rm St}_{2, k})$ 
of $\PGL(4, k)$.

Let $L = L(\nu\ve, \ve\rtimes\nu^{-1/2}\xi)$ be the nontempered quotient of $\nu^{1/2}\ve 1_{2, k}\rtimes\nu^{-1/2}\xi$ (\cite{ST}).
\indexs{L@$L({\rm char., repn.})$}%
The local quasi-packet $\{L, \delta^-\}$ lifts to $I_{(2, 2)}(\xi 1_{2,k},\, \ve\xi 1_{2, k})$ of $\PGL(4, k)$.

There exist intertwining operators $A^L \in \Hom_G(L, \ve L)$, $A^\delta \in \Hom_G(\delta, \ve\delta)$, 
and $A^- \in \Hom_G(\delta^-, \ve\delta^-)$,
such that the following holds for matching functions:
\[
\begin{split}
\la L, f\ra_{A^L} + \la \delta^-, f\ra_{A^-} &= \la \xi 1_{2, k}\otimes_1 1,\, f_1\ra,\\
\la L, f\ra_{A^L} - \la \delta^-, f\ra_{A^-} &= \la (\xi\circ\N_{\K/k})\,1_{2, \K}\otimes_2 1,\, f_2\ra;\\
\quad\\
\la \delta, f\ra_{A^\delta} - \la \delta^-, f\ra_{A^-} &= \la \xi\,{\rm St}_{2, k}\otimes_1 1,\, f_1\ra,\\
\la \delta, f\ra_{A^\delta} + \la \delta^-, f\ra_{A^-} &= \la (\xi \circ \N_{\K/k})\,{\rm St}_{2, \K}\otimes_2 1,\, f_2\ra.
\end{split}
\]
5. Let $\tau$ be an irreducible, cuspidal, non-$\K$-monomial, admissible representation of the group \dt $\GL(2, k)$,
$\chi$ a character of $\K^\times$, such that: There does not exist a character $\ve'$ of $k^\times$ for which
${}^\sigma\chi/\chi = \ve'\circ\N_{\K/k}$ and $\tau \cong \ve'\tau$.
Then, there exists a collection $[\pi]$ of distinct $\ve$-invariant, irreducible, cuspidal, admissible 
\indexs{b@$[{\rm repn.}]$}%
$G$-modules, and an intertwining operator $A(\pi') \in \Hom_G(\pi',\,\ve\pi')$ for each $\pi' \in [\pi]$,
\indexs{A@$A(\pi)$}%
such that the following holds for matching functions:
\[
\sum_{\pi' \in [\pi]}\la \pi', f\ra_{A(\pi')} = \la \tau\otimes_1\chi,\,f_1\ra.
\]
{\sc Remark.}
From the proof of the above twisted character identity, one should expect that $[\pi]$ is a subset of a local packet
which lifts to $\pi(\chi B_{\K/k}\tau)$, the cuspidal $\GL(4, k)$-module which is the twisted endoscopic lift of
\indexs{pi@$\pi({\rm repn.})$}%
the representation $\chi B_{\K/k}\tau$ of $\GL(2, \K)$ (see \cite{AC}).\\
\quad\\
6. Let $\chi$ be a character of $\K^\times$ which is not fixed by $\Gal(\K/k)$.  Let $L = L(\nu\ve, \nu^{-1/2}\pi(\chi))$ be 
the nontempered quotient of the parabolically induced $G$-module $\nu\ve\rtimes\nu^{-1/2}\pi(\chi)$.
Let \dt $\delta = \delta(\nu\ve, \nu^{-1/2}\pi(\chi))$ be the square integrable subrepresentation of $\nu\ve\rtimes\nu^{-1/2}\pi(\chi)$.
\indexs{delta@$\delta({\rm repn., char.})$}%
There exist intertwining operators $A^L$ in $\Hom_G(L, \ve L)$, and $A^\delta$ in  $\Hom_G(\delta, \ve\delta)$ such that the
following holds for matching functions:
\[
\begin{split}
\la L, f\ra_{A^L} &= \la 1_{2, k}\otimes_1 \chi,\, f_1\ra,\\
\la \delta, f\ra_{A^\delta} &= \la {\rm St}_{2, k}\otimes_1 \chi,\,f_1\ra.
\end{split}
\]
7. Let $\tau$ be an irreducible, cuspidal, admissible representation of $\GL(2, \K)$, not fixed by the action of $\Gal(\K/k)$,
such that its central character $\omega_\tau$ is equal to $\mu\circ\N_{\K/k}$ for some character $\mu$ of $k^\times$.
Then, there exists a collection $[\pi]$ of distinct $\ve$-invariant, irreducible, cuspidal, admissible 
\indexs{b@$[{\rm repn.}]$}%
$G$-modules, and an intertwining operator $A(\pi')$ in $\Hom_G(\pi',\,\ve\pi')$ for each $\pi' \in [\pi]$, 
\indexs{A@$A(\pi)$}%
such that the following holds for matching functions:
\[
\sum_{\pi' \in [\pi]}\la \pi', f\ra_{A(\pi')} = \la \tau\otimes_2\mu,\,f_2\ra.
\]
{\sc Remark.}
From the proof of the above twisted character identity, one should expect that $[\pi]$ is a subset of a local packet
which lifts to $\pi(\tau)$, the cuspidal $\GL(4, k)$-module which is the twisted endoscopic lift of $\tau$.\\
\quad\\
8. Let $\tau$ be a cuspidal non-$\K$-monomial, or one dimensional, irreducible, admissible \dt 
representation of $\GL(2, k)$.
There exists an operator $A$, intertwining the induced $G$-module $\ve\rtimes \tau$ with $\ve\rtimes\ve\tau$,
such that
\[
\la \ve\rtimes\tau,\, f\ra_A = \la B_{\K/k}\tau\otimes_2 \omega_\tau\ve,\,f_2\ra
\]
for matching functions.\\
\quad\\
9. Suppose the field extension $\K/k$ is unramified.  Let $\chi$  be a character of $\K^\times$ such that $\chi \neq \sig\chi$ and
$\sig\chi/\chi = \ve'\circ\N_{\K/k}$ for some nontrivial quadratic character $\ve'$ of $k^\times$.  Let $\zeta = \ve'$ or $\ve'\ve$.  
Let $L = L(\ve\zeta\nu, \nu^{-1/2}\pi(\chi))$ and $\delta = \delta(\ve\zeta\nu, \nu^{-1/2}\pi(\chi))$.
There exist intertwining operators $A^L$ in $\Hom_G(L, \ve L)$, and $A^\delta$ in $\Hom_G(\delta, \ve\delta)$, such that
\indexs{A@$A^L$}\indexs{A@$A^\delta$}%
the following holds for matching functions:
\[
\begin{split}
\la L, f\ra_{A^L} &= \la \chi 1_{2, \K}\otimes_2 \chi|_{k^\times}\cdot\zeta,\, f_2\ra,\\
\la \delta, f\ra_{A^\delta} &= \la \chi\, {\rm St}_{2,\K}\otimes_2 \chi|_{k^\times}\cdot\zeta,\, f_2\ra.
\end{split}
\]
Additional character identities involving parabolically induced $G$-modules are stated in
Sections \ref{sec:localH1} and \ref{sec:localH2}.
\section{Assumptions and Restrictions}
Before we continue further, we list here the assumptions on which this work rests and the restrictions 
on our results:
\subsection*{Assumptions}
\begin{enumerate}
\item The multiplicity one theorem for the discrete spectrum of $\GSp(2)$.
As mentioned earlier in this chapter, the multiplicity one theorem has been established for $\PGSp(2)$
by Flicker in \cite{F1}.
\item The results of Arthur's announced in \cite{A}, including the multiplicity one theorem for the automorphic
representations of $\GSp(2)$ which are associated with semisimple Arthur parameters.
\item The existence of the transfer of functions from $\GSp(2)$ to its $\ve$-twisted endoscopic groups,
  and the precise forms which matching functions take (see remarks after the proofs of Theorem \ref{thm:contribH1}
  and Corollary \ref{corollary:contribH2discrete}).  In the nontwisted case, Waldspurger has shown that
the Fundamental Lemma implies the transfer of functions \cite{WA}.  
In all likelihood, the twisted analogue of Waldspurger's theorem also holds.
\end{enumerate}
\subsection*{Restrictions}
We are only able to classify those discrete spectrum, $\ve$-invariant, automorphic representations of $\GSp(2)$ 
which
have at least two elliptic local components.
The restriction ``two elliptic components'' could presumably be reduced to ``one elliptic component,'' through a technique 
developed in \cite{F3} which involves a class of test functions called regular functions.  
\indexi{function!regular}%
It also seems likely that
the work of Arthur's announced in \cite{A} contains a method to do away with any restriction on the automorphic representations.


\chapter{$\ve$-Endoscopy for $\GSp(2)$}\label{chap:endoscopy}
\section{Endoscopic Data}
\indexi{endoscopic!data|(}\indexi{endoscopic!group|(}%
\indexi{splitting|(}\indexi{maximal torus|(}\indexi{root datum|(}%
In this chapter, $F$ is either a local field of characteristic zero or a number field.  
If $F$ is a number field, let $\idc{F}$ denote the \idele class group of $F$.
For a reductive $F$-group $\mb{H}$, let $\hat{H}$ denote the identity component of the 
$L$-group ${}^L H$ of $\mb{H}$ (see \cite{B}).
\indexs{H@$\hat{H}$}%
\indexi{L@$L$-!group}%
Let $\mb{G}$ be the reductive $F$-group $\GSp(2)$.

Kottwitz and Shelstad have defined in \cite{KS} the endoscopic data attached to a
\indexi{endoscopic!data}%
triple $(\mb{G}, \theta, \mb{a})$, where $\theta$ is an automorphism of $\mb{G}$ over $F$.
We now recall this definition in the special case where $\theta$ is trivial.
\begin{define}\label{define:endoscopy}
{\rm 
The quadruple $(\mb{H}, \mathcal{H}, s, \xi)$ is a set of 
\indexs{H@$(\mb{H}, \mathcal{H}, s, \xi)$}%
{\bf $\ve$-endoscopic data} (or in the terminology of \cite{KS}, the
endoscopic data attached to the triple $(\mb{G}, 1, \mb{a})$) if it satisfies each of the following conditions:
\indexs{G@$(\mb{G}, 1, \mb{a})$}%
  \begin{enumerate}
  \item
    $\mb{H}$ is a quasi-split reductive group over $F$.
  \item
    $\mathcal{H}$ is a split extension of $W_F$ by $\hat{H}$ such that the $L$-action
\indexi{L@$L$-!action}\indexi{L@$L$-!group}\indexs{H@$\mathcal{H}$}%
    $\rho_\mathcal{H}$ of $W_F$ on $\hat{H}$ determined by this extension coincides with 
\indexs{rho@$\rho_{\mathcal{H}}$}%
    the Weil group action determined by ${}^L H$ (see \cite{B}).
  \item
    $s$ is a semisimple element in $\hat{G}$.
\indexs{s@$s$}%
  \item
    $\xi : \mathcal{H}\rightarrow {}^L G$ is an $L$-homomorphism satisfying the following two conditions.
\indexi{L@$L$-!homomorphism}\indexs{xi@$\xi$}%
    \begin{enumerate}
  \item\label{condition:endo4a}
    $\text{Int}(s)\circ\xi = a\cdot\xi$;
  \item
    $\xi$ maps $\hat{H}$ isomorphically onto the identity component of $\text{Cent}(s,\hat{G})$,
    the centralizer of $s$ in $\hat{G}$.
    \end{enumerate}
\end{enumerate}
We call $\mb{H}$ an {\bf $\ve$-endoscopic group} of $\mb{G}$.}\\
\indexi{endoscopic!group}%
\begin{remark}
The action $\rho$ of $W_F$ on $\hat{H}$ determined by $\mc{H}$ may not fix any choice of splitting for $\hat{H}$.
\indexi{splitting}%
However, there exists a set $\{h_w \in \hat{H}\}_{w \in W_F}$ such that the action
\indexs{h@$h_w$}%
$\rho_{\mc{H}} : w \mapsto \Int(h_w)\rho(w)$ does fix a choice of splitting for $\hat{H}$.  
As such, $\rho_{\mc{H}}$ constitutes an $L$-action.  For our purpose, $\mc{H}$ may be taken to be
${}^L H$; hence, condition 2 is largely moot.
\end{remark}
\end{define}


\subsection{Elliptic $\ve$-Endoscopic Data of $\mb{G}$}
In general, an $L$-group is a semidirect product $\hat{H} \rtimes W_F$.  
As we shall show, every group considered in this work is split over some quadratic extension of $F$.
For simplicity, we let ${}^L H$ denote the
{finite Galois form} $\hat{H} \rtimes \Gal(K/F)$ of the $L$-group, where
$K$ is the smallest extension of $F$ over which $\mb{H}$ splits.
In particular, since $\mb{G}$ is $F$-split, 
we let ${}^L G$ denote $\hat{G}$, which we identify with $\GSp(2, \mathbb{C})$ (see \cite{B}, \cite{A}).


Let $\Gamma$ be the absolute Galois group  $\Gal(\bar{F}/F)$.  Since $W_F$ is dense
in $\Gamma$, the action of $W_F$ on any group $\hat{H}$ extends
to an action of $\Gamma$ on $\hat{H}$.  
A set of endoscopic data $(H, \mathcal{H}, s, \xi)$
is said to be {\bf elliptic} if $\xi(Z(\hat{H})^\Gamma)^0$ is contained in $Z(\hat{G})$
(Here, upper $0$ denotes identity component).
We may take $\Gamma$ to be $\Gal(K/F)$ if $\mb{H}$ splits over some finite extension $K$ of $F$.

Two sets of $\ve$-endoscopic data $(H, \mathcal{H}, s, \xi)$ and $(H', \mathcal{H}', s', \xi')$
are said to be {\bf equivalent} if there exists $g \in \hat{G}$ such that:
\indexi{endoscopic!data!equivalent}%
\begin{enumerate}
\item
$g\xi(\mathcal{H})g^{-1} = \xi'(\mathcal{H}')$,
\item
$gsg^{-1} \equiv s' \mod Z(\hat{G})$.
\end{enumerate}
With regard to the application of Kottwitz-Shelstad's stabilization of the trace formula,
we will compute elliptic endoscopic data only up to equivalence
(hence the earlier remark that we may take $\mc{H}$ to be ${}^L H$).

Let $\ve$ be a nontrivial quadratic character of $F^\times$ if $F$ is local, or of $\idc{F}$ if
$F$ is global.  Let $W_F$ denote the Weil group of $F$ (see \cite{T}).
As explained in \cite{KS},
corresponding via Langlands correspondence
\indexi{Langlands correspondence}%
to $\ve$ is a cocycle $\mb{a}$ in $\displaystyle H^1(W_F, Z(\hat{G}))$ if
\indexs{a@$\mb{a}$, $a$}%
$F$ is local, or in $\displaystyle H^1(W_F, Z(\hat{G}))/\text{ker}^1(W_F, Z(\hat{G}))$ 
if $F$ is global.

Let $E$ be the quadratic extension of $F$ which corresponds to $\ve$ via local or global
class field theory (see \cite{LANG}).  
Then, $\mb{a}|_{W_E}$ is trivial, and
we identify $\mb{a}$ with a representative in $H^1(\Gal(E/F), Z(\hat{G}))$.

Since $\mb{G}$ is $F$-split, 
the Galois action of $\Gal(E/F)$ on $Z(\hat{G}) = \CC^\times$ is 
\indexi{Galois action}%
trivial.  
Let $\sigma$ be the generator of $\Gal(E/F)$.
\indexs{sigma@$\sigma$}%
The cohomology class of $\mb{a}$ consists of a single character 
$a \in \Hom(\Gal(E/F), \CC^\times)$, defined by $a(\sigma) = -1$.

Suppose $(\mb{H}, \mathcal{H}, s, \xi)$ is a set of $\ve$-endoscopic data.
Then, for any $h_w \in \mathcal{H}$ with image $w$ in $W_F$, we have:
\[
\xi(h_w)^{-1} s \xi(h_w) =
\begin{cases}
-s & \text{ if } w \equiv\sigma \mod W_E,\\
s & \text{ otherwise}.
\end{cases}
\]

Since we shall list the sets of endoscopic data only up to equivalence,
we may assume that the semisimple element $s$ lies in the diagonal torus $\hat{T}$ of $\hat{G}$.
Since ${}^L G = \hat{G}$,
if $\xi(h_w)^{-1} s \xi(h_w) = -s$ for some $h_w \in \mathcal{H}$,
then $\text{Int}(\xi(h_w))$ must correspond to an element of the Weyl group
$W = W(\hat{T}, \hat{G})$ of $\hat{T}$ in $\hat{G}$.
  
Given any $t = \diag(a, b, \lambda/b,\lambda/a)$ in $\hat{T}$,
let the numbers $1, 2, 3, 4$ represent the entries $a, b, \lambda/b, \lambda/a$, respectively.
Then, the actions of $W$  on $t$ are represented by the following set of permutations on
$\{1, 2, 3, 4\}$:
\[
\{1, (12)(34), (23), (3421), (2431),  (42)(31),(23)(41),(14)\}.
\]
Hence, $W$ is isomorphic to $D_4$, the dihedral group of order $8$.

By going through all possible
images of $s$ under $W$, we see that, up to equivalence of endoscopic data,
$s$ must be equal to one of:
\indexs{s@$s$!$s_1$, $s_2$}
\begin{enumerate}
\item
$s_1 = \diag(1, 1, -1, -1)$,
\item
$s_2 = \diag(1, -1, -1, 1)$,
\item
$s_3 = \diag(1, -1, -d, d)$ for some  $d\neq \pm 1$.
\end{enumerate}
\begin{claim}
No elliptic $\ve$-endoscopic data is of the form $(\mb{H}, \mathcal{H}, s_3, \xi)$.
\end{claim}
\begin{proof}
Suppose $(\mb{H}, \mc{H}, s_3, \xi)$ constitutes a set of endoscopic data.
Then, $\hat{H}$ is isomorphic to
\[
\text{Cent}(s, \hat{G})^0
= \hat{T} = \{\diag(a, b, \lambda/b, \lambda/a) : a, b, \lambda \in \CC^\times\}.
\]
By condition \ref{condition:endo4a} in Definition \ref{define:endoscopy}, 
if $w \in W_F$ is equivalent to $\sigma$ modulo $W_E$, then
the action of $w$ on $\hat{H}$ must correspond to the permutation $(12)(34)$.
If $w$ is not equivalent to $\sigma$ modulo $W_E$, then the action of $w$ on $\hat{H}$
is trivial.
Consequently, 
\[
Z(\hat{H})^\Gamma = \hat{H}^\Gamma = 
\{\diag(a, a, \lambda/a, \lambda/a): a,\lambda \in \CC^\times\},
\]
which is connected and not contained in $Z(\hat{G})$, which means
$(H, \mathcal{H}, s_3, \xi)$ is not elliptic.
\end{proof}

In the following sections, we consider the  $\ve$-endoscopic data 
$(\mb{H}, \mc{H}, s, \xi)$ where $s$ is equal to $s_1$ or $s_2$.  
It turns out that there is only one quasi-split $\mb{H}$ 
such that $(\mb{H}, \mc{H}, s, \xi)$ is elliptic for each choice of $s$.
\section{Endoscopic group $\mb{H_1}$}
\indexs{H@$\mb{H_1}$|(}\indexs{H@$\hat{H}$!$\hat{H}_1$|(}\indexs{H@$\mathcal{H}$!$\mathcal{H}_1$|(}%
\indexs{xi@$\xi$!$\xi_1$|(}%
\subsection{The data $(\mb{H_1}, \mathcal{H}_1, s_1, \xi_1)$}
We now construct a set of elliptic $\ve$-endoscopic data $(\mb{H}_1, \mathcal{H}_1, s_1, \xi_1)$, where
\[s_1 = \dmfour{1}{1}{-1}{-1}.\]
Let
\[
\hat{H}_1 = 
\left( \GL(2, \mathbb{C}) \times \mathbb{C}^\times \times \mathbb{C}^\times \right)/
\left\{\lp\diag(z, z)^{-1}, z,z\rp : z \in \mathbb{C}^\times\right\}.
\]
Fix a splitting $\spl_{\hat{H_1}} = \{\hat{B_1},\hat{T_1}, \{X_1\}\}$ of $\hat{H_1}$, where
\indexi{splitting}%
$\displaystyle
\hat{B_1} = \left\{\lp\lp\lsm *&*\\&*\rsm\rp, *, *\rp\right\}
$
is a Borel subgroup of $\hat{H_1}$,
$\hat{T_1}$ is the diagonal torus of $\hat{H_1}$, and $\{X_1\}$ is a singleton consisting of
the image  of the vector 
$X_1 = \lp\lp\lsm 0& 1 \\ 0& 0\rsm\rp, 0,0\rp \in \mathfrak{gl}(2, \CC) \times \CC \times \CC$
in the Lie algebra of $\hat{H_1}$.
Here, lower case, gothic type denotes the Lie algebra of a Lie group.

Let $\sigma$ be the generator of $\Gal(E/F)$.
Define an action of $W_F$ on $\hat{H_1}$ by
\[
{w}(A, x, y) = 
\begin{cases}
(A, y, x)&\text{ if } {w} \equiv \sigma \mod W_E,\\
(A, x, y)&\text{ otherwise.}
\end{cases}
\]
This defines a semi-direct product $\mathcal{H}_1 = \hat{H_1}\rtimes W_F$.
\indexs{H@$\mc{H}_1$}%
\indexs{L@${}^L H$!${}^L H_1$}%
Since the action of $W_F$ fixes $\spl_{\hat{H_1}}$,
$\mc{H}_1$ is an $L$-group with $\hat{H_1}$ as its identity component.

Let $e = \lp\lsm 1 & \\ & -1\rsm\rp$.
\indexs{e@$e$}%
For any $A \in \GL(2, \CC)$, $x, y \in \CC^\times$, let
$d(A, x, y)$ be the element $\lp\lsm x A & \\ & y  eAe \rsm\rp$ in $\hat{G}.
\indexs{d@$d(\;,\;,\;)$}%
$
Then, $d$
induces an isomorphism
\[
d :\hat{H_1} \tilde{\rightarrow}\;\text{Cent}(s_1, \hat{G})^0.
\]
Define an $L$-embedding
\indexi{L@$L$-!embedding}%
$\xi_1:\mathcal{H}_1 \rightarrow \hat{G}$ by $\xi_1|_{\hat{H_1}} = d$ and
\[
\xi_1(1 \rtimes w)
=
\begin{cases}
\lp\lsm &e\\e&\rsm\rp &\text{ if } w \equiv \sigma \mod W_E,\\
1&\text{ otherwise.}
\end{cases} 
\]
Let $\mb{H_1}$ be the unique quasi-split reductive group over $F$ whose $L$-group
${}^L H_1$ is equal to $\mc{H}_1$.
One can check to see that $(\mb{H}_1, \mathcal{H}_1, s_1, \xi_1)$ is a set of elliptic
$\ve$-endoscopic data.
\begin{claim}
Up to equivalence,
$(\mb{H}_1, \mathcal{H}_1, s_1, \xi_1)$ is the unique elliptic \dts $\ve$-endoscopic data attached to $s_1$.
\end{claim}
\begin{proof}
If $(\mb{H}, \mathcal{H}, s_1, \xi)$ is another set of elliptic $\ve$-endoscopic data attached to 
$s_1$, then $\xi(\mathcal{H}) = \xi_1(\mathcal{H}_1)$.
It follows from the definitions
that $(\mb{H}, \mc{H}, s_1, \xi)$ is equivalent to \dt $(\mb{H_1}, \mc{H}_1, s_1, \xi_1)$.
\end{proof}

\subsection{Explicit Description of $\mb{H_1}$}
We now describe explicitly the quasi-split reductive $F$-group $\mb{H}_1$.
Let $\sigma$ be the generator of $\Gal(E/F)$.
Recall that $\mathcal{H}_1 = \hat{H}_1 \rtimes W_F$, 
and the action of $W_F$ factors through
$W_F/W_E = \Gal(E/F)$.  Consequently, $\mb{H}_1$ must split over $E$.

To construct $\mb{H_1}$, we make use of the fact that there exists a maximal torus $\mb{T_1}$ of
\indexi{maximal torus}%
$\mb{H_1}$ dual to the maximal diagonal torus $\hat{T_1}$ of $\hat{H_1}$; namely, 
\indexi{dual group}%
$X_*(\mb{T_1}) = X^*(\hat{T_1})$, where $X_*(\mb{T_1}) := \Hom(\mbb{G}_m, \mb{T_1})$ and
$X^*(\hat{T_1}) := \Hom(\hat{T_1}, \mbb{C}^\times)$.
\indexs{X@$X_*(\;)$}\indexs{X@$X^*(\;)$}%
We have:
\[
X^*(\hat{T_1})  = \{(x, y;z, t) \in \mathbb{Z}^2 \times \mathbb{Z}^2 : x+y = z + t\},
\]
where $(x, y; z, t)(\diag(a, b), c, d)
:= a^x b^y c^z d^t$ for all $(\diag(a, b), c, d) \in \hat{T_1}$.

The Galois action of $\sigma$ on $\hat{T_1}$ is given by
$\sigma(\diag(a, b), c, d) = (\diag(a, b), d, c)$.  This induces an action on
$X^*(\hat{T_1})$ and consequently on $X_*(\mb{T_1})$, namely: 
\[\sigma(x, y;z, t) = (x, y; t, z), \quad \forall (x, y;z, t) \in X_*(\mb{T_1}).\]  
We conclude that
\[
\mb{T_1} = (\Gm \times \Gm \times \rR_{E/F}\Gm)',
\]
where $\rR_{E/F}\Gm$ is the $F$-group obtained from $\Gm$ on restricting scalars from $E$ to $F$,
\indexi{restriction of scalars}%
and the prime indicates that the product of the first two factors
is equal to the norm of the third factor in $\Gm$. 

The algebraic group $\mb{H_1}$ is defined by its root datum, which is the dual of
\indexi{root datum}%
the root datum of $\hat{H_1}$ (see \cite{S}).
Having found the dual of $\hat{T_1}$, it remains to 
\indexi{dual group}%
find the dual of the root lattice and make a choice of splitting.  We skip this routine  procedure.  
The resulting group is
\[
\mb{H_1} = (\GL(2)\times \rR_{E/F}\Gm)',
\]
where the prime indicates that the determinant
of the first factor is equal to the norm of the second factor in $\Gm$.

Since $\mb{H_1}$ is split over $E$ and the group of $E$-points of $\rR_{E/F}\Gm$ is
$E^\times \times E^\times$, the group of $E$-points of $\mb{H_1}$ is $\mb{H_1}(E) =$
\[\lp\GL(2, E)\times E^\times \times E^\times\rp' :=
\left\{
(g, \gamma, \delta) \in \GL(2, E) \times E^\times \times E^\times : \det g = \gamma\delta
\right\}.\]
For any $\gamma \in E$, put $\bar{\gamma} := \sigma \gamma$.
\indexs{bar@$\bar{\;}$}%
For $g = (g_{ij}) \in \GL(2, E)$, put $\bar{g} := \lp \ol{g_{ij}} \rp$.
\indexs{bar@$\bar{\;}$}%
The action of $\sigma$ on $\mb{H}_1(E)$ is induced by
the action of $\sigma$ on $\hat{H_1}$, namely:
\[
\sigma^* : (g, \gamma,\delta) \mapsto \lp\bar{g}, \bar{\delta},\bar{\gamma}\rp,\;\forall 
(g, \gamma, \delta) \in \lp\GL(2, E) \times E^\times \times E^\times\rp'.
\]
\indexs{sigma@$\sigma^*$}%
The $F$-points of $\mb{H}_1$ are the elements in $\mb{H_1}(E)$ which are fixed by $\sigma^*$; thus,
\[
\mb{H_1}(F) = (\GL(2, F) \times E^\times)'\\
:= 
\boxed{
\{(g, \gamma) \in \GL(2, F) \times E^\times : \det g = \N_{E/F} \gamma\}
}\;.
\]
Suppose $F$ is a number field.  The norm mapping $\N_{E/F}: E \rightarrow F$ induces
\indexi{norm mapping}%
a norm mapping $\N_{E/F}:\Ae \rightarrow \Af$ on the \adeles, and
the group of $\Af$-points of $\mb{H_1}$ is
\[
\mb{H_1}(\Af) = (\GL(2, \Af) \times \Ae^\times)'\\
:= 
\boxed{
\{(g, \gamma) \in \GL(2, \Af) \times \Ae^\times : \det g = \N_{E/F} \gamma\}
}\;.
\]
At each place $v$ of $F$, the Galois action on $\mb{H_1}$ as an
$F_v$-group is compatible with the Galois action on $\mb{H_1}$ as an $F$-group.
At a place $v$ which does not split in $E$, $E_v := E\otimes_F F_v$ is a field and
$\Gal(E_v/F_v)$ embeds isomorphically into $\Gal(E/F)$.
Consequently, the $F_v$-group $\mb{H_1}$ splits over $E_v$ and its group of $F_v$-points is as follows:
\[
\mb{H_1}(F_v) =
\lp\GL(2, F_v) \times E_v^\times\rp'\\
:= \left\{(g, \gamma) \in \GL(2, F_v)\times E_v^\times: \det g = \N_{E/F}\gamma\right\}.
\]
If $v$ splits into two places $v_1, v_2$ in $E$, then $E_{v_1} = E_{v_2} = F_v$
and the Galois action on $\mb{H_1}$ as an $F_v$-group is trivial, i.e.
$\mb{H_1}$ is split over $F_v$.  We have $\mb{H_1}(F_v) =$
\[
\lp \GL(2, F_v)\times F_v^\times \times F_v^\times\rp'\\
:= \left\{ (g, a, b) \in \GL(2, F_v) \times F_v^\times \times F_v^\times
: \det g = ab\right\}.
\]
\indexs{H@$\mb{H_1}$|)}\indexs{H@$\hat{H}$!$\hat{H}_1$|)}\indexs{H@$\mathcal{H}$!$\mathcal{H}_1$|)}%
\indexs{xi@$\xi$!$\xi_1$|)}%

\section{Endoscopic group $\mb{H_2}$}
\indexs{H@$\mb{H_2}$|(}\indexs{H@$\hat{H}$!$\hat{H}_2$|(}\indexs{H@$\mathcal{H}$!$\mathcal{H}_2$|(}%
\indexs{xi@$\xi$!$\xi_2$|(}%
\subsection{The Data $(\mb{H_2}, \mathcal{H}_2, s_2, \xi_2)$}
We now construct  a set of elliptic $\ve$-endoscopic data $(\mb{H_2}, \mc{H}_2, s_2, \xi_2)$,
where \[s_2 = \dmfour{1}{-1}{-1}{1}.\]
Let
\[
\hat{H_2} = 
\left\{ (A, B, \lambda)\in \GL(2, \mathbb{C})^2 \times\mathbb{C}^\times  
: \det A = \det B = \lambda\right\}.
\]
Fix a splitting $\spl_{\hat{H_2}} = \{\hat{B_2},\hat{T_2},\{X_a, X_b\}\}$ of
\indexi{splitting}%
$\hat{H_2}$, where
\[
\hat{B_2} = \left\{\lp\lp\lsm *&*\\&*\rsm\rp, \lp\lsm*&*\\&*\rsm\rp, *\rp\right\}
\subset \hat{H_2},
\]
$\hat{T_2}$ is the maximal diagonal torus of $\hat{H_2}$, and
\indexi{maximal torus}%
\[
X_a = \lp\lp\lsm 0&1\\0&0\rsm\rp,\lp\lsm 0&0\\0&0\rsm\rp , 0\rp,\quad
X_b = \lp \lp\lsm 0&0\\0&0\rsm\rp, \lp\lsm 0&1\\0&0\rsm\rp, 0\rp 
\]
are elements in the Lie algebra of $\hat{H_2}$.

Define an action of $W_F$ on $\hat{H_2}$ by
\[
w(A, B, \lambda) = 
\begin{cases}
(B, A, \lambda)&\text{ if } w \equiv \sigma \mod W_E,\\
(A, B, \lambda)&\text{ otherwise}.
\end{cases}
\]
This defines a semi-direct product $\mathcal{H}_2 = \hat{H_2}\rtimes W_F$.
\indexs{H@$\mc{H}_2$}%
\indexs{L@${}^L H$!${}^L H_2$}%
Since the action of $W_F$ fixes $\spl_{\hat{H_2}}$,
$\mc{H}_2$ is an $L$-group, with $\hat{H_2}$ as its identity component.

We define an $L$-group embedding $\xi_2 : \mc{H}_2 \rightarrow \hat{G}$ as follows:
\indexi{L@$L$-!embedding}%
Put 
\[
\left[\lp\lsm a&b\\c&d\rsm\rp, \lp\lsm e&f\\g&h\rsm\rp\right]
:= \lp\lsm a &&&b\\ &\lsm e&f\\g&h\rsm\\c&&&d\rsm\rp.
\indexs{b@$\left[\lp\;\rp,\lp\;\rp\right]$}%
\]
Define the restriction of $\xi_2$ to $\hat{H_2}$ by
\[
\xi_2|_{\hat{H_2}} (A, B, \lambda) = [A, B],\quad \forall (A, B, \lambda) \in \hat{H_2}.
\]
Then, $\hat{H_2}$ is isomorphic to $\text{Cent}(s_2, \hat{G})^0$ via $\xi_2|_{\hat{H_2}}$.
Put $w_0 := \lp\lsm & 1 \\1&\rsm\rp \in \GL(2, \mbb{C})$.
Let
\[
\xi_2(1 \rtimes w) =
\begin{cases}
\lp\lsm w_0 &\\&w_0\rsm\rp &\text{ if } w \equiv \sigma \mod W_E,\\
1&\text{ otherwise.}
\end{cases} 
\]
Let $\mb{H_2}$ be the unique quasi-split reductive group over $F$ whose $L$-group ${}^L H_2$ is $\mc{H}_2$.
One can check to see that $(\mb{H}_2, \mathcal{H}_2, s_2, \xi_2)$ is a set
of elliptic $\ve$-endoscopic data.
\begin{claim}
Up to equivalence, $(\mb{H}_2, \mathcal{H}_2, s_2, \xi_2)$ is the unique elliptic \dts 
$\ve$-endoscopic data attached to $s_2$.
\end{claim}
\begin{proof}
As in the case of $s_1$, 
if $(\mb{H}, \mathcal{H}, s_2, \xi)$ is another set of elliptic $\ve$-endoscopic data attached to 
$s_2$, then $\xi(\mathcal{H}) = \xi_2(\mathcal{H}_2)$.  It follows from the definitions
that $(\mb{H}, \mc{H}, s_2, \xi)$ is equivalent to $(\mb{H_2}, \mc{H}_2, s_2, \xi_2)$.
\end{proof}

\subsection{Explicit Description of $\mb{H_2}$}
We now give an explicit description of the quasi-split, reductive $F$-group $\mb{H_2}$.
The action of $W_F$ on $\hat{H_2}$ factors through $W_F / W_E = \Gal(E/F)$.
Consequently, the group $\mb{H_2}$ is split over the quadratic extension $E$. 

To find $\mb{H_2}$, let $\mb{T_2}$ be the maximal torus of $\mb{H_2}$ which is dual to
\indexi{maximal torus}\indexi{dual group}%
the maximal diagonal torus $\hat{T_2}$ of $\hat{H_2}$.
Thus, $X_*(\mb{T_2})$ is equal to
\[
X^*(\hat{T_2}) = 
\{
(x, y; z, t;w) \in \mathbb{Z}^2 \times \mathbb{Z}^2 \times\mathbb{Z}
\}/\{(n, n; m, m;-n -m) : m,n\in \mathbb{Z}\},
\]
where
\begin{multline*}
(x, y; z, t; w)(\diag(a, b), \diag(c, d), \lambda)
= a^x b^y c^z d^t \lambda^w, \\ \forall (\diag(a, b), \diag(c, d), \lambda) \in \hat{T_2}.
\end{multline*}
The action of $\sigma$ on $\hat{T_2}$ induces an action on $X^*(\hat{T_2})$; namely,
\[
\sigma(x, y; z, t;w) = (z, t; x, y;w),\quad \forall (z, t; x, y;w) \in X^*(\hat{T_2}).
\]  
We conclude that $\mb{T_2}$ is a torus whose group of $E$-points is
\[
\Big\{
\lp\dmtwo{\zeta_1}{\zeta_2}, \dmtwo{\gamma_1}{\gamma_2}, \eta\rp : \zeta_1, \zeta_2, \gamma_1, \gamma_2, \eta\in E^\times
\Big\}
/\{(\alpha I_2, \beta I_2, (\alpha\beta)^{-1}) : \alpha\beta \in E^\times\},
\]
where $I_2$ is the $2 \times 2$ identity matrix.
\indexs{I@$I_2$}%
Put $\bar{\gamma} := \sigma \gamma$ for $\gamma \in E$.
The action of $\sigma$ on $\hat{T_2}$ gives rise to an action of $\sigma$ on $\mb{T_2}(E)$ define by:
\[
\sigma^*(\diag({\zeta_1}, {\zeta_2}), \diag({\gamma_1},{\gamma_2}), \eta)
= (\diag(\bar{\gamma_1},\bar{\gamma_2}),\diag(\bar{\zeta_1}, \bar{\zeta_2}), \eta).
\]
\indexs{sigma@$\sigma^*$}%
Having found the the dual of $\hat{T_2}$, 
\indexi{dual group}%
\indexi{root datum}%
we conclude without going through the details here that:
\[
\mb{H_2}(E) =
\lp\GL(2, E)\times\GL(2, E)\times E^\times\rp 
/ \left\{\lp\alpha I_2, \beta I_2, (\alpha\beta)^{-1}\rp : \alpha, \beta \in E^\times\right\}.
\]

For $g =(g_{ij}) \in \GL(2, E)$, let $\bar{g} = (\ol{g_{ij}})$.
The Galois action on $\hat{H_2}$ induces an action on $\mb{H_2}(E)$ defined by:
\[
\sigma^* : (g_1, g_2, \zeta) \mapsto (\bar{g_2}, \bar{g_1}, \bar{\zeta}),\quad\forall
(g_1, g_2, \zeta) \in \mb{H_2}(E).
\]
\begin{claim}
The group of $F$-points of $\mb{H_2}$ is
\begin{multline*}
\mb{H_2}(F) = \{(g, \bar{g}, c) : g \in \GL(2, E), c \in F^\times\}
/\{(z I_2,\bar{z} I_2,\N_{E/F}z^{-1}): z\in E^\times\}\\
\cong
\boxed{
\lp \GL(2, E)\times F^\times\rp
/\{(z I_2,\N_{E/F}z^{-1}): z\in E^\times\}
}\;.
\end{multline*}
\end{claim}
\begin{proof}
The group of $F$-points of $\mb{H_2}$ consists of the set of equivalence classes represented by
$(g, h, \zeta) \in \GL(2, E)\times\GL(2, E)\times F^\times$ such that
$(g, h, \zeta) \equiv (\bar{h}, \bar{g}, \bar{\zeta})$ modulo
$ \left\{\lp\alpha I_2, \beta I_2, (\alpha\beta)^{-1}\rp : \alpha, \beta \in E^\times\right\}$.

Suppose $(g, h, \zeta) \equiv (\bar{h},\bar{g},\bar{\zeta})$, then there exist
$\alpha, \beta \in E^\times$ such that
\[
\begin{split}
\alpha g &= \bar{h},\\
\beta h &= \bar{g},\\
(\alpha\beta)^{-1}\zeta &= \bar{\zeta}.
\end{split}
\]
This implies that $\alpha = \bar{\beta}^{-1}$ and $\beta^{-1}\zeta \in F^\times$.
Since 
\[
(g, h, \zeta) \equiv 
(g, \bar{g}, b^{-1}\zeta) = (g, h, \zeta) (1, \beta, \beta^{-1}) \mod 
\left\{\lp\alpha I_2, \beta I_2, (\alpha\beta)^{-1}\rp\right\},
\]
the claim follows.
\end{proof}

Suppose $F$ is a number field.  The norm mapping $\N_{E/F} : E \rightarrow F$ induces a mapping
\indexi{norm mapping}%
$\N_{E/F} : \Ae \rightarrow \Af$.  The group of $\Af$-points of $\mb{H_2}$ is:
\[
\mb{H_2}(\Af) = \lp\GL(2, \Ae)\times \Af^\times\rp/\{(zI_2,\N_{E/F}z^{-1}): z\in \Ae^\times\}.
\]
At a place $v$ of $F$ which does not split in $E$, 
the Galois group $\Gal(E_v/F_v)$ is isomorphic to $\Gal(E/F)$.
Consequently, we have:
\[
\mb{H_2}(F_v) = (\GL(2, E_v) \times F_v^\times)/\{(z I_2, \N_{E/F}z^{-1}): z \in E_v^\times\}.
\]
At a place $v$ which splits in $E$, $\mb{H_2}$ is split over $F_v$, and $\mb{H_2}(F_v)$ is equal to
\[
\left\{
(g_1, g_2, c) \in \GL(2, F_v)\times\GL(2, F_v)\times F_v^\times
\right\}/
\left\{
\lp aI_2, bI_2, (ab^{-1})\rp : a, b\in F_v^\times
\right\}.
\]
\subsubsection{Summary}
In summary, we work with the the following quasi-split reductive $F$-groups:
\begin{itemize}
\item
$
\mb{G} = \GSp(2) := \left\{g \in \GL(4): {}^t g J g = \lambda(g) J
\text{ for some } \lambda(g) \in \mathbb{G}_m\right\},
$
where $J = \lp\lsm &&&1\\&&1&\\&-1&&\\-1&&&\rsm\rp$.
\item
$\mb{H_1} = (\GL(2) \times \rR_{E/F}\Gm)'$, whose group of $F$-points is:
\[
\mb{H}_1(F) =\{(g, \gamma) \in \GL(2, F) \times E^\times : \det g = \N_{E/F} \gamma\}.
\]
\item
The reductive group $\mb{H_2}$, whose group of $F$-points is:
\[
\mb{H_2}(F) = \lp \GL(2, E) \times F^\times \rp / \{(zI_2, \N_{E/F}z^{-1}): z \in E^\times\}.
\]
\end{itemize}
\indexi{endoscopic!data|)}\indexi{endoscopic!group|)}%
\indexi{splitting|)}\indexi{maximal torus|)}\indexi{root datum|)}%
\indexs{H@$\mb{H_2}$|)}\indexs{H@$\hat{H}$!$\hat{H}_2$|)}\indexs{H@$\mathcal{H}$!$\mathcal{H}_2$|)}%
\indexs{xi@$\xi$!$\xi_2$|)}%

\section{Norm Correspondence}\label{sec:normmapping}
\indexi{norm correspondence|(}%
Let $F$ be a local field or a number field.

{\bf Terminology: }
\begin{itemize}
\item
For any algebraic $F$-group $\mb{H}$, we say that two elements $h$, $h'$ in $\mb{H}(F)$ 
are {\bf conjugate} (or $F$-conjugate) if there exists an element $g \in \mb{H}(F)$ such that $g^{-1}hg =h'$.
A conjugacy class in $\mb{H}(F)$ consists of all elements in $\mb{H}(F)$ which are conjugate to one
another.  
\indexi{conjugacy}\indexi{conjugacy!F@$F$-conjugacy}%
\item
We say that $h, h' \in \mb{H}(F)$ are {\bf stably conjugate} (or $\bar{F}$-conjugate) if there exists
an element $g \in \mb{H}(\bar{F})$ such that $g^{-1}hg = h'$.  A {\bf stable conjugacy class}
in $\mb{H}(F)$ consists of all elements in $\mb{H}(F)$ which are stably conjugate to a given one.
We say that a conjugacy class is stable if it coincides with the stable conjugacy class which contains it.
\indexi{conjugacy!stable}%
\item
For our purpose, an element in $\mb{H}(F)$ is {\bf semisimple} if it is diagonalizable over $\bar{F}$.
A semisimple element is {\bf regular} if its centralizer in $\mb{H}(F)$ is a maximal torus.
\indexi{element!semisimple}%
\indexi{element!regular}%
\item
A torus in $\mb{H}(F)$ is {\bf elliptic} if it is not contained in any proper parabolic subgroup
of $\mb{H}(F)$.
\indexi{maximal torus!elliptic}%
An element in $\mb{H}(F)$ is {\bf elliptic} if it is contained in an elliptic torus.
\indexi{element!elliptic}%
\end{itemize}

The goal of this work is to relate admissible representations of $\mb{G}(F), \mb{H_1}(F)$, and
$\mb{H_2}(F)$ if $F$ is local, or automorphic representations of
$\mb{G}(\Af), \mb{H_1}(\Af)$, and $\mb{H_2}(\Af)$ if $F$ is a number field.  To do so, we need
to introduce the notion of transfer of stable conjugacy classes 
from $\mb{H_1}, \mb{H_2}$ to $\mb{G}$.


Let $\mb{H}$ be an $\ve$-endoscopic group of $\mb{G}$ with 
$L$-group embedding $\xi : {}^L H \rightarrow {}^L G$. 
\indexi{L@$L$-!embedding}%
Let $\CL{\mb{H}}, \CL{\mb{G}}$ denote the set of stable conjugacy classes of semisimple elements in
\indexs{C@$\CL{\;}$}%
$\mb{H}({F}), \mb{G}({F})$, respectively.

Let $\mb{T} = \left\{\diag(a, b, l/b, l /a):a,b,l\in \mbb{G}_m\right\}$ be the
maximal diagonal torus of $\mb{G}$.
Let $\mb{T_H}$ be a maximal torus of $\mb{H}$ defined over $F$.
Let $\lambda$ be the unique map from $X^*(\mb{T_H})$ to $X_*(\hat{T})$ 
such that the following diagram commutes:
\[
\xymatrix{
X^*(\mb{T_H}) \ar[r]^\lambda \ar[d]& X^*(\mb{T})\ar[d]\\
X_*(\hat{T}_H) \ar[r]^\xi & X_*(\hat{T})
}
\]
\indexs{lambda@$\lambda$}\indexs{xi@$\xi$}%
Here, the vertical isomorphisms are given by the correspondence between root data and dual root data
(see \cite{B}).


Any semisimple $\delta_H \in \mb{H}(F)$ is $\bar{F}$-conjugate to some $t_H \in \mb{T_H}(\bar{F})$.
Likewise, any semisimple $\delta \in \mb{G}(F)$ is $\bar{F}$-conjugate to some
$t \in \mb{T}(\bar{F})$.
Let $\widetilde{\delta_H} \in \CL{\mb{H}}$ 
denote the stable conjugacy class of $\delta_H \in \mb{H}(F)$.  Let $\tilde{\delta}$
\indexs{tilde@$\tilde{\;}$}%
denote the stable conjugacy class of $\delta \in \mb{G}(F)$.

Define a map $\mathcal{A}_{\mb{H}\bs\mb{G}} : \CL{\mb{H}} \rightarrow \CL{\mb{G}}$ as follows:
\indexs{A@$\mathcal{A}_{\;\bs\;}$}%
$\mathcal{A}_{\mb{H}\bs\mb{G}}(\widetilde{\delta_H}) = \tilde{\delta}$ if
there exists $t_H \in \mb{T_H}(\bar{F})$, $t \in \mb{T}(\bar{F})$ such that
\begin{itemize}
\item
$t_H$ is conjugate to $\delta_H$ in $\mb{H}(\bar{F})$,
$t$ is conjugate to $\delta$ in $\mb{G}(\bar{F})$;
\item
$\chi(t_H) = (\lambda(\chi))(t)$ for all $\chi \in X^*(\mb{T_H})$.
\end{itemize}

\begin{define}
{\rm We say that $\delta_H \in \mb{H}(F)$ is a {\bf norm} of $\delta \in \mb{G}(F)$ if
\indexi{norm}%
$\tilde{\delta} = \Acon{\mb{H}}(\widetilde{\delta_H})$.
If $F$ is global, we say that $\delta_H = (\delta_{H, v}) \in \mb{H}(\Af)$ 
is a norm of $\delta =(\delta_v)\in \mb{G}(\Af)$ if $\delta_{H,v} \in \mb{H}(F_v)$ 
is a norm of $\delta_v \in \mb{G}(F_v)$ for every place $v$ of $F$.
}
\end{define}

In general,
$\Acon{\mb{H}}$ may not be one-to-one.  Moreover, a regular element in $\mb{H}$ may be
a norm of a non-regular element in $\mb{G}$.  If $\delta_H \in \mb{H}$ is a norm of a regular element
in $\mb{G}$, we say that $\delta_H$ is \emph{$\mb{G}$-regular}.
\indexi{element!G@$\mb{G}$-regular}%
(Remark: A $\mb{G}$-regular element in $\mb{H}$ is necessarily regular).

Note that $\Acon{\mb{H}}$ is equivalent to a map
\[
\N : \mb{T_H}(\bar{F})/W_H \rightarrow \mb{T}(\bar{F})/W,
\]
\indexs{N@$\N$}%
where $W, W_H$ are the Weyl groups associated with $\mb{T}(\bar{F}), \mb{T_H}(\bar{F})$, respectively.
\indexi{Weyl group}%
More precisely, given any  class $\widetilde{t_{H}} \in  \mb{T_H}(\bar{F})/W_H$ 
represented by $t_{{H}} \in \mb{T_H}(\bar{F})$, 
the class $\N(\widetilde{t_{H}})$ is represented by some $t \in \mb{T}(\bar{F})$ such that  
$(\lambda(\chi))(t) = \chi(t_H)$ for all $\chi \in X^*(\mb{T_H})$.
We shall describe $\N$ explicitly in the next two sections.  

In this work, we identify $\hat{G}$ with $\GSp(2, \mbb{C})$ (see \cite{B}). 
Write:
\[
X^*(\mb{T}) = \{(x, y, z, t) \in \mathbb{Z}^4\}/ \{(m, -m, -m, m) : m\in \mathbb{Z}\},
\]
where $(x, y, z, t)\delta := a^xb^y(l/b)^z(l/a)^t$ for all
$\delta = \diag(a, b, l/b, l/a) \in \mb{T}$.
Write:
\[
X_*(\hat{T}) = \{(x, y, z, t) \in \mathbb{Z}^4 : x + t = y + z \},
\] 
where $(x, y, z, t)\zeta := \diag(\zeta^x, \zeta^y, \zeta^z, \zeta^t) \in \mb{T}$ for all
$\zeta \in \Gm$.

The isomorphism $\iota : X^*(\mb{T}) \rightarrow X_*(\hat{T})$, determined from the identification of
\indexs{iota@$\iota$}%
$\hat{G}$ with $\GSp(2, \mbb{C})$, is defined as follows:
\[
\iota : (x, y, z, t) \mapsto (x + y, x+z, t+y, t+z),
\]
\[
\iota^{-1} : (x, y, z, t) \mapsto (x - z, z, t, 0).
\]

For $i = 1, 2$, let $\mb{T}_i$ denote the maximal diagonal torus of $\mb{H}_i$.
Let $W_i$ be the Weyl group of $\mb{T}_i(\bar{F})$ in $\mb{H}_i(\bar{F})$.
Let $\N_i$ denote the map $\N$ from $\mb{T}_i(\bar{F})/W_i$ to $\mb{T}(\bar{F})/W$
defined above.
\subsection{Norm Correspondence for $\mb{H_1}$}
Recall that 
\[\hat{H_1} = \lp\GL(2, \mathbb{C}) \times \mathbb{C}^\times \times\CC^\times\rp/
\{(\diag(z, z), z^{-1}, z^{-1}) : z \in \CC^\times\}.\]
We have chosen the embedding $\xi_1 : \hat{H_1} \rightarrow \hat{G}$ to be
$\xi_1 : (A, u, v) \mapsto \lp\lsm u A &\\&v eAe\rsm\rp$, where
$e := \lp\lsm 1 & \\&-1\rsm\rp$.
Write:
\[
X^*(\mb{T_1}) = \{(x_1, y_1; z_1, t_1) \in \mathbb{Z}^4\}/\{(z, z, -z, -z) : z \in\mathbb{Z}\},
\]
where 
\[
(x_1, y_1; z_1, t_1)(\diag(a,b),c, d) := a^{x_1}b^{y_1}c^{z_1}d^{t_1},\quad\forall
(\diag(a,b),c, d) \in \mb{T_1}.
\]

\begin{claim}\label{claim:normH1}
The map $\N_1 : \mb{T_1}(\bar{F})/W_1\rightarrow\mb{T}(\bar{F})/W  $ is given by
\indexs{N@$\N$!$\N_1$}%
\[
\N_1 : (\diag(a, l/a);b,l/b) \mapsto
\diag(b,a, l/a, l/b).
\]
\end{claim}
\begin{proof}
The lattice $X^*(\mb{T_1})$ is isomorphic to $X_*(\hat{T}_1)$ in the following way:
\indexi{lattice}%
For any \dt $(x_1, y_1; z_1, t_1) \in X^*(\mb{T_1})$, its image in $X_*(\hat{T}_1)$
is the homomorphism which maps any 
$\zeta \in \CC^\times$ to $\diag((\zeta^{x_1}, \zeta^{y_1}), \zeta^{z_1},\zeta^{t_1})$.
The embedding of $X_*(\hat{T_1})$ in $X_*(\hat{T})$ is given by
\[
\xi_1 : (x_1, y_1; z_1, t_1) \mapsto 
(z_1 + x_1, z_1 + y_1, t_1 + x_1, t_1 + y_1).
\]
The map $\lambda : X^*(\mb{T}_1) \rightarrow X^*(\mb{T})$ is therefore given by
\indexs{lambda@$\lambda$}%
\[
\lambda : (x_1, y_1;z_1, t_1) \mapsto 
\iota^{-1}\circ\xi_1(x_1,y_1;z_1,t_1) = (z_1 - t_1, t_1 + x_1, t_1 + y_1, 0)
\]

Suppose an element $t_1 = (\diag(a_1,l_1/a_1); b_1, l_1/b_1)$ in $\mb{T_1}(\bar{F})$
is a norm of an element $t = (a, b, l/b, l/a)$ in $\mb{T}(\bar{F})$.
By definition, this means that
$\chi(t_1) = (\lambda(\chi))(t')$ for some $t'$ conjugate to $t$ and for all $\chi \in X^*(\mb{T_1})$.  
Without loss of generality, assume $t' = t$.
We have:
\[
{a_1}^{x_1}(l_1/a_1)^{y_1}{b_1}^{z_1}(l_1/b_1)^{t_1} =
{a_1}^{x_1 - y_1}{b}^{z_1 + t_1}(l/b)^{t_1 + y_1}
= {a}^{z_1 - t_1}{b}^{x_1 - y_1}(l/b)^{t_1 + y_1}
\]
for all $(x_1, y_1, z_1, t_1) \in X^*(\mb{T_1})$; hence,
$a = b_1, b = a_1, l = l_1$, which completes the proof.
\end{proof}

\subsection{Norm Correspondence for $\mb{H_2}$}\label{mtH2}
To compute the norm correspondence between $\mb{G}$ and $\mb{H_2}$,
we carry out the same procedure employed for $\mb{H_1}$.
Let $I_2$ be the $2 \times 2$ identity matrix.
Write:
\[
\mb{H_2}(\bar{F}) =
\lp\GL(2, \bar{F})\times\GL(2, \bar{F}) \times \bar{F}^\times\rp/
\{(aI_2, bI_2, (ab)^{-1}) : a, b \in \bar{F}^\times\}.
\]
A typical element in $\mb{T_2}(\bar{F})$ is represented by
\[
(\diag(a, b), \diag(c, d), 1) \in \GL(2, \bar{F})\times\GL(2, \bar{F}) \times \bar{F}^\times.
\]
\begin{claim}\label{claim:normH2}
The map
$\N_2 : \mb{T_2}(\bar{F})/W_2 \rightarrow \mb{T}(\bar{F})/W  $
\indexs{N@$\N$!$\N_2$}%
is given by
\[
\N_2 : (\diag(a, b), \diag(c, d), 1) \mapsto \diag(ac, ad, bc, bd).
\]
\end{claim}
\begin{proof}
We have:
\[
X^*(\mb{T_2}) = \{(x_2, y_2; z_2, t_2) \in \mathbb{Z}^4 : 
x_2+y_2 = z_2 + t_2 \},
\]
where
\begin{multline*}
(x_2, y_2; z_2, t_2)(\diag(a, b), \diag(c, d), 1)
:= a^{x_2} b^{y_2} c^{z_2} d^{t_2},\\
\forall(\diag(a, b), \diag(c, d), 1) \in \mb{T_2}.
\end{multline*}
The embedding of $X_*(\hat{T_2}) = X^*(\mb{T_2})$ in $X_*(\hat{T})$ is given by
\[
\xi_2 : (x_2, y_2; z_2, t_2) \mapsto (x_2, z_2, t_2, y_2).
\]
Thus, $\lambda : X^*(\mb{T_2}) \rightarrow X^*(\mb{T})$ is given by
\indexs{lambda@$\lambda$}%
\[
\lambda : (x_2, y_2; z_2, t_2) \mapsto \iota^{-1}(x_2, z_2, t_2, y_2)
= (x_2 - t_2, t_2, y_2, 0).
\]
Suppose an element $t_2 = (\diag(a_2, b_2), \diag(c_2, d_2), 1)$ in $\mb{T_2}(\bar{F})$ is a norm of 
some $t = (a, b, l/b, l/a) \in \mb{T}(\bar{F})$.  Then by definition
we must have
$\chi(t_2) = (\lambda(\chi))(t')$ for some $t'$ conjugate to $t$ and
for all $\chi \in X^*(\mb{T_2})$.  Without loss of generality, assume that $t' = t$.
We have:
\[
a_2^{x_2} b_2^{y_2}c_2^{z_2} d_2^{t_2} = (a_2 c_2)^{x_2} (b_2 c_2)^{y_2} (d_2/c_2)^{t_2}
= a^{x_2 - t_2} b^{t_2} (l/b)^{y_2} = a^{x_2} (b/a)^{t_2}(l/b)^{y_2}
\]
for arbitrary $(x_2, y_2, z_2, t_2;w) \in X^*(\mb{T_2})$ (The first equality is valid because
of the condition that $x_2 + y_2 = z_2 + t_2$).  We conclude that:
\[
a = a_2 c_2,\; l/b = b_2 c_2,\; b/a = d_2 /c_2;
\]
hence,
$t = \diag(a_2 c_2, a_2 d_2, b_2 c_2, b_2 d_2)$.
\end{proof}
\subsection{Norm Correspondence for $F$-Points}
Let $G = \mb{G}(F)$, $H_i = \mb{H}_i(F)$ ($i = 1, 2$).
We have described the norm correspondence among the $\bar{F}$-points of the groups.
We now describe the norm correspondence among non-elliptic or central elements in $G$, $H_1, H_2$.
Knowing the norm correspondence for these types of elements is a prerequisite to stating certain 
useful character identities in Chapter \ref{chap:local}.
The norm correspondence among elliptic regular elements is addressed in Appendix \ref{chap:fundlemma}.

Recall that $E$ is the quadratic extension of $F$ which corresponds to $\ve$.
Suppose
$E = F(\sqrt{A})$ for some element $A \in F^\times - {F^\times}^2$.
For any element $\gamma = a + b\sqrt{A}$ in $E^\times$, where $a, b \in F$, put
\[
\phi^A(\gamma) := \lp\lmx a & b A\\b&a\rmx\rp \in \GL(2, F).
\]
\indexs{phi@$\phi^A(\;)$}%
The (not necessarily distinct) eigenvalues of $\phi^A(\gamma)$ are ${\gamma, \sigma \gamma}$,
where $\sigma$ is the generator of $\Gal(E/F)$.
Observe that $\det \phi^A(\gamma) = a^2 - b^2 A = \N_{E/F}\gamma$, and $\phi^A$ defines an
embedding of $E^\times$ onto a maximal elliptic torus in $\GL(2, F)$.
Moreover, $\phi^A(\sigma \gamma)$ is equal to $\lp\lsm 1 & \\&-1\rsm\rp^{-1}\phi^A(\gamma)\lp\lsm 1&\\&-1\rsm\rp$.
\begin{claim}
Let $\delta_1 = (\diag(a, b), \gamma)$ be an element in $H_1$.  Then, $\delta_1$ is
a norm of \dt $\displaystyle \lp\lsm a &&\\ &\phi^A(\gamma)&\\&&b\rsm\rp \in G$.
\end{claim}
\begin{proof}
The image of $\delta_1$ in $\mb{T_1}(\bar{F})$ is $(\diag(a, b), \gamma, \sigma{\gamma})$.
The rest follows from Claim \ref{claim:normH1}.
\end{proof}
Note that since norm correspondence is defined up to stable conjugacy, the elements \dt $(\diag(a, b), \gamma)$ 
and $(\diag(a, b), \sigma\gamma)$ in $H_1$ are both norms of the same element in $G$.
\begin{claim}
Let $\lp\diag(\alpha,\beta),c\rp$ be an element in $\GL(2, E)\times F^\times$,
and $\delta_2$ its image in $H_2$.  Then, $\delta_2$ is a norm of
$\displaystyle \lp \lsm c\N_{E/F}\alpha&&\\
&c\phi^A(\alpha\sigma\beta)&\\&&c\N_{E/F}\beta\rsm\rp \in G$.
\end{claim}
\begin{proof}
The image of $\delta_2$ in $\mb{T_2}(\bar{F})$ is represented by 
\[
\lp \lp\lsm \alpha &\\&\beta\rsm\rp,c\lp\lsm\sigma\alpha&\\&\sigma\beta\rsm\rp, 1\rp
\in \GL(2, \bar{F})\times\GL(2, \bar{F}) \times \bar{F}^\times.
\]
The rest follows from Claim \ref{claim:normH2}.
\end{proof}
Note that the elements $(\diag(\alpha, \beta), c)$ and $(\diag(\sigma\alpha, \sigma\beta), c)$ in $H_2$
are both norms of the same element in $G$.
\subsubsection{Central Elements}
The center of $G$ is ${Z} = \{\diag(z, z, z, z) : z \in F^\times\}$.
The maximal $F$-split component of the center of $H_1$ is
\[
Z_0(H_1) = \left\{\lp\diag(z, z), z\rp : z \in F^\times\right\}.
\]
The maximal $F$-split component of the center of $H_2$ is
\[
Z_0(H_2) = \left\{\lp\diag(1, 1), z\rp_*: z \in F^\times\right\},
\]
\indexs{Z@$\mb{Z}_0$}%
where lower * denotes the equivalence class of the element modulo the group
\[\{(\diag(x, x),\N_{E/F}x^{-1}): x \in E^\times\}.\]
\begin{corollary}\label{corollary:normcenter}
The elements $(\diag(z, z), z)$ in $Z_0(H_1)$ and $(\diag(1, 1), z)_*$ in $Z_0(H_2)$ 
are norms of $\diag(z, z, z, z) \in Z$.
\end{corollary}
\subsubsection{Split Case}
Suppose $F$ is a number field.  Let $v$ be a place of $F$ where $\ve_v = 1$.  Then,
$
\mb{H_1}(F_v) = \{ (g, a, b) \in \GL(2, F_v) \times F_v^\times \times F_v^\times : \det g = ab\},
$
and $\mb{H_2}(F_v)$ is 
\[
\{(g_1, g_2, c) \in \GL(2, F_v)\times\GL(2, F_v)\times F_v^\times\}
/\{(a I_2, b I_2, (ab)^{-1}) : a, b \in F_v^\times\}.
\]
\begin{claim}
The element $\delta_1 = (\diag(c, d), a, b) \in \mb{H_1}(F_v)$ is a norm of the element
$\diag(c, a, b, d) \in \mb{G}(F_v)$.
\end{claim}
\begin{proof}
This follows from Claim \ref{claim:normH1}.
\end{proof}
\begin{claim}
The element $\delta_2 = (\diag(a, b), \diag(c, d), 1)_* \in \mb{H_2}(F_v)$ is a norm of \dt
$\diag(ac, ad, bc, bd) \in \mb{G}(F_v)$.
\end{claim}
\begin{proof}
This follows from Claim \ref{claim:normH2}.
\end{proof}
\indexi{norm correspondence|)}%
\section{Matching Functions}\label{sec:matchingfunctions}
\indexi{matching functions}%
Suppose $F$ is a local field.  Let $G = \mb{G}(F)$.
Let $\ve$ be a quadratic character of $F^\times$.
Fix a character $\omega$ of the center of $G$.
For any regular element $t \in G$, let $Z_G(t)$ denote the centralizer of $t$ in
$G$.  In particular, since $t$ is regular, $Z_G(t)$ is a maximal torus.
Let $C(G, \omega)$ denote the space of functions on $G$ which are
smooth, compactly supported modulo center, and transform under the center of $G$ via $\omega^{-1}$.
For $f \in C(G, \omega)$, put
\[
O_G(f, t) := \int_{Z_G(t)\bs G}f(g^{-1}tg)\ve(g)\, dg.
\indexi{orbital integral}%
\]
Let $\mb{H}$ be either $\mb{H_1}$ or $\mb{H_2}$.
Let $H = \mb{H}(F)$.
The character $\omega$ defines a character of $Z_0(H) \cong Z$.
Let $C(H, \omega)$ denote the space of functions on $H$ which are smooth, compactly supported
modulo $Z_0(H)$, and transform under $Z_0(H)$ via $\omega^{-1}$.
For any $\mb{G}$-regular element $t_H \in H$ and function $f_H \in C(H, \omega)$, put
\[
O_H(f_H, t_H) := \int_{Z_H(t_H)\bs H} f_H(h^{-1}t_H h)\,dh
\]
and
\[
\displaystyle SO_H(f_H, t_H) := \sum_{t_H'} O_H(f_H, t_H'),
\indexi{orbital integral!stable}%
\]
where the sum is taken over a set of representatives for the conjugacy classes
of elements $t_H' \in H$ in the stable conjugacy class of $t_H$.

Let $H'$ denote the subgroup of $\mb{G}$-regular elements in $H$.  
A conjecture of Langlands asserts that there exists a function 
$\Delta_{G/H}$ on $H' \times G$ with the following property:
For every $f \in C(G, \omega)$, there exists a function $f_H \in C(H, \omega)$ such that
\begin{equation}\label{eq:matchingdef}
SO_H(f_H, t_H) = \sum_t\Delta_{G/H}(t_H, t)O_G(f, t),\quad\forall\;t_H \in H',
\end{equation}
where the sum is over a set of 
representatives for the $F$-conjugacy classes of regular elements $t \in G$.  
We say that $f, f_H$ are {\bf matching functions} if the above identity holds.
The function $\Delta_{G/H}$, called the {\bf transfer factor},
\indexi{transfer factor}%
is defined in \cite{LS} in the standard case (where $\ve = 1$) 
and in \cite{KS} in a more general twisted context.  In certain special cases, 
one may find a description of it in \cite{H1}, \cite{H2}.
The transfer factor is expected to satisfy a host of properties (see \cite{KS}).
A full discussion of the transfer factor is beyond the scope of this work,
but we do specify it explicitly in the unramified case, which is needed for the Fundamental Lemma (see below).
We list here some properties of $\Delta_{G/H}$ which are of interest to us:
\begin{itemize}
\item
$\Delta_{G/H}(t_H, t) = 0$ if $t_H$ is not a norm of $t$.
\item
For all $g \in G$, $\Delta(t_H, g^{-1}tg) = \ve(g)\Delta(t_H, t)$.
\end{itemize}
Note that, without the second property, the right hand side of \eqref{eq:matchingdef} is
not well defined.

Suppose now that $F$ is a number field.  Let $\ve$ be a quadratic character of $\idc{F}$.
Let $\mb{H}$ be an elliptic $\ve$-endoscopic group of $\mb{G}$ over $F$.
Let $G = \mb{G}(F)$ and $G_v = \mb{G}(F_v)$ for any place $v$ of $F$.  Define $H, H_v$ likewise.
Suppose $f \in C(\mb{G}(\Af), \omega)$ and $f_H \in C(\mb{H}(\Af), \omega)$
are products of local functions $f_v, f_{H,v}$, respectively.
We say that $f$ and $f_H$ are {\bf matching functions } 
if $f_v$ and $f_{H,v}$ are matching for all $v \in V$.
\subsection{The Fundamental Lemma}
\indexi{Fundamental Lemma}%
An ingredient in the conjecture on the existence of matching functions is the
``Fundamental Lemma'', which is actually a conjecture.  
Let $F$ be a local field.  Let $\ve$ be a quadratic character of ${F^\times}$.
\indexs{epsilon@$\ve$}%
Let $\mb{H}$ be an elliptic $\ve$-endoscopic group of $\mb{G}$ over $F$.
Let $\mc{O}$  be the ring of integers in $F$.  
Assume that $\mb{H}$ is defined over $\mc{O}$.
Let $K = \mb{G}(\mc{O})$,
$K_H = \mb{H}(\mc{O})$.
Let $\mc{H}(G)$ be the Hecke algebra of smooth, compactly supported, $K$-biinvariant functions
\indexs{H@$\mc{H}(\;)$}%
on $G$.  
\indexi{Hecke algebra}%
Define $\mc{H}(H)$ likewise, with $K$ replaced with $K_H$.
The $L$-group embedding $\xi : {}^L H \rightarrow \hat{G}$ induces a map
\indexs{xi@$\xi$}%
\indexi{L@$L$-!embedding}%
\[
b_\xi : \mc{H}(G) \rightarrow \mc{H}(H).
\]
\indexs{b@$b_\xi$}%
The Fundamental Lemma asserts that $f$ and $b_\xi(f)$ are matching functions for any $f$ in $\mc{H}(G)$.

Suppose $F$ is a number field.  
Hales shows in \cite{H} that if one can show that the unit elements of
$\mc{H}(G_v)$ and $\mc{H}(H_v)$ have matching orbital integrals for
almost all finite places $v$, then the Fundamental Lemma follows.
Note that even though the title of \cite{H} reads ``standard endoscopy,''
the paper does in fact cover our situation.  
Using this result of Hales' and the results of Flicker's in \cite{F}, we prove in Appendix \ref{chap:fundlemma} the Fundamental
Lemma in the context of $\mb{G}$ and its $\ve$-endoscopic groups.
\chapter{The Trace Formula}\label{chap:tf}





\section{The fine $\chi$-expansion}\label{sec:finechiexp}
\indexi{fine $\chi$-expansion|(}%
\subsection{An Overview}\label{sec:finechiexpoverview}
Let $F$ be a number field.  Let $V$ be the set of places of $F$.
For any finite place $v$, let $\mc{O}_v$ be the ring of integers of $F_v$.
\indexs{O@$\mc{O}_v$}%
Let $\mb{H}$ be a reductive $F$-group.  
Let $\mb{Z}_0$ be the maximal $F$-split component of the center $\mb{Z}$ of $\mb{H}$.
\indexs{Z@$\mb{Z}_0$}%
Let $\ve$ be a character of $\mb{H}(\Af)$ whose restriction to $\mb{Z}(\Af)$ is trivial.
From now on, by a group we mean an $F$-group.

Fix a character $\omega$ of $\mb{Z}_0(F)\bs\mb{Z}_0(\Af)$.
For any place $v$ of $F$, let $H_v = \mb{H}(F_v)$.
Let $C(H_v, \omega_v)$ be the space of smooth functions $f_v$ on $H_v$
\indexs{C@$C(H_v,\omega_v)$}%
such that $f_v$ is compactly supported
modulo $Z_{0,v}$ and $f_v(zh) = \omega_v^{-1}(z)f(h)$ for all $z \in Z_{0,v}, h \in H_v$.

At a finite place $v$ where $\mb{H}$ is defined over $\mc{O}_v$,
let $K_v$ denote the hyperspecial, maximal compact subgroup $\mb{H}(\mathcal{O}_v)$ of
$H_v$.
Let $\mathcal{H}(H_v, \omega_v)$ denote the Hecke algebra of $K_v$-biinvariant functions
in $C(H_v, \omega_v)$.  
If $v$ is archimedean, fix a maximal compact subgroup of $H_v$, and let
$\mathcal{H}(H_v, \omega_v)$  be the set of $K_v$-finite functions in $C(H_v, \omega_v)$.
\indexs{H@$\mathcal{H}(H_v,\omega_v)$}%

Let $C(\HH, \omega)$ denote the span of the smooth, compactly supported functions on $\HH$ which are
\indexs{C@$C(\HH,\omega)$}%
of the form $\otimes_v f_v$, where $f_v \in C(H_v, \omega_v)$
for all $v$ and $f_v$ is a unit in the Hecke algebra $\mc{H}(H_v, \omega_v)$ for almost all finite $v$.
In this work, whenever we mention a function $f \in C(\HH, \omega)$, we assume that
$f$ is the  product of local components $f_v$.  
Since such functions span $C(\HH, \omega)$,
our assumption has no ill effect on our attempt to understand the spectrum of $\HH$.

Fix a minimal ($F$-)parabolic subgroup $\mb{P}_0$ of $\mb{H}$.  Let $\mb{A}_0$ be the maximal $F$-split component
\indexs{P@$\mb{P}_0$}\indexs{A@$\mb{A}_0$}%
of the ($F$-)Levi subgroup of $\mb{P}_0$.
\indexi{Levi subgroup}%
Let $W(A_0, H)$ be the Weyl group of $\mb{A}_0$ in $\mb{H}$.
\indexs{W@$W(A_0,H)$}%
For any Levi subgroup $\mb{M}$ of $\mb{H}$, let $\mb{A}_M$ denote the
\indexs{A@$\mb{A}_M$}%
split component of the center of $\mb{M}$.
Let $X_*(\mb{A}_M) = \Hom(\mbb{G}_m, \mb{A}_M)$.
\indexs{X@$X_*(\mb{A}_M)$}%
Let $\mf{a}_M$ denote $X_*(\mb{A}_M)\otimes_{\mbb{Z}}\mbb{R}$.  Let $\mf{a}_M^*$ denote the dual of $\mf{a}_M$.
\indexs{a@$\mf{a}_M$}\indexs{a@$\mf{a}_M^*$}%
Let $\mc{P}(M)$ denote the set of all ($F$-)parabolic subgroups of $\mb{H}$ containing $\mb{M}$.\indexs{P@$\mc{P}(\;)$}\indexi{parabolic subgroup}%

We fix a maximal compact subgroup $\mbb{K} = \otimes_{v \in V} K_v$ in $\mb{H}(\Af)$, where
\indexs{K@$\mbb{K}$}%
$K_v = \mb{H}(\mc{O}_v)$ for every finite place $v$ such that $\mb{H}$ is defined over $\mc{O}_v$,
and $K_v$ is a fixed maximal compact subgroup in $\mb{H}(F_v)$ for the rest of the places $v$.
\indexi{maximal compact subgroup}%

Let $\mb{M}$  be a Levi subgroup in $\mb{H}$.
For any $\zeta \in \mf{a}_M^*$, let $\chi_\zeta$ denote the character of $\mb{M}(\Af)$ associated with $\zeta$.
\indexs{chi@$\chi_\zeta$}%
For any $m \in \mb{M}(\Af)$, let $\rH_P(m)$ be the element in $\mf{a}_M$ uniquely defined by the condition: 
\indexs{H@$\rH_P(\;)$}%
$e^{\la \zeta, \rH_P(m)\ra} = \chi_\zeta(m)$ for all $\zeta \in \mf{a}_M^*$.

For any parabolic subgroup $\mb{P} \in \mc{P}(M)$, representation $\tau$ of $\mb{M}(\Af)$, and element $\zeta \in \mf{a}_M^*$,
let $I_{P, \tau}(\zeta)$ denote the $\HH$-module normalizedly induced from the $\mb{M}(\Af)$-module
\indexs{I@$I_{P,\tau}(\;)$}\indexi{representation!normalizedly induced}%
which sends an element $m$ in $\mb{M}(\Af)$ to
$\tau(m)e^{\la \zeta, H_P(m)\ra}$.  In other words, $I_{P, \tau}(\zeta)$ is the right regular action on the space
\indexi{representation!right regular}%
of smooth functions $\phi$ on $\HH$ which satisfy:
$\phi(m g) = (\delta^{1/2}\tau)(m)\phi(g)$ for all  
$m \in \mb{M}(\Af)$ and $g \in \HH$.  
Here, $\delta(m)$ is defined as $\abs{\det\lp{\rm Ad}\;m|_{\mf{n}}\rp}$,
where $\mf{n}$ denotes the Lie algebra of the unipotent component of $\mb{P}$.
\indexs{delta@$\delta(\;)$}%
For any function $f$ in $C(\mb{H}(\Af), \omega)$, let
$I_{P, \tau}(\zeta, f)$ denote the convolution operator
\indexs{I@$I_{P,\tau}(\;,\;)$}\indexi{convolution operator}%
\[
\int_{\mb{Z}_0(\Af)\bs\mb{H}(\Af)}\lp I_{P, \tau}(\zeta)\rp(h)f(h)\;dh,
\]
where $dh$ is a fixed Tamagawa measure on $\mb{H}(\Af)$.
\indexi{convolution operator}\indexi{Tamagawa measure}%

Let $W_H(M)$ denote the group of automorphisms of $\mf{a}_M$ obtained by restricting elements
\indexs{W@$W(\;)$}%
$s \in W(A_0, H)$ which satisfy $s(\mb{M}) = \mb{M}$.
We define an action of $W_H(M)$ on the set of all $\mb{M}(\Af)$-modules as follows:
For any element $s$ in $W_H(M)$ and any $\mb{M}(\Af)$-module $\tau$,
let $s\tau$ denote the $\mb{M}(\Af)$-module: 
\[
m \mapsto  \tau((s')^{-1}m), \quad \forall m \in \mb{M}(\Af),
\]
where $s'$ is a preimage of $s$ in $W(A_0, H)$.  Up to equivalence of representations,
$s\tau$ is independent of this choice of preimage.

Suppose $\tau$ is equivalent to $s\ve\tau$ for some $s \in W_G(M)$.
Let $I_{P, \tau}(\ve)$ be the operator on the space of
\indexs{I@$I_{P,\tau}(\ve)$}%
$I_{P, \tau}$ which sends $\phi(g)$ to $\ve(g)\phi(g)$.  Then, $I_{P, \tau}(\ve)$ intertwines
\indexi{intertwining operator}%
$I_{P, \tau}(\zeta)$ with $\ve I_{P, \ve\tau}(\zeta)$ for any $\zeta \in \mf{a}_M$.
Let $M_P(s, \zeta)$ be an operator which intertwines $I_{P, \tau}(\zeta)$ with
\indexs{M@$M_P(\;,\;)$}%
$I_{P, s\tau}(s\zeta)$.  
Then, $I_{P, s\tau}(\ve)M_P(s, \zeta)$ intertwines $I_{P, \tau}(\zeta)$
with $\ve I_{P, s\ve\tau}(s\zeta) = \ve I_{P, \tau}(s\zeta)$.
Put 
\[
I_{P, \tau}(\zeta, f\times\ve) := I_{P, \tau}(\zeta, f)I_{P, \zeta}(\ve).
\indexs{I@$I_{P,\tau}(\;,\;\times\ve)$}%
\]


The $\ve$-twisted trace formula (\cite{CLL}) is the equality
\indexi{trace formula!epsilon-twisted@$\ve$-twisted}%
\begin{equation}\label{eq:twistedTF}
\sum_{\{\mc{O}\}} J_{\mc{O}, \ve}^T(f) = \sum_{\{\chi\}}J_{\chi, \ve}^T(f),
\indexs{J@\JTO}\indexs{J@\JTchi}%
\end{equation}
where the left (resp. right) hand side
of the equation is the integral over the diagonal subgroup of the modified geometric (resp. spectral) kernel of
the operator $\rho(f)\rho(\ve)$ (see Chapter \ref{chap:intro} for the definitions of symbols).  The modification
depends on an element $T$, called a {\bf truncation parameter},  
\indexi{truncation}%
in the Lie algebra of a fixed minimal Levi subgroup of $\mb{H}$.
Each side of \eqref{eq:twistedTF} is defined only when $T$ is sufficiently regular (see \cite{CLL}).
The sum $\sum_{\{\chi\}}J_{\chi, \ve}^T(f)$ is called the (twisted-){\bf fine $\chi$-expansion} of $\mb{H}$.
\indexi{fine $\chi$-expansion}%

It is important to note that the authors of \cite{CLL} have skipped giving a precise form of the fine $\chi$-expansion
for the twisting by $\rho(\ve)$.  We take the liberty here to describe without proof the sum $\sum_{\{\chi\}}J_{\chi, \ve}^T(f)$
(For instance, we can bypass the twisting by considering the untwisted trace formula for the group described in Lemma
\ref{lemma:venondegeneracy}).

For $f$ in $C(\mb{H}(\Af), \omega)$, the fine $\chi$-expansion
is a sum over the set of quadruples
\indexi{quadruple}%
$\{\chi\} = (\mb{M}, \mb{L}, \tau, s)$ consisting of Levi subgroups $\mb{M}, \mb{L}$ of $\mb{H}$, 
an element $s \in W_G(M)$, 
and a discrete spectrum automorphic representation $\tau$ of $\mb{M}(\Af)$ such that:
\indexi{representation!discrete spectrum}%
\begin{itemize}
\item
  $\mb{M} \subset \mb{L}$;
\item
  $\mf{a}_L$ is the subspace of $\mf{a}_M$ fixed pointwise by $s$;
\item
  $\tau$ is trivial on $\exp \mf{a}_M$;
\item
  the restriction of $\tau$ to $\Zo$ is equal to $\omega$;
\item
  $s\tau$ is equivalent to $\ve\tau$.
\end{itemize}

We say that a parabolic subgroup is \emph{standard} if it contains the fixed minimal parabolic subgroup $\mb{P}_0$.
\indexi{parabolic subgroup!standard}%
Let $\mb{P} \in \mc{P}(M)$ be the standard parabolic subgroup with Levi component $\mb{M}$.
The term associated with $(\mb{M}, \mb{L}, \tau, s)$ in the fine $\chi$-expansion is the product of
\begin{equation}\label{eq:finechiexp1}
\frac{\abs{W(A_0, M)}}{\abs{W(A_0, H)}}
\abs{ \det (1 - s)|_{\mf{a}_M/\mf{a}_L}}^{-1}
\end{equation}
and
\begin{equation}\label{eq:finechiexp2}
(2\pi)^{-\dim \lp\mf{a}_L / \mf{a}_H\rp}\int_{i \mf{a}_L^*/ i \mf{a}_H^*}
\tr 
\mc{M}^T_L(P, \zeta) 
I_{P, \tau}(\zeta, f\times\ve) M_P(s, \zeta)d\zeta.
\end{equation}
Here, $\mc{M}^T_L(P, \zeta)$ is the logarithmic derivative of an intertwining operator.  We will not
\indexi{logarithmic derivative}%
reproduce here the defintion of $\mc{M}^T_L(P, \zeta)$.  We only use the fact that $\mc{M}^T_L(P, \zeta) = 1$
when $\mb{L}$ is equal $\mb{H}$.
In particular,
the term in the expansion associated with a quadruple of the form $(\mb{M}, \mb{H}, \tau, s)$ is
\begin{equation}\label{eq:finechidiscrete}
\frac{\abs{W(A_0, M)}}{\abs{W(A_0, H)}}
\abs{ \det (1 - s)|_{\mf{a}_M/\mf{a}_L}}^{-1}
\tr I_{P, \tau}(\zeta, f\times\ve) M_P(s, 0).
\end{equation}
We call the sum of all the terms in the form \eqref{eq:finechidiscrete} the {\bf discrete part} of the fine $\chi$-expansion
of the $\ve$-twisted trace formula.
We denote it by $I^d(H, f,\ve)$.
We call the sum of the rest of the terms the {\bf continuous part} and denote it by $I^c(H, f,\ve)$.
\indexi{trace formula!discrete part}\indexi{trace formula!continuous part}%
If $\ve = 1$, we simply write $I^d(H, f)$ (resp. $I^c(H, f)$) for $I^d(H, f, 1)$ (resp. $I^c(H, f, 1)$).
Observe that the $\ve$-twisted trace of a discrete spectrum representation $\pi$ of $\HH$ corresponds
to the the quadruple $(\mb{H}, \mb{H}, \pi, 1)$.  
For any quadruple of the form $(\mb{M}, \mb{H}, \tau, s)$ and parabolic subgroup $\mb{P}$ in $\mc{P}(M)$,
we say that the (normalizedly) parabolically induced representation 
$I_{P, \tau}$ {\bf occurs $\ve$-discretely} in the spectrum of $\HH$.
\indexi{representation!discretely occuring}%

In this work, we only study the discrete parts of the spectral expansions of the groups.  We now
list the terms which occur in $I^d(\GSp(2), f,\ve)$ and $I^d(H_i, f_i)$ ($i = 1, 2$).

{\bf Notation: }
\begin{itemize}
\item
For any number field $L$, let $C_L$ denote the \idele class group $L^\times\bs\mbb{A}_L^\times$.
\item
For any character $\theta$ of $\idc{E}$, put ${}^\sigma \theta(x) := \theta(\sigma x)$ for all $x \in \Ae^\times$.
\indexs{sigma@${}^\sigma{\rm char.}$}%
We say that $\theta$ is {\bf $\sigma$-invariant} if $\theta = {}^\sigma \theta$.
\indexi{sigma@$\sigma$-invariant}%
\item
Let $\omega_\pi$ denote the central character of a representation $\pi$.
\end{itemize}

{\bf Terminology: }
\begin{itemize}
\item Let $L$ be a quadratic extension of $F$.
We say that an automorphic representation $\pi$ of $\GL(2, \Af)$ is {\bf $L$-monomial} if it is the monomial
representation $\pi(\theta)$ associated with a character $\theta$ of $\idc{L}$.  
\indexi{representation!monomial}%
In particular, if $\mu$ is the class field character of
$L/F$, 
then $\mu\pi(\theta)$ is equivalent to $\pi(\theta)$, and the central character $\omega_{\pi(\theta)}$ is equal to 
$\theta|_{\Af^\times}\cdot\mu$ 
(see \cite[Thm. 4.6]{JL}, \cite{K}).
\item Let ($\ve$-){\bf DOR} stand for ``($\ve$)-discretely occurring representation.''
\end{itemize}
\subsection{$\ve$-DOR of $\GSp(2, \Af)$}
Let $\mb{G} = \GSp(2)$.   Let $\mb{Z}$ be the center of $\mb{G}$.  It is isomorphic to $\mbb{G}_m$.  In particular,
the maximal $F$-split component $\mb{Z}_0$ of the center is $\mb{Z}$ itself.  We fix a character
$\omega$ of $\mb{Z}_0(F)\bs\Zo$.

Recall that we have fixed a quadratic character $\ve$ of $\idc{F}$.
It induces a character of $\GSp(2, \Af)$ defined by $\ve(g) = \ve(\lambda(g))$.
In this  section, we list all quadruples $(\mb{M}, \mb{G}, \tau, s)$ which satisfy the conditions
mentioned in the previous section.

If the quadruple has the form $(\mb{G}, \mb{G}, \tau, s)$, then $s$ is necessarily trivial, and
$\tau$ can be any discrete spectrum automorphic representation of $\mb{G}(\Af)$ with the property
that $\omega_\tau = \omega$.
In this case, the intertwining operator $M_G(1, 0)$ is equal to $1$.
The contribution to the fine $\chi$-expansion corresponding to
$(\mb{G}, \mb{G}, \tau, 1)$ is therefore
\[
\tr \tau(f\times\ve), 
\indexi{convolution operator}%
\]
where $\tau(f \times \ve)$ denotes $\tau(f)\rho(\ve)$.

We now list all the valid quadruples for the proper Levi subgroups $\mb{M}$ of $\mb{G}$.
There are 5 proper Levi subgroups: $\mb{A}_0, \mb{M}_\alpha, \mb{M}_{\alpha}', \mb{M}_\beta, \mb{M}_\beta'$,
which we shall describe shortly.

We choose $\mb{P}_0$ to be the upper triangular subgroup in $\mb{G}$.  Its Levi subgroup is 
$\mb{A}_0 = \{\diag(a, b, \lambda/b, \lambda/a) : a, b, \lambda \in \mathbb{G}_m\}$.  The $F$-split component
\indexs{A@$\mb{A}_0$}%
of the center of $\mb{A}_0$ is $\mb{A}_0$ itself.
A basis of the root system associated with the pair $(\mb{G}, \mb{A}_0)$ is 
\indexi{root system}%
\[
\Delta = \Delta(\mb{G}, \mb{P}_0, \mb{A}_0) = \{\alpha, \beta\},
\] 
where $\alpha(t) := a /b$, $\beta(t) := b^2 / \lambda$
for all $t = \diag(a, b, \lambda/a, \lambda/b)$.
The root system is $R = R^+ \cup -R^+$, where
$R^+ = R^+(\mb{G}, \mb{P}_0, \mb{A}_0) = \{\alpha, \beta, \alpha + \beta, 2\alpha+\beta\}$ 
is a fixed set of positive roots.
Let $X_*(\mb{A}_0) = \Hom(\mbb{G}_m, \mb{A}_0)$.
Let $\mf{a}_0 = X_*(\mb{A}_0)\otimes \mbb{R}$.  We identify $\mf{a}_0$ with the real vector space
$\{(x, y, t -x, t- y) : x, y, t \in \mbb{R}\} \cong \mbb{R}^3$.

Let $\mc{V}$ be the space of weights (\cite{S}) associated with the pair $(\mb{G}, \mb{A}_0)$.
\indexs{V@$\mc{V}$}%
The Weyl group $W(A_0, G)$ of $\mb{A}_0$ in $\mb{G}$
corresponds to the group of automorphisms of $\mc{V}$ which are generated 
by reflections associated with the roots in $R^+$.
For a root $\gamma \in R^+$,
let $s_\gamma$ denote the Weyl group element which corresponds to the reflection over $\gamma$.

The Weyl group $W(A_0, G)$ is isomorphic to $D_4$, the dihedral group of order 8.
Viewed as a permutation group on four letters (corresponding to the four entries of an element in $\mb{A}_0$),
it consists of the following elements:
\begin{itemize}
\item
$s_\alpha = (12)(34)$,
\item
$s_\beta = (23)$,
\item
$s_{2\alpha + \beta}s_\alpha = (1243)$,
\item
$\rho := s_\beta s_\alpha = (1342)$,
\item
$s_{\alpha + \beta} = (24)(13)$,
\item
$s_{2\alpha + \beta} s_\beta = (23)(14)$,
\item
$s_{2 \alpha + \beta} = (14)$.
\end{itemize}
\indexs{s@$s_\alpha$, $s_\beta$,\ldots, etc.}%
Note that $\rho = s_{\beta}s_{\alpha}$ has order $4$, and
$
W(A_0, G) = \la s_\beta, \rho \ra.
$

Let $e$ denote the element $\diag(1, -1)$ in $\GL(2)$.
The proper Levi subgroups of $\mb{G}$ which strictly contain $\mb{A}_0$ are:
\[
\begin{split}
\mb{M}_\alpha = \left\{\blockdiag\lp g, \frac{\lambda}{\det g} e g e\rp : g \in \GL(2), \lambda \in \mbb{G}_m \right\},&\quad
\mb{M}_\alpha' = s_\beta \mb{M}_\alpha s_\beta,\\
\mb{M}_\beta = \{\blockdiag(a, g, (\det g) / a) : a \in \mbb{G}_m, g \in \GL(2)\},&\quad \mb{M}_\beta' = s_\alpha \mb{M}_\beta s_\alpha.
\end{split}
\]
\indexs{M@$\mb{M}_\alpha$, $\mb{M}_\beta$, \ldots , etc.}%

\subsubsection{$\ve$-DOR associated with $\mb{A}_0$}
Suppose $(\mb{A}_0, \mb{G}, \tau, s)$ is a quadruple satisfying the conditions outlined in Section 
\ref{sec:finechiexpoverview}.
The discrete representation $\tau$ of $\mb{A}_0(\Af)$ has the form: 
\[
\tau = \mu_1\otimes\mu_2 \otimes\mu:
\diag(a, b, \lambda /b, \lambda, a) \mapsto \mu_1(a)\mu_2(b)\mu(\lambda)
\indexs{t@$\;\otimes\;\otimes\;$}%
\]
for some characters $\mu_1, \mu_2, \mu$ of $\idc{F}$.
Let $\mu_1\times\mu_2\rtimes\mu$ denote the (normalizedly) parabolically 
\indexs{t@$\;\times\;\rtimes\;$}%
induced representation $I_{P_0}(\mu_1\otimes\mu_2\otimes\mu)$ of $\mb{G}(\Af)$.

For $\mf{a}_0^s$ to equal $\mf{a}_G$, the element $s$ in $W_G(A_0) = W(A_0, G)$ must be one of $\rho, \rho^2, \rho^3$.
By assumption, $s\tau \cong \ve\tau$ and $\omega_\tau = \omega$, which implies that the following holds:
\begin{itemize}
\item
If $s = \rho$ or  $\rho^3$, then $\tau$ must have the form $\ve \otimes \ve \otimes \mu$ for some character
$\mu$ of $\idc{F}$ such that $\mu^2 = \omega$.
\item
If $s = \rho^2$, then $\tau$ must have the form $\ve'\otimes \ve' \ve \otimes \mu$ for some characters
$\ve', \mu$ of $\idc{F}$ such that $\lp\ve'\rp^2 = 1$ and $\ve\mu^2 = \omega$.
\end{itemize}
We have: $\abs{W(A_0, A_0)} = 1$, $\abs{W(A_0, G)} = 8$.
For $s = \rho, \rho^3$, $\det\lp(1 - s)|_{\mf{a}_0/\mf{a}_G}\rp$ is equal to $2$, while
$\abs{\det ( 1 - \rho^2)|_{\mf{a}_0/\mf{a}_G}} = 4$.  

{\bf Notation: }
For any representation $\tau$ of $\mb{M}(\Af)$ and element $w$ in the Weyl group $W(A_0, G)$,
let ${}^w \tau$ denote $\tau\circ w$.

The orbit of $\ve\otimes\ve\otimes\mu$ under $W(A_0, G) = D_4$ is of order $2$.
Consequently,
the set of quadruples $\{ (\mb{A}_0, \mb{G}, {}^w\tau, s) : s = \rho, \rho^2, w \in W(A_0, G)\}$
has $4$ distinct elements, and they all give rise to the same representation
$I_{P_0}\lp\ve\otimes\ve\otimes\mu\rp = \ve\times\ve\rtimes\mu$.
Since $\ve\times\ve\rtimes\mu$ is irreducible, by Schur's lemma
we may normalize the intertwining operators so that $M_{P_0}(\rho, 0) = M_{P_0}(\rho^2, 0) = 1$.

We conclude that the contribution from $\ve\times\ve\rtimes\mu$ to the fine $\chi$-expansion is
\[
4 \cdot \frac{1}{8}\frac{1}{2}\cdot\tr \lp \ve\times\ve\rtimes \mu\rp(f\times\ve).
\]

Similarly, the orbit of $\ve'\otimes\ve'\ve\otimes\mu$ under $W(A_0, G)$ is of order $8$, and
the irreducibility of $\ve'\times\ve'\ve\rtimes\mu$ implies that $M_{P_0}(\rho^2, 0) = 1$.
Consequently, the contribution from $\ve'\times\ve'\ve\rtimes\mu$ to the fine $\chi$-expansion is
\[
8 \cdot \frac{1}{8}\frac{1}{4}\cdot\tr \lp \ve'\times\ve'\ve\rtimes \mu\rp(f\times\ve).
\]

To summarize, the contribution to the fine $\chi$-expansion 
from quadruples of the form $(\mb{A}_0, \mb{G}, \tau, s)$ is:
\begin{multline}
\frac{1}{4}
\sum_{\mu^2 = \omega} \tr\!\!\!\lp \ve\times\ve\rtimes \mu\rp(f\times\ve)
+ \frac{1}{4} 
\sum_{\substack{{\ve'}^2 = 1,\\\ve'\mu^2 = \omega}} \tr\!\!\!\lp \ve'\times\ve'\ve\rtimes \mu\rp(f\times\ve).
\end{multline}

\subsubsection{$\ve$-DOR associated with $\mb{M}_\alpha$}
Suppose $(\mb{M}_\alpha, \mb{G}, \tau, s)$ is a quadruple satisfying the conditions outlined in
Section \ref{sec:finechiexpoverview}.  
The discrete spectrum representation $\tau$ of $\mb{M}_\alpha(\Af)$ has the form:
\[
\tau = \pi_2 \otimes \mu:
\blockdiag\lp g, \frac{\lambda}{\det g} e ge \rp \mapsto \mu(\lambda) \pi_2(g),
\indexs{t@repn.\,$\otimes$ char.}%
\]
where $\pi_2$ is a discrete spectrum representation of $\GL(2, \Af)$ and $\mu$ is a character of
$\idc{F}$.  Let $\mb{P}_\alpha$ be the standard parabolic subgroup of $\mb{G}$ whose Levi component is
$\mb{M}_\alpha$.  Put
\[
\pi_2\rtimes\mu 
:= I_{P_\alpha}(\pi_2\otimes\mu).
\indexs{t@repn.\,$\rtimes$ char.}%
\]

The group $W_G(M_\alpha)$ is represented by the elements 
$\{1, 
(24)(13)\}$ in $W(A_0, G)$, and $\mf{a}_{M_\alpha}$ is
equal to $\{ (x, x, t - x, t - x): x, t \in \mbb{R}\} \cong \mbb{R}^2$.  For $\mf{a}_M^s$ to equal
$\mf{a}_G$, the element $s \in W_G(M_\alpha)$ must be the restriction of $s_{\alpha + \beta} = (24)(13)$.
The following statements hold:
\begin{itemize}
\item
$\pi_2$ is $\ve$-invariant; i.e., $\pi_2$ is equivalent to the representation
\[
\ve\pi_2 : g \mapsto \ve(\det g)\pi_2(g),\quad \forall g \in \GL(2, \Af).
\]
\item
$\omega_{\pi_2} = \ve$ and $\ve\mu^2 = \omega$.
\end{itemize}
By \cite[Thm. 4.6]{JL}, there exists a character $\eta$ of $\idc{E}$ such that 
$\eta\neq {}^\sigma \eta$ and $\pi_2$ is the cuspidal monomial representation $\pi(\eta)$ associated with $\eta$.
Moreover, since $\omega_{\pi(\eta)}$ is equal to $\eta|_{\Af^\times}\cdot\ve$, the character $\eta|_{\Af^\times}$
must be trivial.

We have: $\abs{W(A_0, M_\alpha)} = 2$,
$\abs{\det\lp (1 - s_{\alpha+\beta})|_{\mf{a}_{M_\alpha}/\mf{a}_G}\rp} = 2$.
The orbit of the representation $\tau = \pi(\eta)\otimes\mu$ under $W(A_{M_\alpha}, G)$ is of order $2$, and both 
elements in the orbit give rise to the same induced representation $\pi(\eta) \rtimes \mu$ of
$\mb{G}(\Af)$.
Since $\pi(\eta)\rtimes \mu$ is irreducible (\cite{ST}), we may normalize the
intertwining operator $M_{P_\alpha}(s_{\alpha+\beta}, 0)$ to be $1$.

For each $\{\chi\} = (\mb{M}_\alpha, \mb{G}, \tau, s)$, 
there is a representation $\tau'$ of $\mb{M}_\alpha'(\Af)$ and an element $s'$ in $W_G(M_\alpha')$
such that the quadruple $(\mb{M}_\alpha', \mb{G}, \tau', s')$ contributes a term which is equal to the contribution
from $\{\chi\}$.
The contribution of $\tr\! \lp\pi(\eta)\rtimes\mu\rp(f\times \ve)$ to the fine $\chi$-expansion is therefore:
\begin{equation}\label{eq:finechiexpalpha}
2\cdot 2 \cdot \frac{1}{4} \cdot\frac{1}{2} 
\sum_{\eta|_{\Af^\times} = 1,\;\ve\mu^2 =\omega} 
\tr \lp \pi(\eta)\rtimes \mu\rp(f\times\ve).
\end{equation}

\subsubsection{$\ve$-DOR associated with $\mb{M}_\beta$}\label{sec:finechiexpbeta}
Suppose $(\mb{M}_\beta, \mb{G}, \tau, s)$ is a quadruple satisfying the conditions listed in 
Section \ref{sec:finechiexpoverview}.  
The discrete representation $\tau$ of $\mb{M}_\beta(\Af)$ has the form:
\[
\tau = \mu \otimes \pi_2:
\blockdiag\lp a, g, \frac{\det g}{a}  \rp \mapsto \mu(a) \pi_2(g),
\indexs{t@char.\,$\otimes$ repn.}%
\]
where $\pi_2$ is a discrete spectrum representation of $\GL(2, \Af)$ and $\mu$ is a character of
$\idc{F}$.  Let $\mb{P}_\beta$ be the standard parabolic subgroup of $\mb{G}$ whose Levi component is
$\mb{M}_\beta$.  Put
\[
\mu\rtimes\pi_2 := I_{P_\beta}(\mu\otimes\pi_2).
\indexs{t@char.\,$\rtimes$ repn.}%
\]

The group $W_G(M_\beta)$ is represented by $\{1, (14)\}$ in $W(A_0, G)$, and $\mf{a}_{M_\beta}$ is
equal to $\{ (x, y, y, 2y - x): x, y \in \mbb{R}\} \cong \mbb{R}^2$.  For $\mf{a}_M^s$ to equal
$\mf{a}_G$, the element $s \in W_G(M_\alpha)$ must be the restriction of $s_{2\alpha + \beta} = (14)$.  
The following statements hold:
\begin{itemize}
\item
$\mu^2 = 1.$ 
\item
$\pi_2$ is $\mu\ve$-invariant.
\item
$\mu\omega_{\pi_2} = \omega$.
\end{itemize}
The classification can be refined further:
\begin{itemize}
\item
If $\mu = 1$, then $\pi_2 = \pi(\chi)$ for some character $\chi$ of $\idc{E}$ such that
$\chi \neq {}^\sigma \chi$ and $\chi|_{\Af^\times}\cdot\ve = \omega$ (see \cite{JL}, \cite{K}).
Observe that the orbit of $1\otimes\pi(\chi)$ under $W_G(M_\beta)$ is of
order $1$.
\item
If $\mu \neq 1, \ve$, let $E_{\mu\ve}$ be the quadratic extension of $F$ corresponding to
$\mu\ve$ via global class field theory.  Suppose $\Gal(E_{\mu\ve}/F) = \la \sigma'\ra$.  Then,
$\pi_2 = \pi(\theta)$ for some  character $\theta$ of $\idc{E_{\mu\ve}}$
such that
$\theta \neq {}^{\sigma'}\!\theta$ and $\theta|_{\Af^\times}\cdot\ve = \omega$.

Let $E_\mu$ be the quadratic extension of $F$ corresponding  to $\mu$.
The orbit of $\mu\otimes\pi(\theta)$ under $W_G(M_\beta)$ is of order $2$ if
$\pi(\theta)$ is not $E_\mu$-monomial 
(i.e., not a monomial representation associated with a character of $\idc{E_\mu}$).
It is of order $1$ if $\pi(\theta)$ is $E_\mu$-monomial. 
\item
If $\mu = \ve$, then $\pi_2$ can be any cuspidal or one dimensional automorphic representation
of $\GL(2, \Af)$ with central character $\omega_{\pi_2} = \ve\omega$.
The orbit of $\ve\otimes\pi_2$ under $W_G(M_\beta)$ is of order $1$ if $\pi_2$ is $E$-monomial and
of order $2$ otherwise.
\end{itemize}

We have: $\abs{W(A_0, M_\beta)} = 2$,
$\abs{\det\lp (1 - s_{2\alpha+\beta})|_{\mf{a}_{M_\alpha}/\mf{a}_G}\rp} = 2$.
As in the case of $\mb{M}_\alpha$, for every contribution to the spectral expansion from
a quadruple of the form $(\mb{M}_\beta, \mb{G}, \tau, s)$,
there is an equal contribution from a quadruple of the form $(\mb{M}_\beta', \mb{G}, \tau', s')$.
Moreover, the elements in the orbit of $\tau = \mu\otimes\pi_2$ under $W_G(M_\beta)$
all give rise to the same representation $\mu \rtimes \pi_2$.
Consequently, the contribution from representations of the form $\mu\rtimes \pi_2$ to the fine $\chi$-expansion
is the sum of the following terms:
\begin{enumerate}
\item
$\displaystyle
2\cdot\frac{1}{4} \cdot\frac{1}{2} 
\sum
\tr \lp 1\rtimes \pi(\chi)\rp(f\times\ve)M_{P_\beta}(s_{2\alpha+\beta}, 0)
$,
where the sum is over the characters $\chi$ of $\idc{E}$ such that $\chi|_{\Af^\times} = \ve\omega$.
\item
$\displaystyle
2\cdot 2\cdot\frac{1}{4} \cdot\frac{1}{2} 
\sum
\tr \lp \mu\rtimes \pi(\theta)\rp(f\times\ve)M_{P_\beta}(s_{2\alpha+\beta}, 0)
$.  The sum is over the nontrivial quadratic characters $\mu \neq \ve$ of $\idc{F}$ and cuspidal
$E_{\mu\ve}$-monomial representations
$\pi(\theta)$ of $\GL(2, \Af)$, where $\theta$ is a character of $\idc{E_{\mu\ve}}$,
such that: $\pi(\theta)$ is not $E_\mu$-monomial and $\theta|_{\Af^\times}\cdot\ve = \omega$.
Let $\sigma'$ be the generator of the $\Gal(E_{\mu\ve}/F)$.
It is known that $\pi(\theta)$ is $E_\mu$-monomial if and only if ${{}^{\sigma'}\!\theta}/{\theta} = \mu\circ\N_{E/F}$
(see \cite{LL}).
\item
$\displaystyle
2\cdot\frac{1}{4} \cdot\frac{1}{2} 
\sum
\tr \lp \mu\rtimes \pi(\theta)\rp(f\times\ve)M_{P_\beta}(s_{2\alpha+\beta}, 0)
$.
The sum is over the nontrivial quadratic characters $\mu \neq \ve$ of $\idc{F}$, cuspidal
$E_{\mu\ve}$-monomial representations
$\pi(\theta)$ of $\GL(2, \Af)$, where $\theta$ is a character of $\idc{E_{\mu\ve}}$,
such that: ${{}^{\sigma'}\!\theta}/{\theta}$ is equal to $\mu$ and $\theta|_{\Af^\times}\cdot\ve = \omega$.
\item
$\displaystyle
2\cdot 2\cdot\frac{1}{4} \cdot\frac{1}{2} 
\sum \tr \lp \ve\rtimes \pi\rp(f\times\ve)M_{P_\beta}(s_{2\alpha+\beta}, 0)
$.
The sum is over all cuspidal non-$E$-monomial, or one dimensional, automorphic representations $\pi$ of $\GL(2, \Af)$
with central character $\omega_{\pi} = \ve\omega$.
\item
$\displaystyle
2\cdot\frac{1}{4} \cdot\frac{1}{2} 
\sum \tr \lp \ve\rtimes \pi\rp(f\times\ve)M_{P_\beta}(s_{2\alpha+\beta}, 0).
$
The sum is over all cuspidal $E$-monomial automorphic representations $\pi$ of \dt $\GL(2, \Af)$ with central character
$\omega_\pi = \ve\omega$.
\end{enumerate}

Suppose an automorphic representation $\pi = \otimes_v \pi_v$ of $\mb{G}(\Af)$ is $\ve$-invariant.
Assuming that the multiplicity one theorem holds for $\GSp(2)$, the intertwining operator
$\rho(\ve)$ maps the space of $\pi$ back to itself.  Consequently, 
$\rho(\ve)$ defines at each place $v$ a local operator $\rho(\ve)_v$ which intertwines $\pi_v$ with $\ve_v\pi_v$.
Hence, for $f \in C(\mb{G}(\Af), \omega)$, the twisted character $\tr \pi(f\times\ve)$ may be expressed as the product
$\prod_v \tr \pi_v(f_v\times\ve_v)$ of twisted local characters.  
Here, $\tr \pi_v(f_v\times\ve_v) := \tr \pi_v(f_v)\rho(\ve)_v$.

We now examine the behavior of the intertwining operator $M_{P_\beta}(s_{2\alpha+\beta}, 0)$.
Let $\pi$ be the induced representation $\mu\rtimes\pi_2$, where $\mu$ is a character of $\idc{F}$, and $\pi_2$ is
a cuspidal or one dimensional automorphic representation of $\GL(2, \Af)$.
At any place $v$ of $F$ for which $\mu_v = 1$ and $\pi_{2,v}$ is square integrable or one dimensional,
the local representation $\pi_v$ is of length two, with constituents
$\pi_v^+, \pi_v^-$ (tempered if $\pi_{2,v}$ is square integrable, nontempered otherwise) (see \cite{ST}).  
\indexs{pi@$\pi^+$, $\pi^-$}%
We normalize intertwining operators so that
the local component $M_{P_\beta}(s_{2\alpha +\beta}, 0)_v$ of $M_{P_\beta}(s_{2\alpha +\beta}, 0)$
acts trivially on $\pi_v^+$, and via scalar multiplication by $-1$ on $\pi_v^-$.  
At a place $v$ where $\pi_{2, v}$ is parabolically induced and/or $\mu_v \neq 1$,  the representation 
$\pi_v$ is irreducible.
For convenience, for such $v$ we put $\pi_v^+ := \pi_v$ and
$\pi_v^- := 0$.  Hence, $M_{P_\beta}(s_{2\alpha + \beta}, 0)_v$ acts trivially
on $\pi_v^+$ and via scalar multiplication by $-1$ on $\pi_v^-$, for all $v$.

Using the above notation, we write (recall that we assume $f = \otimes_v f_v$):
\begin{equation}\label{eq:heisenbergpm}
\tr \pi(f\times\ve)M_{P_\beta}(s_{2\alpha+\beta}, 0)
= \prod_{v \in V} \left[ \tr \pi_v^+(f_v\times\ve_v) - \tr\pi_v^-(f_v\times\ve_v)\right].
\end{equation}

\subsection{DOR for $\mb{H_1}(\Af)$}
The norm mapping $\N_{E/F} : E \rightarrow F$ induces the norm mapping
$\N_{E/F} : \Ae \rightarrow \Af$ on the \adeles.
\indexi{norm mapping}%
Let $\GL(2, \Af)^E$ denote the subgroup of $\GL(2, \Af)$ consisting of elements whose
\indexs{G@$\GL(2,\;)^E$}%
determinants lie in $\N_{E/F}\Ae^\times$.
For any representation $\tau$ of $\GL(2, \Af)$ and character $\chi$ of $\idc{E}$,
let $\tau\otimes_1\chi$ denote the representation of $\mb{H_1}(\Af)$ in the
space of $\tau$ defined as follows:
\[
\tau\otimes_1 \chi : (g, x) \mapsto \chi(x)\lp\tau|_{\GL(2, \Af)^E}\rp(g),\quad \forall (g, x) 
\in \mb{H_1}(\Af).
\]

Since any $(g, x) \in \mb{H_1}(\Af)$ satisfies $\det g = \N_{E/F} x$, 
the representations $\tau\otimes_1 \chi$ and 
$\mu^{-1}\tau\otimes_1 (\mu\circ\N_{E/F})\chi$ are equal to each other for any
quasicharacter $\mu$ of $\idc{F}$.  Moreover,
since the representation $\tau\otimes_1\chi$ is defined by $\chi$ and the restriction of $\tau$ to 
the subgroup $\GL(2, \Af)^E$ of $\GL(2, \Af)$, it may be reducible even if $\tau$
is an irreducible representation of $\GL(2, \Af)$ (see \cite{LL}).  


\subsubsection{Discrete Spectrum Representations}
Let $(\mb{H_1}, \mb{H_1}, \pi, 1)$ be a quadruple for
the group $\mb{H_1}$. 
Then, $\pi$ is an irreducible, discrete spectrum, automorphic representation of $\mb{H_1}(\Af)$, and  
it is an irreducible constituent of 
$\tau\otimes_1 \eta$ for some irreducible,  discrete spectrum, automorphic representation $\tau$ of
$\GL(2, \Af)$ and character $\mu$ of $\idc{E}$.  



For any place $v$ of $F$ which remains prime in $E$, let 
\[\GL(2, F_v)^E = \{g \in \GL(2, F_v) : \det g \in \N_{E/F}E_v^\times\}.
\indexs{G@$\GL(2,\;)^E$}\]%
If $v$ splits in $E$, let $\GL(2, F_v)^E = \GL(2, F_v)$.  
The ad\`elic group
$\GL(2,\Af)^E$ is equal to the restricted tensor product
\indexi{restricted tensor product}%
$\otimes_{v \in V} \GL(2, F_v)^E$ (namely, almost all local components lie in $\GL(2, \mc{O}_v) \cap \GL(2, F_v)^E$), 
and $\tau|_{\GL(2, \Af)^E} = \otimes_v \tau_v|_{\GL(2, F_v)^E}$.

Let $v$ be any place of $F$.  From \cite{LL}, the restriction of $\tau_v$ to $\GL(2, F_v)^E$
is reducible if and only if $v$ is prime in $E$, 
and $\tau_v$ is the $E_v$-monomial representation associated with some character $\theta_v$ of $E_v^\times$.
In this case, the representation $\tau_v|_{\GL(2, F_v)^E}$
has length two, with constituents $\tau_v^+$ and $\tau_v^-$.
For convenience, if $\tau_v$ is irreducible, we let $\tau_v^+ = \tau_v$ and $\tau_v^- = 0$.
If $\tau_v|_{\GL(2, F_v)^E}$ is unramified, we let $\tau_v^+$  be the
constituent which contains the (unique up to scalar multiple) $\GL(2, \mc{O}_v)$-fixed vector.
If $\tau_v|_{\GL(2, F_v)^E}$ has length two and is ramified, then which constituent to label
$\pi_v^+$ and which to label $\pi_v^-$ is irrelevant for our purpose.

We let $\pi_v^+ = \tau_v^+\otimes_1 \eta_v$ and $\pi_v^- = \tau_v^-\otimes_1 \eta_v$ for any place $v$ of $F$.
Let $\tau_v\otimes_1 \eta_v$ denote the {\bf local packet} $\{\pi_v^+, \pi_v^-\}$.
\indexi{packet!local}%
We define the {\bf global packet} $\tau\otimes_1 \eta$ to be the \emph{restricted tensor product}
\indexi{packet!global}\indexi{restricted tensor product}%
\[
\otimes_{v \in V} \{\tau_v\otimes_1\eta_v\} := \{\otimes_v \pi_v' : \pi_v' \in \{\pi_v^+, \pi_v^-\},\;
\pi_v' = \pi_v^+ \text{ almost all } v\}.
\]
Recall our assumption that the global test functions are tensor products of local components.
For $f_1 = \otimes_{v \in V}f_{1,v}$ in $C(\mb{H_1}(\Af), \omega)$, put
\[
\tr \lp\tau\otimes_1 \eta\rp(f_1) :=
\prod_{v \in V}\bigg[\tr \pi_v^+(f_v) + \pi_v^-(f_v) \bigg].
\]
The distribution $f_1 \mapsto \tr \lp\tau\otimes_1 \eta\rp(f_1)$ defines a
distribution on $C(\mb{H_1}(\Af), \omega)$ which is invariant under stable conjugacy.

If $\tau$ is not an $E$-monomial representation, then by \cite{LL} each member
of $\tau\otimes_1 \eta$ occurs with multiplicity one in the discrete spectrum of $\mb{H_1}(\Af)$.
Thus,
the contribution of the packet $\tau\otimes_1 \eta$ to the spectral side of the trace formula of 
$\mb{H_1}(\Af)$ is
\[
\tr \lp\tau\otimes_1 \eta\rp(f_1).
\]


If $\tau$ is cuspidal $E$-monomial,
the multiplicity $m(\pi')$ with which each member $\pi'$ of 
\indexi{multiplicity formula}%
$\{\pi\}$ occurs in the discrete spectrum of $\mb{H_1}(\Af)$ is given by the following formula:
\[
m(\pi') = \frac{1}{2}\lp1 + (-1)^{n(\pi')}\rp,
\]
where $n(\pi')$ is the number of places $v$ for which $\pi_v' = \pi_v^-$ (\cite{LL}).
Taking into account the nontrivial multiplicity formula,
the contribution of $\tau\otimes_1\eta$ to the spectral side of the trace formula of $\mb{H_1}(\Af)$
is equal to
$D(\tau\otimes_1 \eta, f_1) :=$
\indexs{D@$D(\;,\;)$}%
\[
\frac{1}{2}\prod_{v \in V} \Big[ \tr \pi_v^+(f_{1,v}) + \tr \pi_v^-(f_{1,v})\Big] + 
\frac{1}{2}\prod_{v \in V} \Big[ \tr \pi_v^+(f_{1,v}) - \tr \pi_v^-(f_{1,v})\Big].
\]
Note that the first term in the above sum is equal to $\frac{1}{2}\cdot\tr\lp\tau\otimes_1\eta\rp(f_1)$.
If $v$ is prime in $E$,
there exists an $h \in \mb{H_1}(\ol{F_v})$ such that conjugation by $h$ swaps $\pi_v^+$ with $\pi_v^-$.
Hence, $\frac{1}{2}\cdot\tr\lp\tau\otimes_1\eta\rp(f_1)$ defines a {stable} distribution on 
$C(\mb{H_1}(\Af), \omega)$, whereas the second term in the sum defines an {unstable} distribution. 
We call $\frac{1}{2}\cdot\tr\lp\tau\otimes_1\eta\rp(f_1)$ the {\bf stable part} of $D(\tau\otimes_1\eta, f_1)$.
\indexi{trace formula!spectral side!stable part}%



\subsubsection{Induced Representations}
There is one proper Levi subgroup in $\mb{H_1}$; namely, the diagonal torus
\[
\mb{M_1} = \{ (\diag(a, b), c) : a, b \in \mbb{G}_m, c \in \rR_{E/F} \mbb{G}_m, ab = \N_{E/F} c\}.
\]
Suppose $(\mb{M_1}, \mb{H_1}, \tau, s)$ is a quadruple in the fine $\chi$-expansion of
$\mb{H_1}(\Af)$.
Any discrete spectrum representation of $\mb{M_1}(\Af)$ has the form:
\begin{multline*}
\tau = \mu_1\otimes\mu_2\otimes_1 \eta:
(\diag(a, b), c) \mapsto \mu_1(a)\mu_2(b)\eta(c),\\\forall\;(\diag(a, b), c) \in \mb{M_1}(\Af),
\end{multline*}
where $\mu_1, \mu_2$ are characters of $\idc{F}$ and $\eta$ is a character of $\idc{E}$.
From the way $\mb{M_1}$ is defined, for any quasi-character $\mu$ of $\idc{F}$, we have:
\begin{equation}\label{eq:equivcharsH1}
\mu_1\otimes \mu_2 \otimes_1 \eta = \mu\mu_1\otimes\mu\mu_2\otimes_1 \lp\mu^{-1}\circ\N_{E/F}\rp\eta.
\end{equation}

The $F$-split component $\mb{A}_{M_1}$ of $\mb{M}_1$ is 
$\{(\diag(a, b), c): a, b, c \in \mbb{G}_m, ab = c^2\}$, and
$\mf{a}_{M_1} = \{ (x, y; t) : x, y, t \in \mbb{R}, x + y = 2t\} \cong \mbb{R}^2$.
The group $W_{H_1}(M_1)$ is equal to $\ZZ/2\ZZ$, generated by
$(12) : (\diag(a, b), c) \mapsto (\diag(b, a), c)$.
The $F$-split component of the center of $\mb{M_1}$ is 
$\mb{A}_{H_1} = \{(\diag(c, c), \pm c) : c\in \mbb{G}_m)\}$, and 
$\mf{a}_{H_1} = \{(z, z; z) : z \in \mbb{R}\}$.

If $\mf{a}_{M_1}^s = \mf{a}_{H_1}$ for some $s \in W_{H_1}(M_1)$, then
$s$ must be equal to $(12)$.  Consequently,
$\tau = \mu_1 \otimes \mu_2 \otimes_1 \eta$ is equivalent to $\mu_2 \otimes \mu_1 \otimes_1 \eta$.
In other words,
\[
\mu_1(a)\mu_2(b)\eta(c)
= \mu_2(a) \mu_1 (b)\eta(c), \quad \forall (\diag(a, b), c) \in \mb{M_1}(\Af).
\]
Since $ab = \N_{E/F} c$, the above equality may be rewritten as
\[
\frac{\mu_1}{\mu_2}(a)\lp\lp\mu_2\circ\N_{E/F}\rp\cdot\eta\rp(c)
= \frac{\mu_2}{\mu_1}(a)\lp\lp\mu_1\circ\N_{E/F}\rp\cdot\eta\rp(c),\;\; \forall a \in \Af^\times,\;
c \in \Ae^\times.
\]
Thus, ${\mu_1}/{\mu_2}$ is equal to either $1$ or $\ve$.
By \eqref{eq:equivcharsH1}, we conclude that $\tau$ must have one of the following forms:
\begin{itemize}
\item
$1\otimes 1 \otimes_1 \eta$,
\item
$1\otimes \ve \otimes_1 \eta$.
\end{itemize}

We have: $\abs{W(M_1, M_1)} = 1, \abs{W(M_1, H_1)} = 2$, and 
$\abs{\det (1 -(12))|_{\mf{a}_{M_1}/\mf{a}_{H_1}}} $ is equal to $2$.

Let $\mb{P_1}$ be the upper triangular parabolic subgroup of $\mb{H_1}$ whose
Levi component is $\mb{M_1}$.
By \eqref{eq:finechidiscrete}, with $\ve = 1$,
the contribution to the fine $\chi$-expansion of $\mb{H_1}(\Af)$ from $\mb{M_1}$ is
\begin{multline*}
\frac{1}{2}\cdot\frac{1}{2} \sum_{\eta|_{\Af^\times} = \omega} 
\bigg[
\tr \lp I(1, 1)\otimes_1\eta\rp(f_1)M_{P_1}((12), 0)\\ +
\tr \lp I(1, \ve)\otimes_1\eta\rp(f_1)M_{P_1}((12), 0)\bigg].
\end{multline*}

The representation $I(1, 1)\otimes_1 \eta$ is irreducible, and we normalize 
the intertwining operator $M_{P_1}((12), 0)$ to be $1$.  
For any place $v$, the distribution on the space
$C(\mb{H_1}(F_v),\omega_v)$ defined by
$f_{1,v} \mapsto \tr\lp I(1, 1)\otimes_1 \eta\rp_v(f_{1, v})$ is stable.

In the case of $I(1, \ve)\otimes_1 \eta$,
the restriction of $I(1,\ve)_v$ to $\GL(2, F_v)^E$
is irreducible if the place $v$ splits in $E$, and it is reducible of length two, with constituents
$I_v^+, I_v^-$, if $v$ is prime in $E$.
\indexs{I@$I^+$, $I^-$}%
We let $I_v^+$ denote the unramified constituent of
$I(1, \ve)_v$ if $I(1, \ve)_v$ is unramified, and we let $I_v^+ = I(1,\ve)_v$, $I_v^- = 0$ if
$v$ splits in $E$.  Write the intertwining operator $M_{P_1}((12), 0)$ for $I(1, \ve)\otimes_1 \eta$ 
as the tensor product  $\otimes_{v \in V}M_{P_1}((12), 0)_v$ of local intertwining operators.
\indexi{intertwining operator}%
If $v$ splits in $E$, we normalize $M_{P_1}((12), 0)_v$ to be $1$.
If $v$ is prime in $E$, we may normalize $M_{P_1}((12), 0)_v$ so that it acts trivially
on $I_v^+$ and it acts via multiplication with $-1$ on $I_v^-$ (see \cite[Sect. 5]{LL}).  
Hence, for $f_1 \in C(\mb{H_1}(\Af), \omega)$, we have
\[
\tr \lp I(1, \ve)\otimes_1\eta\rp(f_1)M_{P_1}((12), 0)
= \prod_{v \in V} \lp \tr I_v^+(f_{1,v}) - \tr I_v^-(f_{1,v})\rp.
\]
If $v$ is prime in $E$,
there exists an element of $h \in \mb{H_1}(\ol{F_v})$ such that the action of conjugation by $h$ swaps
$I_v^+$ and $I_v^-$ (\cite{LL}), multiplying the trace by a factor of $-1$.
The distribution 
$f_{1, v} \mapsto \tr \lp I(1, \ve)\otimes_1\eta\rp_v(f_{1, v})M_{P_1}((12, 0)_v$ is therefore \emph{unstable}.

\subsubsection{Stable Spectrum}\label{sec:stablespectrum}
Recall that $I^d(H_1, f_1)$ is the discrete part of the fine $\chi$-expansion for $\mb{H_1}(\Af)$.
We let $SI^d(H_1, f_1)$ denote the stable part of $I^d(H_1, f_1)$; namely, $SI^d(H_1, f_1)$ is the sum of the following:
\begin{itemize}
\item
$\displaystyle \sum \tr \pi(f_1)$.  

The sum is over representations $\pi$ of the form $\tau\otimes_1 \eta$, where  
$\tau$ is a cuspidal non-$E$-monomial or one dimensional representation of
$\GL(2, \Af)$, and $\eta$ is a character of $\idc{E}$, such that $\eta|_{\Af^\times}\cdot\omega_\tau = \omega$.
\item
$\displaystyle \frac{1}{2}\sum \tr \pi(f_1)$.  

The sum is over representations $\pi$ of the form
$\pi = \pi(\theta)\otimes_1 \eta$, where $\theta, \eta$ are characters of $\idc{E}$, $\theta \neq {}^\sigma \theta$,
$\lp\theta\eta\rp|_{\Af^\times}\cdot\ve = \omega$,
and $\pi(\theta)$ is the cuspidal, $E$-monomial automorphic representation of $\GL(2, \Af)$ associated with $\theta$.
\item
$\displaystyle \frac{1}{4}\sum \tr\lp I(1, 1)\otimes_1 \eta\rp(f_1)$.  

The sum is over characters $\eta$ of $\idc{E}$ such that $\eta|_{\Af^\times} = \omega$.
\end{itemize}

\subsection{DOR for $\mb{H_2}(\Af)$}
For any representation $\pi$ of $\GL(2, \Ae)$ and character $\mu$ of $\idc{F}$ such that
$\omega_{\pi} = \mu\circ\N_{E/F}$, let $\pi\otimes_2 \mu$ denote the representation of
$\mb{H_2}(\Af)$ defined as follows:
\[
\pi\otimes_2 \mu: (g, c) \mapsto \mu(c)\pi(g),\quad \forall (g, c) \in \mb{H_2}(\Af).
\]

There is one proper Levi subgroup in $\mb{H_2}$; namely, the diagonal torus whose
group of $F$-points is:
\[
\mb{M_2}(F) = \{ (\diag(a, b), c)_* : a, b \in E^\times, c \in F^\times\},
\]
where lower star denotes the equivalence class of the element in $\GL(2, E)\times F^\times$ modulo
the subgroup
$\left\{\lp\diag(z, z), \N_{E/F}z^{-1}\rp : z \in E^\times\right\}$.

Suppose $(\mb{M_2}, \mb{H_2}, \tau, s)$ is a quadruple in the fine $\chi$-expansion of
$\mb{H_2}$.
Any discrete spectrum representation of $\mb{M_2}(\Af)$ has the form:
\[
\tau = \eta_1\otimes\eta_2\otimes_2 \mu:
(\diag(a, b), c)_* \mapsto \eta_1(a)\eta_2(b)\mu(c),\;\forall\;(\diag(a, b), c)_*\;,
\]
where $\eta_1, \eta_2$ are characters of $\idc{E}$ and $\mu$ is a character of $\idc{F}$.
Because of the equivalence relation on $\mb{H_2}(\Af)$,
the characters must satisfy $\eta_1\eta_2 = \mu\circ\N_{E/F}$.

The $F$-points of the split component $\mb{A}_{M_2}$ of $\mb{M}_2$ is 
$\{(\diag(a, b), c)_*: a, b, c \in F^\times\}$, and
\[
\mf{a}_{M_2} = \{ (x, y; t) : x, y, t \in \mbb{R}\}/\{ (z, z, -2z): z \in \mbb{R}\}\cong \mbb{R}^2.
\]
The group $W_{H_2}(M_2) = W(H_2)$ is equal to $\ZZ/2\ZZ$, generated by the permutation
$(12) : (\diag(a, b), c)_* \mapsto (\diag(b, a), c)_*$.
The $F$-points of the split component of the center of $\mb{H_2}$ is 
$\mb{A}_{H_2}(F) = \{(\diag(1, 1), c)_* : c\in F^\times)\}$, and 
$\mf{a}_{H_2} = \{(0, 0; z) : z \in \mbb{R}\}$.

If $\mf{a}_{M_2}^s = \mf{a}_{H_2}$ for some $s \in W_{H_2}(M_2)$, then
$s$ must be $(12)$.  Consequently,
\[
\tau = \chi \otimes \chi \otimes_1 \omega
\]
for some character $\chi$ of $\idc{E}$ such that $\chi^2 = \omega\circ\N_{E/F}$.

We have: $\abs{W(M_2, M_2)} = 1, \abs{W(M_2, H_2)} = 2$, and 
$\abs{\det (1 -(12))|_{\mf{a}_{M_2}/\mf{a}_{H_2}}}$ is equal to $2$.

Let $\mb{P_2}$ be the upper triangular parabolic subgroup of $\mb{H_2}$ whose
Levi component is $\mb{M_2}$.
Then, the contribution to the fine $\chi$-expansion of $\mb{H_2}(\Af)$ from $\mb{M_2}$ is:
\[
\frac{1}{2}\cdot\frac{1}{2} \sum_{\chi^2 = \omega\circ\N_{E/F}} 
\tr \lp I(\chi, \chi)\otimes_2\omega\rp(f_2)M_{P_2}((12), 0),
\]
where $f_2$ is a test function in $C(\mb{H_2}(\Af), \omega)$.
Since $I(\chi, \chi)\otimes_2 \omega$ is irreducible, we normalize $M_{P_2}((12), 0)$ to be $1$.
\subsection{Summary}\label{sec:tfsummary}
In summary, $I^d(G, f,\ve)$ is equal to the sum of following:
\begin{itemize}
\item
$\displaystyle \sum_{\pi \in L_d(\mb{G}(\Af), \omega)} \tr \pi(f \times \ve).$
\item
$\displaystyle
\frac{1}{4}
\sum_{\mu^2 = \omega} \tr \lp \ve\times\ve\rtimes \mu\rp(f\times\ve)$.
\item
$\displaystyle
\frac{1}{4} 
\sum_{{\ve'}^2 = 1,\;\ve'\mu^2 = \omega} \tr \lp \ve'\times\ve'\ve\rtimes \mu\rp(f\times\ve)
$
\item
$\displaystyle \frac{1}{4}
\sum_{\chi|_{\Af^\times} = 1,\;\ve\mu^2 =\omega} 
\tr \lp \pi(\chi)\rtimes \mu\rp(f\times\ve)M_{P_\alpha}(s_{\alpha+\beta}, 0).$
\item
$\displaystyle
\frac{1}{4} 
\sum_{\chi \in \widehat{\idc{E}};\; \chi|_{\Af^\times}\cdot\ve = \omega} 
\tr \lp 1\rtimes \pi(\chi)\rp(f\times\ve)M_{P_\beta}(s_{s\alpha+\beta}, 0).
$
\item
$\displaystyle
\frac{1}{2} 
\sum
\tr \lp \mu\rtimes \pi(\theta)\rp(f\times\ve)M_{P_\beta}(s_{2\alpha+\beta}, 0).
$

The sum is over the nontrivial quadratic characters $\mu$ of $\idc{F}$ and
cuspidal monomial representations
$\pi(\theta)$ of $\GL(2, \Af)$, where $\theta$ is a character of $\idc{E_{\mu\ve}}$,
such that: $\pi(\theta)$ is not $E_\mu$-monomial and $\theta|_{\Af^\times}\cdot\ve = \omega$.
\item
$\displaystyle
\frac{1}{4}
\sum
\tr \lp \mu\rtimes \pi(\theta)\rp(f\times\ve)M_{P_\beta}(s_{2\alpha+\beta}, 0).
$

The sum is over the nontrivial quadratic characters $\mu$ of $\idc{F}$ and cuspidal monomial representations
$\pi(\theta)$ of $\GL(2, \Af)$, where $\theta$ is a character of $\idc{E_{\mu\ve}}$,
such that $\theta|_{\Af^\times}\cdot\ve = \omega$ and $\pi(\theta)$ is $E$-monomial.
In particular, we have
$\frac{{}^{\sigma'}\theta}{\theta} = \mu$, where $\sigma'$ is the generator of $\Gal(E_{\mu\ve}/F)$.
\item
$\displaystyle
\frac{1}{2} 
\sum \tr \lp \ve\rtimes \pi\rp(f\times\ve)M_{P_\beta}(s_{2\alpha+\beta}, 0).
$

The sum is over all cuspidal non-$E$-monomial, or one dimensional, automorphic representations $\pi$ of $\GL(2, \Af)$
such that $\ve \omega_{\pi} = \omega$.
\item
$\displaystyle
\frac{1}{4}
\sum \tr \lp \ve\rtimes \pi\rp(f\times\ve)M_{P_\beta}(s_{2\alpha+\beta}, 0).
$

The sum is over all cuspidal $E$-monomial automorphic representations $\pi$ of \dt $\GL(2, \Af)$ such that
$\ve\omega_\pi = \omega$.
\end{itemize}
The stable spectral expansion $SI^d(H_1, f_1)$ is equal to:
\begin{multline*}
\sum_{\substack{\tau \in L_d(\GL(2, \Af)),\\ \tau\neq \ve\tau,\;\eta|_{\Af^\times}\cdot\omega_\tau = \omega}}
\tr \lp \tau\otimes_1 \eta\rp(f_1) 
+ \frac{1}{2}\sum_{\substack{\theta\neq {}^\sigma \theta,\\\lp\theta\eta\rp|_{\Af^\times}\cdot\ve = \omega}}
\tr \lp\pi(\theta)\otimes_1 \eta\rp(f_1) \\
+ \frac{1}{4} \sum_{\eta|_{\Af^\times} = \omega}
\tr \lp I(1, 1)\otimes_1 \eta \rp(f_1),
\end{multline*}
and $I^d(H_2, f_2)$ is equal to:
\begin{equation*}
\sum_{\pi \in L_d(\mb{H_2}(\Af), \omega)} \tr \pi(f_2) +
\frac{1}{4}\sum_{\chi^2 = \omega\circ\N_{E/F}} \tr \lp I(\chi, \chi)\otimes_2 \omega\rp(f_2).
\end{equation*}
\indexi{fine $\chi$-expansion|)}%
\section{Comparison of Geometric Sides of Trace Formulas}
\subsection{Weighted Orbital Integrals}\label{sec:simplifyO}
Let $\mb{H}$ be a reductive $F$-group.  Let $\mb{Z}_0$ be the split component of 
the center $\mb{Z}$ of $\mb{H}$. Let $\ve$ be a character of $\mb{H}$ trivial on $\mb{Z}$.
Fix a minimal parabolic subgroup $\mb{P}_0$ of $\mb{H}$.  
Let $\mb{A}_0$ be the split component of the Levi subgroup of $\mb{P}_0$.  
Put $\mf{a}_0 := X_*(\mb{A}_0)\otimes_{\mbb{Z}} \mbb{R}$.

We say that a regular element in $\mb{H}(F)$ is {\bf parabolic} if it is not elliptic; 
\indexi{element!parabolic}%
in other words, its centralizer in $\mb{H}(F)$ is a maximal torus which is contained in a 
proper ($F$-)parabolic subgroup.
The terms in the geometric side of the ($\ve$-twisted) trace formula associated with
the semisimple conjugacy classes of parabolic elements are  
$\ve$-twisted-{\bf weighted orbital integrals}, defined as follows:
\indexi{orbital integral!weighted}%

Suppose $\gamma$ is a parabolic element in $\mb{H}(F)$.  
Let $\mb{M}$ be the minimal Levi subgroup of $\mb{H}$ containing $\gamma$.
Let $\mb{H}(\Af)_\gamma$ be the centralizer of $\gamma$ in $\mb{H}(\Af)$.
The weighted orbital integral associated with the semisimple conjugacy class of $\gamma$ 
has the following form (\cite[Lecture 5]{CLL}):
\[
O^T(\gamma, f) = 
c(\gamma) \int_{\mb{H}(\Af)_\gamma\bs\mb{H}(\Af)}f(h^{-1}\gamma h)\ve(h)v(h, T)\;dh.
\]
Here, $c(\gamma)$ is some constant dependent on $\gamma$, $T$ is the truncation parameter in $\mf{a}_0$,
and $v(h, T)$ is the volume of a certain convex hull in $\mf{a}_H\bs\mf{a}_M$ determined by
\indexs{v@$v(\;,\;)$}%
the projections of $T$ and $\rH_{P_0}(h)$ onto $\mf{a}_H\bs\mf{a}_M$.

Let $n$ be the real dimension of $\mf{a}_H\bs\mf{a}_M$.
For each $T$, the volume function $v(h, T)$ is a degree $n$ polynomial in the coordinates of $\rH_{P_0}(h)$.
Let $\{h_1, h_2, \dots\}$ be the eigenvalues of $h$ .
The coordinates of $\rH_{P_0}(h) \in \mf{a}_M$ are expressed in terms of $\log \abs{h_i'}$,
where $h_i' = \N_{L_i/F}h_i$, $L_i$ is the smallest extension field containing $h_i$, and
$\abs{\cdot}$ is the absolute value function on $\Af$.
For any place $v$ of $F$, let $\abs{\cdot}_v$ be the $p$-adic absolute value on $F_v$.
Then, $\abs{h_i} = \prod_{v} \abs{h_i}_v$,  and 
$\log \abs{h_i'} = \sum_{v}\log \abs{h_i'}_v$.
We rewrite $O^T(\gamma, f)$ as follows:
\begin{multline}
c(\gamma) \sum_{S}\Bigg[\int_{H_{S,\gamma}\bs H_{S}}
f_{S}(h_{S}^{-1}\gamma h_{S})\ve_{S}(h_{S})v(h, T, S)\;dh_{S}\\
\times \prod_{w \notin S}\int_{H_{w,\gamma}\bs H_{w}} f_{w}(h_{w}^{-1}\gamma h_{w})\ve_w(h_w)\;dh_{w}\Bigg].
\end{multline}
Here, the sum is over the sets $S$ of $n$ distinct places of $F$.
The subscript $S$ in $H_S, f_S, \ve_S$ signifies the product of local components 
\indexs{H@$H_S$}\indexs{f@$f_S$}\indexs{epsilon@$\ve_S$}%
$H_v, f_v, \ve_v$ over the places $v$ in $S$.  
The term $v(h, T, S)$ is a polynomial in the local components of $v(h, T)$ at the places in $S$.
Observe that the product over $w \notin S$ consists of non-weighted local $\ve_w$-twisted orbital integrals.

We say that a local function $f_v$ on $\mb{H}(F_v)$ is {\bf elliptic} if the orbital integral
\indexi{function!elliptic}%
of $f_v$ vanishes at the parabolic elements in $\mb{H}(F_v)$.  Matrix coefficients of cuspidal representations
\indexi{matrix coefficient}%
of $\mb{H}(F_v)$ are examples of such functions.
By choosing test functions $f$ whose local components at $n + 1$ distinct places are elliptic, 
the weighted orbital integral $O^T(\gamma, f)$ vanishes.

We say that an element in $\mb{H}(F)$ is {\bf singular} if it has two equal eigenvalues.
\indexi{element!singular}%
From \cite{AWO}, the weighted orbital integral at a singular element may be expressed
as a limit of a linear combination of weighted orbital integrals at parabolic elements.  Hence,
by choosing $f$ with $n + 1$ elliptic components with $n$ sufficiently large, 
the weighted orbital integrals at singular elements also vanish.  
Here, 
$n$ is sufficiently large if it is greater than or equal to the semisimple rank of $\mb{H}$.

Consider the case of $\mb{G} = \GSp(2)$.  Suppose $\gamma$ is a parabolic element in $\mb{G}(F)$ for which the 
$\ve$-twisted orbital integral is nonzero.   Then, the centralizer of $\gamma$ in $\mb{G}(\Af)$ must
necessarily lie in the kernel of $\ve$.
The centralizer of $\gamma$ in $\mb{G}(\Af)$ is a maximal non-elliptic torus,
The following is a list of representatives of the conjugacy classes of the maximal non-elliptic
tori in $\mb{G}(\Af)$:
\begin{enumerate}
\item
$
\mb{T}(\Af)
= \{\diag(a, b, \lambda/b, \lambda/a) : a, b, \lambda \in \Af^\times\};
$
\item
$
\mb{T}(\Af)
= \left\{
\lp\lsm t_2 &\\ &\frac{\lambda}{\det t_2} e t_2 e \rsm\rp : t_2 \in \mb{T_2}(\Af), 
\lambda \in \Af^\times
\right\},
$

where $\mb{T_2}$ is a maximal elliptic torus in $\GL(2)$ and $e = \diag(1, -1)$;
\item
$
\mb{T}(\Af) =
\left\{\lp\lsm a &&\\ &t_2& \\&&b \rsm\rp : a, b \in \Af^\times, t_2 \in \mb{T_2}(\Af),\;ab = \det t_2\right\},
$

where $\mb{T_2}$ is a maximal elliptic torus in $\GL(2)$.
\end{enumerate}
In the first two cases, the similitude factors of the elements in the torus range over all of
$\Af^\times$; hence these tori do not lie in the kernel of $\ve$.  In the third case,
the similitude factor is $\det t_2$.  An elliptic maximal torus in
$\GL(2, \Af)$ which lies in the kernel of $\ve\circ \det$ is conjugate to
$\left\{\lp\lsm c & d A\\d&c\rsm\rp : c, d \in \Af, (c, d) \neq (0, 0) \right\}$,
where $A$ is an element in $F^\times - {F^\times}^2$ such that $E = F(\sqrt{A})$.
Consequently, if the centralizer of $\gamma$ must be conjugate to
\[
\mb{T}(\Af) =
\left\{\lp\lmx a &&\\&\lsm c&d A\\d&a\rsm&\\&&b \rmx\rp : 
a, b \in \Af^\times,\; c, d \in \Af,\; ab = c^2 - d^2 A\right\}.
\]
Since the semisimple rank of
$\mb{T}$ is $1$, using test functions with two elliptic components suffices to eliminate the weighted
orbital integrals from the geometric side of the $\ve$-twisted trace formula for $\mb{G}$.

As for the $\ve$-endoscopic groups $\mb{H_1}$, $\mb{H_2}$, their semisimple $F$-ranks are equal to $1$.  Hence,
using test functions with two elliptic components suffices to eliminate the weighted orbital integrals.

Fix two places $u_1, u_2$ of $F$ which are prime in $E$.  Let $\EL{\mb{G}(\Af), \omega}$ denote the
\indexs{E@$\EL{\;,\;}$}%
space of functions in $C(\mb{G}(\Af), \omega)$ whose local components at $u_1, u_2$ are elliptic.
For $i = 1, 2$, define $\EL{\mb{H}_i(\Af), \omega}$ likewise for $\mb{H}_i$.
Since elliptic regular elements in $\mb{H}_i(F_v)$ are norms of elliptic regular elements 
in $\mb{G}(F_v)$ for any place $v$ which is prime in $E$, it is reasonable to expect that,
for any function $f$ in $\EL{\mb{G}(\Af), \omega}$, 
there exists matching functions $f_i \in \EL{\mb{H}_i(\Af), \omega}$.

Let $\mb{H} = \mb{G}, \mb{H}_1$ or $\mb{H_2}$.
Let $\{\mco\}_e$ denote the set of $F$-conjugacy classes of elliptic regular elements in  $\mb{H}(F)$.
\indexs{O@$\{\mco\}_e$}%
For each $\mc{O} \in \{\mc{O}\}_e$, put
\[
J_{\mco}(f) := O_\gamma(f) = \int_{\mb{H}(\Af)_\gamma\bs\mb{H}(\Af)}f(h^{-1}\gamma h)\ve_H(h)\;dh/du,
\indexs{J@$J_{\mco}(\;)$}%
\]
where $\gamma$ is a representative in $\mb{H}(F)$ of the conjugacy class $\mc{O}$, and $\ve_H$ is $\ve$ if $\mb{H} = \mb{G}$,
$1$ if $\mb{H} = \mb{H}_i$ ($i = 1,2$).
Here, $dh$, $du$ are Tamagawa measures on $\mb{H}(\Af)_\gamma$ and $\mb{H}(\Af)$, respectively.
\begin{prop}\label{corollary:ellipticOexp}
Let $f$ be a function in $\EL{\mb{H}(\Af), \omega}$.  
Then, the geometric side of the $\ve_H$-twisted trace formula for $\mb{H}$ is
equal to
\[
\sum_{\mco \in \{\mco\}_e}J_{\mco}(f).
\]
\end{prop}
\begin{proof}
This follows from \cite[Lecture 5]{CLL} and the previous discussion.
\end{proof}
\subsection{Kottwitz-Shelstad's Formula}
\indexi{Kottwitz-Shelstad's trace formula|(}%
Let $i = 1$ or $2$.  For any strongly $\mb{G}$-regular element $\delta \in \mb{H}_i(F)$
and test function $f_i \in C(\mb{H}_i(\Af), \omega)$,
let
\[
SO_{\delta}(f_i) = \sum_{\{\delta'\}} O_{\delta'}(f_i),
\indexi{orbital integral!stable}%
\]
where the sum is over a set of representatives for the $\mb{H}_i(\Af)$-conjugacy
classes of elements $\delta' \in \mb{H}_i(\Af)$ in the $\mb{H}_i(\ol{\Af})$-conjugacy class
of $\delta$.

If we have an equation relating the geometric sides of the trace formulas for $\mb{G}$ and
\indexi{trace formula!geometric side}%
its $\ve$-endoscopic groups $\mb{H_1}$, $\mb{H_2}$,
\indexi{endoscopic!group}%
then we may deduce a relation of the spectral sides.
\indexi{trace formula!spectral side}%
The tool with which we use to relate the geometric sides of the trace formulas is
the following formula of Kottwitz-Shelstad's in \cite[Sect. 7.4]{KS} (applied to our situation):
\begin{equation}\label{eq:KStrace}
T_e(G, f) = \iH{1} ST_e^{**}(H_1, f_1) + \iH{2} ST_e^{**}(H_2, f_2),
\end{equation}
where:
\begin{itemize}
\item
$f$ and $f_i$ ($i = 1, 2$) are matching functions in $C(\mb{G}(\Af), \omega)$ and
$C(\mb{H}_i(\Af), \omega)$, 
respectively;
\item
$T_e(G, f) = \sum_{\mc{O} \in \{\mc{O}\}_e}J_{\mc{O}}(f)$;
\indexs{T@$T_e(\;,\;)$}%
\item
$ST_e^{**}(H_i, f_i) := 
a_{\mb{H}_i} \sum_{\{\delta_{i}\}}SO_{\delta_{i}}(f_i),$ where
\indexs{S@$ST_e^{**}(\;,\;)$}%
the sum is over a set of representatives for the $\mb{H}_i(\bar{F})$-conjugacy classes
of elliptic
$\mb{G}$-regular elements $\delta_{i} \in \mb{H}_i(F)$;
\item
For any algebraic $F$-group $\mb{H}$,
\[
a_{\mb{H}} := \abs{\pi_0\lp Z(\hat{H})^\Gamma\rp}\cdot \abs{\text{ker}^1(F, Z(\hat{H}))}^{-1};
\indexs{a@$a_{\mb{H}}$}%
\]
\item
For $i = 1, 2$,
\[
\iH{i} :=
a_{\mb{G}} \abs{\text{Out}(\mb{H}_i, \mathcal{H}_i,s_i, \xi_i)}^{-1}\cdot
a_{\mb{H}_i}^{-1},
\]
where $\text{Out}(\mb{H}_i, \mathcal{H}_i,s_i, \zeta_i)$ is the group of outer automorphisms of
\indexs{O@${\rm Out}(\;,\;,\;,\;)$}%
$\mb{H}_i$ determined by automorphisms of the endoscopic data $(\mb{H}_i, \mc{H}_i, s_i, \xi_i)$
(see \cite[Sect. 2.1]{KS}).
\end{itemize}

The group $Z(\hat{G}) = \mathbb{C}^\times$ is connected and 
the action of $\Gal(\bar{F}/F)$ on $Z(\hat{G})$ is trivial; hence,
$a_{\mb{G}} = 1$.  
For both endoscopic groups, $\abs{\pi_0\lp Z(\hat{H_i})^\Gamma\rp} = 1$;
hence,
$
a_{\mb{H}_i} = \abs{\text{ker}^1(F, Z(\hat{H_i}))}^{-1}.
$
\begin{claim}\label{claim:tateshaf}
For $i = 1, 2$, $\ker^1( F, Z(\hat{H_i}))$ is trivial.
\end{claim}
\begin{proof}

By definition, for $i = 1, 2$, 
$\ker^1( F, Z(\hat{H_i}))$ is the kernel of the natural embedding
\[
H^1(\Gamma_F, Z(\hat{H_i})) \rightarrow \prod_v H^1\lp\Gamma_{F_v}, Z(\hat{H_i})\rp,
\]
where $\Gamma_F = \Gal(\bar{F}/F)$ and
$\Gamma_{F_v}$ is the decomposition group of $v$ in $\Gal(\bar{F}/F)$.
Observe that the action of $\Gamma_F$ on $Z(\hat{H_i})$ factors through $\Gal(E/F) = \la \sigma \ra$.
Let $\Gamma_E = \Gal(\bar{F}/E) \subset \Gamma_F$, then 
$H^1\lp\Gamma_E, Z(\hat{H_i})\rp = \Hom_{{\rm cts}}(\Gamma_E,  Z(\hat{H_i}))$ because $\Gamma_E$ acts
trivially on $Z(\hat{H_i})$.

For all but finitely many places $w$ of $E$, $w$ either
lies above a place $v$ of $F$ which is prime and unramified in $E$, or it lies above a place $v$ of $F$
which splits into two places $w, w'$ of $E$.
In the former case, we have $\Gal(E_w/F_v) = \Gal(E/F)$.  In the later, we have
$\Gal(E_{w}/F_v) = \Gal(E_{w'}/F_v)= 1$.


Given any $x \in \ker^1( F, Z(\hat{H_i})) \subset H^1(\Gamma_F, Z(\hat{H_i}))$, the restriction
$x|_{\Gamma_{E}}$ of $x$ to $\Gamma_E$ is a (continuous) homomorphism from $\Gamma_E$ to $Z(\hat{H_i})$.
By assumption, the image of $x$ in $\prod_v H^1\lp\Gamma_{F_v}, Z(\hat{H_i})\rp$ is trivial, which means,  
in particular, that the image of $x|_{\Gamma_E}$ in  $H^1\lp\Gamma_{E_w}, Z(\hat{H_i})\rp$ is trivial for all
but finitely many places $w$ of $E$.  Since $x|_{\Gamma_E}$ is a continuous homomorphism, 
we conclude that $x|_{\Gamma_E} = 1$.

Fix an element $\sigma_0 \in \Gamma_F - \Gamma_E$.
Then, any $x \in \ker^1(F, Z(\hat{H_i}))$ produces  
a cocycle $y \in H^1(\Gal(E/F), Z(\hat{H_i}))$ defined by
$y_\sigma := x_{\sigma_0}$.  Since $x|_{\Gamma_E} = 1$, $y$ is independent of the choice of $\sigma_0$.
Pick a place $v$ of $F$ which remains prime in $E$.
Since $x$ is an element of
$\ker^1(F, Z(\hat{H_i}))$, the image of $y$ in  $H^1(\Gal(E_v/F_v), Z(\hat{H_i}))$ is trivial. 
But, $H^1(\Gal(E_v/F_v), Z(\hat{H_i}))$ is equal to  $H^1(\Gal(E/F), Z(\hat{H_i}))$ because $v$ is prime in $E$, so
$y$ itself is trivial.  The claim follows.
\end{proof}
\begin{claim}
$\abs{{\rm Out}(\mb{H}_i, \mathcal{H}_i,s_i, \xi_i)} = 2$ for $i = 1, 2$.
\end{claim}
\begin{proof}
For $i = 1, 2$, the only nontrivial element of $\text{Out}(\mb{H}_i, \mathcal{H}_is_i, \xi_i)$
comes from the action of $\sigma \in \Gal(E/F)$ on $\hat{H_i}$.
\end{proof}
\begin{corollary}\label{corollary:KStracerefined}
The following holds for matching functions:
\begin{equation}\label{eq:KStracerefined}
T_e(G, f) = 
\frac{1}{2}\sum_{\{\delta_{1}\}}SO_{\delta_{1}}(f_1) +
\frac{1}{2}\sum_{\{\delta_{2}\}}SO_{\delta_{2}}(f_2)
\end{equation}
where the sums are over stable conjugacy classes of elliptic 
$\mb{G}$-regular elements in the $F$-points of the 
$\ve$-endoscopic groups.
\end{corollary}
\section{Application of Kottwitz-Shelstad's Formula}\label{sec:appKS}
Recall, for $f \in C(\mb{G}(\Af), \omega)$, $I^d(G, f,\ve)$ denotes 
the $\ve$-discrete part of the fine $\chi$-expansion
of $\mb{G}$, $I^c(G, f,\ve)$ denotes the continuous part of the expansion, and
$I^d(H_i, f_i)$, $I^c(H_i, f_i)$ ($i = 1, 2$) denote analogous objects for $\mb{H}_i$.

Fix a two places $u_1, u_2$ of $F$ which are prime in $E$.
By Proposition \ref{corollary:ellipticOexp} and the trace formula \eqref{eq:twistedTF}, for any function $f$ in $\EL{\mb{G}(\Af), \omega}$,
we have
\[
T_e(G, f) = I^d(G, f, \ve) + I^c(G, f, \ve).
\]
Since $\mb{H_2}$ is a quotient of $\rR_{E/F}\GL(2) \times \mbb{G}_m$, there is no distinction between
stable and ordinary conjugacy classes.
For any function $f_2$ in $\EL{\mb{H_2}(\Af), \omega}$ which matches a function in $\EL{\mb{G}(\Af), \omega}$, 
the orbital integral of $f_2$ is zero at any element in $\mb{H_2}(F)$
which is not elliptic $G$-regular; hence, for such $f_2$ we have
\[
ST^{**}_e(H_2, f_2) = T_e(H_2, f_2) = I^d(H_2, f_2) + I_{H_2}^c(H_2, f_2).
\]
For $\mb{H_1}$, where there is a distinction between stable and ordinary conjugacy, the
situation is more complicated.

As in the case of $\mb{H_2}$, if the function $f_1 \in \EL{\mb{H_1}(\Af), \omega}$  
matches some $f \in \EL{\mb{G}(\Af), \omega}$, then
the orbital integral of $f_1$ is zero at all
elliptic regular elements which are not $\mb{G}$-regular.
Hence, $ST^{**}_e(H_1, f_1) = ST_e(H_1, f_1)$, where $ST_e(H_1, f_1)$ is the sum of stable orbital integrals
over the stable conjugacy classes of elliptic regular elements in $\mb{H_1}(F)$.

Applying Kottwitz-Shelstad's formula in the context of standard endoscopy for $\mb{H_1}$,
we have:
\begin{equation}\label{eq:stablekappa}
I^d(H_1, f_1) + I^c(H_1, f_1)
= T_e(H_1, f_1) = ST_e(H_1, f_1) + \frac{1}{2} T_e^\kappa(H_1, f_1),
\end{equation}
where $T_e^\kappa(H_1, f_1)$ is
a sum over stable conjugacy classes in $\mb{H}(F)$
of $\kappa$-orbital integrals $O_{\delta_{H_1}}^\kappa(f_1)$.
In this case, there is only one nontrivial $\kappa$
(in fact, $\kappa = \ve$), and it is associated with a
unique (up to equivalence) standard endoscopic group $\mb{T_1}$ of $\mb{H_1}$.
In particular, $\mb{T_1}$ is the $F$-torus whose $F$-points are:
\[
\mb{T_1}(F) = \{(x, y) \in E^\times \times E^\times : \N_{E/F} x = \N_{E/F} y\}.
\]

We rewrite \eqref{eq:stablekappa} as:
\[
I^d(H_1, f_1) - \frac{1}{2}T_e^\kappa(H_1, f_1) + I^c(H_1, f_1) = ST_e(H_1, f_1).
\]
Recall from Section \ref{sec:stablespectrum} that $SI^d(H_1, f_1)$ denotes the stable part of 
the discrete spectral expansion $I^d(H_1, f_1)$.
From the theory of \cite{LL}, $\frac{1}{2}T_e^\kappa(H_1, f_1)$ is equal to the unstable part of
$I^d(H_1, f_1)$; thus,
$I^d(H_1, f_1) - \frac{1}{2}T_e^\kappa(H_1, f_1)$ is equal to $SI^d(H_1, f_1)$. 

We conclude that the following equation holds for any $f_1$ in 
\dt $\EL{\mb{H_1}(\Af), \omega}$ matching some function in $\EL{\mb{G}(\Af), \omega}$:
\begin{equation}\label{eq:stablespectrum}
ST_e(H_1, f_1) = SI^d(H_1, f_1) + I^c(H_1, f_1).
\end{equation}

In conclusion, the following equations hold for $f \in \EL{\mb{G}(\Af), \omega}$ and
matching  functions $f_i$ in $\EL{\mb{H}_i, \omega}$ ($i = 1, 2$):
\begin{multline}\label{eq:prediscretecont}
I^d(G, f,\ve) + I^c(G, f,\ve) \\
= \frac{1}{2}\Big( SI^d(H_1, f_1) + I^c(H_1, f_1) + I^d(H_2, f_2) + I^c(H_2, f_2) \Big).
\end{multline}
\indexi{Kottwitz-Shelstad's trace formula|)}%

\chapter{Global Lifting}\label{chap:global}
\section{The $\ve$-Trace Identity}\label{sec:vetraceidentity}
\subsection{Some Notation and Terminology}
\subsubsection{Global Data}\label{sec:globaldata}
\indexi{global datum|(}%
Let $F$ be a number field.
Let $\mb{H}$ be a reductive group over $F$ with $L$-group ${}^L H$.
Let $V$ be the set of places of $F$.
Let $S$ be a finite set of places of $V$.  Let $\{c_v : v \in V - S\}$ be a set of
conjugacy classes in ${}^L H$.

Define a (countably-)infinite-tuple as follows:
\[
\mathcal{C}(S, {}^L H) := \left[ c_v\right]_{v \in V - S} \subset \prod_{v \in V - S}{}^L H.
\]
We call $\mathcal{C}(S, {}^L H)$ a {\bf global datum} in ${}^L H$.  For simplicity, we sometimes
replace ${}^L H$ by a finite Galois form.  For example, if $\mb{H}$ is split over $F$, then we use
$\hat{H}$ instead.

For any automorphic representation $\pi$ of $\mb{H}(\Af)$, there exists a finite set of places $S$
such that $\pi_v$ is unramified for all $v \notin S$.
For each $v \notin S$, denote by $c(\pi_v)$ the Frobenius-Hecke (abbrev. F-H) class in ${}^L H$ 
\indexi{Frobenius-Hecke class}%
which parametrizes $\pi_v$ (see \cite{B}).
Then, $\pi$ defines the global datum
\[
\mathcal{C}(\pi, S, {}^L H) := \left[ c(\pi_v)\right]_{v \in V - S} \subset \prod_{v \in V - S}{}^L H.
\]
We call $\mathcal{C}(\pi, S, {}^L H)$ a {\bf global datum} (or the {\bf $S$-global datum}) of $\pi$.
The global datum $\mc{C}(\pi, S, {}^L H)$ is {\bf defined} if and only if $\pi_v$ is unramified for all $v \notin S$.
Observe that if $\mc{C}(\pi, S, {}^L H)$ is defined, then so is $\mc{C}(\pi, S', {}^L H)$
for any finite subset $S'$ of $V$ containing $S$.


\subsubsection{Global Packets}
\indexi{packet!global}%




Let $\{\{\Pi_v\} : v \in V\}$ be a collection of local (quasi-)packets $\{\Pi_v\}$ of $\mb{G}(F_v)$
such that, for almost all (finite) $v$, $\{\Pi_v\}$ contains a unique unramified representation $\Pi_v^0$.
We define the {\bf (quasi-)packet} $\{\Pi\}$ to be the restricted tensor product
\indexi{restricted tensor product}%
\[
\otimes_{v\in V}\{\Pi_v\}
:=
\left\{\otimes_{v\in V} \Pi'_v
:
\Pi'_v \in \{\Pi_v\}\text{ for all } v,
\Pi'_v = \Pi_v^0 \text{ for almost all } v.
\right\}.
\]
For any automorphic representation $\Pi$ of $\mb{G}(\Af)$,
we denote by $\{\Pi\}$ the global packet to which $\Pi$ belongs.  
Under this notation, $\{\Pi\} = \{\Pi'\}$ if $\Pi' \in \{\Pi\}$.  

We say that a (quasi-)packet of representations of $\mb{G}(\Af)$ is a {\bf discrete spectrum (quasi-)packet} 
if it contains a discrete spectrum representation.
We say that a discrete spectrum packet of $\mb{G}(\Af)$ is {\bf $\ve$-invariant} 
\indexi{packet!epsilon-invariant@$\ve$-invariant}%
if it contains an $\ve$-invariant discrete spectrum representation.
We say that a discrete spectrum packet $\{\Pi\}$ is {\bf stable} if every member of 
$\{\Pi\}$ occurs with the same multiplicity in the discrete spectrum of $\mb{G}(\Af)$.
Otherwise, we say that $\{\Pi\}$ is {\bf unstable}.

Let $\{\Pi\} = \otimes_{v \in V}\{\Pi_v\}$ be a (quasi-)packet of representations of $\mb{G}(\Af)$.
Let $S$ be a finite set of places
such that, for all $v \notin S$, the local packet $\{\Pi_v\}$ contains a unique
unramified representation $\Pi_v^0$.
For $v \notin S$, let $c(\{\Pi_v\})$ denote the F-H class in $\hat{G}$ which parametrizes $\Pi_v^0$.
Let $\{\Pi\}_S$ denote the subset of $\{\Pi\}$ of representations 
$\Pi$ such that $\Pi_v = \Pi_v^0$ for all $v \notin S$.  
We call $\{\Pi\}_S$ the {\bf $S$-part} of $\{\Pi\}$.
\indexi{packet!global!$S$-part}%
In general, the cardinality of a global packet may be infinite, but the cardinality of $\{\Pi\}_S$ is finite.
Put $\mc{C}(\{\Pi\}, S, \hat{G}) := \left[ c(\{\Pi_v\})\right]_{v \in V - S} \subset \prod_{v \notin S}\hat{G}$.


Let $\Vr$ denote the set of finite places of $F$ which are unramified in $E$.
\indexs{V@$\Vr$}%
Let $\{\pi_i\} = \otimes_{v \in V}\{\pi_{i,v}\}$ be a global (quasi-)packet of $\mb{H}_i(\Af)$, where $i = 1$ or $2$.
Let $S$ be a finite set of places containing $V - \Vr$ such that 
the $S$-part $\{\pi_i\}_S$ of $\{\pi_i\}$ is nonempty.
Let ${}^L H_i$ denote the finite Galois form $\hat{H_i}\rtimes \Gal(E/F)$ of the
$L$-group of $\mb{H}_i$.  
Then,
$\{\pi_i\}$ defines the global datum 
\[\mathcal{C}(\{\pi_i\}, S, {}^L H_i) = [c(\{\pi_{i,v}\})]_{v \notin S} \subset \prod_{v \notin S} {}^L H_i.\]
As we shall discuss in Section \ref{FHclass}, each $c(\{\pi_{i,v}\}) \subset {}^L H_i$ lifts via 
the $L$-group embedding $\xi_i : {}^L H_i \hookrightarrow \hat{G}$ 
to a conjugacy class $\xi_i(c(\{\pi_{i,v}\}))$ in $\hat{G}$.  We obtain the following global datum:
\[
\mathcal{C}(\xi_i(\{\pi_i\}), S, \hat{G}) := 
\left[ \xi_i\lp c(\{\pi_{i,v}\})\rp \right]_{v \notin S} \subset \prod_{v \notin S}\hat{G}.
\]
If $\mathcal{C}(\xi_i(\{\pi\}), S, \hat{G})$ coincides with the global datum $\mc{C}(\{\Pi\}, S, \hat{G})$
for some global (quasi-)packet $\{\Pi\}$ of $\mb{G}(\Af)$,
we say that the packet $\{\pi_i\}$ {\bf lifts} to $\{\Pi\}$ and write $\{\Pi\} = \xi_i(\{\pi\})$.
\indexi{lift}%
\indexi{global datum|)}%
\subsection{Langlands' Trick}
For $f \in C(\mb{G}(\Af))$,
let $I^d(G, f, \ve)$ denote the discrete part of the $\ve$-twisted fine $\chi$-expansion for $\mb{G}$.
Let  $I^c(G, f, \ve)$ denote the continuous part of the expansion.
Define $I^d(H_i, f_i)$, $I^c(H_i, f_i)$ ($i = 1, 2$) likewise (without twisting) for $\mb{H}_i$.

The following proposition depends on the validity of the Fundamental Lemma (see Chapter \ref{chap:fundlemma} in the Appendix).
\indexi{Fundamental Lemma}%
\begin{prop}\label{prop:discreteeq}
Fix two finite places $u_1$, $u_2$ of $F$ which are prime in $E$.
The following holds for matching functions $f, f_1, f_2$ whose components at $u_1, u_2$ are elliptic:
\indexi{elliptic function}%
\begin{equation}\label{eq:discretediff}
I^d(G, f, \ve) - \frac{1}{2}SI^d(H_1, f_1) - \frac{1}{2}I^d(H_2, f_2) = 0
\end{equation}
\end{prop}
We content ourselves with giving only a sketch of the proof, which is essentially
a retelling of an argument of Langlands' in \cite[Chap. 11]{L1}.  

Recall that $\Vr$ is the set of finite places of $F$ which are unramified in $E$, and that $\omega$ is a
fixed character of the center of $\mb{G}$.
Let $S$ be a finite set of places of $F$ containing $\lp V - \Vr\rp \cup \{u_1, u_2\}$ such that
$\omega_v$ is unramified, and $f_v, f_{1, v}, f_{2,v}$ are spherical, for all $v \notin S$.

Rewrite equation \eqref{eq:prediscretecont} as follows:
\begin{multline}\label{eq:discretecont}
I^d(G, f,\ve) - \frac{1}{2}SI^d(H_1, f_1) - \frac{1}{2}I^d(H_2, f_2)\\
= -\lp I^c(G, f,\ve) - \frac{1}{2} I^c(H_1, f_1) - \frac{1}{2}I^c(H_2, f_2)\rp.
\end{multline}

Fix a finite place $u \notin S$ which splits in $E$.  Regard $f_v, f_{1,v}, f_{2,v}$ as fixed for $v \notin S\cup\{u\}$, and
let $f_u$ vary in $\mc{H}(G_u, \omega_u)$ (hence $f_{1, u}$ and $f_{2, u}$ also vary by the matching condition).
Fix the element $\diag(1, 1, \lambda, \lambda) \in \hat{G} = \GSp(2, \CC)$, where $\lambda = \omega_u(\vp_u)$.
Let $\mc{T}(G_u, \omega_u)$ be the set of conjugacy classes in $\hat{G}$ whose similitude
factors are equal to $\omega_u(\vp_u)$.  Then, any class $\mc{C}$ in $\mc{T}(G_u, \omega_u)$
is represented by
\[
\dmfour{x}{y}{\lambda y^{-1}}{\lambda x^{-1}}
\]
for some $x, y \in \CC^\times$.
If $\mc{C}$ corresponds to a unitarizable unramified representation of
$G_u$, then $(x, y)$ must lie in a compact subset $X$ of $\CC^2$.  
More precisely, 
using the results of \cite{ST}, it can be shown that $X = T \cup U$, where
\[
T = \left\{(\abs{\vp_u}^s, \abs{\vp_u}^t) : (s, t) \in i\mf{a}_0/i\mf{a}_G = i\mbb{R} \times i\mbb{R}\right\}
\subset \mbb{C}^2
\]
(recall that $\mf{a}_0 = X_*(\mb{A_0})\otimes_{\mbb{Z}} \mbb{R}$, where $\mb{A}_0$ is the maximal
diagonal torus of $\mb{G}$), and $U$ is a compact nondiscrete subset of a real plane (dependent on the choice
of $a$, $b$) in $\mbb{C}^2$.

Conversely, any element $(x, y)$ in $X$ gives rise to a conjugacy class in $\hat{G}$ which corresponds
to a unitarizable unramified representation.

The Hecke algebra $\mc{H}(G_u, \omega_u)$ yields via Satake transform
\indexi{Satake transform}%
the algebra $\mc{H}^{\vee}$ of Laurent series $\phi$ in $x, y$ with the property that 
\indexs{H@$\mc{H}^{\vee}$}%
\begin{equation}\label{eq:conditiononCX}
\phi(x, y) = \phi(x,\lambda y^{-1}) = \phi(y, x).
\end{equation}
The above property holds
because the Satake transforms of Hecke functions are defined on the set of conjugacy classes in $\hat{G}$.
\indexi{Satake transform}%

Let $C({X})'$ denote the space of complex valued continuous functions $\phi$ on ${X}$ satisfying \eqref{eq:conditiononCX}.
By the Stone-Weierstrass theorem, $\mc{H}^\vee$ is dense in $C(X)'$ under the sup-norm on $C(X)'$.

To prove the proposition, we use each side of equation \eqref{eq:discretecont} to define a functional on $C({X})'$.  
We then argue that these two functionals must be identically zero.
\begin{claim}\label{claim:discretefunctional}
The left hand side of equation \eqref{eq:discretecont} defines a continuous linear functional on $C({X})'$ 
which is given by an atomic measure on ${X}$.
\indexi{atomic measure}%
\end{claim}
\begin{proof}
A nonzero term in $I^d(G, f, \ve)$, with corresponding quadruple $(\mb{M}, \mb{G}, \tau, s)$, has the form:
\begin{multline}\label{eq:discretetermprod}
c_{(\mb{M}, \mb{G}, \tau, s)} \tr I_{P, \tau}(f)M^\ve_P(s, 0)\\
=
c_{(\mb{M}, \mb{G}, \tau, s)} \prod_{v \notin S} \tr I_{P, \tau_v} (f_v)M^\ve_P(s, 0)_v
\cdot \prod_{v \in S} \tr I_{P, \tau_v}(f_v)M_P^\ve(s, 0)_v,
\end{multline}
where $c_{(\mb{M}, \mb{G},\tau,s)}$ is a constant dependant on $(\mb{M}, \mb{G}, \tau, s)$.
Since by assumption $f_v$ is spherical for all $v \notin S$, the representations $I_{P, \tau_v}$ must be 
unramified, or else \eqref{eq:discretetermprod} is equal to zero.
For $v \notin S$, the space of $\mb{G}(\mc{O}_v)$-fixed vectors in $I_{P, \tau_v}$ is
one dimensional; hence, we can (and do) normalize $M_P^\ve(s, 0)_v$ to be $1$.
Let $\mc{C}_{\tau_v}$ ($v \notin S$) be the F-H class in $\hat{G}$ which parametrizes $I_{P, \tau_v}$.  Then
\eqref{eq:discretetermprod} may be written as
\[
c_{(\mb{M}, \mb{G}, \tau, s)}\prod_{v \notin S} \fhat_v(\mc{C}_{\tau_v})\cdot 
\prod_{v \in S}\tr I_{P, \tau_v}(f_v) M^\ve_P(s, 0)_v,
\]
where $\fhat$ is the Satake transform of $f$.
\indexi{Satake transform}%
Analogous statements hold for $\mb{H_1}$ and $\mb{H_2}$, with $\ve$ replaced with $1$.

Through the embedding of $L$-groups, every conjugacy class $\mc{C}_i$ in ${}^L H_i$
($i = 1, 2$) lifts to a conjugacy class $\mc{C}$ in $\hat{G}$.
Moreover, if $\mc{C}_i$ corresponds to a unitarizable unramified representation of $H_{i,v}$,
then the unramified representation of $G_v$ parametrized by $\mc{C}$ must also be
unitarizable (see Sect. \ref{FHclass}, \cite{ST}).

Let $v$ be a place in $\Vr$.
Let $f_{i, v}$ ($i = 1, 2$), $f_v$ be 
functions in $\mc{H}(H_{i,v}, \omega_v)$, $\mc{H}(G_v, \omega_v)$, respectively,
such that 
\[
\fhat_{i, v}(\mc{C}_i) = \fhat_{v}(\mc{C})
\]
for every conjugacy class $\mc{C}_i$ in ${}^L H_i$ which lifts to a class $\mc{C}$ in $\hat{G}$.
By the Fundamental Lemma, such Hecke functions have matching orbital integrals.

Given any global datum $\mc{G} = \mc{C}(S, \hat{G})$, let
$Q(G, \mc{G})$ denote the set of equivalence classes of $\ve$-invariant,
\indexs{Q@$Q(\;,\;)$}%
$\ve$-discretely occurring automorphic representations of $\mb{G}(\Af)$
whose components at each $v \notin S$ are unramified and are parametrized by
the local component $\mc{G}_v$ of $\mc{G}$.


Let $Q(H_2, \mc{G})$ denote the set of equivalence classes of the discretely occurring automorphic representations $\pi_2$ of 
$\mb{H_2}(\Af)$ such that, for all $v \notin S$, $\pi_{2,v}$ is unramified, and the F-H class in ${}^L H_2$
parametrizing $\pi_{2,v}$ lifts to $\mc{G}_v$.



Let $\{\pi_1\}$ be a global (quasi-)packet of representations of $\mb{H_1}(\Af)$.
For any function $f_1$ in $C(\mb{H_1}(\Af), \omega)$,
put 
\[
\displaystyle \tr \{\pi_1\}(f_1) := \sum_{\pi_1'\in \{\pi_1\}}m(\pi_1')\tr \pi_1'(f_1),
\indexs{t@$\tr\{{\rm packet}\}(\;)$}%
\] 
where $m(\pi_1')$ is the multiplicity of $\pi_1'$ in the discrete spectrum of $\mb{H_1}(\Af)$.
Let \dt \[\str\!\{\pi_1\}(f_1)\] denote the stable part of $\tr \{\pi_1\}(f_1)$ (see Section \ref{sec:stablespectrum}).
\indexs{t@$\str$}%
Assume that $f_1$ is the tensor product of local components $f_{1,v}$.
Put \[\tr\! \{\pi_{1,v}\}(f_{1,v}) := \sum_{\pi_{1,v}' \in \{\pi_{1,v}\}}\tr\pi_{1,v}'(f_{1,v}).\]
Suppose the members of $\{\pi_1\}$ are the constituents of $\tau\otimes_1 \eta$, where
$\tau$ is an automorphic representation of $\GL(2, \Af)$, and $\eta$ is a character
of $\idc{E}$.
By \cite{LL}, we have
\[
\str \{\pi_1\}(f_1) = a_{\{\pi_1\}} \prod_{v \in V}\tr \{\pi_{1,v}\}(f_{1,v}),
\]
where 
$a_{\{\pi_1\}} = \begin{cases} \frac{1}{2} &\text{if } \tau \text{ is } E\text{-monomial},\\
1 &\text{ otherwise}.\end{cases}$
\indexs{a@$a_{\{{\rm packet}\}}$}%

Let $\bar{Q}(H_1, \mc{G})$ denote the collection of equivalence classes of
\indexs{Q@$\bar{Q}(\;,\;)$}%
global packets of $\mb{H_1}$ whose global data lift to $\mc{G}$.

Let
\begin{multline*}
\alpha_{\mc{G}} 
=
\sum_{\pi \in Q(G, \mc{G})} 
m(\pi)c_\pi \prod_{v \in S} \tr \pi_v(f_v\times \ve_v)M(\pi)_v \\
- \frac{1}{2} 
\sum_{\{\pi_1\} \in \bar{Q}(H_1, \mc{G})} 
m(\{\pi_1\})c_{\{\pi_1\}} a_{\{\pi_1\}}\prod_{v \in S} \tr\{\pi_{1,v}\}(f_{1,v}) M(\{\pi_1\})_v\\
- \frac{1}{2}
\sum_{\pi_2 \in Q(H_2, \mc{G})} 
m(\pi_2)c_{\pi_2} \prod_{v \in S} \tr \pi_2(f_{2,v})M(\pi_2)_v.
\end{multline*}
Here, $c_\pi$ is the constant defined by \eqref{eq:finechiexp1},
$M(\pi)$ is the intertwining operator $M^\ve_P(s, 0)$.  
The symbols $c_{\pi_i}, M(\pi_i)$ ($i = 1, 2$) are defined likewise.
\indexs{M@$M(\;)$}%
The factor $m(\cdot)$ denotes the multiplicity of the representation/packet.

By the Fundamental Lemma,
the left hand side of equation \eqref{eq:discretecont} may be written as the sum
\begin{equation}\label{eq:discretegd}
\sum_{\mc{G} \in \{\mc{G}\}}\alpha_{\mc{G}}\prod_{v \notin S}\fhat_v(\mc{G}_v)
\end{equation}
over the set $\{\mc{G}\}$ consisting of the
global data $\mc{G} = \mc{C}(S, \hat{G})$ for which $\alpha_{\mc{G}}$ is nonzero.  Note that the set $\{\mc{G}\}$
is countable by the definition of the discrete part of the fine $\chi$-expansion 
(see Sect. \ref{sec:finechiexpoverview}).



Recall that we have fixed a finite place $u$ which splits in $E$ such that the test functions have spherical local
components at $u$.
Let $\{r_{u, j}\}_{j = 1, 2, \dots}$ be the set of distinct $\mc{G}_u$'s for $\mc{G} \in \{\mc{G}\}$.  Let 
\[
c_j = \sum_{\{\mc{G} \in \{\mc{G}\}\;:\;\mc{G}_u = r_{u, j}\}} 
\alpha_{\mc{G}}\prod_{v \notin S \cup \{u\}}\fhat_v(\mc{G}_v).
\]
Then, \eqref{eq:discretegd} is equal to
\begin{equation}\label{eq:atomicmeasure}
\sum_{j = 1, 2, \dots}c_j \fhat_u(r_{u, j}).
\end{equation}

The representations which appear in the ($\ve$-)discrete spectra for the groups have unitarizable
local components at every place.  Hence, the conjugacy classes $r_{u, j}$ in 
\eqref{eq:atomicmeasure} correpond to elements $(x_j, y_j)$ in the compact set $X$.
Since $\mc{H}^\vee$ is dense in $C({X})'$, \eqref{eq:atomicmeasure} extends to a continuous linear functional
$\mc{D}_d$ on $C({X})'$ defined by:
\[
\mc{D}_d(\phi) = \sum_{j = 1, 2, \dots} c_j\, \phi(x_j, y_j),\quad \forall
\phi \in C({X})'.
\]
Note that the correspondence between $r_j$ and $(x_i, y_i)$ involves a choice of representative in $\hat{G}$ of the conjugacy class $r_j$.  However, that choice does not affect the definition of the functional.
The claim follows.
\end{proof}

\begin{claim}\label{claim:contfunctional}
The right hand side of equation \eqref{eq:discretecont} defines a sum of continuous linear functionals on $C({X})'$, 
each of which is defined by a measure which is absolutely continuous with respect to the Lebesque measure on a subtorus of $T$.
\indexi{Lebesque measure}%
\end{claim}
\begin{proof}
The right hand side of \eqref{eq:discretecont} consists of terms in the continuous parts
of the ($\ve$-twisted) fine $\chi$-expansions for $\mb{G}$, $\mb{H_1}$, and $\mb{H_2}$.
Let $\mb{H}$ be either $\mb{G}$, $\mb{H_1}$, or $\mb{H_2}$.
Let $\{\chi\} = (\mb{M}_H, \mb{L}_H, \tau, s)$ be the quadruple associated with a
term in the continuous part of the fine $\chi$-expansion for $\mb{H}$.  In particular, $\mb{L}_H$
is a proper Levi subgroup of $\mb{H}$.
The term associated with $\{\chi\}$ has the form:
\begin{equation}\label{eq:contterm}
c_{\{\chi\}}\int_{i\mf{a}_{L_H}^*/i \mf{a}_H^*}
\tr \mc{M}^T_{L_H}(P_H, \zeta) I_{P_H, \tau}(\zeta, f_H\times\ve) M_{P_H}(s, \zeta)\,d\zeta,
\end{equation}
where $c_{\{\chi\}}$ is a constant, $\mb{P}_H$ is the standard parabolic with Levi component $\mb{M}_H$, and
$f_H$ is a test function in $C(\mb{H}(\Af), \omega)$ satisfying the hypothesis of
Proposition \ref{prop:discreteeq}.  Here, we take $\ve$ to be $1$ if $\mb{H}$ is one of the twisted endoscopic groups.  
We assume that \eqref{eq:contterm} is nonzero.  In particular, since $f_u$ is spherical, 
$\tau_u$ is necessarily unramified.
For $\zeta \in i\mf{a}_{L_H}^*/i \mf{a}_H^*$,
let $\mc{C}_{H_u, \tau_u}(\zeta)$ denote the conjugacy class in ${}^L H$ which parametrizes $I_{P_{H,u}, \tau_u}(\zeta)$.
Put $\mc{C}_{H_u, \tau_u} := \mc{C}_{H_u, \tau_u}(0)$.

Let $\tau_u^0$ be the unramified representation of the diagonal torus $A_{0, u}$ of $G_u$ such that
$\mc{C}_0 := \mc{C}_{G_u, \tau^0_u}$ is represented by $\diag(1, 1, \lambda, \lambda)$ in $\hat{G}$.
For each quadruple $(\mb{M}, \mb{L}, \tau, s)$ which appears in the $\ve$-twisted fine $\chi$-expansion for $\mb{G}$
such that $\tau_u$ is unramified, there is an element $\zeta_{\tau}$ in $i\mf{a}_0/i\mf{a}_G$ such that 
$I_{P_{M, u}, \tau_u}$ ($\mb{P}_M$ being the standard parabolic with Levi component $\mb{M}$) is equal to
$I_{P_{0, u}, \tau^0_u}(\zeta_\tau)$ ($\mb{P}_0$ being the standard parabolic with Levi component $\mb{A}_0$).
Note that, for each $\zeta_L \in i\mf{a}_L/i\mf{a}_G$, we have 
$I_{P_{M, u}, \tau_u}(\zeta_L) = I_{P_{0, u}, \tau^0_u}(\zeta_\tau + \zeta_L)$.

For each proper Levi subgroup $\mb{L}$ in $\mb{G}$, let $Q_G(L)$ denote the set of quadruples (associated with
terms in the $\ve$-twisted fine $\chi$-expansion for $\mb{G}$) of the form
$(\mb{M}, \mb{L}, \tau, s)$.  As in \cite[Chap. 11]{L1}, it may be shown that
there exists a function $d_{L}$ on $Q_G(L) \times i\mf{a}_{L}/i\mf{a}_G$, integrable with respect to the second
argument, such that the sum of all the terms in the $\ve$-twisted fine $\chi$-expansion for $\mb{G}$ associated with 
the quadruples in $Q_G(L)$ is  equal to
\[
\sum_{\{\chi\} = (\mb{M}, \mb{L}, \tau, s) \in Q_G(L)} \int_{i\mf{a}_{L}^*/i \mf{a}_G^*} 
 \fhat_u(\mc{C}_0(\zeta_\tau + \zeta))\,d_{L}(\{\chi\}, \zeta)
\,d\zeta.
\]
Here, $\mc{C}_0(\zeta_\tau + \zeta_{L}) := \mc{C}_{G_u, \tau^0_u}(\zeta_\tau + \zeta_{L})$.

For $j = 1, 2$, let $\mb{A}_j$ be the diagonal torus of $\mb{H}_j$.  
Let $\mc{C}_j$ be a conjugacy class in ${}^L H_j$ which lifts to $\mc{C}_0$.  Let $\tau_{j, u}^0$ 
be the unramified representation of $A_{j, u}$ which is parametrized by  $\mc{C}_j$.
For each quadruple $(\mb{M}_j, \mb{L}_j, \tau_j, s)$ which appears in the fine $\chi$-expansion for $\mb{H}_j$ such that
$\tau_{j, u}$ is unramified, there is an element $\zeta_{\tau_j}$ in $i\mf{a}_j^*/i\mf{a}_{H_j}^*$ such that
$I_{P_{M_j, u}, \tau_{j, u}} = I_{P_{j, u}, \tau_{j, u}^0}(\zeta_{\tau_j})$.  Here, $\mb{P}_{M_j}$ is the standard parabolic subgroup
of $\mb{H}_j$ with Levi component $\mb{M}_j$, and $\mb{P}_j$ is the minimal standard parabolic subgroup of $\mb{H}_j$.
For each proper Levi subgroup $\mb{L}_j$ in $\mb{H}_j$, let $Q_{H_j}(L_j)$ denote the set of quadruples (associated with
terms in the fine $\chi$-expansion for $\mb{H}_j$) of the form $(\mb{M}_j, \mb{L}_j, \tau_j, s)$.
As in the case of $\mb{G}$, there is a function $d_{j, L_j}$ on $Q_{H_j}(L_j)\times i\mf{a}_{j}/i\mf{a}_{H_j}$,
integrable with respect to the second argument, such that the sum of all the terms in the fine $\chi$-expansion for $\mb{H}_j$ 
associated with the quadruples in $Q_{H_j}(L_j)$ is equal to
\[
\sum_{\{\chi_j\} = (\mb{M}_j, \mb{L}_j, \tau_j, s) \in Q_{H_j}(L_j)}
\int_{i\mf{a}_{L_j}^*/i \mf{a}_{H_j}^*}
\fhat_{j, u}(\mc{C}_j(\zeta_{\tau_j} + \zeta))\, d_{j, L_j}(\{\chi_j\}, \zeta)
\,d\zeta.
\]

For $\zeta_j \in i\mf{a}_j^*/i\mf{a}_{H_j}^*$ ($j = 1, 2$), let $\zeta_j^G$ denote the element in $i\mf{a}_0^*/i\mf{a}_{G}^*$
which is the lift of $\zeta_j$ via the $L$-group embedding $\xi_j : {}^L H_j \hookrightarrow \hat{G}$.
For the image $f_j \in \mc{H}(H_{j, u}, \omega_u)$ of $f \in \mc{H}(G_u, \omega_u)$
under the map between the Hecke algebras which is dual to $\xi_j$, we have
\[
\fhat(\mc{C}_0(\zeta_j^G)) = \fhat_j(\mc{C}_j(\zeta_j)),\quad \forall \zeta_j \in  i\mf{a}_j^*/i\mf{a}_{H_j}^*.
\]
By the Fundamental Lemma, $f_j$ and $f$ have matching orbital integrals.

Collecting terms on the right hand side of \eqref{eq:discretecont} and factoring, 
we obtain for each (up to conjugation) proper Levi subgroup $\mb{L}$ of $\mb{G}$, 
a function $B_L$ on $Q_G(L) \times i\mf{a}_L^*/i \mf{a}_G^*$, integrable with respect to the second argument,
 such that the right hand side of 
\eqref{eq:discretecont} is equal to
\begin{equation}\label{eq:contintegralG}
\sum_{\mb{L}} \sum_{\{\chi\} = (\mb{M}, \mb{L}, \tau, s) \in Q_L(G)}
\int_{i\mf{a}_L^*/i \mf{a}_G^*} \fhat_{u}(\mc{C}_0(\zeta_\tau + \zeta))\,B_L(\{\chi\}, \zeta)\, d\zeta.
\end{equation}
The sum is over, up to conjugation, the set of proper Levi subgroups $\mb{L}$ of $\mb{G}$.

Since by the Stone-Weierstrass Theorem $\mc{H}^\vee$ is dense in $C(X)'$,
\indexi{Stone-Weierstrass Theorem}%
the expression \eqref{eq:contintegralG} extends to a continuous linear functional $\mc{D}_c$ 
on $C({X})'$.
Moreover, $\mc{D}_c$ is a sum of functionals, each of which given by a measure which is absolutely continuous with respect to the 
Lebesque measure on a subtorus of $T$.
\end{proof}

\begin{proof}[Proof of Proposition \ref{prop:discreteeq}]
By equation \eqref{eq:discretecont}, we have an equality of functionals: $\mc{D}_d = \mc{D}_c$.
By Claims \ref{claim:discretefunctional}, \ref{claim:contfunctional}, and the Riesz Representation Theorem (see \cite[Chap. 7]{Fo}), 
\indexi{Riesz Representation Theorem}%
we conclude that both $\mc{D}_d$ and $\mc{D}_c$ are zero.
\end{proof}

We have shown that \eqref{eq:atomicmeasure} is zero.  By the generalized linear independence
of characters (see \cite{FK1}), we conclude that each 
\[
c_j = \sum_{\{\mc{G} : \mc{G}_u = r_{u, j}\}} \alpha_{\mc{G}}\prod_{v \notin S \cup \{u\}}\fhat_v(\mc{G}_v)
\] 
is zero.  By applying the same argument repeatedly, it can be shown that if $U$ is any finite set of places disjoint from $S$,
then
\[
\sum \alpha_{\mc{G}}\prod_{v \notin S\cup U}\fhat_v(\mc{G}_v) = 0,
\]
where the sum is over those global data $\mc{G}$ with the property that 
$\mc{G}_v$ is equal to a fixed conjugacy class $r_{v}$ in $\hat{G}$ for each $v \in U$.
Using yet another argument of Langlands' in \cite[Chap. 11]{L1}, it can then be shown
that each $\alpha_{\mc{G}}$ in \eqref{eq:discretegd} is equal to zero.  
Equivalently, using the notation introduced earlier, the following corollary to Proposition
\ref{prop:discreteeq} holds:
\begin{corollary}\label{corollary:finediscreteeq}
Fix two places $u_1, u_2$ of $F$ which are prime in $E$.
For any 
\begin{itemize}
\item
finite set of places $S$ containing $(V - \Vr) \cup \{u_1, u_2\}$, 
\item 
global datum $\mc{G} = \mc{C}(S, \hat{G}) = [c_v]_{v \in V -S}$, where the
$c_v$'s are conjugacy classes in $\hat{G}$,
\item
matching functions $f \in \EL{\mb{G}(\Af), \omega}$, $f_i \in \EL{\mb{H}_i(\Af), \omega}$ $(i = 1, 2)$
whose local components at all $v\notin S$ are spherical,
\end{itemize}
the following holds:
\begin{multline}\label{eq:traceeq}
\sum_{\pi \in Q(G, \mc{G})} 
m(\pi)c_\pi \prod_{v \in S} \tr \pi_v(f_v\times \ve_v)M(\pi)_v\\
= \frac{1}{2} 
\sum_{\{\pi_1\} \in \bar{Q}(H_1, \mc{G})} 
m(\{\pi_1\})c_{\{\pi_1\}} a_{\{\pi_1\}}\prod_{v \in S}  \tr\{\pi_{1,v}\}(f_{1,v}) M(\{\pi_1\})_v\\
+ \frac{1}{2}
\sum_{\pi_2 \in Q(H_2, \mc{G})} 
m(\pi_2)c_{\pi_2} \prod_{v \in S} \tr \pi_{2, v}(f_{2,v})M(\pi_2)_v.
\end{multline}
\begin{remark}
Note that the above identity is of the trivial form $0 = 0$ if the global datum $\mc{G}$
correponds to representations whose central characters are not equal to the character $\omega$.
\end{remark}
\end{corollary}

We call equation \eqref{eq:traceeq} an {\bf $\ve$-trace identity} with respect to $S$. 
\indexi{epsilon-trace identity@$\ve$-trace identity}%
A representation which appears nontrivially on either side of the equation is said to
{\bf contribute} to the $\ve$-trace identity.
\indexi{contribution}%
By definition, the representations which appear nontrivially on the left hand side of \eqref{eq:traceeq}
form the  $S$-part of a global (quasi-)packet of $\mb{G}(\Af)$.
On the other hand,
the representations $\pi_i$ which contribute to the right hand side may comprise the $S$-parts of more
than one global packet of $\mb{H}_i(\Af)$.  For each group $\mb{H}_i$, the set of all
contributing representations form what we call a {\bf multi-packet} of $\mb{H}_i(\Af)$.
\indexi{packet!multi-}%

Let $\{\pi\} = \bigcup_{1\leq j \leq n}\{\pi_{j}\}$ be a multi-packet of $\mb{H}_i(\Af)$ ($i = 1$ or $2$),
where $\{\pi_{j}\}$ ($1 \leq j \leq n$) is a global (quasi-)packet of $\mb{H}_i(\Af)$.
Let $\{\pi\}_S = \bigcup_{1 \leq j \leq n}\{\pi_j\}_S$, the union of the $S$-parts of the global packets
$\{\pi_j\}$.
We call $\{\pi\}_S$ the {\bf $S$-part} of the multi-packet $\{\pi\}$.


If the $S$-part of a global packet $\{\Pi\}$ of $\mb{G}(\Af)$ 
and the $S$-parts of multi-packets $\{\pi_i\}$ of $\mb{H}_i(\Af)$ ($i = 1, 2$)
contribute to an $\ve$-trace identity, we express the fact in the following table:
\[
\begin{tabular}{|c|c|c|}
\hline
$G$ & $H_1$ & $H_2$\\
\hline
$\{\Pi\}_S$ & $\{\pi_1\}_S$ & $\{\pi_2\}_S$\\
\hline
\end{tabular}
\]
Sometimes we drop the subscript $S$ for brevity.

In the next section, we examine the correspondence between global data and automorphic
representations.
In addition, we demonstrate how the global data of the $\ve$-endoscopic groups lift to those of $\mb{G}$.
\section{Frobenius-Hecke Classes}\label{FHclass}
\indexi{Frobenius-Hecke class|(}%
Let $V$ be the set of places of $F$.  
Let $\Vr$ be the set of places of $F$ which are unramified in $E$.
For each finite place $v$ of $F$, fix a uniformizer $\vp_v$ of $F_v$.  We sometimes drop the subscript
$v$ when the context is clear.  If $v \in \Vr$ does not split in $E$, let $\vp_v$ also be the
uniformizer of $E_v$.  If a place $v$ of $F$ splits into two distinct places $v_1, v_2$ of $E$,
then $E_{v_i} = F_v$ ($i = 1, 2$), and we let $\vp_v$ be the uniformizer of $E_{v_i}$.
\subsection{Frobenius-Hecke Classes for $\mb{G}$}
We now describe the global data for parabolically induced representations of
$\mb{G}(\Af)$.  The results are well known (see \cite{F1}) and we state them without proof.
\begin{itemize}
\item
Let $\alpha, \beta, \mu$ be characters of $\idc{F}$.
Let $\mb{T}$ be the maximal diagonal torus of $\mb{G}$.  Let $\mb{P}_0$ be the
minimal upper triangular parabolic subgroup of $\mb{G}$ containing $\mb{T}$.
Define a representation of $\mb{T}(\Af)$ as follows:
\[
\alpha\otimes\beta\otimes\mu : \diag(a, b, \lambda/b, \lambda/a)
\mapsto \alpha(a)\beta(b)\mu(\lambda).
\]
We extend $\alpha\otimes\beta\otimes\mu$ to a representation of $\mb{P}_0(\Af)$
by setting it to be $1$ on the unipotent component of $\mb{P}_0(\Af)$.

Let $\alpha \times \beta \rtimes \mu$ denote the (normalizedly) parabolically induced
representation \dt $I_{\mb{P_0}(\Af)}^{\mb{G}(\Af)}(\alpha\otimes\beta\otimes\mu)$.
Let $S$ be a finite set of places such that $\alpha_v, \beta_v, \mu_v$ are unramified
for all $v \notin S$.  Then,
\[
\mc{C}(\alpha\times\beta\rtimes\mu, S, \hat{G})
= \left[\dmfour{\alpha_v\beta_v\mu_v}{\alpha_v\mu_v}{\beta_v\mu_v}{\mu_v}\right]_{v \notin S}.
\]
Here, for any unramified character $\chi_v$ of $F_v^\times$, we also let
$\chi_v$ denote the value of the character at $\vp_v$.
\item
Let $e = \diag(1, -1)$.
Let $\mb{P}_\alpha$ be the upper triangular Siegel parabolic subgroup of $\mb{G}$.
Its Levi component is
\[
\mb{M}_\alpha := \left\{\lp\lsm g &\\
&\frac{\lambda}{\det g} e g e\rsm\rp:
\lambda \in \mbb{G}_m,\; g \in \GL(2)\right\}.
\]
Let $\pi$ be an automorphic representation
of $\GL(2, \Af)$.  Let $\mu$ be a character of $\idc{F}$.
Define a representation of $\mb{M}_\alpha(\Af)$ as follows:
\[
\pi\otimes\mu : \lp\lsm g &\\
&\frac{\lambda}{\det g} e g e\rsm\rp
\mapsto \mu(\lambda)\pi(g).
\]
The representation $\pi\otimes\mu$ extends to a representation of $\mb{P}_\alpha(\Af)$.
Let $\pi\rtimes \mu$ denote the normalizedly induced representation
$I_{\mb{P}_{\alpha}(\Af)}^{\mb{G}(\Af)}(\pi\otimes\mu)$.
Let $S$ be a finite set of places such that $\mu_v, \pi_v$ are unramified
for all $v \notin S$.
Let $\displaystyle \mc{C}(\pi, S, \GL(2, \mbb{C})) = \left[ c(\pi_v)\right]_{v \notin S}$ be
the global datum in $\GL(2, \mbb{C})$ which parametrizes $\pi$.
Then,
\[
\mc{C}(\pi\rtimes\mu, S, \hat{G})
= \left[ \lp\lsm \mu_v\omega_{\pi_v}&&\\ &\mu_vc(\pi_v)&\\&& \mu_v\rsm\rp\right]_{v\notin S}.
\]
Here, $\omega_{\pi_v}$ is the central character of $\pi_v$.
\item
Let $\mb{P}_\beta$ be the upper triangular Heisenberg parabolic subgroup of $\mb{G}$.
Its Levi component is
\[
\mb{M}_\beta := \left\{\lp\lsm a &&\\&g&\\&&\frac{\det g}{a}\rsm\rp:
a \in \mbb{G}_m,\; g \in \GL(2)\right\}.
\]
Let $\mu$ be a character of $\idc{F}$.
Let $\pi$ be an automorphic $\GL(2, \Af)$-module.  
Define a representation of $\mb{M}_\beta(\Af)$ as follows:
\[
\mu\otimes\pi : \lp\lsm a&&\\&g&\\&&\frac{\det g}{a}\rsm\rp
\mapsto \mu(a)\pi(g).
\]
The representation $\pi\otimes\mu$ extends to a representation of $\mb{P}_\beta(\Af)$.
Let $\mu\rtimes \pi$ denote the normalizedly induced representation
$I_{\mb{P}_{\beta}(\Af)}^{\mb{G}(\Af)}(\pi\otimes\mu)$.
Let $S$ be a finite set of places such that $\mu_v, \pi_v$ are unramified for all $v \notin S$.
Let $\displaystyle \mc{C}(\pi, S, \GL(2, \mbb{C})) = \left[ c(\pi_v)\right]_{v \notin S}$ 
be the global datum parametrizing $\pi$.  
Suppose for each $v \notin S$ the conjugacy class $c(\pi_v)$ in $\GL(2, \CC)$ is represented by
a diagonal matrix $t(\pi_v) = \diag(c_1(\pi_v), c_2(\pi_v)) \in \GL(2, \CC)$.
\indexs{t@$t(\;)$}%
Then,
$\mc{C}(\pi\rtimes\mu, S, \hat{G})$ is represented by
\[
\left[ \lp\lsm \mu_v t(\pi_v)&\\ &t(\pi_v)\rsm\rp\right]_{v\notin S} \in \prod_{v \notin S}\hat{G}.
\]
\end{itemize}
As for the F-H classes parametrizing discrete spectrum automorphic representations
of $\mb{G}(\Af)$, we postpone their description to Section \ref{sec:arthurunstable}.

We now turn to describing the global data parametrizing the automorphic representations
of the $\ve$-endoscopic groups.  The groups are nonsplit over $F$, and as a result
the description of their global data is less straightforward.
As a warmup, we first review the correspondence between F-H classes and
automorphic representations in the case of $\rR_{E/F}\GL(1)$,
the group  obtained from $\GL(1)$ upon restriction of scalars from $E$ to $F$.
\subsection{The Case of $\rR_{E/F} \GL(1)$}
\indexs{R@$\rR_{E/F}$!GL1@$\GL(1)$}%
Recall that $E$ is the quadratic extension of $F$ which corresponds via global class field
theory to $\ve$.  Let $\sigma$ be the generator of $\Gal(E/F)$.
The group
$E^\times$ may be identified with the group of $F$-points of $\mb{H} = \rR_{E/F} \GL(1)$.  
In particular,
$\mb{H}$ is the unique elliptic $\ve$-endoscopic group of $\GL(2)$
(see \cite{K}).  Recall that endoscopic groups are by definition quasi-split.
The $L$-group of $\mb{H}$ is 
${}^L H = \lp\mathbb{C}^\times \times \mathbb{C}^\times\rp \rtimes \Gal(E/F)$,
and the Galois action is given by $\sigma(a, b) = (b, a)$ for $(a, b) \in \mathbb{C}^\times \times
\mathbb{C}^\times$.  

If a place $v \in \Vr$ does not split in $E$, then $\mb{H}(F_v) = E_v^\times$.  
If $v$ splits into two places $v_1, v_2$ of $E$, then
$E_{v_1} = E_{v_2} = F_v$, and $\mb{H}(F_v) = F_v^\times \times F_v^\times$.

Let $\chi$ be a character of $E^\times\bs\Ae^\times$.  Consider $\chi$ as
an automorphic representation of $\mb{H}(\Af)$.  We would like to compute the
F-H classes in ${}^L H$ which parametrize $\chi$.
The representation $\chi_v$ is unramified for almost all $v \in \Vr$.
This statement is equivalent to the following:
\begin{enumerate}
\item
If $v$ does not split in $E$, then $\chi_v = \chi'_v\circ N_{E/F}$
for some unramified character $\chi'_v$ of $F_v^\times$.
Here, to simplify the notation, we let $\N_{E/F}$ denote the norm mapping $\N_{E_v/F_v}:E_v \rightarrow F_v$.
\indexi{norm mapping}%
\item
If $v$ splits into distinct places $v_1, v_2$ of $E$, then $\chi_v = \chi_{v_1} \otimes \chi_{v_2}$,
where $\chi_{v_i}$ ($i = 1, 2$) is an unramified character of $E_{v_i}^\times = F_v^\times$.  
\end{enumerate}

Consider case 1.
To maintain consistency with the expression of its $L$-group, we write: 
$\mb{H}(\bar{F}) = \{(a, b) \in \GL(1, \bar{F})\times\GL(1, \bar{F})\}$
and $\mb{H}(F_v) = \{(x, \sigma x):x \in E_v^\times\}$.
An unramified character $\chi_v$ of $\mb{H}(F_v)$ has the form
\[
\chi_v = \mu_1\otimes\mu_2 : (x, \sigma x) \mapsto \mu_1(x)\mu_2(\sigma x),
\]
where $\mu_{1,v}, \mu_{2,v}$ are unramified characters of $E_v^\times$.
Let $\vp_v \in F_v$ be the uniformizer of both $F_v$ and $E_v$; this is possible because 
$v$ is unramified in $E$.
The F-H class in ${}^L H$ parametrizing $\chi_v$ is represented by
(see \cite{B}):
\[
(\mu_{1}(\vp_v),\mu_{2}(\sigma\vp_v))\rtimes \sigma = 
(\mu_{1}(\vp_v),\mu_{2}(\vp_v))\rtimes \sigma.
\]  
The element $(\mu_{1}(\vp_v),\mu_{2}(\vp_v))\rtimes \sigma$ is conjugate in ${}^L H$ to
\[
(\lambda \mu_{1}(\vp_v),\lambda^{-1}\mu_{2}(\vp_v))\rtimes \sigma
\]
for any  
$\lambda \in \mathbb{C}^\times$.  Letting
$\lambda = \mu_{2}(\vp_v)$, we see that the F-H class is represented by 
\[
t(\chi_v) = (\mu_{1}\mu_{2}(\vp_v), 1) \rtimes \sigma = (\chi_v(\vp_v), 1)\rtimes \sigma.
\]
In particular, the F-H class parametrizing
$\chi_v$ depends only on the restriction of the character
to the maximal $F_v$-split torus $T_d = \{(x, x) : x \in F_v^\times\}$ of $\mb{H}(F_v)$.  

In case 2, where $v$ of $F$ splits into $v_1, v_2$ of $E$, the group $\mb{H}$ is split over
$F_v$, and the decomposition group of $v$ in $\Gal(E/F)$ is trivial.
The group of $F_v$-points of $\mb{H}$ is $\mb{H}(F_v) = F_v^\times \times F_v^\times$.
If $\chi_v = \chi_{v_1}\otimes\chi_{v_2}$ for some unramified characters $\chi_{v_1},\chi_{v_2}$ of 
$F_v^\times$, then the F-H class of $\chi_v$ is represented by
\[
t\lp\chi_v\rp = \lp \chi_{v_1}(\vp_v), \chi_{v_2}(\vp_v) \rp \rtimes 1.
\]

\subsection{Classes in ${}^L H_1$}
Recall from Chapter \ref{chap:endoscopy} that ${}^L H_1 = \hat{H_1} \rtimes \Gal(E/F)$, where
\[
\hat{H_1} = \Big[\GL(2, \mathbb{C}) \times \mathbb{C}^\times \times \mathbb{C}^\times\Big]/
\{(\diag(z,z), z^{-1}, z^{-1}) : z \in \mathbb{C}^\times\}.
\]
Let $\tau$ be an automorphic representation of $\GL(2, \Af)$.  Let $\chi$ be a character of $\idc{E}$. 
Let $S$ be a finite set of places of $F$ containing $V - \Vr$ 
such that $\tau_v$ is unramified for all $v \notin S$, and
$\chi_w$ is unramified whenever the place $w$ of $E$ lies above a place $v \notin S$.

Suppose $\tau$ is parametrized by the
the global datum \[\mc{C}\lp \tau, S, \GL(2, \CC)\rp = \left[ c(\tau_v)\right]_{v \notin S},\]
where, for $v \notin S$,
$c(\tau_v)$ is a conjugacy class in $\GL(2, \CC)$ represented by some diagonal matrix
$t(\tau_v) = \dmtwo{c_{1,v}}{c_{2,v}}$.


If a place $v \in \Vr$ remains prime in $E$, let $v$ also denote the (unique) place of $E$ 
which lies above $v$.
If $\chi_v$ is unramified, then $\chi_v = \chi'_v\circ\N_{E/F}$
for some unramified character $\chi'_v$ of $F_v^\times$.
If $v \in \Vr$ splits into two distinct places $v_1, v_2$ of $E$,
then $E\otimes_F F_v = F_v \oplus F_v = E_{v_1} \oplus E_{v_2}$.  At such a place $v$, 
let $\chi_v$ denote the character $\chi_{v_1}\otimes\chi_{v_2}$ of 
$E_{v_1}^\times\oplus E_{v_2}^\times$.
 

\begin{claim}
The local component at $v \notin S$ of the global datum 
$\mc{C}\lp\tau\otimes_1 \chi, S, {}^L H_1\rp$ \dt parametrizing the global packet $\tau\otimes_1\chi$ of 
$\mb{H_1}(\Af)$ is as follows:
\[
\mc{C}\lp\tau\otimes_1 \chi, S, {}^L H_1\rp_v
\ni \begin{cases}
\lp t(\tau_v), 1, \chi_v'\rp\rtimes\sigma &\text{if } v \text{ is prime in } E,\\
\lp t(\tau_v), \chi_{v_1}, \chi_{v_2}\rp\rtimes 1 &\text{if } v \text{ splits in } E.
\end{cases}
\]
\\{\rm
Here, for an unramified character $\eta_v$ of  $F_v^\times$, we also let $\eta_v$ denote the value of
the character at $\vp_v$.
}
\end{claim}
\begin{proof}
At a place $v \notin S$ which does not split in $E$, we write $\mb{H_1}(F_v)$ in the following form:
\[
\{ (g, x, \sigma x) \in \GL(2, F_v) \times E_v^\times\times E_v^\times : \det g = \N_{E/F} x\}.
\]
The character $\chi_v$ has the form $\chi_v(x) = \mu_{1,v}(x)\mu_{2,v}(\sigma x)$ for all 
$x \in E_v^\times$, where
$\mu_{1,v}, \mu_{2,v}$ are unramified character of $E_v^\times$.
The F-H class in ${}^L H_1$ parametrizing the unramified representation $\tau_v\otimes_1 \chi_v$ 
is represented by
$
\lp t(\tau_v), \mu_{1, v}, \mu_{2,v}\rp \rtimes \sigma,
$
which is conjugate in ${}^L H_1$ to 
\[
\lp t(\tau_v), \lambda^{-1}\mu_{1,v}, \lambda\mu_{2,v}\rp \rtimes \sigma
\]
for any $\lambda \in \mathbb{C}^\times$.  Letting $\lambda = \mu_1(\vp_v)$, we see that
the F-H class of $\tau_v\otimes_1\chi_v$ is represented by
\[
\lp t(\tau_v), 1, \chi_v\rp \rtimes \sigma.
\]

The proof for the split case is similar, and we skip it.
\end{proof}
Let $e = \lp\lsm 1 & \\& -1\rsm\rp \in \GL(2, \CC)$.
Recall from Chapter \ref{chap:endoscopy} that the $L$-group
embedding $\xi_{1} : {}^L H_1 \rightarrow \hat{G}$ sends
$1 \rtimes \sigma$ to $\lp\lsm & e\\e &\rsm\rp$  and
$(g, a, b) \in \hat{H}_1$ to $d(g, a, b) := \lp\lsm  a g & \\ & b ege\rsm\rp$.
\begin{corollary}\label{corollary:FHliftH1}
The global datum $\mc{C}\lp\tau\otimes_1 \chi, S, {}^L H_1\rp$ parametrizing
$\tau\otimes_1 \chi$ lifts to
the global datum  $\mc{C}\lp S, \GSp(2, \CC)\rp$ described as follows:
\begin{itemize}
\item
At a place $v \notin S$ which is prime in $E$,
the local component at $v$ of the global datum $\mc{C}\lp S, \GSp(2, \CC)\rp$ is represented by
\[
\lp\lsm  & t(\tau_v) e\\ {\chi'}^2 t(\tau_v) e&\rsm\rp
\sim
\lp\lsm \chi_v'c_{1,v} &&& \\ & -\chi_v'c_{2,v}&& \\ &&-\chi_v' c_{1,v} & 
\\&&&\chi_v' c_{2,v}\rsm\rp \in \GSp(2, \CC)
\]
$($Recall: $t(\tau_v) = \lp\lsm c_{1,v} & \\ & c_{2,v}\rsm\rp \in \GL(2, \CC)$$)$.
\item
At a place $v$ which splits into two places $v_1, v_2$ of $E$, 
the local component at $v$ of $\mc{C}\lp S, \GSp(2, \CC)\rp$ is represented by 
\[
\lp\lsm
\chi_{v_1}t(\tau_v) & \\
& \chi_{v_2} t(\tau_v)
\rsm\rp \in \GSp(2, \CC).
\]
\end{itemize}
\end{corollary}
Note that a priori all we know about $\mc{C}(S, \GSp(2, \CC))$ is that it is an infinite-tuple
of conjugacy classes in $\GSp(2, \CC)$.  At this stage we make no prediction on whether
$\mc{C}(S, \GSp(\CC))$ is the global datum of any (quasi-)packet of $\GSp(2, \Af)$.
\subsection{Classes in ${}^L H_2$}
Recall from Chapter \ref{chap:endoscopy} that ${}^L H_2 = \hat{H_2} \rtimes \Gal(E/F)$, where
\[
\hat{H_2} = 
\left\{
(g_1, g_2, \lambda) \in \GL(2,\CC)^2 \times \CC^\times : \det g_1 = \det g_2 = \lambda
\right\},
\]
and the generator $\sigma$ of $\Gal(E/F)$ swaps the two $\GL(2, \CC)$-factors.

Let $\pi_2 = \pi\otimes_2\mu$ be an automorphic representation of $\mb{H_2}(\Af)$,
where $\pi$ is an automorphic representation of $\GL(2, \Ae)$ and $\mu$ is a character 
of $\idc{F}$ such that
the central character $\omega_\pi$ of $\pi$ is equal to $\mu\circ\N_{E/F}$.  

If $v \in \Vr$ is prime in $E$, then $\pi_v$ is a representation of
$\GL(2, E_v)$ which satisfies $\omega_{\pi_v} = \mu_v\circ\N_{E/F}$.
If $v$ splits into
two distinct places $v_1, v_2$ of $E$, then
\[
\mb{H_2}(F_v) = 
\lp\GL(2, F_v) \times \GL(2, F_v) \times F_v^\times\rp/\{(a I_2, b I_2, (ab)^{-1}):a,b\in F_v^\times\}, 
\]
and $\lp\pi\otimes_2\mu\rp_v$ takes the form $\pi_{v_1} \otimes \pi_{v_2} \otimes_2 \mu_v$, 
where $\pi_{v_i}$ ($i = 1, 2$) is a representation of $\GL(2, E_{v_i}) = \GL(2, F_v)$, and
$\omega_{\pi_{v_1}} = \omega_{\pi_{v_2}} = \mu_v$.

Identify $\GL(2, \Ae)$ with the group of $\Af$-points of $\rR_{E/F}\GL(2)$, the $F$-group obtained from
$\GL(2)$ via the restriction of scalars from $E$ to $F$.
The $L$-group (actually, a finite Galois form of the $L$-group) of $\rR_{E/F}\GL(2)$ 
is the semidirect product 
${}^L \rR_{E/F}\GL(2) = \GL(2, \CC) \times \GL(2, \CC) \rtimes \Gal(E/F)$, 
where the action of $\sigma \in \Gal(E/F)$ swaps the two $\GL(2, \CC)$-factors.

Let $S$ be a finite set of places of $F$ containing $V - \Vr$ such that
$\pi_{2,v}$ is unramified for all $v \notin S$.  
Consider $\pi$ as an automorphic representation of $\lp\rR_{E/F}\GL(2)\rp(\Af)$.
Then, $\pi_{v}$ is unramified for all $v \notin S$.
Suppose the F-H classes $\{c(\pi_v)\}_{v\notin S}$ in ${}^L \rR_{E/F}\GL(2)$ 
parametrizing $\pi$ are as follows:
\[
c(\pi_v) \ni
\begin{cases}
\lp
\lp\lsm
t_{1}(\pi_v) &\\& t_{2}(\pi_v)
\rsm\rp
,
\lp\lsm
1 &\\& 1
\rsm\rp 
\rp\rtimes \sigma
&
\text{if } v \text{ is prime in } E,\\
\lp d(\pi_{v_1}), d(\pi_{v_2})\rp
\rtimes 1
&
\text{if } v \text{ splits into } v_1,v_2 \text{ in } E,
\end{cases}
\]
where $t_{1}(\pi_v)$, $t_{2}(\pi_v)$ are elements in $\CC^\times$, and
$d(\pi_{v_i})$ ($i = 1, 2$) is a diagonal matrix in $\GL(2, \CC)$ whose
conjugacy class parametrizes the unramified representation $\pi_{v_i}$ of \dt $\GL(2, E_{v_i})$.
\begin{claim}
The global datum $\mc{C}(\pi_2, S, {}^L H_2)$ parametrizing
the automorphic representation $\pi_2 = \pi\otimes_2 \mu$ of $\mb{H_2}(\Af)$ is as follows:
\begin{itemize}
\item
If $v \notin S$ is prime in $E$, then the local component $\mc{C}(\pi_2, S, {}^L H_2)_v$
of the global datum is represented by
\[
 \lp \lp\lsm \mu_v^{1/2}&\\& \mu_v^{1/2}\rsm\rp,
\lp\lsm\mu_v^{-1/2}t_{1}(\pi_v)&\\
& \mu^{-1/2}t_{2}(\pi_v)\rsm\rp, 
\mu_v\rp\rtimes\sigma.
\]
\item
If $v \notin S$ splits into $v_1, v_2$ in $E$, then $\mc{C}(\pi_2, S, {}^L H_2)_v$ is
represented by
\[
\lp d(\pi_{v_1}) , d(\pi_{v_2}),
\mu_v
\rp \rtimes 1.
\]
\end{itemize}
\end{claim}
\begin{proof}
First, we deal with the case where $v$ is prime in $E$.
Write $\mb{H_2}(F_v)$ in the form
\[
\{(g, \sigma g, x) \in \GL(2, E_v)^2 \times F_v^\times\}
/
\{(z I_2, z I_2, \N_{E/F}z^{-1}) : z \in E_v^\times\}.
\]
If $\pi_{2,v} = \pi_v\otimes_2\mu_v$ is unramified, then $\pi_v$ and $\mu_v$
are unramified.
In particular, $\pi_{v}$ is the induced representation $I(\alpha, \beta)$
of $\GL(2, E_v)$
for some unramified characters $\alpha, \beta$ of $E_v^\times$.
The condition $\omega_\pi = \mu\circ\N_{E/F}$ implies that $\omega_{\pi_v}|_{F_v^\times} = 
\lp\alpha\beta\rp|_{F_v^\times} = \mu_v^2$.  The F-H class $c(\pi_v)$
in ${}^L \rR_{E/F}\GL(2)$ parametrizing $\pi$ is represented by
\[
t(\pi_v) = \lp 
\lp\lsm \alpha & \\ & \beta\rsm\rp,
\lp\lsm 1&\\&1\rsm\rp
\rp \rtimes
\sigma.
\]
In particular, we may take $t_1(\pi_v) = \alpha(\vp_v)$ and $t_2(\pi_v) = \beta(\vp_v)$.
Suppose
the F-H class $c(\pi_{2, v})$ in ${}^L H_2$ parametrizing $\pi_{2,v}$ is represented by
\[
\lp \diag(\alpha_1, \beta_1),\diag(\alpha_2, \beta_2), \mu_v\rp\rtimes\sigma 
\in {}^L H_2,
\]
where $\alpha_i, \beta_i$ ($i = 1, 2$) are complex numbers.
The following conditions must be satisfied:
\begin{itemize}
\item
$\alpha_1\alpha_2 = \alpha(\vp_v)$,
\item
$\beta_1\beta_2 = \beta(\vp_v)$,
\item
$\alpha_1\alpha_2\beta_1\beta_2 = \mu_v^2(\vp_v)$.
\end{itemize}
The element
$
\lp \diag(\alpha_1, \beta_1),\diag(\alpha_2, \beta_2), \mu_v\rp\rtimes\sigma
$
is conjugate in ${}^L H_2$ to
\[
\lp \lp\lsm \lambda_1\alpha_1&\\& \lambda_2\beta_1\rsm\rp,
\lp\lsm \lambda_1^{-1}\alpha_2&\\& \lambda_2^{-1}\beta_2\rsm\rp, \mu_v\rp\rtimes\sigma
\]
for any $\lambda_1, \lambda_2 \in \CC^\times$.  
Fix a choice of $\mu_v(\vp_v)^{1/2}$.
Letting $\lambda_1 = \mu_v(\vp_v)^{1/2}\alpha_1^{-1}$ and $\lambda_2 = \mu_v(\vp_v)^{1/2}\beta_1^{-1}$,
we see that $c(\pi_{2,v})$ is represented by
\begin{multline*}
\lp \lp\lsm \mu_v^{1/2}&\\& \mu_v^{1/2}\rsm\rp,
\lp\lsm\mu_v^{-1/2}\alpha_1\alpha_2&\\& \mu_v^{-1/2}\beta_1\beta_2\rsm\rp, 
\mu_v\rp\rtimes\sigma\\
=
\lp \lp\lsm \mu_v^{1/2}&\\& \mu_v^{1/2}\rsm\rp,
\lp\lsm\mu_v^{-1/2}\alpha&\\
& \mu^{-1/2}\beta\rsm\rp, 
\mu_v\rp\rtimes\sigma.
\end{multline*}
The nonsplit case of the claim follows.
Notice that the conjugacy class of the above element does not change if
$\mu_v(\vp_v)^{1/2}$ is replaced with $-\mu_v(\vp_v)^{1/2}$.

The proof for the split case is similar, and we skip it.
\end{proof}
Recall the notation $[(a_{ij}), (b_{ij})] : = 
\lp\lsm a_{11} &&& a_{12}\\& b_{11}&b_{12}&\\&b_{21}&b_{22}&\\a_{21}&&&a_{22}\rsm\rp$.
The $L$-group embedding $\xi_{2,v} : {}^L H_2 \rightarrow \hat{G}$ sends
$1 \rtimes \sigma$ to $\lp\lsm &1&&\\1&&&\\&&&1\\&&1&\rsm\rp$ and
$(g_1, g_2, c)\rtimes 1$ to $[g_1, g_2]$.
\begin{corollary}\label{corollary:FHliftH2}
The global datum $\mc{C}\lp\pi\otimes_2 \mu, S, {}^L H_2\rp$ parametrizing $\pi\otimes_2 \mu$
lifts to the global datum  $\mc{C}\lp S, \GSp(2, \CC)\rp$ described as follows:
\begin{itemize}
\item
At a place $v \notin S$ which is prime in $E$,
$\mc{C}\!\lp S, \GSp(2, \CC)\rp_v$ is represented by
\[
\lp\lsm &\mu_v^{1/2}&& \\  \mu_{v}^{-1/2}t_{1}(\pi_v) &&& \\ 
&&& \mu_v^{-1/2}t_{2}(\pi_v) \\&&\mu_v^{1/2}&\rsm\rp.
\]
\item
At a place $v$ which splits into two places $v_1, v_2$ of $E$, 
$\mc{C}\!\lp S, \GSp(2, \CC)\rp_v$ is represented by 
\[
[d(\pi_{v_1}), d(\pi_{v_2})].
\]
\end{itemize}
{\rm
As remarked earlier, at this stage we make no prediction on whether 
the global datum $\mc{C}(S, \GSp(2, \CC))$ corresponds to a global (quasi-)packet of $ \GSp(2, \Af)$.
}
\end{corollary}
\indexi{Frobenius-Hecke class|)}%

\section{Packets}\label{sec:globalcomparison}
\indexi{packet|(}%
In this section, we describe the multi-packets for each $\ve$-endoscopic group; in other words,
we describe the inequivalent packets which simultaneously contribute to an $\ve$-trace identity.



\subsection{(Quasi-)Packets of $\GSp(2, \Af)$}\label{sec:arthurunstable}
\indexi{packet!quasi-|(}
For any place $v$, we let $\nu_v$ be the normalized absolute value function on $F_v$.
\indexs{nu@$\nu$}%
In other words, if $v$ is finite, $\nu_v(x) := q_v^{-\ord x}$ for all $x \in F_v$, 
where $q_v$ is the cardinality of the residue field of $F_v$, and $\ord x$ is the $p$-adic valuation of $x$.
We let $\nu = \otimes_v \nu_v$ be the normalized absolute value function on $\Af$.
For any finite place $v$ of $F$ and unramified representation $\tau_v$ of $\GL(2, F_v)$, let $t(\tau_v)$ be
a diagonal matrix in $\GL(2, \CC)$ whose conjugacy class parametrizes $\tau_v$.

Let $1_2$ denote the trivial representation $1_{\GL(2, \Af)}$ of $\GL(2, \Af)$.
\indexs{1@$1_{\GL(2,\Af)}$}%
For any place $v$, let $1_{2, v}$ and $\St{v}$ denote the trivial and Steinberg representations of 
\indexs{S@${\rm St}_2$}%
$\GL(2, F_v)$, respectively.
\subsubsection{Unstable Packets}
\indexi{packet!unstable|(}%
From \cite{A}, the unstable (quasi-)packets of $\GSp(2, \Af)$ fall into three types.
They are described as follows:

{\bf 1.}
Let $\tau_1, \tau_2$ be distinct cuspidal automorphic representations of the group $\GL(2, \Af)$ such that
$\omega_{\tau_1} = \omega_{\tau_2}$.   Suppose $\tau_{1,v}$ and $\tau_{2,v}$ are 
unramified for all $v$ outside of some finite set of places $S$.
The packet $[\tau_1, \tau_2]$ is parametrized by the global datum (product of conjugacy classes):
\indexs{b@$[\;,\;]$}%
\[
\mc{C}\lp [\tau_1,\tau_2], S, \hat{G}\rp \ni \Big[\left[t(\tau_{1,v}), t(\tau_{2,v})\right]\Big]_{v \notin S}.
\]
In particular, $[\tau_1,\tau_2]$ lifts via the $L$-group
embedding $\GSp(2, \mbb{C}) \hookrightarrow \GL(4, \mbb{C})$ to the induced representation
$I_{(2, 2)}(\tau_1, \tau_2)$ of $\GL(4, \Af)$.

If the central characters of $\tau_1$ and $\tau_2$ are trivial, then we know from 
\cite[V. 10]{F1} that $[\tau_1, \tau_2]$ is the restricted tensor product
$\{\Pi\} = \otimes_v \{\Pi_v\}$.  Here, for all $v$ where both $\tau_1$ and $\tau_2$ are square integrable,
the local packet $\{\Pi_v\}$ consists of two square integrable representations $\Pi_v^+, \Pi_v^-$.
\indexs{Pi@$\Pi^+$, $\Pi^-$}%
At the rest of the places, $\{\Pi_v\}$ consists of a single representation, which we denote by
$\Pi_v^+$.

The multiplicity formula for $\Pi' \in \{\Pi\}$ is $m(\Pi') = \frac{1}{2}(1 + (-1)^{n(\Pi')})$, where
\indexi{multiplicity formula}%
$n(\Pi')$ is the number of places $v$ for which $\Pi'_v = \Pi_v^-$.
\\
{\sc Remark.}  In \cite{A}, $[\tau_1, \tau_2]$ is denoted by $(\tau_1\boxtimes 1)\boxplus(\tau_2\boxtimes 1)$
and is called a \emph{Yoshida type} packet.

{\bf 2.}  Let $\tau$ be a cuspidal automorphic representation of $\GL(2, \Af)$,
$\mu$ be a unitary character of $\idc{F}$, such that $\omega_\tau = \mu^2$.
Suppose $\tau_v, \mu_v$ are unramified for all $v$ outside of some finite set of places $S$.
The quasi-packet $[\tau, \mu 1_2]$ is parametrized by the global datum
\begin{multline*}
  \mc{C}\lp[\tau, \mu 1_2], S, \hat{G}\rp
  \ni\Big[ [t(\tau_v), t(\mu_v \one{v})]\Big]_{v \notin S}\\
  =\Big[\blockdiag\lp \mu_v q_v^{1/2}, t(\tau_v), \mu_vq_v^{-1/2}\rp\Big]_{v \notin S}.
\end{multline*}
Here, $\mu_v := \mu_v(\vp_v)$.

In particular, the quasi-packet $[\tau, \mu 1_2]$ lifts to the induced representation $I_{(2, 2)}(\tau, \mu 1_2)$
of \dt $\GL(4, \Af)$.

If $\mu^2 = \omega_\tau = 1$, then from \cite[V. 10]{F1}
we know that $[\tau, \mu 1_2]$ is the global quasi-packet $\{\Pi\} = \otimes_v \{\Pi_v\}$,
where $\{\Pi_v\} = \{\Pi_v^\times, \Pi_v^-\}$ for all $v$ where $\tau_v$ is square integrable.  Here,
\indexs{Pi@$\Pi^\times$, $\Pi^-$}%
$\Pi_v^\times$ is the nontempered quotient $L(\nu_v^{1/2}\mu_v\tau_v, \mu_v\nu_v^{-1/2})$ of $\nu_v^{1/2}\mu_v\tau_v\rtimes \mu_v\nu_v^{-1/2}$
\indexs{L@$L({\rm repn., char.})$}%
and $\Pi_v^-$ is the cuspidal representation denoted by $\delta^-(\nu_v^{1/2}\mu_v\tau_v, \mu_v\nu_v^{-1/2})$ in \cite[V. 8]{F1}.
\indexs{delta@$\delta^-({\rm repn., char.})$}%
At each place $v$ where $\tau_v$ is the induced representation $I(\eta_v, \eta_v^{-1})$ for some character $\eta_v$ of $F_v^\times$,
$\{\Pi_v\}$ is the singleton consisting of the irreducible induced representation 
$\Pi_v^\times := \eta_v\mu_v \one{v}\rtimes \eta_v^{-1}$.

The multiplicity formula for $\Pi' \in \{\Pi\}$ is 
\[m(\Pi') = \frac{1}{2}\lp 1 + \ve(\mu\tau, 1/2) (-1)^{n(\Pi')}\rp,\]
\indexi{multiplicity formula}%
where $n(\Pi')$ is the number of places $v$ for which $\Pi'_v = \Pi_v^-$, and $\ve\!\lp\mu\tau, \frac{1}{2}\rp$ is 
the value at $1/2$ of the epsilon factor in the functional equation of the $L$-function of $\mu\tau$.
\\{\sc Remark.}  In \cite{A}, $[\tau, \mu 1_2]$ is denoted by $(\mu\boxtimes\nu(2))\boxplus(\tau\boxtimes 1)$
and is called a \emph{Saito, Kurokawa type} packet.
\indexi{packet!Saito, Kurokawa}%

{\bf 3.} 
Let $\mu_1, \mu_2$ be distinct unitary characters of  $\idc{F}$ such that $\mu_1^2 = \mu_2^2$.
Suppose $\mu_{i,v}$ ($i = 1, 2$)  is unramified for all $v$ lying outside of some finite set of places $S$.
The quasi-packet $[\mu_1 1_2, \mu_2 1_2]$ is parametrized by the global datum
\begin{multline*}
  \mc{C}\lp[\mu_1 1_2, \mu_2 1_2], S, \hat{G}\rp
  \ni
  \Big[\left[t\lp\mu_{1,v} \one{v}\rp, t\lp \mu_{2,v} \one{v}\rp\right]\Big]_{v \notin S}\\
  =
  \Big[\diag\lp\mu_{1,v}q_v^{1/2},\mu_{2,v}q_v^{1/2},\mu_{2,v}q_v^{-1/2},\mu_{1,v}q_v^{-1/2}\rp\Big]_{v \notin S}.
\end{multline*}
In particular, the quasi-packet $[\mu_1 1_2, \mu_2 1_2]$ lifts to the representation 
\dt $I_{(2, 2)}(\mu_1 1_2,\;\mu_2 1_2)$ of $\GL(4, \Af)$.

If $\mu_1^2 = \mu_2^2 = 1$, let $\zeta$ be the nontrivial quadratic character 
$\frac{\mu_1}{\mu_2}$ of $\idc{F}$.  In this case, we know from \cite[V. 10]{F1} that
$[\mu_1 1_2, \mu_2 1_2]$ is the quasi-packet $\{\Pi\} = \otimes_v \{\Pi_v^\times, \Pi_v^-\}$.
\indexs{Pi@$\Pi^\times$, $\Pi^-$}%
Here, $\Pi_v^\times$ is the nontempered quotient 
\[L(\nu_v\zeta_v, \zeta_v\rtimes \nu_v^{-1/2}\mu_{2,v})\]
\indexs{L@$L({\rm char., repn.})$}%
of the induced representation
$I_v = \nu_v\zeta_v\times\zeta_v \rtimes \nu_v^{-1/2}\mu_{2,v}$ of $\GSp(2, F_v)$, and
\[\Pi_v^- = \begin{cases}
  \delta^-(\zeta_v\nu_v^{1/2}\St{v}, \zeta_v\mu_{2,v}\nu_v^{-1/2})&\text{ if } \zeta_v \neq 1,\\
  L(\zeta_v\nu_v^{1/2}\St{v}, \zeta_v\mu_{2,v}\nu_v^{-1/2})&\text{ if } \zeta_v = 1,      
\end{cases}
\]
where $L(\zeta_v\nu_v^{1/2}\St{v}, \zeta_v\mu_{2,v}\nu_v^{-1/2})$ is the unique nontempered
\indexs{L@$L({\rm repn., char.})$}%
subquotient of the induced representation $\zeta_v\nu_v^{1/2}\St{v}\rtimes \zeta_v\mu_{2,v}\nu_v^{-1/2}$, and
$\delta^-(\zeta_v\nu_v^{1/2}\St{v}, \zeta_v\mu_{2,v}\nu_v^{-1/2})$
\indexs{delta@$\delta^-({\rm repn., char.})$}%
is the cuspidal member of the local packet which contains the unique square integrable subrepresentation
$\delta(\zeta_v\nu_v^{1/2}\St{v}, \zeta_v\mu_{2,v}\nu_v^{-1/2})$  of $I_v$.
\indexs{delta@$\delta({\rm repn., char.})$}%
At almost every place $v$, $\Pi_v^+$ is unramified.  

The multiplicity formula for
$\Pi' \in \{\Pi\}$ is $\frac{1}{2}(1 + (-1)^{n(\Pi')})$, where $n(\Pi')$ is the number
of places $v$ for which $\Pi'_v = \Pi_v^-$.
\\{\sc Remark.}  In \cite{A}, $[\mu_1 1_2, \mu_2 1_2]$ is denoted 
by $(\mu_1\boxtimes\nu(2))\boxplus(\mu_2\boxtimes\nu(2))$ and is called a \emph{Howe, Piatetskii-Shapiro type} packet.
\indexi{packet!unstable|)}%
\subsubsection{A Stable Packet}
\indexi{packet!stable|(}%
Also contained in \cite{A} is a classification of
the stable (quasi-)packets of $\GSp(2, \Af)$.  We now describe one which is of interest to us.
\indexi{Arthur, James}%
Let $\zeta$ be a nontrivial quadratic character of $\idc{F}$.  Let $\tau$ be a $\zeta$-invariant, 
cuspidal automorphic representation of $\GL(2, \Af)$.
Suppose $\zeta_v$, $\tau_v$  are unramified for all $v$ lying outside
of some finite set of places $S$.
The stable quasi-packet $\{L(\nu\zeta, \nu^{-1/2}\tau)\}$ is parametized by the global datum
\[  
\mc{C}\lp\{L(\nu\zeta, \nu^{-1/2}\tau)\}, S, \hat{G}\rp
\ni \Big[\blockdiag\!\lp q_v^{1/2}\zeta_v t(\tau_v), q_v^{-1/2} t(\tau_v)\rp\Big]_{v \notin S}.
\]

In particular, $\{L(\nu\zeta, \nu^{-1/2}\tau)\}$ lifts to the Langlands quotient
$J(\nu^{1/2}\tau, \nu^{-1/2}\tau)$ of the induced representation 
\indexi{Langlands quotient}%
$I_{(2, 2)}(\nu^{1/2}\tau, \nu^{-1/2}\tau)$ of $\GL(4, \Af)$.

If $\omega_\tau = \zeta$, then $\{L(\nu\zeta, \nu^{-1/2}\tau)\}$ is the quasi-packet
defined in \cite[V. 10]{F1}.  It is the restricted tensor product $\otimes_v \{L_v\}$, where
$\{L_v\}$ is the local quasi-packet defined as follows:
\begin{itemize}
\item
If $\zeta_v \neq 1$ and $\tau_v$ is cuspidal, $\{L_v\}$ is a singleton consisting of
the unique nontempered quotient $L(\nu_v\zeta_v, \nu_v^{-1/2}\tau_v)$ of
$\nu_v\zeta_v \rtimes \nu_v^{-1/2}\tau_v$.
\item
Suppose $\zeta \neq 1$ and $\tau_v$ is the induced representation $I(\mu_v, \zeta_v \mu_v)$
for some (possibly trivial) quadratic character $\mu_v$ of $F_v^\times$.  Then,
\[
\{L_v\} = \{L(\nu_v\zeta_v, \zeta_v \rtimes \nu_v^{-1/2}\tau_v),\;
\delta^-(\zeta_v\nu_v^{1/2}\St{v}, \zeta_v\mu_{v}\nu_v^{-1/2})\}.
\]
\item
Suppose $\zeta_v = 1$, then by the theory of monomial representations $\tau_v$ must be induced.
Suppose $\tau_v = I(\mu_v, \mu_v^{-1})$, where $\mu_v$ is a unitary character
of $F_v^\times$.
If $\mu_v^2 \neq 1$,
$\{L_v\}$ is the singleton consisting of the irreducible representation
$\mu_v^{-2} \rtimes \mu_v \one{v}$.
If $\mu_v^2 = 1$, $\{L_v\}$ consists of the nontempered subquotients 
$L(\nu_v, \one{v}\rtimes \nu_v^{-1/2}\mu_v)$,
$L(\nu_v^{1/2}\St{v},\nu_v^{-1/2}\mu_v)$ of
$\nu_v \times 1 \rtimes \nu_v^{-1/2}\mu_v$ (see \cite{ST}).
\end{itemize}
Every member of $\{L(\nu\zeta, \nu^{-1/2}\tau)\}$ occurs with multiplicity one in the
discrete spectrum of $\GSp(2)$.
\\{\sc Remark.}  In \cite{A}, $\{L(\nu\zeta, \nu^{-1/2}\tau)\}$ 
is denoted by $\tau\boxtimes\nu(2)$ and is called a \emph{Soudry type} packet.
\indexi{packet!stable|)}%
\subsubsection{$\ve$-Invariant Unstable Packets}
\indexi{packet!unstable|(}\indexi{packet!epsilon-invariant@$\ve$-invariant|(}%
For any character $\chi$ of $\idc{E}$, put ${}^\sigma\chi(x) := \chi(\sigma x)$ for all $x \in \idc{E}$.
\begin{lemma}\label{lemma:listunstable}
  Suppose $\{\Pi\}$ is an $\ve$-invariant, unstable, discrete spectrum global \dts $($quasi-$)$packet of $\GSp(2, \Af)$.  Then,
  \begin{enumerate}
  \item
    $\{\Pi\} \neq [\tau, \mu 1_2]$ for any cuspidal automorphic representation $\tau$ of the group \dt $\GL(2, \Af)$,
    unitary character $\mu$ of $\idc{F}$, such that $\omega_\tau = \mu^2$.
  \item
    If $\{\Pi\} = [\tau_1, \tau_2]$  for
    two distinct cuspidal automorphic representations $\tau_1, \tau_2$ of $\GL(2, \Af)$, 
    then one of the following conditions is satisfied:
    \begin{itemize}
    \item
      The representation $\tau_1$ is not $E$-monomial, and $\tau_2 = \ve\tau_1$.
    \item
      There exist unitary characters $\zeta, \chi$ of $\idc{E}$
      such that $\tau_1, \tau_2$ are equal to the $E$-monomial representations $\pi(\zeta\chi)$, $\pi(\zeta\;{}^\sigma\!\chi)$,
      respectively.
    \end{itemize}
  \item
    If $\{\Pi\} = [\mu_1 1_2, \mu_2 1_2]$, where $\mu_1, \mu_2$ are distinct unitary characters of $\idc{F}$, then
    $\mu_2 = \ve\mu_1$.
  \end{enumerate}
\end{lemma}
\begin{proof}
Suppose $\{\Pi\} = [\tau, \mu 1_2]$ 
for some cuspidal automorphic representation $\tau$ of \dt $\GL(2, \Af)$ and a unitary character $\mu$ of $\idc{F}$ such 
that $\mu^2 = \omega_{\tau}$.  
Let $\Pi$ be an $\ve$-invariant discrete spectrum representation in $\{\Pi\}$.
For a sufficiently large finite set of places $S$ such that $\Pi_v$ is unramified for all $v\notin S$, we have:
\[
\mc{C}\lp\Pi, S, \hat{G}\rp
\ni
\Big[\blockdiag\lp t(\tau_v), t(\mu_v \one{v})\rp\Big]_{v \notin S}.
\]
Since $\Pi$ is $\ve$-invariant,
$\blockdiag\lp t(\tau_v), t(\mu_v \one{v})\rp$ must be conjugate in $\hat{G} = \GSp(2, \CC)$ to
$\blockdiag\lp \ve_vt(\tau_v), \ve_vt(\mu_v \one{v})\rp$ for each $v \notin S$.
In particular, at each place $v \notin S$ such that $\ve_v \neq 1$, 
$\tau_v$ is equivalent to either $\ve_v \tau_v$ or $\ve_v\mu_v \one{v}$.
If $\tau_v \cong \ve_v\tau_v$, then
$\mu_v$ must be equal to $\ve_v\mu_v$, which cannot hold because $\ve_v \neq 1$.  
Hence, $\tau_v$ is equivalent to $\ve_v\mu_v \one{v}$ at infinitely many places $v$.
In particular, $\tau$ should be one dimensional, but by assumption $\tau$ is cuspidal.
The first statement of the lemma follows.

Suppose $\{\Pi\} = [\tau_1, \tau_2]$,
where $\tau_1, \tau_2$ are two distinct cuspidal automorphic representations  of $\GL(2, \Af)$
whose central characters satisfy $\omega_{\tau_1} = \omega_{\tau_2}$.  
Let $\Pi$ be an $\ve$-invariant discrete spectrum representation in $\{\Pi\}$.
For a finite set of places $S$ such that $\Pi_v$ is unramified for all $v \notin S$, we have:
\[
\mc{C}\lp\Pi, S, \hat{G}\rp
\ni
\Big[\left[t(\tau_{1,v}), t(\tau_{2,v})\right]\Big]_{v \notin S}.
\]
View $\GSp(2)$ as a twisted endoscopic group with respect to the triple $\lp\GL(4), \theta, 1\rp$ 
(see \cite[Chap. 2]{KS}), where $\theta$ is the automorphism
$g \mapsto {}^t g ^{-1}$ of $\GL(4)$.
The F-H classes in $\GSp(2, \CC)$ parametrizing $\Pi$, viewed as classes in $\GL(4, \CC)$, parametrize
the induced representation 
$I_{(2,2)}(\tau_1,\tau_2)$ of $\GL(4, \Af)$.

Since $\Pi$ is an $\ve$-invariant representation, the element
$\left[t(\tau_{1,v}), t(\tau_{2,v})\right]$ is conjugate to
$\left[\ve_vt(\tau_{1,v}), \ve_vt(\tau_{2,v})\right]$ in $\GSp(2, \CC)$, and hence in $\GL(4, \CC)$, 
for each $v \notin S$.
Consequently, $I_{(2,2)}(\tau_1, \tau_2)$ is equivalent to $I_{(2,2)}(\ve\tau_1, \ve\tau_2)$, which implies that
the set $\{\tau_1, \tau_2\}$ is equal to $\{\ve\tau_1, \ve\tau_2\}$ (see \cite[Thm. 4.4]{JS1}).




If $\tau_2 \neq \ve\tau_1$, then $\tau_1$ and $\tau_2$ must both be $E$-monomial.  Consequently, there are characters
$\mu_1, \mu_2$ of $\idc{E}$ such that $\tau_1 = \pi(\mu_1), \tau_2 = \pi(\mu_2)$.  
The central character of the monomial representation $\pi(\mu_i)$ ($i = 1, 2$) is $\mu_i|_{\Af^\times}\cdot\ve$; hence,
the condition $\omega_{\tau_1} = \omega_{\tau_2}$ implies that $\mu_1|_{\Af^\times} = \mu_2|_{\Af^\times}$.
Since $\frac{\mu_2}{\mu_1}|_{\Af^\times}$ is trivial, there exists a character
$\chi$ of $\idc{E}$ such that $\frac{\mu_2}{\mu_1} = \frac{{}^\sigma\chi}{\chi}$.  Let
$\zeta = \frac{\mu_1}{\chi}$, then $\tau_1 = \pi(\zeta\chi)$ and 
$\tau_2 = \pi(\zeta\;{}^\sigma\!\chi)$.
The second statement of the lemma follows.

Suppose  $\{\Pi\} = [\mu_1 1_2, \mu_2 1_2]$,
where $\mu_1, \mu_2$ are two distinct unitary characters  of $\idc{F}$ such that $\mu_1^2 = \mu_2^2$.
Let $\Pi$ be an $\ve$-invariant discrete spectrum representation in $\{\Pi\}$.
It lifts to the induced representation $I_{(2,2)}\lp\mu_1 1_2, \mu_2 1_2\rp$ of $\GL(4, \Af)$.
By the same reasoning used in the second case, $\Pi \cong \ve\Pi$ implies that the induced representations
$I_{(2,2)}\lp\mu_1 1_2, \mu_2 1_2\rp$ and
$I_{(2,2)}\lp\ve\mu_1 1_2, \ve\mu_2 1_2\rp$ are equivalent.
Hence, $\mu_2 = \ve\mu_1$, and the final statement of the lemma follows.
\end{proof}
\indexi{packet!quasi-|)}
\indexi{packet!unstable|)}\indexi{packet!epsilon-invariant@$\ve$-invariant|)}%
\subsection{Multi-Packets of $\mb{H_1}(\Af)$}
\indexi{packet!multi-|(}%
Suppose $\tau$ is an automorphic representation of $\GL(2, \Af)$, and $\chi$ is a character
of $\idc{E}$.
Recall from Section \ref{sec:stablespectrum} that we let $\tau\otimes_1 \chi$ denote the
global packet of $\mb{H_1}(\Af)$ consisting of those irreducible constituents of 
$\tau|_{\GL(2, \Af)^E} \otimes_1 \chi$ which are unramified at almost all places.

The following theorem rests on the assumption that a result of Flicker's (\cite{F1}) extends to an
analoguous result in an announced but unpublished work of Arthur's (\cite{A1}).
\begin{thm}\label{thm:contribH1}
Let
$\tau$ be an automorphic representation of $\GL(2, \Af)$.  Let $\chi$ be a character
of $\idc{E}$.  
Suppose the packet
$\tau\otimes_1 \chi$ of $\mb{H_1}(\Af)$ contributes to an $\ve$-trace identity.  The following holds:
\begin{itemize}
\item
  If $\chi \neq {}^\sigma \chi$, and $\tau$ is not $E$-monomial,
  then the only other contribution to the $\ve$-trace identity from $\mb{H_1}$
  is $\tau \otimes_1 {}^\sigma \chi$.
 \item
  If $\chi \neq {}^\sigma \chi$, and  $\tau$ is the cuspidal monomial representation $\pi(\theta)$ for 
  some character $\theta$ of $\idc{E}$, then the following
  packets of $\mb{H_1}(\Af)$ contribute to the $\ve$-trace identity provided that they are distinct from
  $\pi(\theta)\otimes_1 \chi$:
  \begin{itemize}
  \item
    $\pi(\theta)\otimes_1 {}^\sigma \chi$,
  \item
    $\pi(\chi)\otimes_1 \theta$,
  \item
    $\pi(\chi)\otimes_1\! {}^\sigma \theta$.
  \end{itemize}
\item
  If $\chi = {}^\sigma \chi$, then $\tau\otimes_1 \chi$ is the only
  packet of $\mb{H_1}(\Af)$ which contributes to the $\ve$-trace identity.
\end{itemize}
\begin{remark}
The theorem only comments on contribution from $\mb{H_1}$, it says nothing regarding contribution from
$\mb{H_2}$.
\end{remark}
\end{thm}

\begin{proof}[Sketch of proof.]
The proof requires twisted endoscopic lifting results for 
the triple \dt $(\GL(4), \theta, 1)$, where $\theta$ is the automorphism 
$g \mapsto {}^t g^{-1}$ of $\GL(4)$ (see \cite{KS}, \cite{A}, \cite{F}).
Consider the $F$-split $\theta$-twisted endoscopic group  $\mb{C}$
of $\GL(4)$ (see \cite[p. 22--23]{F}).
It is equal to $(\GL(2)\times\GL(2))'$, where the prime indicates that the two $\GL(2)$ components
have equal determinants.
Any automorphic representation of $\mb{C}(\Af)$ has the form 
\[
\tau_1 \times \tau_2 : (g, h) \mapsto \tau_1(g)\otimes\tau_2(h),\quad\forall (g, h) \in \mb{C}(\Af),
\] 
where $\tau_1, \tau_2$ are automorphic representations of $\GL(2, \Af)$.  Because of the condition
on the determinants of the two $\GL(2)$ components, $\tau_1 \times \tau_2$ is the same as 
$\mu \tau_1\times \mu^{-1}\tau_2$ for any quasi-character $\mu$ of $\idc{F}$.
Suppose $\tau_1 \times \tau_2$ contributes to a $\theta$-twisted trace identity derived from Kottwitz-Shelstad's formula for
$(\GL(4), \theta, 1)$.
Expecting the results of lifting from ${\rm SO}(4)$ to $\PGL(4)$ to be a guide (see proof of \cite[Prop. 4.1]{F1}) ,
the only other contribution from $\mb{C}$ to the trace identity should be $\tau_2 \times \tau_1$, provided
$\tau_1 \neq \tau_2$.
Although the work of \cite{F1} has not yet been extended to the case of nontrivial central characters,
Arthur has announced in \cite{A} the lifting results for the triple $(\GL(4), \theta, 1)$.  It is 
likely that a proof of the statement which we assume here will be contained in a complete treatise on the problem.
For now, we shall take for granted the above hypothesis regarding the fibers of the lifting from $\mb{C}$ to $\GL(4)$.


Recall the hypothesis of our theorem:  A global packet $\tau \otimes_1 \chi$ of
$\mb{H_1}(\Af)$ contributes to an $\ve$-trace identity.
If $\pi_1' = \tau' \otimes_1 \chi'$ is another global packet of $\mb{H_1}(\Af)$ which contributes to
the $\ve$-trace identity, then the global data parametrizing the representations $\tau'\times\pi(\chi')$ 
and $\tau \times \pi(\chi)$ of  $\mb{C}(\Af)$
lift to the same datum in $\GL(4, \CC)$.
Hence, by our hypothesis
$\tau'\times \pi(\chi')$ is equivalent to $\tau\times \pi(\chi)$ and/or $\pi(\chi)\times \tau$.

If two representations $\tau_1\times\tau_2$, $\tau_1'\times\tau_2'$ of $\mb{C}(\Af)$ are equivalent
to each other, then there exists a quasi-character $\mu$ of $\idc{F}$ such that
$\tau_1' = \mu^{-1} \tau_1$ and $\tau_2' = \mu\tau_2$.
Therefore, if $\tau'\times\pi(\chi')$ is equivalent to $\pi(\chi)\times\tau$, then
$\tau' = \mu^{-1}\pi(\chi)$ and $\pi(\chi') = \mu\tau$ for some quasi-character $\mu$ of $\idc{F}$.  
Since twisting by a character does not change whether a representation of $\GL(2, \Af)$ is $E$-monomial,
$\tau'\times\pi(\chi')$ cannot be equivalent to $\pi(\chi)\times \tau$ if $\tau$ is not $E$-monomial.

Suppose $\tau'\times\pi(\chi') = \tau\times\pi(\chi)$.
If $\pi(\chi') = \mu\pi(\chi)$, then $\chi'$ is equal to either $\lp \mu\circ\N_{E/F}\rp\cdot\chi$ 
or $\lp \mu\circ\N_{E/F}\rp\cdot{}^\sigma \chi$ (\cite[Chap. 7]{L1}).
Consequently, if $\tau$ is not $E$-monomial and $\chi \neq {}^\sigma \chi$,
we must have $\tau' = \mu^{-1}\tau$, and $\chi' =\lp\mu\circ\N_{E/F}\rp{}^\sigma \chi$.
For any quasi-character $\mu$ of $\idc{F}$, $\mu^{-1}\tau\otimes_1\lp\mu\circ\N_{E/F}\rp\chi$ defines the
same representation.  We conclude that if $\tau$ is non-$E$-monomial and $\chi \neq {}^\sigma \chi$, then
$\tau \otimes_1 {}^\sigma\chi$ is the only other packet of $\mb{H_1}(\Af)$ which contributes to the 
$\ve$-trace identity.

If $\tau$ is non-$E$-monomial and $\chi = {}^\sigma \chi$, then there is no other contribution from $\mb{H_1}(\Af)$.

Suppose $\tau$ is the cuspidal monomial representation 
$\pi(\theta)$ associated with some character $\theta \neq \chi$ of $\idc{E}$ such that $\theta \neq {}^\sigma \theta$.
Then, the representation
$\tau'\times \pi(\chi')$ is equivalent to either: (i) $\pi(\theta)\times \pi(\chi)$,
or (ii) $\pi(\chi)\times \pi(\theta)$.  In the first case, we have 
$\tau'\otimes_1 \chi' = \pi(\theta)\otimes_1{}^\sigma \chi$ if $\chi \neq {}^\sigma \chi$.
In the second case, $\tau'\otimes_1 \chi'$ must be equal to either $\pi(\chi)\otimes_1 \theta$
or $\pi(\chi)\otimes_1 {}^\sigma \theta$ provided that they are distinct from $\pi(\theta)\otimes_1 \chi$.
If $\chi = {}^\sigma\chi = \mu\circ\N_{E/F}$ for some
character $\mu$ of $\idc{F}$, then $\pi(\chi) = I(\mu, \mu\ve)$.  From the discussion in 
Section \ref{sec:stablespectrum}, neither $I(\mu,\mu\ve)\otimes_1\theta$ nor
$I(\mu,\mu\ve)\otimes_1 {}^\sigma \theta$ contributes to the stable spectrum of $\mb{H_1}$.

Contingent upon published results on the $\theta$-twisted endoscopic lifting for $\GL(4)$,
the theorem follows.
\end{proof}
Suppose a global packet of $\mb{H_1}(\Af)$ has the form $\tau\otimes_1 \chi$, where
$\tau$ is a representation of $\GL(2, \Af)$ and $\chi$ is a character of $\idc{E}$
which is not equal to ${}^\sigma \chi$.
If $\tr\!\lp \tau\otimes_1 \chi\rp(f_1)$ appears in an $\ve$-trace identity, 
then so does $\tr\! \lp\tau\otimes_1 {}^\sigma\chi\rp(f_1)$.
If an element $(h, x) \in \mb{H_1}(\Af)$ is a norm of some element $g \in \GSp(2, \Af)$,
then $(h, \sigma x)$ is also a norm of $g$.  Therefore, we expect that
if $f \in C(\mb{G}(\Af), \omega)$ and $f_1 \in C(\mb{H_1}(\Af), \omega)$ are matching functions, then
$f_1(h, x) = f_1(h, \sigma x)$ for all $(h, x) \in \mb{H_1}(\Af)$.
Consequently, the terms $\tr\! \lp \tau\otimes_1 \chi \rp(f_1)$ and
$\tr\!\lp\tau\otimes_1 {}^\sigma\chi\rp(f_1)$ in the trace identity combine into a single term
$2\cdot \tr\!\lp \tau\otimes_1 \chi\rp(f_1)$.
\subsection{Multi-Packets of $\mb{H_2}(\Af)$}
For any representation $\tau$ of $\GL(2, \Ae)$, 
let ${}^\sigma \tau$ denote the $\GL(2, \Ae)$-module which sends $(g_{ij}) \in \GL(2, \Ae)$
\indexs{sigma@${}^\sigma{\rm repn.}$}%
to $\tau((\sigma g_{ij}))$.  
For any representation $\pi_2 = \tau\otimes_2 \mu$ of $\mb{H_2}(\Af)$,
where $\omega_\tau = \mu\circ\N_{E/F}$, put ${}^\sigma\pi_2 := {}^\sigma \tau\otimes_2 \mu$.



\begin{lemma}\label{lemma:contribH2}
Suppose a cuspidal automorphic representation $\pi_2 = \tau\otimes_2 \mu$ of $\mb{H_2}(\Af)$
contributes to an $\ve$-trace identity.  
\begin{itemize}
\item
If $\pi_2 \neq {}^\sigma \pi_2$, then the only other 
cuspidal automorphic representation of $\mb{H_2}(\Af)$ which 
contributes to the $\ve$-trace identity is ${}^\sigma \pi_2$.
\item
If $\pi_2 = {}^\sigma\pi_2$, then no other cuspidal automorphic representation of $\mb{H_2}(\Af)$
contributes to the $\ve$-trace identity.
\end{itemize}
\begin{remark}
The condition ``cuspidal'' is crucial.  The lemma says nothing regarding non-cuspidal
contributions from $\mb{H_2}(\Af)$.
\end{remark}
\end{lemma}

Suppose $\pi_2'$ is another cuspidal automorphic representation of $\mb{H_2}(\Af)$ which 
contributes to the $\ve$-trace identity.  
Let $v$ be any finite place for which the representations $\pi_{2,v}, \pi_{2,v}'$ are unramified.
From Section \ref{sec:globalcomparison}, 
it follows that the F-H classes parametrizing $\pi_{2,v}, \pi_{2,v}'$
lift via $\xi_2 : {}^L H_2 \rightarrow \hat{G}$ to the same F-H class in $\hat{G}$.

To examine what the above statement entails, we need to consider separately the case
where the place $v$ is prime in $E$ and the case where $v$ splits in $E$.

\begin{claim}
Let $v$ be a finite place of $F$ which is unramified, prime in $E$.
If $\pi_{2,v}$ and $\pi_{2,v}'$ are unramified and their Frobenius-Hecke classes
lift to the same conjugacy class in $\hat{G}$, then $\pi_{2,v}$ and $\pi_{2,v}'$ 
are equivalent.
\end{claim}
\begin{proof}
Since $\pi_{2}$ is cuspidal, the unramified representation $\pi_{2,v}$ is an 
irreducible induced representation $I(\chi_1, \chi_2)\otimes_2 {\mu_v}$, where
$\mu_v$ is an unramified character of $F_v^\times$, and
$\chi_1, \chi_2$ are unramified characters of $E_v^\times$.
Since ${\mu_v} \circ \N_{E/F} = \chi_1\chi_2$, $\pi_{2,v}$ has the form 
$I\!\lp\chi_1, \frac{{\mu_v}\circ\N_{E/F}}{\chi_1}\rp\otimes_2 {\mu_v}$.
Since $v$ is prime in $E$, and $\chi_1$ is unramified, we have $\chi_1 = \chi_1'\circ\N_{E/F}$ 
for some unramified character $\chi_1'$ of $F_v^\times$.
Likewise,
$\pi_{2,v}' = I\lp\theta_1, \frac{{\mu'_v}\circ\N_{E/F}}{\theta_1}\rp\otimes_2 {\mu'_v}$,
where $\theta_1 = \theta_1'\circ\N_{E/F}$ for some unramified character $\theta_1'$
of $F_v^\times$, and $\mu_v'$ is an unramified character of $F_v^\times$.


From Section \ref{FHclass}, the lift of the F-H class parametrizing $\pi_{2,v}$ is represented by
\[
\hat{g} =
\ddmfour{\chi_1'}{-\chi_1'}{\frac{-{\mu_v}}{\chi_1'}}{\frac{{\mu_v}}{\chi_1'}}
\in \GSp(2, \mbb{C}),
\]
while the lift of the class parametrizing $\pi_{2,v}'$ is represented by
\[
\hat{g}' =
\ddmfour{\theta_1'}{-\theta_1'}{\frac{-{\mu'_v}}{\theta_1'}}{\frac{{\mu'_v}}{\theta_1'}}
\in \GSp(2, \CC).
\]

If $\hat{g}$ and $\hat{g}'$ belong to the same conjugacy class in $\hat{G}$,
then there must exist an element $w$ in the Weyl group $W$ of $\hat{G}$ such that
$w \hat{g} = \hat{g}'$.
Observe that ${\mu_v}(\vp_v)$ and $\mu'_v(\vp_v)$ are the respective similitude factors
of $\hat{g}$ and $\hat{g}'$.  Since the similitude factor
is invariant under conjugation, we conclude that ${\mu_v} = \mu'_v$.
Since $E_v/F_v$ is unramified, $\ve_v(\vp_v)$ is equal to $-1$.
By permuting the entries of $\hat{g}$ under $W$, 
we conclude that $\theta'$ is equal to one of the following characters:
\begin{enumerate}
\item
$\frac{{\mu_v}}{\chi_1'}$,
\item
$\ve_v\chi_1'$,
\item
$\frac{{\ve_v}{\mu_v}}{\chi_1'}$,
\item
$\chi_1'$.
\end{enumerate}
Since $\ve_v\circ\N_{E/F} = 1$,
in cases 2 and 4 we have $\theta_1 = \theta_1'\circ\N_{E/F}
= \chi_1'\circ\N_{E/F} = \chi_1$.  Hence, $\pi_{2,v}' = \pi_{2,v}$.
In cases 1 and 3, we have
$\pi_{2,v}' = I\lp\frac{{\mu_v}\circ\N_{E/F}}{\chi_1},\chi_1\rp\otimes_2 {\mu_v}$, which
is equivalent to $\pi_{2,v}$.
Hence, if the F-H classes parametrizing the unramified representations $\pi_{2,v}$ and $\pi_{2,v}'$
lift to the same conjugacy class in $\hat{G}$, then $\pi_{2,v}$ and $\pi_{2,v}'$ 
are equivalent.
\end{proof}

\begin{claim}
Let $v$ be a finite place of $F$ which splits in $E$.
If  $\pi_{2,v}$ and $\pi_{2,v}'$ are unramified, and their
Frobenius-Hecke classes lift to the same conjugacy
class in $\hat{G}$, then $\pi_{2,v}'$ is equivalent to either $\pi_{2,v}$ 
or ${}^\sigma \pi_{2,v}$.
\end{claim}
\begin{proof}
If $v$ is a finite place of $F$ which splits in $E$,
then $H_{2,v} = \mb{H_2}(F_v)$ has the form
\[
\left[\GL(2, F_v) \times \GL(2, F_v) \times F_v^\times\right]/\{(a I_2, b I_2, (ab)^{-1}) : a,b \in F_v^\times\}.
\]
Suppose $v_1, v_2$ are the two places of $E$ which lie above $v$.
We have
$\pi_{2,v} = \tau_1 \times \tau_2 \times \mu_v$, where $\tau_i = \tau_{v_i}$ ($i = 1, 2$), and
${}^\sigma \pi_{2,v} = \tau_2 \times \tau_1 \times \mu_v$.
Since
\[
\omega_{\tau_1}(a)\omega_{\tau_2}(b)\mu_v(a)^{-1}\mu_v(b)^{-1} = 1, \quad \forall a, b \in F_v^\times,
\]
we must have
$\omega_{\tau_1} = \omega_{\tau_2} = \mu_v$.  In particular, $\mu_v$ is uniquely determined by
either $\tau_1$ or $\tau_2$.

Observe that $H_{2,v}$ is isomorphic via $(g, h, c) \mapsto (g, ch)$ to
\[
\left[\GL(2, F_v) \times \GL(2, F_v)\right]/\{(z I_2,z^{-1} I_2) : z \in F_v^\times\}.
\]
We identify representations of $H_{2,v}$ with representations of the above group, and 
write $\pi_{2,v} = \tau_1 \times \tau_2$, ${}^\sigma\pi_{2,v} = \tau_2 \times \tau_1$.


If $\pi_{2,v} = \tau_1 \times \tau_2$ is unramified, then $\tau_i = I(\alpha_i,\beta_i)$ ($i=1,2$),
where $\alpha_i, \beta_i$ are unramified characters of $F_v^\times$.
Likewise, $\pi_{2,v}' =
I(\alpha_1',\beta_1') \times I(\alpha_2',\beta_2')$ for some unramified characters $\alpha_i', \beta_i'$ of
$F_v^\times$.


The F-H class parametrizing $\pi_{2,v}$ lifts to the conjugacy class in $\hat{G}$ represented by
\[\hat{g} = [\diag(\alpha_1,\beta_1),\diag(\alpha_2,\beta_2)],\]
while the F-H class parametrizing $\pi_{2,v}'$ lifts to the conjugacy class in $\hat{G}$ represented by
\[\hat{g}' = [\diag(\alpha_1',\beta_1'),\diag(\alpha_2',\beta_2')].\]
By assumption, there exists an element $w \in W$ such that $w\hat{g} = \hat{g}'$.
If we label the entries of the maximal diagonal torus of $\hat{G}$ by
$1, 2, 3, 4$, then $W = D_4$ is generated by the permutations $(3421)$ and $(23)$.

Suppose $(23)\hat{g} = \hat{g}'$.
The action of $(23)$ sends
$\hat{g} = [\diag(\alpha_1,\beta_1),\diag(\alpha_2,\beta_2)]$ to
$[\diag(\alpha_1,\beta_1),\diag(\beta_2,\alpha_2)]$.  Since
$I(\beta_2, \alpha_2)$ is equivalent to $I(\alpha_2,\beta_2)$, 
$\pi_{2,v}'$ is equivalent to $\pi_{2,v}$.

The action of $(3421)$ takes $\hat{g}$ to
$[\diag(\alpha_2, \beta_2), \diag(\beta_1, \alpha_1)]$.  Hence,
if $\hat{g}' = (3421)\hat{g}$, then $\pi_{2,v}' = \tau_2 \times \tau_1 = {}^\sigma\pi_{2,v}$.

We conclude that if the F-H classes parametrizing the unramified representations $\pi_{2,v}$ and
$\pi_{2,v}'$ lift to the same conjugacy class in $\hat{G}$, then
$\pi_{2,v}'$ is equivalent to either $\pi_{2,v}$ or ${}^\sigma \pi_{2,v}$.
In particular, $\mu_v'$ is equal to $\mu_v$.
\end{proof}
From the previous two claims,
we see that $\mu_v = \mu'_v$ for almost all $v$, which implies that $\mu = \mu'$ by the
weak approximation theorem for $\Af$ (The weak approximation theorem says that $F$ is dense in
$\prod_{v \in S}F_v$ for any finite set of places $S$.  It follows that
$F^\times$ is dense in $\prod_{v \in S}F_v^\times$.  Using that $\mu$ and $\mu'$ are continuous
characters of $\Af^\times$ whose restrictions to $F^\times$ are trivial, our assertion follows).

Suppose $\pi_2' = \tau'\otimes \mu$, where $\tau'$ is a cuspidal automorphic representation of $\GL(2, \Ae)$.
To prove Lemma \ref{lemma:contribH2}, it remains to show that:
\begin{claim}\label{claim:contribH2JS}
$\tau'$ is equivalent to either
$\tau$ or $\;{}^\sigma \tau$.
\end{claim}
Before we prove the claim, we need a digression:
Let $k$ be a number field.
Let $V_k$ be the set of places of $k$.  Let $n, m$ be positive integers.
\indexs{V@$V_k$}%
Let $\pi$ be an automorphic representation of $\GL(n, \adl_k)$ which is unramified at places 
outside of a finite set $S \subset V_k$.  For each $v \notin S$, 
let $t(\pi_v)$ be a matrix in $\GL(2, \CC)$ whose conjugacy class parametrizes $\pi_v$.
Let $\pi'$ be an automorphic representation of $\GL(m, \adl_k)$ which is unramified at places outside
of $S$. Define $t({\pi'_v}) \in \GL(m, \CC)$ likewise for $\pi'$.

Consider $t(\pi_v), t(\pi'_v)$ as endomorphisms of $\CC^n$, $\CC^m$, respectively.
Their tensor product $t(\pi_v)\otimes t(\pi'_v)$, being an endomorphism of $\CC^{mn}$, is represented
by an $mn \times mn$ complex matrix.  
Define a partial $L$-function in $s \in \CC$ as follows:
\indexi{L@$L$-function!partial}%
\[
L^S(s, \pi\otimes{\pi'})
:= \prod_{v \notin S}\det\lp 1 - t_v(\pi)\otimes t_v(\pi')q_v^{-s}\rp^{-1},
\]
where $q_v$ is the cardinality of the residue field of $k_v$.
If $\pi$ and $\pi'$ are cuspidal, 
it is known (\cite{JS}) that: (i) $L^S(s, \pi\otimes\pi')$ does not vanish on the line ${\rm Re}\;s = 1$;
(ii) $L^S(s, \pi\otimes\pi')$ has a pole of order $1$ at $s = 1$ if and only if 
$n = m$ and $\pi'$ is equivalent to $\tilde{\pi}$, the contragradient of $\pi$.
\indexi{representation!contragratient}%

\begin{proof}[Proof of Claim \ref{claim:contribH2JS}.]
Let $S$ be a finite set of places of $F$ such that 
$\pi_{2,v}, \pi_{2,v}'$ are unramified for all $v \notin S$.
We know, for $v \notin S$, that:
\[
\pi_{2,v}' = \begin{cases}
\pi_{2,v} & \text{ if } v \text{ is unramified, prime in }E,\\
\pi_{2,v}\;\text{ or }\;{}^\sigma\pi_{2,v} & \text{ if } v \text{ splits in } E.
\end{cases}
\]
Consider $\tau, \tau'$ as automorphic representations of $\GL(2, \Ae)$, 
with $E$ as the base field.
Let $S_E$ be the set of places of $E$ which lie above the places in $S$.  In particular, $S_E$ is finite.
For any place $w$ of $E$ outside of $S_E$, $\tau_w$ and $\tau_w'$ are unramified.  We have:
\[
\tau'_w = \tau_w \text{ or } \lp{}^\sigma \tau\rp_w.
\]
Notice the order of operation: We first apply $\sigma$ to the global $\tau$, and then
we take the component at the place $w$.
We have the following equality of products of partial $L$-functions:
\[
L^{S_E}(s, \tau'\otimes \tilde{\tau})L^{S_E}(s, {}^\sigma \tau'\otimes \tilde{\tau})
=
L^{S_E}(s, \tau\otimes \tilde{\tau})L^{S_E}(s, {}^\sigma \tau\otimes \tilde{\tau}).
\]
Since $\pi_2$ and $\pi_2'$ are cuspidal automorphic representations of $\mb{H_2}(\Af)$, 
$\tau_w$ and $\tau_w'$ are cuspidal automorphic representations of $\GL(2, \Ae)$.
By \cite{JS} and \cite{Sh}, 
the product on the right has a single pole at $s = 1$;  therefore, the product on the left
must have a pole, which implies that $\tau'$ is equivalent to either
$\tau$ or ${}^\sigma \tau$.
\end{proof}

Let $1_{\GL(2, \Ae)}$ denote the trivial representation of $\GL(2, \Ae)$.
\indexs{1@$1_{\GL(2,\Af)}$}%
\begin{lemma}\label{lemma:contribH21dim}
Let $\chi$ be a character of $\idc{E}$, $\mu$ be a character of $\idc{F}$, 
such that $\chi^2 = \mu \circ \N_{E/F}$.  Suppose the automorphic representation 
\[
\pi_2 = \lp\chi\circ\det\rp\otimes_2 \mu = \chi 1_{\GL(2, \Ae)}\otimes_2\mu
\]
of $\mb{H_2}(\Af)$ contributes to an $\ve$-trace identity.
\begin{itemize}
\item
If $\chi \neq {}^\sigma \chi$, then the only
other contribution from $\mb{H_2}$ to the $\ve$-trace identity is 
${}^\sigma \pi_2 = {}^\sigma\chi 1_{\GL(2, \Ae)}\otimes_2 \mu$. 
\item
If $\chi = {}^\sigma \chi$, then no other automorphic representation of $\mb{H_2}(\Af)$
contributes to the $\ve$-trace identity.
\end{itemize}
\end{lemma}
\begin{proof}
Suppose an automorphic representation $\pi_2' = \tau\otimes_2 \mu' \neq \pi_2$ 
of $\mb{H_2}(\Af)$ contributes to the trace 
identity.
As in the case where $\pi_2$ is cuspidal, we consider the places $v$, unramified in $E$, 
where both $\pi_{2,v}$ and $\pi_{2,v}'$ are unramified.  Of these places, we examine separately
those which split in $E$ and those which do not.

If $v$ does not split in $E$, then by the same argument used in the proof of
Lemma \ref{lemma:contribH2}, we conclude that $\pi_{2,v} = \pi_{2,v}'$.
Suppose $v$ splits into two places $v_1, v_2$ of $E$.  For $i = 1, 2$, let $\chi_i = \chi_{v_i}$.  Since
$\chi^2 = \mu\circ\N_{E/F}$, we have $\chi_1^2 = \chi_2^2 = \mu_v$.
Let $\nu_v$ be the normalized absolute value function on $F_v$.
Then,
$\pi_{2,v}' = \tau_{v_1}'\times \tau_{v_2}'$ is parametrized by the conjugacy class of
one of the following two elements in ${}^L{H_2}$:
\begin{enumerate}
\item
$
\lp\dmtwo{\chi_1\nu_v^{1/2}}{\chi_1\nu_v^{-1/2}}, \dmtwo{\chi_2\nu_v^{1/2}}{\chi_2\nu_v^{-1/2}},
\mu_v\rp \rtimes 1,
$
\item
$
\lp\dmtwo{\chi_2\nu_v^{1/2}}{\chi_2\nu_v^{-1/2}},\dmtwo{\chi_1\nu_v^{1/2}}{\chi_1\nu_v^{-1/2}},
\mu_v\rp \rtimes 1.
$
\end{enumerate}
In the first case, we have $\pi_{2,v}' = \pi_{2,v}$.  
In the second case, we have $\pi_{2,v}' = {}^\sigma \pi_{2,v}$.
Moreover, as in the cuspidal case, we have $\mu_v = \mu'_v$ for almost all $v$.

Hence, $\pi_{2}' = \chi' 1_{\GL(2, \Af)}\otimes_2 \mu$, where $\chi'$ is a character of $\idc{E}$ such that
$\chi'_v$ is equal to $\chi_v$ or ${}^\sigma \chi_v$ for almost all $v$.  Applying the partial $L$-function argument
(explained in the proof of Lemma \ref{lemma:contribH2})
to $\chi$ as a representation of $\GL(1, \Ae)$, we conclude that $\chi'$ is equal to either
$\chi$ or ${}^\sigma \chi$.  The lemma follows.
\end{proof}
\begin{corollary}\label{corollary:contribH2discrete}
Let $\pi_2$ be a discrete spectrum automorphic representation of $\mb{H_2}(\Af)$ which
contributes to an $\ve$-trace identity. 
Then, a discrete spectrum representation $\pi_2' \neq \pi_2$ of $\mb{H_2}(\Af)$ contributes 
to the trace identity if and only if
$\pi_2 \neq {}^\sigma \pi_2$, in which case $\pi_2'$ is equal to ${}^\sigma \pi_2$.
\end{corollary}
\begin{proof}
The group $\mb{H_2}(\Af)$ is a quotient of $\GL(2, \Ae)\times\Af^\times$; hence,
its discrete spectrum representations are either cuspidal or one dimensional.  The corollary
follows from Lemmas \ref{lemma:contribH2} and \ref{lemma:contribH21dim}.
\end{proof}
If $(h, c) \in \mb{H_2}(\Af)$ is a norm of $g \in \mb{G}(\Af)$, then $(\sigma h, c)$ is also a norm
of $g$.  Hence, we expect that
if $f \in C(\mb{G}(\Af), \omega)$ and $f_2 \in C(\mb{H_2}(\Af), \omega)$ are matching functions, then
$f_2(h, c) = f_2(\sigma h, c)$ for all $(h, c) \in \mb{H_2}(\Af)$.
Consequently, if $\pi$ is an automorphic  representation of $\GL(2, \Ae)$ which is not equal to ${}^\sigma \pi$ and
$\mu$ is a character of $\idc{F}$ such that $\omega_\pi = \mu\circ\N_{E/F}$,
the contributions from  $\tr\!\lp\pi\otimes_2\mu\rp(f_2)$ and $\tr\!\lp{}^\sigma\pi\otimes_2 \mu\rp(f_2)$ to an
$\ve$-trace identity combine into a single term $2 \cdot\tr\!\lp\pi\otimes_2 \mu\rp(f_2)$.
\indexi{packet|)}%
\indexi{packet!multi-|)}%
\section{Contributions}\label{sec:contributions}
We now address the question:  Under what circumstances are there contributions from
both $\ve$-endoscopic groups to an $\ve$-trace identity.

{\bf Notation: }
\begin{itemize}
\item
Let $\Pi$ be an automorphic representation of $\mb{G}(\Af)$.  Let $f$ be a function in \dt $C(\mb{G}(\Af), \omega)$.
Recall that $\Pi(f\times\ve)$ is defined in Section \ref{sec:finechiexp}.
Put
\[
\la \Pi, f \ra_\ve := \tr \Pi(f \times \ve).
\indexs{b@$\la\;,\;\ra_\ve$}%
\]
If $\Pi$ is $\ve$-invariant, $\la \Pi, f\ra_\ve$ is equal to a product $\prod_v \tr \Pi_v(f_v\times\ve_v)$ of
twisted local characters.  Put $\la \Pi_v, f_v\ra_{\ve_v} := \tr \Pi_v(f_v\times\ve_v)$.
\item
Let $i = 1$ or $2$.
For a representation $\pi_i$ of $\mb{H}_i(\Af)$ and test function $f_i$ in \dt $C(\mb{H}_i(\Af), \omega)$,
put
\[
\la \pi_i, f_i\ra := \tr \pi_i(f_i).
\indexs{b@$\la\;,\;\ra$}%
\]
\item
For any automorphic representation $\tau$ of $\GL(2, \Af)$, let $B_{E/F}\tau$ denote the
\indexi{base change}%
automorphic representation of $\GL(2, \Ae)$ which is obtained by base change from $\tau$ (see \cite{L1}, \cite{F4}).
\item
Let $V$ be the set of places of $F$.
Let $\Vr$ denote the set of finite places of $F$ which are unramified in $E$.
\item
Let $S$ be a set of places of $F$.
For any ad\`elic object (automorphic representation, test function, trace, \dots, etc.), let
subscript $S$ denote the tensor product of local components over the places in $S$.
\end{itemize}
Let $\pi$ be an irreducible admissible representation of a $p$-adic group $H$.
In \cite{HC1}, Harish-Chandra proves the existence of a locally integrable function $\chi_\pi$
\indexi{Harish-Chandra}%
on $H$ such that \dt 
\[
\int_H \chi_\pi(h)f(h)\;dh = \tr\! \pi(f),\quad\forall f \in C_c^\infty(H).
\]  
We call $\chi_\pi$ the {\bf Harish-Chandra character} of $\pi$.
\indexi{Harish-Chandra character}%
This notion extends to the $\ve$-twisted characters as follows:
Let $k$ be a local field.  Let $\ve_k$ be a quadratic character of $\mb{G}(k)$.
Let $\pi$ be an irreducible, admissible, $\ve_k$-invariant $\mb{G}(k)$-module,
i.e. there exists a nontrivial intertwining operator $A$ in $\Hom_{\mb{G}(k)}(\pi, \ve_k\pi)$.
Note that the space of $\ve_k \pi$ is equal to that of $\pi$.
The {\bf $\ve_k$-twisted Harish-Chandra character} $\chi^A_{\pi}$ of $\pi$ is a 
\indexi{Harish-Chandra character!twisted}%
function on $\mb{G}(k)$ such that
\[\int_{\mb{G}(k)} \chi^A_{\pi}(g)f(g)\;dg = \tr \pi(f)A,\quad \forall f \in C_c^\infty(\mb{G}(k)).\]
It can be shown that $\chi_\pi^A$ is locally constant on the regular set, locally integrable, and unique.
It satisfies $\chi_\pi^A(h^{-1}gh) = \ve_k(h)\chi_\pi^A(g)$.

\quad\\
{\bf Definition: }
We say that an irreducible admissible representation of a $p$-adic group $H$ is {\bf elliptic} 
\indexi{representation!elliptic}%
if its Harish-Chandra character is not identically zero on the elliptic regular set of $H$.
Likewise, we say that an $\ve_k$-invariant, irreducible, admissible
 $\mb{G}(k)$-module $\pi$ is {\bf $\ve_k$-elliptic} 
\indexi{representation!$\ve$-elliptic}%
if $\chi_\pi^A$ is not identically zero on the elliptic regular set of $\mb{G}(k)$.  
It follows from the orthogonality relations that each square integrable admissible representation of
a connected reductive $p$-adic group is elliptic.  
In the twisted case,
the \emph{$\ve_k$-invariant} square integrable representations are similarly elliptic.
This too follows from the orthogonality relations of characters.

In order to apply Corollary \ref{corollary:finediscreteeq} to obtain global lifting results,
we fix once and for all two distinct finite places $w_1, w_2$ of $F$ which are prime in $E$, and
we work with automorphic representations whose local components at $w_i$ ($i = 1, 2$) are elliptic.
For $\mb{H} = \mb{G}$, $\mb{H_1}$, or $\mb{H_2}$,
recall that $\ELw{\mb{H}(\Af), \omega}$ denotes the set of functions $f = \otimes_v f_v$ 
in $C(\mb{H}(\Af), \omega)$ whose local components $f_{w_1}, f_{w_2}$ are elliptic.

For a reductive $F$-group $\mb{H}$ and an automorphic representation $\pi$ of $\mb{H}(\Af)$,
we define the {\bf bad places} for $\pi$ to be the (finite) set of places of $F$ which is the union of
\indexi{bad places}%
$V - \Vr$, $\{w_1, w_2\}$, and the set of places $v$ where $\pi_v$ is ramified (i.e. not unramified).

\subsection{Base Change}\label{sec:globalbc}
\indexi{base change}%
Let $\tau$ be a cuspidal non-$E$-monomial, or one dimensional, automorphic representation of 
$\GL(2, \Af)$ such that $\tau_{w_i}$ ($i = 1, 2$) is elliptic.
Let $\{\Pi\}$ be the unstable (quasi-)packet $[\tau, \ve\tau]$ of $\mb{G}(\Af)$.
This is the packet which lifts to the induced representation $I_{(2, 2)}(\tau, \ve\tau)$ of
$\GL(4, \Af)$.
The table notation was introduced in Section \ref{sec:vetraceidentity}.
\begin{prop}\label{prop:globalbc}
Let $S$ be the set of bad places for $\tau$.
The following table holds for matching test functions with elliptic components at $w_1, w_2$
and spherical components at all $v \notin S$:
\begin{equation}\label{basechange0}
\begin{tabular}{|c|c|c|}
\hline
$G$ & $H_1$ & $H_2$\\
\hline
$\{\Pi\}$&
$\tau\otimes_1 1$ & 
$B_{E/F}\tau\otimes_2\omega_{\tau}$\\
\hline
\end{tabular}
\end{equation}
\begin{remark}
The representation $\tau \otimes_1 \mu\circ\N_{E/F}$ is equivalent to $\mu\tau\otimes_1 1$ for any character
$\mu$ of $\idc{F}$; hence, by Hilbert 90,
the proposition covers all $\tau\otimes_1 \chi$ where $\chi = {}^\sigma\chi$.
\end{remark}
\end{prop}

\begin{proof}
For $v \notin S$, let $t(\tau_v)$ be a diagonal matrix in $\GL(2, \CC)$ whose conjugacy class parametrizes the
unramified representation $\tau_v$.
From Section \ref{FHclass}, the global datum \dt $\mc{C}(\{\Pi\}, S, \hat{G})$ parametrizing $\{\Pi\}$ is
represented by
\[
\Big[\left[t(\tau_v), \ve_v(\vp_v)t(\tau_v)\right]\Big]_{v \notin S} \in \prod_{v\notin S}\hat{G},
\]
where $\vp_v$ is a fixed uniformizer of $F_v$ for each $v \notin S$.
By comparing F-H classes, we see that $\tau\otimes_1 1$ and $B_{E/F}\tau\otimes_2 \omega_\tau$
contribute to the $\ve$-trace identity defined by the global datum $\mc{C}(\{\Pi\}, S, \hat{G})$.
By Theorem \ref{thm:contribH1}, $\tau\otimes_1 1$ is the only contribution from $\mb{H_1}$.
By Corollary \ref{corollary:contribH2discrete}, $B_{E/F}\tau\otimes_2\omega_\tau$ is the only
discrete spectrum contribution from $\mb{H_2}$.  Since the test function $f_2$ on $\mb{H_2}(\Af)$
has elliptic components, no discretely
occuring induced representation of $\mb{H_2}(\Af)$ contributes to the trace identity.
The proposition follows.
\end{proof}
\begin{corollary}\label{corollary:globalbc}
Let
$f \in \ELw{\mb{G}(\Af), \omega_\tau}$ and $f_i \in \ELw{\mb{H}_i, \omega_\tau}$
$(i = 1, 2)$ be matching functions
whose components at all $v \notin S$ are spherical,
the following character identity holds:
\begin{equation}\label{eq:globalbc}
\sum_{\Pi' \in \{\Pi\}} m(\Pi') \la \Pi', f\ra_{\ve, S}
= \frac{1}{2}\la \tau\otimes_1 1, f_1\ra_S + \frac{1}{2}\la B_{E/F}\tau\otimes_2 \omega_\tau, f_2\ra_S.
\end{equation}
Here $m(\Pi')$ is the multiplicity of $\Pi'$ in the discrete spectrum of $\mb{G}(\Af)$.
\begin{remark}
If $\omega_\tau = 1$, then by \cite[V. 10]{F1} the packet $\{\Pi\}$ is
the restricted tensor product $\otimes_{v \in V}\{\Pi_v^+, \Pi_v^-\}$,
where $\Pi_v^- = 0$ unless $\tau_v$ is square integrable.
The multiplicity formula is
$m(\Pi') = \frac{1}{2}(1 + (-1)^{n(\Pi')})$, where $n(\Pi')$ is the number of 
places $v$ for which $\Pi'_v = \Pi_v^-$.
\end{remark}
\end{corollary}
\subsection{Exclusive Contributions}\label{sec:exclcontrib}
\indexi{contribution!exclusive}%


We now address under what circumstances do the representations of only one of the
$\ve$-endoscopic groups contribute to an $\ve$-trace identity.  

\label{begin:unbase}
\subsubsection{Case of $\mb{H_1}$}
\begin{prop}\label{prop:nonbcH1}
Let $\pi_1 = \tau\otimes_1\chi$ be a global $($quasi-$)$packet of representations of 
$\mb{H_1}(\Af)$ such that:
\begin{itemize}
\item
$\tau$ is a non-$E$-monomial automorphic representation of $\GL(2, \Af)$;
\item
$\tau_{w_i}$ $(i = 1, 2)$ is elliptic;
\item
$\chi \neq {}^\sigma \chi$.
\end{itemize}
If $\pi_1$ contributes to an $\ve$-trace identity 
for matching functions $f$ in \dt $\ELw{\mb{G}(\Af), \omega}$ and
$f_i$ in $\ELw{\mb{H}_i(\Af), \omega}$ $(i = 1, 2)$,
then no representation of $\mb{H_2}(\Af)$ contributes to the trace identity.
\\\rm{
The proof is given after Lemma \ref{lemma:nbcH1ell1'}.
}
\end{prop}

Suppose an automorphic representation $\pi_2$ of $\mb{H_2}(\Af)$ contributes to the $\ve$-trace identity.
Let $S$ be the union of the respective sets of bad places for $\pi_1$, $\pi_2$.
That is, $\pi_{1 , v}$, $\pi_{2,v}$ are unramified for all $v \notin S$.
Let $\{\Pi\}$ be the global (quasi-)packet of representations of $\mb{G}(\Af)$
whose global datum is the lift of $\mc{C}(\pi_1, S, {}^L H_1)$. Put ${}^\sigma \pi_1 := \tau\otimes_1 {}^\sigma \chi$.
By Corollary \ref{corollary:finediscreteeq}, Theorem \ref{thm:contribH1}, 
and Corollary \ref{corollary:contribH2discrete},
we have:
\begin{multline}\label{eq:stabaux}
\sum_{\Pi'\in \{\Pi\}} m(\Pi')\la \Pi', f \ra_{\ve, S}
= \frac{1}{2}\la \pi_{1}, f_{1} \ra_S 
+ \frac{1}{2}\la {}^\sigma\pi_{1}, f_{1} \ra_S\\
+ \frac{1}{2} \la \pi_{2}, f_{2} \ra_S
+ \frac{d}{2} \la {}^\sigma \pi_{2}, f_{2} \ra_S,
\end{multline}
where $m(\Pi')$ is the multiplicity
of $\Pi'$, and
$
d = \begin{cases} 1 &\text{if } \pi_2 \neq {}^\sigma \pi_2,\\0&\text{otherwise}.\end{cases}
$
\begin{lemma}\label{lemma:nbcH1stable}
The $($quasi-$)$packet $\{\Pi\}$ is stable.
\end{lemma}
\begin{proof}
Suppose $\{\Pi\}$ is unstable.  Then,
according to Lemma \ref{lemma:listunstable}, one of the following holds:
\begin{enumerate}
\item
$\{\Pi\} = [\pi, \ve\pi]$,
where $\pi$ is a cuspidal non-$E$-monomial automorphic representation of $\GL(2, \Af)$.
\item
$\{\Pi\} = [\mu 1_2, \ve\mu 1_2]$ 
for some character $\mu$ of $\idc{E}$.
\item
$\{\Pi\} = [\pi(\theta\chi), \pi(\theta\;{}^\sigma\!\chi)]$,
where $\theta, \chi$ are characters of $\idc{E}$ such that 
none of $\theta$, $\chi$, $\theta\chi$, $\theta\,{}^\sigma\!\chi$ is fixed by the action of $\Gal(E/F)$.
\end{enumerate}

In the first (resp. second) case, we know from Proposition \ref{prop:globalbc}
that $\{\Pi\}$ is the lift of the packet $\pi_1' = \pi\otimes_1 1$
(resp. $\mu 1_2\otimes_1 1$) of $\mb{H_1}(\Af)$.
Since $\chi \neq {}^\sigma \chi$, by Theorem \ref{thm:contribH1},
$\pi_1'$ and $\pi_1$ cannot both contribute to the trace identity, whence a contradiction.

In the third case, one can see via a comparison of F-H classes that
$\{\Pi\}$ is the lift of $\pi(\theta)\otimes_1\chi$.
Since we assume that $\tau$ is not $E$-monomial, we have a contradiction
by Theorem \ref{thm:contribH1}.
\end{proof}
\begin{corollary}\label{corollary:nbcH1stable}
Each member of $\{\Pi\}$
occurs with multiplicity one in the discrete spectrum of $\GSp(2, \adl_{F})$.
\end{corollary}

Let $\{\{\Pi_v\}\}_{v \in V}$ be the collection of local packets such that $\{\Pi\}$ is the 
restricted tensor product $\otimes_{v \in V}\{\Pi_v\}$.
For each member $\Pi' = \otimes_v \Pi'_v$ of $\{\Pi\}$,
$\la \Pi', f\ra_\ve$ is equal to a product $\prod_v \la \Pi'_v, f_v\ra_{\ve_v}$ of  twisted traces.
For simplicity, we often drop the subscript $v$ from $\ve_v$ if the meaning is clear.
By Lemma \ref{lemma:nbcH1stable}, we have:
\[
\sum_{\Pi' \in \{\Pi\}}m(\Pi')\la \Pi', f \ra_{\ve, S} =
\prod_{v \in S}\la \{\Pi_v\}, f_v\ra_{\ve},
\]
where $\la \{\Pi_v\}, f_v\ra_\ve := \sum_{\Pi'_v \in \{\Pi'_v\}}\la \Pi'_v, f_v\ra_\ve$.

Let $S_0 = S - \{w_1, w_2\}$.
Let
\[C_{S_0} = \la \{\Pi\}, f \ra_{\ve, S_0},\;
D_{1,S_0} = \la \pi_{1}, f_{1}\ra_{S_0},\;
D_{2, S_0} = \la \pi_{2}, f_{2} \ra_{S_0}.\]

From the matching condition on $f$,$f_1, f_2$, we have
$\la \pi_1, f_1 \ra = \la {}^\sigma \pi_1, f_1 \ra$ and
$\la \pi_2, f_2 \ra = \la {}^\sigma \pi_2, f_2\ra$.
We rewrite \eqref{eq:stabaux} as follows:
\begin{multline}\label{eq:nbcbeforetype1}
C_{S_0} \la \{\Pi_{w_1}\}, f_{w_1}\ra_\ve \la \{\Pi_{w_2}\}, f_{w_2} \ra_\ve\\
=  D_{1, S_0}\la \pi_{1, w_1}, f_{1, w_1} \ra \la \pi_{1, w_2}, f_{1, w_2} \ra
+ n\cdot D_{2, S_0} \la \pi_{2, w_1}, f_{2, w_1} \ra \la \pi_{2, w_2}, f_{2, w_2} \ra,
\end{multline}
where $n = \frac{1}{2}$ or $1$, depending on whether $\pi_2$ is equal to ${}^\sigma \pi_2$ or not.
Here, for each place $v$, the symbol
$\pi_{1,v}$ denotes the local packet of $H_{1,v}$ such that $\pi_1 = \otimes_v \pi_{1,v}$.

\begin{define}
{\rm 
Let $v$ be a place of $F$ which is prime in $E$.
We say that an elliptic regular element in $G_v$ is of {\bf type $i$} ($i = 1, 2)$
\indexi{element!type}%
if its norms lie only in $H_{i, v}$.
We say that an elliptic function $f_v$ on $G_v$ is of {\bf type $i$} if
\indexi{function!type}%
the orbital integral of $f_v$ is nonzero only at the elliptic regular elements of type $i$.
}
\end{define}
From the norm correspondence among elliptic regular elements (computed in Appendix \ref{chap:fundlemma}), 
type $i$ ($i = 1, 2$) elements exist.  Consequently, type $i$ functions exist.

\begin{lemma}\label{lemma:nbcH1ell1'}
There exists a type 1 elliptic function $f_{{w_2}}$ on $G_{w_2}$ such that 
the trace \dt $\la \{\Pi\}_{{w_2}}, f_{{w_2}} \ra_\ve$ is nonzero.
\end{lemma}
\begin{proof}
Let $f_{w_1}$ be a type 1 elliptic function such that
$\la \pi_{1,w_1}, f_{1,w_1} \ra \neq 0$.  This is possible because by hypothesis
$\pi_{1,w_1}$ is elliptic.
Then, we may choose $f_{2, w_1}$ to be zero and obtain the following equation:
\[
C_{S_0} \la \{\Pi_{w_1}\}, f_{w_1}\ra_\ve \la \{\Pi_{{w_2}}\}, f_{{w_2}} \ra_\ve
=  D_{1, S_0}\la \pi_{1, w_1}, f_{1, w_1} \ra \la \pi_{1, {w_2}}, f_{1, {w_2}} \ra.
\]
For each $v \in S_0$, choose $f_v$ and matching $f_{1,v}$ such that $\la \pi_{1,v}, f_{1,v}\ra \neq 0$.
Since by assumption $\pi_{1,w_2}$ is elliptic, there exists a type 1 elliptic 
$f_{w_2}$ and matching $f_{1,w_2}$ such that $\la \pi_{1,w_2}, f_{1,w_2}\ra \neq 0$.  The lemma follows.

\end{proof}
\begin{proof}[Proof of Proposition \ref{prop:nonbcH1}]
Let $f_{w_1}$ be a type 2 elliptic function on $G_{w_1}$.  Let 
$f_{1,w_1}$ be zero.  Then, \eqref{eq:nbcbeforetype1} becomes
\begin{equation}\label{eq:fv12,fv2arb}
C_{S_0} \la \{\Pi_{w_1}\}, f_{w_1}\ra_\ve \la \{\Pi_{{w_2}}\}, f_{{w_2}} \ra_\ve
=
n D_{2, S_0} \la \pi_{2, w_1}, f_{2, w_1} \ra \la \pi_{2, {w_2}}, f_{2, {w_2}} \ra.
\end{equation}
Let $f_{w_2}$ be a type 1 function on $G_{w_2}$
such that $\la \{\Pi\}_{{w_2}}, f_{{w_2}} \ra_\ve$ is nonzero.
It exists by Lemma \ref{lemma:nbcH1ell1'}.
Let $f_{2, w_2}$ be zero.  Then,
\[
C_{S_0} \la \{\Pi_{w_1}\}, f_{w_1}\ra_\ve = 0
\]
for all elliptic functions $f_{w_1}$ of type 2.

Now let $f_{{w_2}}$ be an arbitrary test function on $G_{w_2}$.  
The equation \eqref{eq:fv12,fv2arb}, where $f_{w_1}$ is of type 2, implies that
\begin{equation}
0 = C_{S_0} \la \{\Pi_{w_1}\}, f_{w_1}\ra_\ve \la \{\Pi_{{w_2}}\}, f_{{w_2}} \ra_\ve
=
n D_{2, S_0} \la \pi_{2, w_1}, f_{2, w_1} \ra \la \pi_{2, {w_2}}, f_{2, {w_2}} \ra.
\end{equation}
There are two possibilities:  (i) $\la \pi_{2,w_2}, f_{2, w_2}\ra = 0$ for every
$f_{2,w_2}$ matching some elliptic function $f_{w_2}$ on $G_{w_2}$; (ii)
$\la \pi_{2,w_1}, f_{2,w_1}\ra = 0$ for every $f_{2,w_1}$ matching some type 2 elliptic function on $G_{w_1}$.
If case (i) holds, we are done.  Suppose case (ii) holds.
For any $h_2 = (g, c) \in H_{2,w_1}$, put $\sigma h_2 := (\sigma g, c)$.
Since $f_{2,w_1}$ matches a function on $G_{w_1}$, we have $f_{2,w_1}(h_2) = f_{2,w_1}(\sigma h_2)$
for all $h_2 \in H_{2,w_1}$.  Hence, 
\[
0 = \la \pi_{2,w_1}, f_{2,w_1}\ra
= \int_{\la \sigma \ra \bs H_{2,w_1}}\lp\chi_{\pi_{2,w_1}}(h_2) + \chi_{\pi_{2,w_1}}(\sigma h_2)\rp f_{2,w_1}(h_2)\;dh_2.
\]
In particular,
$\chi_{\pi_{2,w_1}}(h_2) + \chi_{\pi_{2,w_1}}(\sigma h_2) = 0$ for all elliptic elements $h_2$ in $H_{2,w_1}$.
It follows that if $f_{w_1}$ is \emph{any} elliptic function on $G_{w_1}$, then $\la \pi_{2,w_1}, f_{2,w_1}\ra = 0$
for any function $f_{2,w_1}$ on $H_{2,w_1}$ which matches $f_{w_1}$.  The proposition follows.
\end{proof}
\begin{lemma}\label{lemma:exclH11dim}
Let $\chi$ be a character of $\idc{E}$ such that $\chi \neq {}^\sigma \chi$.
Suppose the one dimensional representation $\pi_1 = 1_2\otimes_1 \chi$ of $\mb{H_1}(\Af)$
contributes to an $\ve$-trace identity.  Then, no 
automorphic representation of $\mb{H_2}(\Af)$ contributes to the trace identity. 
\end{lemma}
\begin{proof}
Let $S$ be the finite set of bad places for $\pi_1$.
By the Chebotarev Density Theorem, there exists a place $v \notin S$ which splits 
into two places $v_1, v_2$ of $E$ such that $\chi_{v_1} \neq \chi_{v_2}$.  
The F-H class in $\hat{H_1}$ parametrizing $\pi_{1,v}$ lifts to the conjugacy class in
$\GSp(2, \CC)$ represented by
\[
\hat{g} = \dmfour{\chi_{v_1}\nu_v^{1/2}}{\chi_{v_1}\nu_v^{-1/2}}
{\chi_{v_2}\nu_v^{1/2}}{\chi_{v_2}\nu_v^{-1/2}}.
\]
Suppose a representation $\pi_2 = \tau\otimes_2 \mu$ of $\mb{H_2}(\Af)$ contributes to the trace identity,
where $\tau$ is an automorphic representation of $\GL(2, \Ae)$ and $\mu$ is a character of
$\idc{F}$ such that $\omega_{\tau} = \mu\circ\N_{E/F}$.  Then,
$\tau_{v_1}$ must be parametrized by the conjugacy class in $\GL(2, \CC)$ containing
either $\diag(\chi_{v_1}\nu_v^{1/2},\chi_{v_2}\nu_v^{-1/2})$ or
$\diag(\chi_{v_2}\nu_v^{1/2},\chi_{v_1}\nu_v^{-1/2})$.
In either case, since $\chi_{v_1} \neq \chi_{v_2}$,
$\tau_{v_1}$ is not unitarizable, which is a contradiction because all local components
of a unitary automorphic representation of $\GL(2, \Ae)$ are unitarizable.
\end{proof}
Let $\rR_{E/F}\GL(2)$ be the $F$-group obtained from $\GL(2)$ on restricting scalars from $E$ to $F$.
\indexs{R@$\rR_{E/F}$!GL2@$\GL(2)$}%
It is an $\ve$-twisted endoscopic group of $\GL(4)$, and its group of $\Af$-points is $\GL(2, \Ae)$.
For an automorphic  representation $\tau_E$ of $\GL(2, \Ae)$, let $\pi(\tau_E)$ denote the automorphic
representation of $\GL(4, \Af)$ which is obtained from $\tau_E$ via the twisted endoscopic lifting from
$\rR_{E/F}\GL(2)$ to $\GL(4)$ (see \cite[Sect. 3.6]{AC}).  The representation $\pi(\tau_E)$ is cuspidal if and only
if $\tau_E$ is cuspidal and $\tau_E\ncong {}^\sigma\tau_E$.

In the case of $\mb{H_1}$, we in general denote a (quasi-)packet by $\pi_1$, i.e. without braces.
For a (quasi-)packet $\pi_1$ of $\mb{H_1}(\Af)$, let $\xi_1^*(\pi_1)$ denote the (quasi-)packet of $\mb{G}(\Af)$
which is the lift of $\pi_1$.
\begin{lemma}\label{lemma:H1GL4}
Let $\tau$ be a cuspidal or one dimensional automorphic representation of $\GL(2, \Af)$.  Let
$\chi$ be a character of $\idc{E}$.  Then, the $($quasi-$)$packet $\xi_1^*(\tau\otimes_1 \chi)$ of $\mb{G}(\Af)$ 
lifts to the automorphic representation $\pi(\chi B_{E/F}\tau)$ of $\GL(4, \Af)$.
\end{lemma}
\begin{proof}
Let $S$ be the union of the sets of bad places for $\tau, \chi$.  Hence, $\tau_v, \chi_v$ are unramified
for all $v \notin S$.
For each $v \notin S$, let $t(\tau_v) = \diag(a_v, b_v)$ 
be a diagonal matrix in $\GL(2, \CC)$ whose conjugacy class parametrizes $\tau_v$.
The global datum in $\widehat{\GL(4)} = \GL(4, \CC)$ parametrizing $\pi(\chi B_{E/F}\tau)$ is as follows:
\begin{itemize}
\item
Suppose $v \notin S$ is prime in $E$.  Let $\chi_v'$ be an unramified character of $F_v^\times$ such that 
$\chi_v = \chi'_v\circ\N_{E/F}$.  Then,
\[
\mc{C}\lp\pi(\chi B_{E/F}\tau), S, \GL(4, \CC)\rp_v
\ni \dmfour{\chi'_va_v}{\chi'_vb_v}{-\chi'_va_v}{-\chi'_vb_v}.
\]
\item
Suppose $v \notin S$ splits into two places $v_1, v_2$ of $E$.  Then,
\[
\mc{C}\lp\pi(\chi B_{E/F}\tau), S, \GL(4, \CC)\rp_v
\ni \dmfour{\chi_{v_1}a_v}{\chi_{v_1}b_v}{\chi_{v_2}a_v}{\chi_{v_2}b_v}.
\]
\end{itemize}
The global datum $\mc{C}\lp\xi_1^*(\tau\otimes_1 \chi), S, \GSp(2, \CC)\rp$ is described by
Corollary \ref{corollary:FHliftH1}.
The proposition follows on comparing the datum $\mc{C}\lp\pi(\chi B_{E/F}\tau), S, \GL(4, \CC)\rp$ with the natural image of
\dt $\mc{C}\lp\xi_1^*(\tau\otimes_1 \chi), S, \GSp(2, \CC)\rp$ in $\prod_{v \notin S}\GL(4, \CC)$.
\end{proof}
\begin{lemma}\label{lemma:chiBtau}
Let $\tau$ be an automorphic representation of $\GL(2, \Af)$.  Let $\chi$ be a character of $\idc{E}$.
The automorphic representation $\pi = \chi B_{E/F}\tau$ of $\GL(2, \Ae)$ is invariant under $\Gal(E/F)$ if and only
if there exists a $($possibly trivial$)$ quadratic character $\ve'$ of $\idc{F}$ such that:
{\rm (i)} ${}^\sigma \chi / \chi = \ve' \circ \N_{E/F}$, and {\rm (ii)} $\tau$ is equivalent to $\ve'\tau$ or $\ve'\ve\tau$.
\end{lemma}
\begin{proof}
Suppose $\pi \cong {}^\sigma \pi$.  Then, the central character $\omega_\pi = \chi^2\cdot(\omega_\tau\circ\N_{E/F})$
of $\pi$ must be equal to 
$\omega_{{}^\sigma\! \pi} = {}^\sigma\!\chi^2\cdot(\omega_{\tau} \circ \N_{E/F})$.
Hence,
$\chi^2 = {}^\sigma\!\chi^2$, which by Hilbert 90 implies that ${}^\sigma \chi/\chi = \ve'\circ\N_{E/F}$ for some (possibly trivial)
quadratic character $\ve'$ of $\idc{F}$.  We have:
\[
\chi B_{E/F}(\ve'\tau) = \chi\cdot\lp \ve'\circ\N_{E/F}\rp B_{E/F}\tau = {}^\sigma\chi\;B_{E/F}\tau \cong \chi B_{E/F}\tau.
\]
By \cite{L1}, \cite{F4}, $B_{E/F}(\ve'\tau) \cong B_{E/F}\tau$ implies that $\tau$ is equivalent to $\ve'\tau$ or $\ve'\ve\tau$.

Suppose ${}^\sigma \chi/\chi = \ve'\circ \N_{E/F}$ for some quadratic character $\ve'$ of $\idc{F}$, and $\tau$ is
equivalent to $\ve'\tau$ or $\ve'\ve\tau$.
Then, ${}^\sigma \pi = \chi\cdot\lp \ve'\circ\N_{E/F}\rp B_{E/F}\tau = \chi B_{E/F}(\ve'\tau)$, which is equivalent to either
$\chi B_{E/F}\tau$ or $\chi B_{E/F}(\ve\tau)$.  Since $B_{E/F}(\ve\tau) \cong B_{E/F}\tau$, the lemma follows.
\end{proof}
\begin{prop}\label{prop:H1allcuspidal}
Let $\tau$ be a cuspidal, non-$E$-monomial, automorphic representation of $\GL(2, \Af)$,
$\chi$ be a character of $\idc{F}$, such that: 
There does not exist a character $\ve'$ of $\idc{F}$ such that
${}^\sigma\chi/\chi = \ve'\circ\N_{E/F}$ and $\tau \cong \ve'\tau$.
Then, each member of $\{\Pi\} = \xi_1^*(\tau\otimes_1 \chi)$ is a cuspidal automorphic representation of $\mb{G}(\Af)$.
\end{prop}
\begin{proof}
By Lemma \ref{lemma:H1GL4}, $\{\Pi\}$ lifts to $\pi(\chi B_{E/F}\tau)$ of $\GL(4, \Af)$.
Since $\tau$ is not \dt $E$-monomial, $\chi B_{E/F}\tau$ is cuspidal.
By Lemma \ref{lemma:chiBtau} and \cite[Sect. 3.6]{AC}, we conclude that 
the automorphic representation \dt $\pi(\chi B_{E/F}\tau)$ is cuspidal.  

Suppose a representation $\Pi$ in $\{\Pi\}$ is not cuspidal.  Then, $\Pi$ is an irreducible constituent of a parabolically
induced representation $I$ of $\mb{G}(\Af)$.  By the results in \cite[Sect. V. 1]{F1}, $I$ lifts to a
parabolically  induced representation of $\GL(4, \Af)$.
By the strong multiplicity one (or rigidity) theorem for $\GL(4)$, we conclude that the cuspidal representation 
$\pi(\chi B_{E/F}\tau)$ is an irreducible constituent of a parabolically induced representation of $\GL(4, \Af)$, 
which is a contradiction.
\end{proof}
The case where ${}^\sigma \chi/\chi = \ve'\circ\N_{E/F}$ and $\tau\cong\ve'\tau$ will be addressed by
Proposition \ref{prop:H1E'monomial} in Section \ref{sec:somegloballifts}.

\subsubsection{Case of $\mb{H_2}$}
For any character $\theta$ of $\idc{E}$, recall that $\pi(\theta)$ denotes the $E$-monomial representation 
of $\GL(2, \Af)$ associated with $\theta$.
\begin{claim}\label{claim:nonbcH2noI}
Let $\tau$ be an automorphic representation of $\GL(2, \Ae)$ such that $\tau \neq {}^\sigma \tau$
and $\omega_{\tau} = \mu\circ\N_{E/F}$ for some character $\mu$ of $\idc{F}$.
Suppose the representation $\pi_2 = \tau\otimes_2 \mu$ of $\mb{H_2}(\Af)$ contributes to an $\ve$-trace
identity.  Then, no discrete spectrum packet
of $\mb{H_1}(\Af)$ of the form $\pi(\theta)\otimes_1 \chi$, where $\theta, \chi$ are characters
of $\idc{E}$,
contributes to the trace identity.
\end{claim}
\begin{proof}
Let $S$ be the set of bad places for $\pi_2$.
For a place $v \notin S$ which does not split in $E$, $\mu_v$ is an unramified character of $F_v^\times$, and
$\tau_v$ is an induced representation of $\GL(2, E_v)$ of the form
$I(\eta_v\circ\N_{E/F}, (\mu_v\eta_v^{-1})\circ\N_{E/F})$, where $\eta_v$ is an unramified character of $F_v^\times$.
From Section \ref{FHclass}, we see that the F-H class in ${}^LH_2$
parametrizing $\pi_{2,v}$ lifts to the conjugacy class in $\hat{G}$
represented by 
\[
\hat{g}_v = \diag(\eta_v, -\eta_v, -\mu_v\eta_v^{-1}, \mu_v\eta_v^{-1}).
\]
If $v \notin S$ splits into two places $v_1, v_2$ of $E$, then the F-H class parametrizing $\pi_{2,v}$ lifts to
the conjugacy class in $\hat{G}$ represented by
\[
\hat{g}_v = \left[ t(\tau_{v_1}), t(\tau_{v_2})\right],
\]
where $t(\tau_{v_i})$ ($i = 1, 2$) is a diagonal matrix whose conjugacy class in $\GL(2, \CC)$ 
\hyphenation{paramet-rizes}
parametrizes the unramified representation $\tau_{v_i}$ of $\GL(2, E_{v_i}) = \GL(2, F_v)$.

A lifting of automorphic representations of
$\GL(2, \Ae)$ to $\ve$-invariant automorphic representations of $\GL(4, \Af)$ is
described in \cite[Sect. 3.6]{AC}.
In particular, if an automorphic representation $\pi$ of $\GL(2, \Ae)$ satisfies $\pi \neq {}^\sigma \pi$, then
it lifts to a cuspidal automorphic representation of $\GL(4, \Af)$.

Viewing 
$\GSp(2, \CC)$ as a subgroup of $\GL(4, \CC)$, the 
global datum  in \dt $\prod_{v \notin S}\GL(4, \CC)$ represented by $\left[\hat{g}_v\right]_{v \notin S}$
parametrizes the automorphic $\GL(4, \Af)$-module $\Pi$ which is the lift of $\tau$.  
Since $\tau \neq {}^\sigma \tau$, $\Pi$ is cuspidal.

On the other hand, $\pi(\theta)\otimes_1 \chi$ lifts to a packet of $\mb{G}(\Af)$
which in turn lifts to the induced representation 
$I_{(2,2)}(\pi(\theta\chi),\pi(\theta\;{}^\sigma\!\chi))$ of $\GL(4, \Af)$.
By the rigidity theorem for $\GL(4)$, the claim follows.
\end{proof}
Recall that we have fixed two finite places $w_1, w_2$ of $F$ which are prime in $E$.
\begin{lemma}\label{lemma:bcH2ve}
Let $\tau$ be a discrete spectrum representation of $\GL(2, \Af)$ with the property
that 
$\tau_{w_1}$ and $\tau_{w_2}$ are elliptic.
Suppose the discrete spectrum representation $B_{E/F}\tau\otimes_2 \omega_\tau\ve$ of
$\mb{H_2}(\Af)$ contributes to an $\ve$-trace identity 
with respect to matching 
$f, f_i$ $(i = 1, 2)$ in
$\ELw{\mb{G}(\Af), \omega}$, $\ELw{\mb{H}_i(\Af), \omega}$, respectively.
Then, no discrete spectrum packet of $\mb{H_1}(\Af)$ contributes to the trace identity.
\end{lemma}
\begin{proof}
Suppose a packet $\pi_1 = \pi\otimes_1 \chi$ of $\mb{H_1}(\Af)$ 
contributes to the trace identity, where $\pi$
is an automorphic representation of $\GL(2, \Af)$, and $\chi$ is a character of 
$\idc{E}$.  
There are two cases to consider:  The case where $\chi = {}^\sigma \chi$
and the case where $\chi \neq {}^\sigma \chi$.

Suppose $\chi = {}^\sigma \chi$.  Then, $\chi = \mu\circ\N_{E/F}$ for some character
$\mu$ of $\idc{F}$, and $\pi_1 = \mu\pi\otimes_1 1$.  By Proposition \ref{prop:globalbc},
the representation
$B_{E/F}\lp\mu\pi\rp\otimes_2 \omega_{\mu\pi}$ of $\mb{H_2}(\Af)$
must also contribute to the trace identity.
By Corollary \ref{corollary:contribH2discrete}, $B_{E/F}\tau\otimes_1 \omega_\tau\ve$ must be
equivalent to $B_{E/F}\lp\mu\pi\rp\otimes_1\omega_{\mu\pi}$.
The theory of base change for $\GL(2)$ (see \cite{L1}, \cite{F4}) implies that
 $\tau$ is equivalent to either $\mu\pi$ or $\ve\mu\pi$.  In either case,
$\omega_\tau\ve$ is equal to $\omega_{\mu\pi}\ve \neq \omega_{\mu\pi}$, whence a contradiction.

Suppose  $\chi \neq {}^\sigma \chi$.  If $\pi_1 = \pi\otimes_1 \chi$ contributes to the trace
identity, then $\pi_{w_1}, \pi_{w_2}$ must be elliptic, or else the condition that the functions
$f_{1, w_1}, f_{1, w_2}$ are elliptic would imply that $\la \pi_1, f_1 \ra = 0$.
If $\pi_{w_1}, \pi_{w_2}$ are elliptic,
then by Proposition \ref{prop:nonbcH1} and Claim \ref{claim:nonbcH2noI}
no automorphic representation of $\mb{H_2}(\Af)$ may contribute to the trace identity.
The lemma follows.
\end{proof}
\begin{prop}\label{prop:nonbcH2}
Let $\mu$ be a character of $\idc{F}$.
Let $\tau$ be a discrete spectrum automorphic representation of $\GL(2, \Ae)$ 
with the following properties:
\begin{itemize}
\item
$\tau \neq {}^\sigma \tau$;
\item
$\omega_\tau = \mu\circ\N_{E/F}$;
\item
$\tau_{w_1}$ and $\tau_{{w_2}}$ are elliptic.
\end{itemize}
Suppose the automorphic representation $\pi_2 = \tau\otimes_2 \mu$ of $\mb{H_2}(\Af)$
contributes to an $\ve$-trace identity with respect to matching test functions
with elliptic local components at $w_1, w_2$.
Then, no discrete spectrum packet of $\mb{H_1}(\Af)$ contributes to the trace identity.
\end{prop}
\begin{proof}
Suppose a packet $\pi_1$ of $\mb{H_1}(\Af)$ contributes to the trace identity.
Since $\tau \neq {}^\sigma\tau$, by Claim \ref{claim:nonbcH2noI} the packet
$\pi_1$ is not of the form $\pi(\theta)\otimes_1 \chi$ for any characters $\theta, \chi$ of $\idc{E}$.
By Proposition \ref{prop:globalbc} and Lemma \ref{lemma:contribH2}, $\pi_1$ is not of the form
$\pi\otimes_1 1$, where $\pi$ is an automorphic representation of $\GL(2, \Af)$.

Consequently, by Lemma \ref{lemma:nbcH1stable}, 
the packet $\{\Pi\} = \otimes_v \{\Pi_v\}$ of $\mb{G}(\Af)$ which contributes to the trace
identity must be stable; namely, every member of $\{\Pi\}$ occurs
with multiplicity one in the discrete spectrum of $\mb{G}(\Af)$.

The argument we then use mirrors completely the one used in the proof of Proposition \ref{prop:nonbcH1}.
Thus, we shall content ourselves with giving only a sketch:
First, we show that there exists a test function $f_{w_2}$ on $G_{w_2}$ such that $\la \{\Pi_{w_2}\}, f_{w_2}\ra_\ve$
is nonzero.
Then, we show that either
$\pi_1$ does not exist or $\pi_{1,w_1}$ is non-elliptic.  Since by assumption the test function $f_{1, w_1}$ on $H_{1, w_1}$
is elliptic, the proposition follows.
\end{proof}
\begin{lemma}\label{lemma:exclH21dim}
Let $\chi$ be a character of $\idc{E}$ such that $\chi \neq {}^\sigma \chi$ but $\chi^2 = \mu\circ\N_{E/F}$
for some character $\mu$ of $\idc{F}$.  Suppose the 
one dimensional representation $\pi_2 = \lp\chi\circ\det\rp \otimes_2 \mu = \chi 1_{\GL(2, \Ae)}\otimes_2 \mu$ contributes
to an $\ve$-trace identity.  Then, no global packet of $\mb{H_1}(\Af)$ contributes to the trace identity.
\end{lemma}

\begin{proof}
Recall from Chapter \ref{chap:endoscopy} that
\[
\hat{H_1} = 
\lp\GL(2, \CC)\times\CC^\times \times \CC^\times\rp/\{(\diag(z, z)^{-1}, z, z) : z \in \CC^\times\},
\]
and ${}^L H_1 = \hat{H_1}\rtimes \Gal(E/F)$, where the action of $\sigma$ swaps the two $\CC^\times$-components.
Suppose a packet $\pi_1 = \tau \otimes_1 \theta$ of  $\mb{H_1}(\Af)$ contributes to the
trace identity.
By the Chebotarev Density Theorem, there exists a place $v$ of $F$ which splits in $E$ such that
$\pi_{1,v}$ and $\pi_{2,v}$ are unramified.  Suppose $v$ splits into two places $v_1, v_2$ of $E$.
Since $\chi^2$ is invariant under the action of $\Gal(E/F)$, the character $\kappa_v := \frac{\chi_{v_1}}{\chi_{v_2}}$ is 
quadratic.  The F-H class parametrizing $\pi_{2,v}$ lifts to the conjugacy class in $\hat{G}$
represented by
\begin{multline*}
\hat{g}_v = \diag(\chi_{v_1}\nu_v^{1/2}, \chi_{v_2}\nu_v^{1/2}, \chi_{v_2}\nu_v^{-1/2}, \chi_{v_1}\nu_v^{-1/2})\\
\begin{split}
& =
\dmfour{\chi_{v_1}\nu_v^{1/2}}{\chi_{v_2}\nu_v^{1/2}}
{\nu_v^{-1}\kappa_v\chi_{v_1}\nu_v^{1/2}}{\nu_v^{-1}\kappa_v\chi_{v_2}\nu_v^{1/2}}\\
& \substack{(23)\\\sim}
\dmfour{\chi_{v_1}\nu_v^{1/2}}{\chi_{v_2}\nu_v^{-1/2}}
{\chi_{v_2}\nu_v^{1/2}}{\chi_{v_1}\nu_v^{-1/2}}.
\end{split}
\end{multline*}

From Section \ref{FHclass}, the F-H class in ${}^L H_1$ parametrizing $\pi_{1,v}$ is represented by
\[
\hat{h}_v = (t(\tau_v), \theta_{v_1}, \theta_{v_2})\rtimes 1,
\]
where $t(\tau_v)$ is a diagonal matrix in $\GL(2, \CC)$ whose conjugacy class parametrizes $\tau_v$.
Note that the Galois component of $\hat{h}_v$ is trivial because $\mb{H_1}$ is split over $F_v$.
Since $\pi_1$ contributes to the trace identity, $\hat{h}_v$
lifts via the embedding $\xi_1 : {}^L H_1 \rightarrow \hat{G}$ (see Chapter \ref{chap:endoscopy}) to
the conjugacy class of $\hat{g}_v$ in $\hat{G}$.  In other words, the element
\[
\xi_1(\hat{h}_v) = \lp\lmx \theta_{v_1}t(\tau_v) &\\ & \theta_{v_2}t(\tau_v)\rmx\rp \in \hat{G}
\]
is conjugate to $\hat{g}_v$.
Hence, $\pi_{1,v}$ is parametrized by the conjugacy class of one of the following types of
elements in ${}^L {H_1}$:
\begin{enumerate}
\item
$\lp\diag(\chi_{v_1}\nu_v^{1/2}, \chi_{v_2}\nu_v^{1/2}), \eta_{1}, \eta_{2}\rp\rtimes 1$,
where the unordered set $\{\eta_{1}, \eta_{2}\}$ is equal to $\{1, \kappa_v\nu_v^{-1}\}$
(Notice that there is no minus sign in front of the $1/2$ in $\chi_{v_2}\nu_v^{1/2}$);
\item
$\lp \diag(\zeta_{1}\nu_v^{1/2}, \zeta_{2}\nu_v^{-1/2}), \eta_{1}, \eta_{2}\rp\rtimes 1$, 
where the unordered set $\{\zeta_1, \zeta_2\}$ is equal to \dt $\{\chi_{v_1}, \chi_{v_2}\}$, and $\{\eta_1, \eta_2\} = \{1, \kappa_v\}$.
\end{enumerate}

In the first case, assume without loss of generality that $\eta_1 = 1$ and $\eta_2 = \kappa_v\nu_v^{-1}$.
Then, there exists a quasi-character $\lambda_v$ of $F_v^\times$ such that 
$\tau_{v} = \lambda_v \nu_v^{1/2}I(\chi_{v_1}, \chi_{v_2})$,
$\theta_{v_1} = \lambda_v^{-1}$, and $\theta_{v_2} = \lambda_v^{-1}\kappa_v\nu_v^{-1}$.
Since $\tau$ is a unitary automorphic representation of the group $\GL(2, \Af)$, $\tau_v$ must be unitarizable, which implies
that $\lambda_v$ is non-unitary.  
Since $\theta_{v_1} = \lambda_v^{-1}$ and $\theta$ is unitary, we have a contradiction.

In the second case, assume without loss of generality that $\eta_1 = 1$ and $\eta_2 = \kappa_v$.
Then, there exists a quasi-character $\lambda_v$ of $\idc{F_v}$ such that $\tau_{v} = \lambda_v I(\zeta_1\nu_v^{1/2}, \zeta_2\nu_v^{-1/2})$,
$\theta_{v_1} = \lambda_v^{-1}$, and $\theta_{v_2} = \lambda_v^{-1}\kappa_v$.
The condition $\chi \neq {}^\sigma \chi$ implies that $\zeta_1 \neq \zeta_2$;
hence, $\lambda_v$ must be non-unitary for $\tau_v$ to be unitarizable.  But then $\theta_{v_1} = \lambda_v^{-1}$ is
non-unitary, which is a contradiction.
\end{proof}

For an automorphic representation $\pi_2$ of the group $\mb{H_2}(\Af)$, let $\xi_2^*(\pi_2)$ denote the (quasi-)packet of
$\mb{G}(\Af)$ which is the lift of $\pi_2$.
Recall that, for an automorphic representation $\tau$ of $\GL(2, \Ae)$, $\pi(\tau)$ denotes the automorphic
representation of $\GL(4, \Af)$ obtained from $\tau$ via the twisted endoscopic lifting from $\rR_{E/F}\GL(2)$ to $\GL(4)$.
\begin{lemma}\label{lemma:H2GL4}
Let $\tau$ be an automorphic representation of $\GL(2, \Ae)$, $\mu$ be a character of $\idc{F}$, such that
$\tau \ncong {}^\sigma \tau$ and $\omega_\tau = \mu\circ\N_{E/F}$.
The $($quasi-$)$packet $\{\Pi\} = \xi_2^*(\tau\otimes_2\mu)$ lifts to the automorphic representation $\pi(\tau)$ of $\GL(4, \Af)$.
\begin{remark}
The (quasi-)packets $\xi_2^*(\tau\otimes_2\mu)$ and $\xi_2^*(\tau\otimes_2\mu\ve)$
are inequivalent, but they both lift to $\pi(\tau)$.
\end{remark}
\end{lemma}
\begin{proof}
Let $S$ be the set of bad places for $\tau\otimes_2 \mu$.
At a place $v \notin S$ which is prime in $E$, $\tau_v$ is the induced representation $I(\alpha_v, \beta_v)$ for some
unramified characters $\alpha_v, \beta_v$ of $E_v^\times$.  Since $\alpha_v, \beta_v$ are unramified, there exist characters
$a_v, b_v$ of $F_v^\times$ such that $\alpha_v = a_v\circ\N_{E/F}$ and $\beta_v = b_v\circ\N_{E/F}$.
The condition $\omega_\tau = \mu\circ\N_{E/F}$ implies that $a_v^2 b_v^2 = \mu_v^2$.  Consequently,
$a_vb_v(\vp_v) = \pm \mu_v(\vp_v)$ (Here, $\vp_v$ is a fixed uniformizer of $F_v$).
The F-H class $\mc{C}(\pi(\tau), S, \GL(4, \CC))_v$ parametrizing $\pi(\tau)_v$ is represented by
\[
\hat{h}_v = \dmfour{a_v}{b_v}{-a_v}{-b_v}  = \dmfour{a_v}{\pm\mu_v/a_v}{-a_v}{\mp \mu_v/a_v}\in \GL(4, \CC).
\]
From Corollary \ref{corollary:FHliftH2}, the conjugacy class in $\GSp(2, \CC)$ parametrizing the unramified component of $\{\Pi_v\}$
is represented by
\[
\hat{g}_v = \lp\lsm
&\mu_v^{1/2}&&\\ \mu_v^{-1/2}a_v^2&&&\\&&&\mu_v^{-1/2} b_v^2\\&&\mu_v^{1/2}&
\rsm\rp
 \sim \dmfour{\mu_v/a_v}{-\mu_v/a_v}{-a_v}{a_v}.
\]
Hence, the image of $\hat{g}_v$ in $\GL(4, \CC)$ is conjugate to $\hat{h}_v$.

Suppose $v \notin S$ splits into two places $v_1, v_2$ of $E$.  Then, $\GL(2, E_{v_i}) = \GL(2, F_v)$ ($i = 1, 2$).
Let $\tau_i = \tau_{v_i}$.  Let $t(\tau_i)$ be a diagonal matrix in $\GL(2, \CC)$ parametrizing the unramified representation $\tau_i$
of $\GL(2, F_v)$.
The conjugacy class $\mc{C}(\pi(\tau), S, \GL(4, \CC))_v$, which parametrizes $\pi(\tau)_v$, is represented by the element
$\hat{h}_v = \blockdiag(t(\tau_1), t(\tau_2)),$ while the conjugacy class in $\GSp(2, \CC)$ parametrizing the unramified
component of $\{\Pi_v\}$ is represented by $\hat{g}_v = [t(\tau_1), t(\tau_2)]$.  The image of $\hat{g}_v$ in
$\GL(4, \CC)$ is conjugate to $\hat{h}_v$.

The lemma follows.
\end{proof}
\begin{prop}\label{prop:H2allcuspidal}
Let $\tau$ be an automorphic representation of $\GL(2, \Ae)$, $\mu$ be a character of $\idc{F}$, such that
$\tau \ncong {}^\sigma \tau$ and $\omega_\tau = \mu\circ\N_{E/F}$.  Then, every automorphic representation in
$\{\Pi\} = \xi_2^*(\tau\otimes_2\mu)$ is cuspidal.
\end{prop}
\begin{proof}
By lemma \ref{lemma:H2GL4}, the (quasi-)packet $\{\Pi\}$ lifts to the automorphic representation $\pi(\tau)$ of $\GL(4, \Af)$.
By the assumption that $\tau \ncong {}^\sigma \tau$, the representation $\pi(\tau)$ is cuspidal.

Suppose a representation $\Pi$ in $\{\Pi\}$ is not cuspidal.  Then, it is an irreducible constituent of a parabolically  induced
representation, which by \cite[Sect. V. 1]{F1} lifts to a parabolically induced representation of $\GL(4, \Af)$.
By the rigidity theorem for $\GL(4)$, we obtain the contradiction that the cuspidal $\pi(\tau)$ is an irreducible constituent
of a parabolically induced representation.
\end{proof}
\section{Some Global Lifting Results}\label{sec:somegloballifts}
We now list several additional global lifting results.  
All test functions are assumed to be matching and have elliptic local components at two fixed places $w_1, w_2$
which are prime in $E$.

The stable quasi-packet $\{L(\zeta\nu, \nu^{-1/2}\pi)\}$ is defined in
Section \ref{sec:arthurunstable},
where $\zeta$ is a nontrivial quadratic character of $\idc{F}$, and $\pi$ is a cuspidal, 
$\zeta$-invariant, automorphic representation of $\GL(2, \Af)$.
Recall that $1_2$ denotes the trivial representation of $\GL(2, \Af)$.
\begin{claim}\label{claim:globalsoudryH1}
Let $\chi$ be a character of $\idc{E}$ such that $\chi \neq {}^\sigma \chi$.
Let $\mu$ be a character of $\idc{F}$.
The following table holds:
\begin{equation}
\begin{tabular}{|c|c|c|}
\hline
$G$ & $H_1$ & $H_2$\\
\hline
$\left\{L\lp\ve\nu, \nu^{-1/2}\mu\pi(\chi)\rp\right\}$&
$\mu 1_2\otimes_1 \chi$ & 
\\
\hline
\end{tabular}
\end{equation}
\end{claim}
\begin{proof}
That there is no contribution from $\mb{H_2}$ follows from Lemma \ref{lemma:exclH11dim}.
The rest follows on comparing the lift of the global datum of $\mu 1_2\otimes_1 \chi$
with the global datum of $\left\{L\!\lp\ve\nu, \nu^{-1/2}\mu\pi(\chi)\rp\right\}$.
\end{proof}

\begin{claim}\label{claim:globalsoudryH2}
Let $\chi$ be a character of $\idc{E}$ such that $\chi \neq {}^\sigma \chi$ but
$\chi^2 = {}^\sigma \chi^2$.  Then, ${{}^\sigma\chi}/{\chi} = \ve'\circ\N_{E/F}$ for some quadratic
character $\ve'$ of $\idc{F}$.  
Let $\zeta = \ve'$ or $\ve'\ve$.  In particular, $\chi^2 = (\chi|_{\Af^\times}\cdot\zeta)\circ \N_{E/F}$.
The following table holds:
\begin{equation}
\begin{tabular}{|c|c|c|}
\hline
$G$ & $H_1$ & $H_2$\\
\hline
$\left\{L\lp\ve\zeta\nu, \nu^{-1/2}\pi(\chi)\rp\right\}$& & $\chi 1_{\GL(2, \Ae)}\otimes_2 \chi|_{\Af^\times}\cdot\zeta$
\\
\hline
\end{tabular}
\end{equation}
\end{claim}
\begin{proof}
That there is no contribution from $\mb{H_1}$ follows from Lemma \ref{lemma:exclH21dim}.
The rest follows from the comparison of global data.
\end{proof}
In what follows, the terms which appear in the $\ve$-trace identity
often have nontrivial coefficients (as summarized in Section \ref{sec:tfsummary}).
So, instead of using tables to describe the global lifting, we record explicitly the $\ve$-trace identities.
There is often possible cancellation among the coefficients in these trace identities, but we retain them in order to emphasize how 
they arise from the trace formula.

Note that if the representations involved do not have elliptic components at $w_1, w_2$, then
the $\ve$-trace identity is in the trivial form $0 = 0$ due to the restriction on the test functions.

In each of the cases examined, we let $S$ be a finite set of places
containing all the bad places for the representations involved.
That is, $S$ is a finite set containing the union of $V - \Vr$ and $\{{w_1, w_2}\}$
such that the local components of the representations are unramified at all $v \notin S$.
The components at $v \notin S$ of the test functions are assumed to be spherical.
For any ad\`elic object (automorphic representation, test function, trace, \dots, etc),
subscript $S$ denotes the tensor product of local components over the places in $S$.

First, we consider the cases where the packet $\pi_1$ of $\mb{H_1}(\Af)$ has the form
$\pi(\theta)\otimes_1\chi$, where $\theta, \chi$ are characters of $\idc{E}$.  
The following claims follow directly from the results in Section \ref{sec:tfsummary} and Claim \ref{claim:nonbcH2noI}, and
we skip their proofs.

\subsubsection{Monomial Representations of $\mb{H_1}(\Af)$}\label{sec:monomialglobal}
For a character $\theta$ of $\idc{E}$, let $\bar{\theta}$ denote ${}^\sigma \theta$.  
For the $S$-part $\{\Pi\}_S$ of a  (quasi-)packet $\{\Pi\}$ of $\mb{G}(\Af)$, put 
$\la \{\Pi\}_S, f\ra := \sum_{\Pi' \in \{\Pi\}_S} m(\Pi')\la \Pi', f\ra_\ve$, 
where $m(\Pi')$ is the multiplicity of $\Pi'$ in the discrete spectrum of $\mb{G}$.

Let $\theta, \chi$ be characters of $\idc{E}$.
\begin{enumerate}
\item
  Suppose $\theta \neq \bar{\theta}$, $\chi \neq \bar{\chi}$, and
  $\frac{\bar{\theta}}{\theta} \neq \frac{\bar{\chi}}{\chi}$, $\frac{\chi}{\bar{\chi}}$.
  Let $\{\Pi\}$ be the packet $[\pi(\theta\chi), \pi(\theta\bar{\chi})]$ of $\mb{G}(\Af)$.  It lifts to
  the induced representation $I_{(2, 2)}(\pi(\theta\chi), \pi(\theta\bar{\chi}))$ of $\GL(4, \Af)$.
  Those members of $\{\Pi\}$ which are unramified at places outside of $S$ form
  $\{\Pi\}_S$, which is equal to the tensor product
 \[\lp\otimes_{v \in S}\{\Pi_v^+, \Pi_v^-\}\rp\otimes\lp\otimes_{v \notin S} \Pi_v^+\rp.\]
  The multiplicity in the discrete spectrum of $\mb{G}$ of each $\Pi' \in \{\Pi\}_S$ is 
  given by $m(\Pi') = \frac{1}{2}\lp 1 + (-1)^{n(\Pi')}\rp$, 
  where $n(\Pi')$ is the number of places $v$ for which $\Pi'_v = \Pi_v^-$.
  The following identity holds for matching functions
  $f$ in $\ELw{\mb{G}(\Af), \omega}$ and $f_1$ in $\ELw{\mb{H_1}(\Af), \omega}$ whose components
  at all $v \notin S$ are spherical:
  \begin{multline}\label{eq:monomial1}
    \sum_{\Pi' \in \{\Pi\}_S} m(\Pi') \la \Pi', f\ra_{\ve, S}\\
    \shoveleft{ =\frac{1}{2}\cdot\frac{1}{2}
      \Big[ 
      \la \pi(\theta) \otimes_1 \chi, f_{1}\ra_S + 
      \la \pi(\theta) \otimes_1 \bar{\chi}, f_{1} \ra_S}\\
      + \la \pi(\chi)\otimes_1 \theta, f_{1}\ra_S +
      \la \pi(\chi)\otimes_1 \bar{\theta}, f_{1}\ra_S
      \Big]\\
      =\frac{1}{2} \la \pi(\theta)\otimes_1 \chi, f_{1} \ra_S
      + \frac{1}{2} \la \pi(\chi) \otimes_1 \theta, f_{1} \ra_S.
  \end{multline}
\item
  Suppose $\theta \neq \bar{\theta}, \chi \neq \bar{\chi}$, but $\theta\chi = \mu\circ\N_{E/F}$ for some character
  $\mu$ of $\idc{F}$.
  Then, $\pi(\theta), \pi(\chi)$ are cuspidal, and 
  \[\pi(\theta)\otimes_1 \chi = \pi(\chi)\otimes_1 \bar{\theta}, 
  \quad\pi(\chi)\otimes_1\theta = \pi(\theta)\otimes_1\bar{\chi}.\]
  \begin{enumerate}
    \item
      Suppose $\chi^2 \neq \ol{\chi^2}$.  Then, $\pi\!\lp{\bar{\chi}}/{\chi}\rp$ is cuspidal, and
      the induced representation $\pi\!\lp{\bar{\chi}}/{\chi}\rp\rtimes \mu$ is irreducible.
      The following holds:
      \begin{multline}
        \frac{1}{2}\la \pi\!\lp \frac{\bar{\chi}}{\chi}\rp \rtimes \mu, f\ra_{\ve, S}\\
        \begin{split}&= \frac{1}{2}
        \Bigg[ \frac{1}{2} \la\pi(\theta)\otimes_1 \chi, f_1\ra_S + 
          \frac{1}{2}\la \pi(\theta)\otimes_1\bar{\chi}, f_1\ra_S\Bigg]\\
        &= \frac{1}{2}\la \pi(\theta)\otimes_1\chi, f_1\ra_S.\end{split}
      \end{multline}
      At a place $v$ of $F$, prime in $E$, where $\theta_v \neq \bar{\theta}_v$, the representation $\pi(\theta_v)$ is cuspidal.
      Consequently, the $\ve_v$-twisted Harish-Chandra character of the induced representation 
      $\pi\!\lp{\bar{\chi}_v}/{\chi_v}\rp \rtimes \mu_v$ is not identically zero on the elliptic regular set in $G_v$.
    \item
      Suppose $\chi^2 = \ol{\chi^2}$.  Then, ${\bar{\chi}}/{{\chi}} = \ve'\circ\N_{E/F}$
      for some quadratic character $\ve'$ of $\idc{F}$.  Since by assumption $\chi \neq \bar{\chi}$,
      the character $\ve'$ is nontrivial and different from $\ve$.
      Let $E'$ be the quadratic extension of $F$ associated with $\ve'$.  The condition $\theta\chi = \mu \circ\N_{E/F}$
      implies that ${\bar{\theta}}/{{\theta}} = \ve'\circ\N_{E/F}$, from which it follows that
      the $E$-monomial $\pi(\theta)$ is also $E'$-monomial (\cite{LL}).  Then, $\pi(\theta)\otimes_1 \bar{\chi}$
      is equal to $\ve'\pi(\theta)\otimes_1 \chi = \pi(\theta)\otimes_1 \chi$.  Hence, if
      $\pi(\theta)\otimes_1 \chi$ contributes to an $\ve$-trace identity, it is the only contribution from
      $\mb{H_1}$.  
      Since $\pi\!\lp{\bar{\chi}}/{\chi}\rp$ is equivalent to the induced representation $I(\ve',\ve'\ve)$, we have:
      \begin{equation}
      \frac{1}{4} \la \ve'\times\ve'\ve\rtimes\mu, f\ra_{\ve, S}
      = \frac{1}{2}\cdot\frac{1}{2}\la \pi(\theta)\otimes_1 \chi, f_1\ra_S.
      \end{equation}
      In particular, at a place $v$ of $F$, prime in $E$, where $\theta_v \neq \bar{\theta}_v$,
      the $\ve_v$-twisted Harish-Chandra character of $\ve'\times\ve'\ve\rtimes\mu$ is not identically zero
      on the elliptic regular set in $G_v$.
  \end{enumerate}
\item
  Suppose $\theta =\mu\circ\N_{E/F}$ for some character $\mu$ of $\idc{F}$, and $\chi \neq \bar{\chi}$.
  Since \[\pi(\theta)\otimes_1 \chi = \mu\pi(1)\otimes_1 \chi = \pi(1)\otimes_1 \lp\mu\circ\N_{E/F}\rp\cdot\chi,\]
  we assume without loss of generality that $\mu = 1$.  The $E$-monomial representation $\pi(1)$ is equivalent to $I(1, \ve)$.
  From Section \ref{sec:stablespectrum},
  the representations $I(1, \ve)\otimes_1 \chi$ and $I(1, \ve)\otimes_1 \bar{\chi}$ do not contribute to the
  stable discrete part of the spectral expansion for $\mb{H_1}(\Af)$.
  The (normalizedly) induced representation $1\rtimes \pi(\chi)$ of $\mb{G}(\Af)$
  is reducible, and its irreducible constituents form the global packet
  \dt $\{\Pi\} = \otimes_{v \in V}\{\Pi_v^+, \Pi_v^-\}$ (see Section \ref{sec:finechiexpbeta}).
  The intertwining operator which appears in the term associated with $1\rtimes\pi(\chi)$ in the spectral expansion for
  $\mb{H_1}(\Af)$ acts trivially on $\Pi^+_v$ and through the scalar multiplication by $-1$ on $\Pi^-_v$.
  The following holds:
  \begin{equation}\label{eq:globalmonomialinduced}
  \frac{1}{4}\prod_{v \in S}\lp\la \Pi_v^+, f_v\ra_\ve - \la\Pi_v^-, f_v\ra_\ve\rp = \frac{1}{2}\cdot \frac{1}{2} 
  \la \pi(\chi)\otimes_1 1, f_1\ra_S.
  \end{equation}
  
\item
  Suppose $\theta = 1$, and $\chi = \mu\circ\N_{E/F}$ for some character $\mu$ of $\idc{F}$.  Then,
  \[\pi(\theta)\otimes_1 \chi = \pi(\chi)\otimes_1 \theta = I(\mu,\mu\ve)\otimes_1 1,\] which 
  does not contribute to the stable discrete part of the spectral expansion for $\mb{H_1}(\Af)$.
  If we set up an $\ve$-trace identity using the global datum parametrizing $\pi(\theta)\otimes_1 \chi$, then
  there is a nontrivial contribution from $\mb{H_2}$, namely,
  the discretely occuring representation 
  \[B_{E/F}\;\pi(\chi)\otimes_2 \mu^2\ve = \lp\mu\circ\N_{E/F}\rp I(1, 1)\otimes_2 \mu^2\ve.\]
  The following holds for matching functions:
  \begin{equation}\label{eq:1vemu}
  \frac{1}{8}\la 1 \times \ve \rtimes \mu, f\ra_{\ve, S} = 
  \frac{1}{2}\cdot\frac{1}{4}\la \lp\mu\circ\N_{E/F}\rp I(1, 1)\otimes_2 \mu^2\ve, f_2\ra_S.
  \end{equation}
  The identity \eqref{eq:1vemu}  cannot be proven using the global techniques employed thus far, for the right hand side comes from
  an induced representation, and our test functions are required to have two elliptic components.
  Rather, \eqref{eq:1vemu} follows from the local character identities deduced in Section \ref{sec:induced}.
\end{enumerate}
\subsubsection{A Special Case of Non-$E$-Monomial Representations}
Let $E'$ be a quadratic extension of $F$ different from $E$.
Let $EE'$ be the compositum of $E$ and $E'$.  It is a biquadratic extension of $F$.
Let $\sigma, \sigma'$ be the generators of $\Gal(EE'/F)$, such that $E$ is the fixed field of $\sigma'$,
$E'$ is the fixed field of $\sigma$, and $\Gal(E/F) = \la\sigma\ra$, $\Gal(E'/F)= \la \sigma'\ra$.
Let $\ve'$ be the quadratic character of $\idc{F}$ associated with $E'/F$ via class field theory.
Let $\theta$ be a character of $\idc{E'}$ such that $\theta \neq {}^{\sigma'}\!\theta$.
Let $\pi_{E'}(\theta)$ be the $E'$-monomial representation of $\GL(2, \Af)$ associated with $\theta$.
Let $\chi$ be a character of $\idc{E}$ such that $\chi \neq {}^\sigma \chi$ and ${{}^\sigma\chi}/{\chi} = \ve'\circ\N_{E/F}$.

Let $\mu$ be a character of $\idc{E'}$ which satisfies $\mu\circ\N_{EE'/E'} = \chi\circ\N_{EE'/E}$.  Then,
${{}^{\sigma'}\! \mu}/{\mu} = \ve \circ \N_{E'/F}$.
Suppose $\mu\theta \neq {}^{\sigma'}\!\!\lp\mu\theta\rp$.  This implies in particular that
${{}^{\sigma'}\!\theta}/{\theta} \neq \ve\circ\N_{E'/F}$.
Let $\pi_{E'}(\mu\theta)$ be the cuspidal $E'$-monomial representation of $\GL(2, \Af)$ associated with
$\mu\theta$.

Let $\{\Pi\}$ be the global packet of $\mb{G}(\Af)$ which consists of the irreducible constituents of
$\ve'\ve\rtimes\pi_{E'}(\mu\theta)$.
For any place $v$, let $\Pi_v^+, \Pi_v^-$ denote the two constituents of $\Pi_v$ if it is
reducible (This occurs if and only if $\ve'_v\ve_v = 1$ and $\pi_{E'}(\mu\theta)_v$ is square integrable).  
If $\Pi_v$ is irreducible, put $\Pi_v^+ := \Pi_v$ and $\Pi_v^- := 0$.
Then, $\{\Pi\} = \otimes_v \{\Pi_v^+, \Pi_v^-\}$.

Let $\{\pi_{E'}(\theta)\otimes_1 \chi\}$ denote the multi-packet of $\mb{H_1}(\Af)$
containing $\pi_{E'}(\theta)\otimes_1\chi$.
\begin{prop}\label{prop:H1E'monomial}
The following table holds:
\begin{equation}
  \begin{tabular}{|c|c|c|}
    \hline
    $G$ & $H_1$ & $H_2$\\
    \hline
    $\{\Pi\}$ &
    $\{\pi_{E'}(\theta)\otimes_1 \chi\}$ & 
    \\
    \hline
  \end{tabular}
\end{equation}
\end{prop}
\begin{proof}
We first show that $\pi_{E'}(\theta)\otimes_1 \chi$ contributes to the $\ve$-trace identity by
comparing F-H classes.
Then, we show that no automorphic representation of $\mb{H_2}(\Af)$ contributes to the trace identity.

We consider the all but finitely many odd finite places $v$ of $F$ such that: (i) the place $v$ is unramified in
both $E$ and $E'$; (ii) the characters $\theta_v, \mu_v, \chi_v$ are unramified.
For any such place $v$, the local packet $\{\Pi_v\}$
is a singleton consisting of a fully induced unramified representation, which we denote by $\Pi_v^+$.

Since any $p$-adic local field of odd residual characteristic has only one unramified quadratic
extension, there are four cases to consider:
\begin{enumerate}
\item
$\ve_v \neq 1$.
  \begin{enumerate}
  \item
    $\ve'_v = \ve_v$, i.e. $E_v = E'_v$, and $(EE')_v = E_v\oplus E_v$.
    
    We have $\Pi_v^+ = 1 \rtimes \pi_E(\mu_v\theta_v)$.  
    Since $\theta_v$ factors through $\N_{E'/F}$ via some character
    $\theta'_v$ of $F_v^\times$, the local packet
    $\pi(\theta_v)\otimes_1 \chi_v$ is equal to $I(\theta'_v,\theta'_v\ve_v)\otimes_1 \chi_v$,
    which lifts to $1\rtimes \theta'_v\pi(\chi_v) = 1 \rtimes \pi(\theta_v\chi_v)$.
    The condition $\chi\circ\N_{EE'/E} = \mu\circ\N_{EE'/E'}$ implies that,
    for any $(z_1, z_2) \in (EE')_v^\times$, $\chi_v(z_1z_2)$ is equal to $\mu_v(z_1z_2)$.  Hence,
    $\chi_v = \mu_v$, and $\pi_{1,v}$ lifts to $\Pi_v^+$.
  \item
    $\ve'_v = 1$, i.e. $E'_v = F_v \oplus F_v$, and $(EE')_v = E_v\oplus E_v$.
    
    In this case, $v$ splits into two places $v_1, v_2$ of $E'$.  For $i = 1, 2$, 
    let $\mu_i, \theta_i$ denote $\mu_{v_i}, \theta_{v_i}$, respectively.
    Then, 
    $\Pi_v^+ = \ve_v \rtimes I(\mu_1\theta_1,\mu_2\theta_2)$.

    Let $\chi_v'$ be an unramified  character of $F_v^\times$ such that $\chi_v = \chi'_v\circ\N_{E/F}$.
    Then,
    $\pi_{1,v} = I(\theta_1, \theta_2)\otimes_1 \chi_v$ lifts to
    $\ve_v \rtimes I(\chi'_v\theta_1, \ve_v\chi'_v\theta_2)$.

    The condition
    $\chi\circ\N_{EE'/E} = \mu\circ\N_{EE'/E'}$ implies that,
    for all $(z_1, z_2)$ in $(EE')_v^\times$, we have the equality
    \[\chi'_v(\N_{E/F}z_1)\chi'_v(\N_{E/F}z_2)= \mu_1(\N_{E/F}z_1)\mu_2(\N_{E/F}z_2);\]
    hence,
    $\mu_i\circ\N_{E/F} = \chi'_v\circ\N_{E/F}$ for $i = 1, 2$.  
    The condition ${{}^{\sigma'}\!\mu}/{\mu} = \ve \circ \N_{E'/F}$ implies that
    ${\mu_1}/{\mu_2} = \ve_v$.  
    Consequently,
    the set $\{\mu_1,\mu_2\}$ is equal to $\{\chi'_v, \ve_v\chi'_v\}$.
    Thus, $\pi_{1,v}$ lifts to $\ve_v \rtimes I(\mu_1\theta_1, \mu_2\theta_2)$
    or $\ve_v \rtimes I(\mu_2\theta_1, \mu_1\theta_2)$.  
    The two representations are equivalent because $\mu_2 = \ve_v\mu_1$,
    and the representation $\ve_v \rtimes I(\chi'_v \theta_1, \ve_v\chi'_v\theta_2)$ is equivalent to
    \[
    \ve_v^{-1}\rtimes \ve_v I(\chi'_v \theta_1, \ve_v\chi'_v\theta_2)
    = \ve_v \rtimes I(\ve_v\chi'_v \theta_1, \chi'_v\theta_2).
    \]
    Thus, $\pi_{1,v}$ lifts to $\Pi_v^+$.
  \end{enumerate}
\item
    $\ve_v = 1$, i.e. $E_v = F_v \oplus F_v$.
  \begin{enumerate}
  \item
    $\ve'_v \neq 1$, i.e. $[E'_v, F_v] = 2$, and $(EE')_v = E'_v \oplus E'_v$.
    
    We have $\Pi_v^+ = \ve_{v}\rtimes I(\mu'_v\theta'_v, \mu'_v\theta'_v\ve'_v)$, where
    $\mu'_v, \theta'_v$ are unramified characters of $F_{v}^\times$ such that 
    $\mu_v = \mu'_v\circ\N_{E'/F}$ and $\theta_v = \theta'_v\circ\N_{E'/F}$.
    
    By assumption, $v$ splits into two places $v_1, v_2$ of $E_v$, and
    $\pi_{1,v}$ is equal to \dt $I(\theta'_v,\theta'_v\ve'_v) \otimes_1 (\chi_1\otimes \chi_2)$, where
    $\chi_i$ ($i = 1, 2$) denotes $\chi_{v_i}$.  The unramified representation $\pi_{1,v}$ 
    lifts to $\ve'_v\rtimes\theta'_vI(\chi_1, \chi_2)$ of $G_v$.

    The condition
    ${{}^\sigma \chi}/{\chi} = \ve'\circ\N_{E/F}$ implies that
    \[\frac{\chi_1}{\chi_2}(x)\frac{\chi_2}{\chi_1}(y) = \ve'_v(xy), \quad\forall
    (x,y) \in E_v^\times = F_v^\times \oplus F_v^\times;\] hence, ${\chi_1}/{\chi_2} = \ve'_v$.

    By the condition $\mu\circ\N_{EE'/E'} = \chi\circ\N_{EE'/E}$, we have
    \begin{multline*}
      \chi_1(\N_{E'/F}z_1)\chi_2(\N_{E'/F}z_2) = \mu_v(z_1 z_2),\\
      \forall (z_1, z_2) \in (EE')_v^\times = {E'}_v^\times \oplus {E'}_v^\times.
    \end{multline*}
    Consequently,
    $\chi_1\circ\N_{E'/F} = \chi_2 \circ\N_{E'/F} = \mu_v = \mu'_v\circ\N_{E'/F}$,
    and $\{\chi_1, \chi_2\} = \{\mu'_v, \ve'_v\mu'_v\}$.
    It follows that $\pi_{1,v}$ lifts to $\Pi_v^+$.
  \item
    $\ve'_v = 1$, i.e. 
    $E'_v = E_v = F_v\oplus F_v$, and $(EE')_v = F_v\oplus F_v \oplus F_v \oplus F_v$.

    Suppose $v$ splits into $v_1, v_2$ in $E$ and into $v_1', v_2'$ in $E'$.  
    For $i = 1, 2$, let
    $\mu_i, \theta_i, \chi_i$ denote $\mu_{v_i'}, \theta_{v_i'}, \chi_{v_i}$, respectively.
    Then, $\Pi_v^+$ is equal to $1 \rtimes I(\mu_1\theta_1, \mu_2\theta_2)$.  
    The representation $\pi_{1,v} = I(\theta_1,\theta_2)\otimes_1 (\chi_1\otimes \chi_2)$ lifts to
    ${\chi_1}/{\chi_2}\rtimes \chi_2 I(\theta_1, \theta_2)$.
    The condition ${{}^\sigma \chi}/{\chi} = \ve'\circ\N_{E/F}$ implies that
    $\chi_1 = \chi_2$.    
    Likewise, $\mu_1 = \mu_2$.

    The condition $\mu\circ\N_{EE'/E'} = \chi\circ\N_{EE'/E}$ implies that
    \[
    \mu_1(xy)\mu_2(zt) = \chi_1(xz)\chi_2(yt),\quad
    \forall (x, y, z, t) \in (EE')_v^\times;
    \]
    hence, $\mu_1 = \mu_2 = \chi_1 = \chi_2$.  We conclude that $\pi_{1,v}$ lifts to
    $\Pi_v^+$.
  \end{enumerate}
\end{enumerate}
The packet $\{\Pi\}$ lifts to
$I(\pi_{E'}(\mu\theta), \ve\pi_{E'}(\mu\theta))$ of $\GL(4, \Af)$.  By the same argument used in
the proof of Claim \ref{claim:nonbcH2noI}, we conclude that if an automorphic representation 
$\pi_2 = \pi\otimes_2 \eta$ of $\mb{H_2}(\Af)$ contributes to the $\ve$-trace identity, then 
$\pi = {}^\sigma \pi$.  Consequently, $\pi = B_{E/F}\tau$ for some automorphic representation
$\tau$ of $\GL(2, \Af)$, and $\eta$ is equal to either $\omega_{\tau}$ or $\ve\omega_{\tau}$.
In the first case, we know from Proposition \ref{prop:globalbc} that $\tau\otimes_1 1$ must
contribute to the trace identity, which is a contradiction by Theorem \ref{thm:contribH1}.
The second case is ruled out by Lemma \ref{lemma:bcH2ve}.  The proposition follows.
\end{proof}
\begin{claim}\label{claim:monomialE'}
  The following holds for matching functions:
  \begin{equation}\label{eq:monomialE'}
    \frac{1}{2}\prod_{v \in S} \lp\la \Pi_v^+, f_v\ra - \la\Pi_v^-, f_v\ra_\ve\rp = 
    \frac{1}{2}\la \pi_{E'}(\theta)\otimes_1 \chi, f_1\ra_S.
  \end{equation}
\end{claim}
\begin{proof}
  This follows from Proposition \ref{prop:H1E'monomial} and the equation \eqref{eq:heisenbergpm}.
  Since  ${{}^\sigma \chi}/{\chi} = \ve'\circ\N_{E/F}$, the packet
  $\pi_{E'}(\theta)\otimes_1 {}^\sigma \chi$ of $\mb{H_1}(\Af)$ is equal to 
  $\ve'\pi_{E'}(\theta)\otimes_1 \chi = \pi_{E'}(\theta)\otimes_1 \chi$.  Consequently,
  $\pi_{E'}(\theta)\otimes_1 \chi$ is the only contribution from $\mb{H_1}$ to the trace identity.
  Since by assumption $\theta \neq {{}^{\sigma'}\!\theta}$,
  the cuspidal $E'$-monomial representation $\pi_{E'}(\mu\theta)$ is not $E$-monomial, whence
  the coefficient ${1}/{2}$.
\end{proof}
\subsubsection{Induced Automorphic Representations of $\GSp(2, \Af)$}\label{sec:globalinduced}
\indexi{representation!induced}%
Let $\tau$ be a cuspidal non-$E$-monomial, or one dimensional, automorphic representation of $\GL(2, \Af)$.
Let $\Pi$ be the (normalizedly) induced representation $\ve \rtimes \tau$ of $\mb{G}(\Af)$.
For any place $v$ of $F$, let $\Pi_v^+, \Pi_v^-$ denote the two constituents of $\Pi_v$ if it is
reducible.
If $\Pi_v$ is irreducible, put $\Pi_v^+ := \Pi_v$ and $\Pi_v^- := 0$.
From Section \ref{sec:finechiexpbeta},
the intertwining operator $M_{P_\beta}(s_{2\alpha+\beta}, 0)_v$ acts trivially on $\Pi_v^+$ and as the 
scalar multiplication by $-1$ on $\Pi_v^-$. 
\begin{claim}\label{claim:bcH2veeq}
  The following twisted trace identity holds for matching functions
  $f$ in \dt $\ELw{\mb{G}(\Af), \omega}$ and $f_2$ in $\ELw{\mb{H_2}(\Af), \omega}$:
  \begin{equation}
    \frac{1}{2}\prod_{v \in S}\lp \la \Pi_v^+, f_v\ra_\ve - \la \Pi_v^-, f_v\ra_{\ve} \rp
    = \frac{1}{2}\la B_{E/F} \tau \otimes_1 \omega_\tau \ve, f_2 \ra_S.
  \end{equation}
\end{claim}
\begin{proof}
This follows from Lemma \ref{lemma:bcH2ve} and a comparison of F-H classes.
\end{proof}

The following global character identities follow from local identities which we shall prove in
Section \ref{sec:induced}.  We record them here for the completeness of the global picture.
These global identities hold for all matching test functions, being products of local identities.

For matching functions,
\begin{enumerate}

\item
  $
    \la \ve'\ve\rtimes\pi(\chi), f\ra_{\ve}
    = 
    \la \chi I(1, 1)\otimes_2 \chi|_{\Af^\times}\cdot\ve', f_2\ra,
  $
    
  where $\chi$ is a character of $\idc{E}$ such that ${}^\sigma \chi/\chi = \ve'\circ\N_{E/F}$ for some
  nontrivial quadratic character $\ve'$ of $\idc{F}$;
\item
  $
  \la \ve\times\ve\rtimes\mu, f\ra_\ve =
  \la \mu I(1, 1)\otimes_1 1, f_1 \ra = 
  \la \lp\mu\circ\N_{E/F}\rp I(1, 1)\otimes_2 \mu^2, f_2 \ra,
  $

  where $\ve, \mu$ are charaters of $\idc{F}$;
\item
  $
  \la \ve\rtimes\pi(\chi), f\ra_\ve =
  \la I(1, 1)\otimes_1 \chi, f_1 \ra,
  $

  where $\chi$ is a character of $\idc{E}$.
\end{enumerate}
 
\section{Final Words}\label{sec:finalwords}
Let $\sigma$ be the generator of $\Gal(E/F)$. 
The Galois group $\Gal(E/F)$ acts on the representations $\tau$ of $\GL(2, \Ae)$ via
${}^\sigma\tau(g) = \tau(\sigma g)$ for all $g \in \GL(2, \Ae)$.

For an automorphic representation $\tau$ of $\GL(2, \Ae)$, 
let $\pi(\tau)$ denote the automorphic representation of $\GL(4, \Af)$ which is obtained from
$\tau$ via the twisted endoscopic lifting from $\rR_{E/F}\GL(2)$ to $\GL(4)$.

For an automorphic representation $\pi$ of $\GL(2, \Af)$, 
let $B_{E/F}\pi$ denote the automorphic representation of $\GL(2, \Ae)$
obtained from $\tau$ via base change (\cite{L1}, \cite{F4}).

For an automorphic representation/packet $\pi_i$ of $\mb{H}_i(\Af)$ ($i = 1, 2$),
let $\xi_i^*(\pi_i)$ denote the (quasi-)packet of $\GSp(2, \Af)$ which is the lift of $\pi_i$.

We say that a global (quasi-)packet is {\bf discrete spectrum} if it contains a discrete spectrum automorphic representation.
\indexi{packet!discrete spectrum}%
A discrete spectrum \mbox{(quasi-)packet} of the group \dt $\GSp(2, \Af)$ is said to be {\bf $\ve$-invariant} if 
\indexi{packet!epsilon-invariant@$\ve$-invariant}%
it contains a \emph{discrete spectrum}, $\ve$-invariant, automorphic representation.

We now give the list of (quasi-)packets which comprise the image of the
$\ve$-endoscopic lifting in the discrete spectrum of $\GSp(2)$.  These packets
are necessarily $\ve$-invariant.
\begin{enumerate}
\item {\bf Unstable, cuspidal:}
\indexi{packet!stable}\indexi{packet!unstable}\indexi{packet!quasi-}%
  \begin{enumerate}
  \item
    $
    [\pi, \ve\pi]
    = \xi_1^*(\pi\otimes_1 1) = \xi_2^*(B_{E/F}\pi\otimes_2 \omega_{\pi}), 
    $

    where $\pi$ is a cuspidal, non-$E$-monomial, 
    automorphic representation of \dt $\GL(2, \Af)$.  The packet lifts to the induced representation 
    $I_{(2, 2)}(\pi, \ve\pi)$ of $\GL(4, \Af)$.
  \item
    $
    [\pi(\theta\chi), \pi(\theta\;{}^\sigma\!\chi)]
    = \xi_1^*(\pi(\theta)\otimes_1 \chi), 
    $

    where $\theta, \chi$ are characters of $\idc{E}$ such that none of $\theta, \chi, \theta\chi, \theta\;{}^\sigma\!\chi$
    is invariant under the action of $\sigma$.
    The packet lifts to 
    the induced representation $I_{(2, 2)}(\pi(\theta\chi), \pi(\theta\;{}^\sigma\!\chi))$ of
    $\GL(4, \Af)$.
  \end{enumerate}
\item {\bf Unstable, residual:}

  $  [\mu 1_{2}, \ve\mu 1_{2}]
    = \xi_1^*(\mu 1_{2}\otimes_1 1) 
    = \xi_2^*\lp\lp\mu\circ\N_{E/F}\rp 1_{\GL(2, \Ae)}\otimes_2 \mu^2\rp,$

  where $\mu$ is a character of $\idc{F}$, and $1_2$ denotes the trivial representation of $\GL(2, \Af)$.
  The packet lifts to $I_{(2, 2)}(\mu 1_{2},\;\ve\mu 1_{2})$ of $\GL(4, \Af)$.
\item {\bf Stable, cuspidal:}
  \begin{enumerate}
  \item
  $
  \xi_1^*(\tau\otimes_1 \chi),$

  where $\tau$ is a cuspidal, non-$E$-monomial, automorphic representation of
  \dt $\GL(2, \Af)$, $\chi$ a character of $\idc{E}$, such that there is no quadratic character 
  (trivial or otherwise) $\ve'$ of $\idc{F}$
  for which ${}^\sigma\chi/\chi = \ve'\circ\N_{E/F}$ and $\tau$ is equivalent to $\ve'\tau$ or $\ve'\ve\tau$.

  The packet lifts to the $\ve$-invariant, cuspidal, automorphic representation \dt $\pi(\chi B_{E/F}\tau)$ of $\GL(4, \Af)$.
  \item
    $\xi_2^*(\tau_E\otimes_2 \mu),$
    
    where $\tau_E$ is a cuspidal automorphic representation of $\GL(2, \Ae)$,
    $\mu$ a character of $\idc{F}$, such that $\tau_E \neq {}^\sigma \tau_E$ and $\omega_{\tau_E} = \mu\circ\N_{E/F}$.
    
    The packet lifts to the $\ve$-invariant, cuspidal, automorphic representation $\pi(\tau_E)$ of $\GL(4, \Af)$. 
  \end{enumerate}
\item {\bf Stable, residual:}
  \begin{enumerate}
    \item
$
\{L(\ve\nu, \nu^{-1/2}\pi(\chi))\}
= \xi_1^*(1_{2}\otimes_1 \chi), 
$

where $\chi$ is a character of $\idc{E}$ not fixed by $\Gal(E/F)$.
The quasi-packet lifts to $J(\nu^{1/2}\pi(\chi), \nu^{-1/2}\pi(\chi))$, 
the Langlands quotient of the induced representation $I_{(2,2)}(\nu^{1/2}\pi(\chi), \nu^{-1/2}\pi(\chi))$ of $\GL(4, \Af)$.
\item
$
\{L(\ve\zeta\nu, \nu^{-1/2}\pi(\chi))\}
=\xi_2^*(\chi 1_{\GL(2, \Ae)}\otimes_2 \chi|_{\Af^\times}\cdot\zeta),
$

where $\zeta \neq \ve$ is a nontrivial quadratic character of $\idc{F}$, $\chi$
a character of $\idc{E}$, such that ${}^\sigma \chi /\chi = \zeta \circ \N_{E/F}$.
The quasi-packet lifts to the Langlands quotient $J(\nu^{1/2}\pi(\chi), \nu^{-1/2}\pi(\chi))$ of $\GL(4, \Af)$.  
Note that the packet
\dt $\{L(\zeta\nu, \nu^{-1/2}\pi(\chi))\}$ is disjoint from $\{L(\ve\zeta\nu, \nu^{-1/2}\pi(\chi))\}$, but they
both lift to $J(\nu^{1/2}\pi(\chi), \nu^{-1/2}\pi(\chi))$.
\end{enumerate}
\end{enumerate}
\subsubsection{Surjectivity of the Twisted Endoscopic Lifting}
Recall the definition of the operator $\rho_\omega(\ve)$, or $\rho(\ve)$, on $L(\mb{G}(\Af), \omega)$:
\[
\rho_\omega(\ve) : \phi \mapsto \ve\phi, \quad \forall \phi \in L(\mb{G}(\Af), \omega).
\]
\begin{lemma}\label{lemma:venondegeneracy}
If an irreducible, discrete spectrum,
 automorphic representation $\Pi$ of \dt $\mb{G}(\Af)$, with central character $\omega$, is $\ve$-invariant, 
then the distribution $f \mapsto \tr \Pi(f)\rho_\omega(\ve)$ on $C(\mb{G}(\Af), \omega)$ is nonzero.
\end{lemma}
\begin{proof}
Let $\widetilde{G} = \mb{G}(\Af)\times \{\pm 1\} \rtimes \la \ve \ra$, where $\ve^2 = 1$, and
\indexs{G@$\widetilde{G}$}%
$(1, 1, \ve)(g, 1, 1) = (g, \ve(g), \ve)$.

By assumption, $\Pi$ is a subrepresentation of the right-regular representation $\rho_\omega$ on \dt $L(\mb{G}(\Af), \omega)$.
Hence, $\rho_\omega(\ve)\Pi$ is a subrepresentation and consequently belongs to the discrete spectrum.
By the assumption that the multiplicity one theorem holds for $\mb{G}$,
\indexi{multiplicity one property}%
the operator $\rho_\omega(\ve)$ maps the space of $\Pi$ to itself, with
\begin{equation}\label{eq:rhoveprop}
\rho_\omega(\ve)^{-1}\Pi(g)\rho_\omega(\ve) = \ve(g) \Pi(g),\quad \forall g \in \mb{G}(\Af).
\end{equation}

We extend $\Pi$ to a representation $\Pi_1$ of $\widetilde{G}$, on the space of $\Pi$, as follows:
\[
\Pi_1(g, 1, 1) := \Pi(g),\quad
\Pi_1(1, \pm 1, 1) := \pm 1,\quad
\Pi_1(1, 1, \ve) := \rho_\omega(\ve).
\]
For all $g \in \mb{G}(\Af)$, we have:
\[
\Pi_1(\,\underbrace{(1,1,\ve)(g, 1,1)}_{=\, (g,\, \ve(g),\, \ve)}\,) 
= \rho_\omega(\ve)\pi(g), \quad \Pi_1((g, \ve(g), \ve)) = \ve(g)\pi(g)\rho_\omega(\ve).
\]
It follows from \eqref{eq:rhoveprop} that $\Pi_1$ is indeed a representation of $\widetilde{G}$.

We extend $\Pi$ also to a representation $\Pi_2$ of $\widetilde{G}$ such that $\Pi_2$ is the same as $\Pi_1$ on the subgroup $\mb{G}(\Af)\times\{\pm 1\}$ of $\widetilde{G}$, 
but $\Pi_2(1, 1, \ve) := -\rho_\omega(\ve)$.  Since $\Pi$ is irreducible, it follows from Schur's lemma that $\Pi_1$ and $\Pi_2$ are inequivalent.

The center of $\widetilde{G}$ is $\widetilde{Z} = \{(z, 1, \ve) : z \in \mb{Z}(\Af)\}$, where $\mb{Z}$ is the center of $\mb{G}$.
For $i = 1, 2$, the character of $\Pi_i$ is 
$\langle \Pi_i, \tilde{f}\,\rangle = \tr\!\! \int_{\widetilde{Z}\bs\widetilde{G}} \Pi_i(\tilde{g})\tilde{f}(\tilde{g})\,d\tilde{g}$, which is equal to
\[
\tr \int_{\{\tilde{g} = (g, \pm 1, 1)\}}\!\!\pm\Pi(g)\tilde{f}(\tilde{g})\,d\tilde{g} 
+ \tr \int_{\{\tilde{g} = (g, \pm 1, \ve)\}}\!\!\pm\Pi(g)A_i\tilde{f}(\tilde{g})\,d\tilde{g},
\]
where $A_1 = \rho_\omega(\ve)$ and $A_2 = -\rho_\omega(\ve)$.
Hence, the difference $\langle \Pi_1, \tilde{f}\,\rangle - \langle \Pi_2, \tilde{f}\,\rangle$ is equal to
\[
\begin{split}
&2\; \tr \int_{\{\tilde{g} = (g, \pm 1, \ve)\}}\pm\Pi(g)\rho_\omega(\ve) \tilde{f}(\tilde{g})\,d\tilde{g}\\
=&\, 2 \; \tr\! \!\int_{\{(g, 1, \ve)\}}\!\!\!\!\Pi(g)\rho_\omega(\ve) \tilde{f}((g, 1, \ve))\,d\tilde{g}
  - 2 \; \tr\!\! \int_{\{(g, -1, \ve)\}}\!\!\!\!\Pi(g)\rho_\omega(\ve) \tilde{f}((g, -1, \ve))\,d\tilde{g}.
\end{split}
\]

Suppose the twisted character $\la \Pi, f\ra_{\!\rho_\omega(\ve)}$ is identically zero.  Then, both terms in the sum above are zero,
which implies that $\langle \Pi_1, \tilde{f}\,\rangle - \langle \Pi_2, \tilde{f}\,\rangle = 0$ for all test functions $\tilde{f}$ on $\widetilde{G}$.
Since $\Pi_1 \ncong \Pi_2$, we have a contradiction by the linear independence of characters.
\end{proof}
\begin{lemma}\label{lemma:veindependence}
Let $\{\Pi_1, \Pi_2, \dots, \Pi_n\}$  be a finite set of inequivalent automorphic representations of $\mb{G}(\Af)$
which satisfy the hypothesis of Lemma \ref{lemma:venondegeneracy}.
The twisted characters $\{\la \Pi_i, f\ra_{\rho_\omega(\ve)}\}_{1\leq i\leq n}$ are linearly independent.
\end{lemma}
\begin{proof}
For each $1 \leq i \leq n$, let $\Pi_{i, 1}$, $\Pi_{i, 2}$ be the representations of $\widetilde{G}$ defined as in the proof of
Lemma \ref{lemma:venondegeneracy}.  Suppose there exist complex numbers $c_1, c_2,\dots, c_n$, not all zero, such that
\[
\sum_{i = 1}^n c_i \la \Pi_i, f\ra_{\rho_\omega(\ve)} = 0.
\]
Then, by the same argument used in the proof of Lemma \ref{lemma:venondegeneracy}, we have:
\[
\sum_{i = 1}^n c_i\lp \langle \Pi_{i, 1}\,, \tilde{f}\rangle - \langle \Pi_{i, 2}\,, \tilde{f}\rangle\rp = 0
\]
for all test functions $\tilde{f}$ on $\widetilde{G}$.  Hence, we have a contradiction
by the linear independence of the characters $\{\langle \Pi_{i,j}\,, \tilde{f}\rangle\}_{\substack{1\leq i \leq n\\ j = 1, 2}}\,$.
\end{proof}
\begin{corollary}
Each $\ve$-invariant, discrete spectrum, automorphic representation of \dt $\mb{G}(\Af)$ with two elliptic components
\indexi{representation!$\ve$-invariant}\indexi{representation!discrete spectrum}%
\indexi{representation!automorphic}\indexi{representation!elliptic}%
belongs to a packet listed in this section.
\end{corollary}
\begin{proof}
Suppose a discrete spectrum automorphic representation $\Pi$ of $\mb{G}(\Af)$, with central character $\omega$ 
and two elliptic components, is $\ve$-invariant.  
If we set up an $\ve$-trace identity using the global datum of $\Pi$, 
by Lemma \ref{lemma:veindependence} the $\mb{G}$ side of the identity is nonzero.
Consequently, there must be nonzero contribution from the endoscopic groups side of the trace identity.  
The image of the lifting of the discrete spectrum automorphic representations
of the endoscopic groups consists precisely of those representations which appear in the list, whence the claim.
\end{proof}

\chapter{The Local Picture}\label{chap:local}
Let $k$ be a $p$-adic field with ring of integers $\mc{O}$.

For an algebraic $k$-group $\mb{H}$, let $H = \mb{H}(k)$.
Let $H^\reg$ denote the set of regular elements in $H$.
\indexs{g@${\rm group}^\reg$}%
Let ${Z}_0({H})$ be the maximal $k$-split component of the center of $H$.  
\indexs{Z@$\mb{Z}_0$}%
Let $\bar{H}$ denote the quotient $H / {Z}_0({H})$.
Let $C(H)$ denote the space of smooth, compactly supported modulo ${Z}_0({H})$ functions on $H$.
\indexs{C@$C(\;)$}%
For a character $\omega$ of $Z_0({H})$,
let $C(H, \omega)$ denote the space of functions $f$ in $C(H)$ which 
\indexs{C@$C(\;,\omega)$}%
satisfy
\[f(zh) = \omega(z)^{-1}f(h), \quad\forall z \in Z_0({H}),\;h \in H.\]

Let $\displaystyle J = \lp\lsm &&&1\\&&1&\\&-1&&\\-1&&&\rsm\rp \in \GL(4)$.
Let $\mb{G}$ be the reductive $k$-group
\indexs{G@$\GSp(2)$}%
\[
\displaystyle \GSp(2) = \left\{g \in \GL(4): {}^t g J g = \lambda(g) J
\text{ for some } \lambda(g) \in \mbb{G}_m\right\}.
\]
We call $\lambda(g)$ the similitude factor of $g$.
\indexi{similitude factor}\indexs{lambda@$\lambda(\;)$}%
The character $\lambda \in \Hom(\mb{G}, \mbb{G}_m)$ descends to a character in
$\Hom(G, k^\times)$.
Let $\mb{Z} = \{\diag(z, z, z, z):z\in \mbb{G}_m\}$ be the center of $\mb{G}$.
Since $\mb{Z}$ is $k$-split, the group $Z_0(G)$ is equal to $Z$.

Let $\ve$ be a quadratic character of $k^\times$.
We let $\ve$ denote also the character \[g \mapsto \ve(\lambda(g))\] on $G$.  
\indexs{epsilon@$\ve$}%
The similitude factor of any $\diag(z, z, z, z) \in Z$ is $z^2$; hence, $\ve$ is 
trivial on the center of $G$.


{\bf Notation:}
\begin{itemize}
\item
For a representation $\pi$ of $G$, let $\ve\pi$ denote the representation of $G$ on the
space of $\pi$ defined by
\[
\ve\pi  = \ve\otimes\pi: g \mapsto \ve(g)\pi(g),\quad\forall g \in G.
\]
\item
For $g \in \GL(2, k)$, put $\ve(g) := \ve(\det g)$.
For a representation $\tau$ of $\GL(2,k)$,
let $\ve\tau$ denote the $\GL(2, k)$-module on the space of $\tau$ where $g$ in $\GL(2, k)$ acts by
$\ve(g)\pi(g)$.
\end{itemize}

Let $G'$ be either $\GSp(2, k)$ or $\GL(2, k)$.  
For an admissible (\cite{BZ}) representation $(\pi, V)$ of $G'$,
where $V$ is the complex vector space on which $G'$ acts via $\pi$,
put
\[
\Hom_{G'}(\pi, \ve\pi) := \{A \in \Hom_{\mbb{C}}(V, V) : \ve(g)\pi(g)A =  A\pi(g),\;\forall g \in G'\}.
\indexs{H@$\Hom_{\rm group}({\rm repn., repn.})$}%
\]
We call $A \in \Hom_G(\pi, \ve\pi)$ an {\bf intertwining operator}.
\indexi{intertwining operator}%
\begin{define}
{\rm
We say that an admissible representation $\pi$ of $G'$ is {\bf $\ve$-invariant}
\indexi{representation!$\ve$-invariant}%
if $\pi$ is equivalent to $\ve\pi$, i.e. there exists a nontrivial intertwining operator $A$ in
$\Hom_{G'}(\pi, \ve\pi)$.
}
\end{define}



Let $(\pi, V)$ be an irreducible, admissible, $\ve$-invariant representation of $G$.
Let $A$ be a nontrivial operator in $\Hom_G(\pi, \ve\pi)$.
Since $\ve$ is quadratic,  $A^2$ intertwines the irreducible $\pi$ with itself.
Hence, by Schur's lemma $A^2$ is a scalar multiplication on $V$.  
Multiplying $A$ by the scalar $\lp A^2\rp^{-1/2}$ if necessary,
we assume $A^2 = 1$.   Then, $A$ is unique up to a sign.

{Fix once and for all a Haar measure (\cite{BZ}) $dg$ on $\bar{G} = G/Z$.}
\indexi{Haar measure}%
Let $\omega_\pi$ denote the central character of $\pi$.
For $f \in C(G, \omega_\pi)$, define the convolution operator $\pi(f)A$ on $V$ as follows:
\indexi{convolution operator}%
\[
\pi(f)A v := \int_{\bar{G}}  f(g)\pi(g)Av\; dg,\quad\forall v \in V.
\indexi{convolution operator}%
\]
Since $f$ is locally constant, compactly supported modulo $Z$, and $\pi$ is admissible, 
$\pi(f) A$ has finite rank and its trace $\tr \pi(f)A$ is finite.
We often drop the $\omega_\pi$ from $C(G, \omega_\pi)$ and write simply $C(G)$.
It should be understood implicitly that the test function $f$ transforms under $Z$ via $\omega_\pi^{-1}$.


Recall from chapter  \ref{chap:endoscopy} the definitions of the $\ve$-endoscopic groups
$\mb{H_1}, \mb{H_2}$ of $\mb{G}$ over $k$.
\indexs{H@$\mb{H_1}$}\indexs{H@$\mb{H_2}$}\indexi{endoscopic!group}%
From local class field theory, there is a quadratic extension $\K$ of $k$ which corresponds to $\ve$.  
In particular, $\ve$ is the unique nontrivial quadratic character of $k^\times$ whose restriction
to $\N_{\K/k}\K^\times$ is trivial.
The group of $k$-points of $\mb{H}_i$ ($i = 1, 2$) is as follows:
\begin{enumerate}
\item
$H_1 = \lp \GL(2, k) \times \K^\times\rp' := \{(g, x) \in \GL(2, k)\times \K^\times : \det g = \N_{\K/k} x\}$;
\item
$H_2 = \lp \GL(2, \K) \times k^\times \rp / \K^\times$,\\
where $\K^\times$ embeds into $\GL(2, \K)\times k^\times$ via $x \mapsto (\diag(x, x), \N_{\K/k}x^{-1})$.
\end{enumerate}
We have:
\[
Z_0(H_1) = \{(\diag(z, z), z): z\in k^\times\},
\indexs{Z@$\mb{Z}_0$!$Z_0(H_1)$, $Z_0(H_2)$}%
\]
\[
Z_0(H_2) = \{\lp\diag(1, 1), z\rp_*: z \in k^\times\},
\] 
where the lower * signifies the image in $H_2$ of an element in $\GL(2, \K)\times k^\times$.

Let $i = 1, 2$.
{Fix once and for all a Haar measure $dh_i$ on $\bar{H_i} =  H_i/ Z_0(H_i)$.}
Let $\pi_i$ be an irreducible admissible representation of $H_i$.
For $z \in Z_0(H_i)$, $\pi_i(z)$ acts as a scalar $\omega_{\pi_i}(z)$ on the space of $\pi_i$.
We call $\omega_{\pi_i}$ the {\bf split central character} of $\pi_i$.  
\indexi{central character!split}%
For $f_i \in C(H_i, \omega_{\pi_i})$, let $\pi(f_i)$ denote the convolution operator
\[
\int_{\bar{H_i}} f_i(h)\pi_i(h)\;dh_i.
\]
It has finite rank; hence, its trace $\tr\! \pi_i(f_i)$ is finite.  
As in the case of $G$, we often write $C(H_i)$ instead of $C(H_i, \omega_{\pi_i})$.  It should be understood that
the test function $f_i$ transforms under $Z_0(H_i)$ via $\omega_{\pi_i}^{-1}$.

{\bf Notation/Terminology:}
\begin{itemize}
\item
Put $\GL(2, k)^\K := \{ g \in \GL(2, k) : \det g \in \N_{\K/k}\K^\times\}$.
\item
For any quadratic extension $\mc{L}$ of $k$, we say that an admissible representation $\pi$
of $\GL(2, k)$ is {\bf $\mc{L}$-monomial} if it is the monomial representation $\pi(\theta)$ associated
\indexi{representation!monomial}%
with a character $\theta$ of $\mc{L}^\times$ (see \cite[Thm. 4.6]{JL}, \cite{K}).
\item
We say that a representation is {\bf fully induced} if it is an irreducible parabolically induced representation.
\indexi{representation!induced}%
\item
For a representation $\tau$ of $\GL(2, k)$ and a character $\chi$ of $\K^\times$, let
$\tau\otimes_1\chi$ denote the following representation of $H_1$ on the space of $\tau$:
\indexs{t@$\otimes_1$}%
\[
\tau\otimes_1 \chi : (g, x) \mapsto \chi(x)\tau(g),\quad\forall (g, x) \in H_1.
\]
Note that $\tau\otimes_1 \chi$ is reducible of length two if $\tau$ is $\K$-monomial (\cite{LL}).
\item
For a representation $\pi$ of $\GL(2, \K)$, a character $\mu$ of $k^\times$, such that
$\omega_\pi = \mu \circ \N_{\K/k}$, let $\pi\otimes_2 \mu$ denote the following representation
\indexs{t@$\otimes_2$}%
of $H_2$ on the space of $\pi$:
\[
\pi\otimes_2 \mu : (g, c)_* \mapsto \mu(c)\pi(g),\quad\forall (g, c)_* \in H_2.
\]
The lower $*$ signifies the image of $(g, c) \in \GL(2, \K)\times k^\times$ in $H_2$.
Note that $\pi\otimes_2 \mu$ is well defined because of the condition on $\omega_\pi$ and $\mu$.
\item
For an admissible $\ve$-invariant $G$-module $\pi$, a nontrivial intertwining
operator $A \in \Hom_G(\pi, \ve\pi)$, and a function $f$ in $C(G)$, put
\[
\la \pi, f \ra_A := \tr \pi(f) A.
\indexs{b@$\la\;,\;\ra_A$}%
\]
It depends implicitly on the fixed Haar measure $dg$ on $\bar{G}$.
\item
Let $i = 1$ or $2$.  For an admissible $H_i$-module $\pi_i$ and a function $f_i \in C(H_i)$, put
\[
\la \pi_i, f_i \ra := \tr \pi_i(f_i).
\]
It depends on the fixed Haar measure $dh_i$ on $\bar{H_i}$.
\item
Let $\sigma$ be the generator of $\Gal(\K/k)$.  
For $\gamma \in \K$, put $\bar{\gamma} := \sigma\gamma$.  For a character
$\chi$ of $\K^\times$, put ${}^\sigma \chi(\gamma) := \chi(\bar{\gamma})$ for all $\gamma \in \K^\times$.

For a representation $\pi$ of $\GL(2, \K)$, let ${}^\sigma \pi$ denote the representation
\[{}^\sigma\pi((g_{ij})) = \pi((\sigma g_{ij})), \quad\forall\; (g_{ij}) \in \GL(2, \K).\]
\indexs{sigma@${}^\sigma{\rm repn.}$}%

\end{itemize}
\section{Induced Representations}\label{sec:induced}
\subsection{The Intertwining Operator}\label{sec:intop}
\indexi{intertwining operator|(}%
Let $P = MN$ be the upper triangular Heisenberg parabolic subgroup of $G$ with Levi component
\indexi{parabolic subgroup!Heisenberg}%
\[
M = \left\{ \lp\lsm a &\\& g_2 &\\&&\frac{\det{g_2}}{a}\rsm\rp : a \in k^\times, g_2 \in GL(2,k)\right\}
\]
and unipotent component
\[
N = \left\{\lp\lsm 1 &* &* &*\\&1 &0&*\\&&1&*\\&&&1\rsm\rp \right\}.
\]
We identify $M$ with $k^\times \times \GL(2, k)$ via
$\lp\lsm a &\\& g_2 &\\&&(\det{g_2})/a\rsm\rp \mapsto (a, g_2)$.

Let $\mu$ be a character of $k^\times$. Let $\chi$ be a character of $\K^\times$.  Let
$(\tau_2, V_2)$ be the monomial representation $\pi(\chi)$ of $\GL(2, k)$ associated with $\chi$.
Let $(\mu\otimes\tau_2, V_2)$ be the representation of $P$ defined by
\[
\mu\otimes\tau_2 : (a, g_2)n \mapsto \mu(a)\tau_2(g_2),\quad\forall\; (a, g_2) \in M,\; n\in N.
\indexs{t@${\rm char.}\otimes{\rm repn.}$}%
\]
Let $(\pi, V) =\mu\rtimes \tau_2 := I_P^G(\mu\otimes\tau_2)$ be the representation of $G$ normalizedly induced (\cite{BZ})
\indexs{t@${\rm char.}\rtimes{\rm repn.}$}%
\indexi{representation!normalizedly induced}%
from $\mu\otimes\tau_2$.

Let $\tau = \mu\otimes\tau_2$.
For $m \in M$, let $\delta_M(m) = \abs{\det({\rm Ad}\;m|_{\mf{n}})}$,
\indexs{delta@$\delta_M(\;)$}%
where $\mf{n}$ denotes the Lie algebra of $N$.  
By the definition of normalizedly induced representations, the vector space $V$
consists of the smooth (\cite{BZ})
functions $\fee :  G \rightarrow V_2$ which satisfy 
\[\fee\lp m n g\rp
= (\delta_M^{1/2}\tau)(m) \fee(g),\quad \forall m \in M,\; n \in N,\; g \in G.\]
The group $G$ acts on $V$ via right translation $\lp\pi(g)\fee\rp(h) = \fee(h g)$.

The similitude factor of 
$(a, g_2) \in M$ is equal to $\det g_2$; hence,
$\ve\otimes(\mu\rtimes \tau_2) = \mu\rtimes\ve\tau_2$.   Since $\tau_2 = \pi(\chi)$ is $\ve$-invariant
(in the context of $\GL(2, k)$),
the induced representation $\pi = \mu\rtimes\tau_2$ is $\ve$-invariant (in the context of $\GSp(2, k)$).
In what follows,
we define a nontrivial operator $A$ in $\Hom_G(\pi, \ve\pi)$, and then we
compute explicitly the twisted character $\la \pi, f\ra_A$ for $f \in C(G)$.

Let $\GL(2, k)^\K = \{g \in \GL(2, k) : \det g \in \N_{\K/k}\K^\times\}$.  It is a subgroup of index $2$ in $\GL(2, k)$.
\indexs{G@$\GL(2, k)^\K$}%
From \cite[Sect. 2]{LL}, there exist subspaces $V^+_2, V^-_2$ in the space $V_2$ of $\tau_2$ such that
\indexs{V@$V_2^+$, $V_2^-$}%
$V_2 = V^+_2 \oplus V^-_2$ as a $\GL(2, k)^\K$-module.  Moreover, $V^-_2$ is equal to $\tau_2(w)V^+_2$ for any element
$w \in \GL(2, k) - \GL(2, k)^\K$.
Let $A_2$ be the operator on $V_2$ defined by
\indexs{A@$A_2$}%
\[
A_2(v^+ + v^-) := v^+ - v^-, \quad \forall v^+ \in V^+_2,\; v^- \in V^-_2.
\]
It is an intertwining operator in $\Hom_{\GL(2, k)}(\tau_2, \ve\tau_2)$ whose square is the identity.
Note that by the symmetry between $V_2^+$ and $V_2^-$, a choice is involved in the definition of $A_2$.

Define a map $A$ from $V$ to the set of smooth $V_2$-valued functions on $G$ as follows:
\begin{equation}\label{eq:defineA}
(A\fee)(g) := \ve(g) A_2(\fee(g)),\quad \forall g \in G,\;\fee \in V.
\end{equation}
\begin{claim}
The map $A$ is an automorphism of $V$ such that $A^2 = 1$, and
\[
\ve\pi(g) A = A \pi(g), \quad g \in G.  
\]
In other words, $A$ is a nontrivial \intw
in \dts $\Hom_G(\pi, \ve\pi)$.
\end{claim}
\begin{proof}
First, we show that $A$ is an endomorphism of $V$.  It then follows trivially that 
$A^2 = 1$ and hence $A$ is invertible.

For every $\fee \in V$, we need to show that
$(A\fee) (mng) = (\delta_M^{1/2}\tau)(m)\fee(g)$ for all $m = (a, g_2) \in M$, $n \in N$, $g \in G$.
The restriction of $\ve$ to $N$ is trivial, for the similitude factor of a unipotent element is $1$.
We have:
\[
\begin{split}
(A\fee)(m n g)
  &= \ve(m)\ve(n)\ve(g)A_2(\fee(mng))\\
  &= \ve(\det g_2) \ve(g) A_2\!\lp \delta_M^{1/2}(m)\mu(a)\tau_2(g_2)\fee(g)\rp\\
  &= \delta_M^{1/2}(m)\mu(a)\ve(\det g_2) \ve(\det g_2)\tau_2(g_2)\ve(g)A_2(\fee(g))\\
  &= (\delta_M^{1/2}\tau)(m)\ve(g)A_2(\fee(g))\\
  &= (\delta_M^{1/2}\tau)(m) (A\fee)(g).
\end{split}
\]
Hence, $A\fee$ lies in $V$.

Next, we show that
\begin{equation}\label{eq:A}
\ve(g) \pi(g) A\fee = A \pi(g)\fee,\quad \forall g \in G,\;\fee \in V.
\end{equation} 
For all $h \in G$, on the left-hand side of \eqref{eq:A} we have:
\[
(\ve(g)\pi(g)A\fee)(h) = \ve(g)(A\fee)(hg) = \ve(g)\ve(hg)A_2(\fee(hg))
= \ve(h)A_2(\fee(hg)).
\]
On the right-hand side, we have:
\[
(A \pi(g) \fee)(h) = \ve(h)A_2((\pi(g)\fee)(h)) = \ve(h)A_2(\fee(hg)).
\]
Hence, left-hand side equals right-hand side, which completes the proof.
\end{proof}

\indexi{intertwining operator|)}%
\subsection{A Trace Identity between $G$ and $H_1$}

Suppose the field $\K$ is equal to $k(\sqrt{B})$ for some $B \in k^\times - {k^\times}^2$.  We define an
embedding $\phi : \K^\times \hookrightarrow \GL(2, k)$ as follows:
\indexs{phi@$\phi$}%
\[
\phi : a + b\sqrt{B} \mapsto \lp\lsm a  & b B \\ b & a\rsm\rp, \quad\forall a, b \in k;\;(a, b)\neq (0, 0).
\] 
Note that, for $\gamma \in \K^\times$, the (possibly non-distinct) eigenvalues of $\phi(\gamma)$ are
$\gamma, \bar{\gamma}$, and $\det \phi(\gamma) = \N_{\K/k}\gamma$.
Let
\[\displaystyle
T_\K = \left\{\phi(\gamma) : \gamma \in \K^\times\right\}.
\indexs{T@$T_\K$}%
\]
For an element $t_2 \in \GL(2, k)$ which is conjugate to $\phi(\gamma)$ for some $\gamma \in \K^\times$,
put 
\[
D_{T_\K}(t_2) := \frac{\abs{(\gamma - \bar{\gamma})^2}}{\abs{\N_{\K/k}\gamma}}.
\indexs{D@$D_{T_\K}(\;)$}%
\]

Recall that $\tau_2 = \pi(\chi)$ is the $\K$-monomial representation associated with a character $\chi$ of $\K^\times$.   
We have defined in the previous section an intertwining operator $A_2$ in $\Hom_{\GL(2, k)}(\tau_2, \ve\tau_2)$.
Let $\chi_{\tau_2}^\ve$ be the $\ve$-twisted Harish-Chandra character associated with the distribution 
\indexi{Harish-Chandra character!twisted}%
$f \mapsto \tr\!\tau_2(f)A_2$ on $C(\GL(2, k))$.  In particular, for $g, h \in \GL(2, k)$, 
we have $\chi_{\tau_2}^\ve(h^{-1}gh) = \ve(h)\chi_{\tau_2}^\ve(g)$.  

Let $t_2$ be a regular element in $\GL(2, k)$.
If $\chi_{\tau_2}^\ve(t_2) \neq 0$, then the centralizer of $t_2$ must be a torus which lies in the kernel of $\ve$.
All such tori are conjugate to $T_\K$;
hence, $\chi_{\tau_2}^\ve(t_2)$ is zero unless $t_2$ is conjugate to an element in $T_\K$.
By \cite[Lemma 7.19]{L1}, for each regular $t_2 = \phi(\gamma) \in T_\K$, $\gamma \in \K^\times - k^\times$, we have:
\begin{equation}\label{eq:twistedchi2}
{D^{1/2}_{T_\K}(t_2)}\chi_{\pi(\chi)}^\ve(t_2) = d\cdot \lambda(\K/k, \psi) \cdot \ve(\gamma - \bar{\gamma})
\indexs{lambda@$\lambda(\;,\psi)$}%
{\Big[\chi(\gamma) + \chi(\bar{\gamma})\Big]}.
\end{equation}
Here, $\lambda(\K/k, \psi)$ is a constant independent of $t_2$, and $d$ is $\pm 1$ depending on the choice of sign involved in the
definition of $A_2$.  For simplicity, we assume $A_2$ is chosen such that $d = 1$.

Let 
$\displaystyle
T = \left\{ \lp\lsm a &&\\&t_2&\\&&(\det t_2) / a\rsm\rp:a\in k^\times, t_2 \in T_\K\right\}.
$
It is an elliptic maximal torus in $M$.  
Observe that, for $\gamma \in \K^\times$,
\[
\ve\lp\lp\lsm a &&\\&\phi(\gamma)&\\&&(\det \phi(\gamma))/a\rsm\rp\rp = \ve(\det\phi(\gamma)) = \ve(\N_{\K/k}\gamma) = 1;
\]
hence, $T$ lies in the kernel of $\ve$.

Fix a Haar measure $dt$ on $\bar{T} = T/Z$.  Let $d\bar{g} = dg / dt$.
For $t \in {T}$ and $f \in C(G)$, put
\[
O_G^{\ve}(f, t) := \int_{T \bs G} \ve(\bar{g})f(\bar{g}^{-1}t \bar{g})\;d\bar{g}.
\indexi{orbital integral!twisted}%
\]
Note that $\ve(\bar{g})$ is well defined because $T \subset \ker \ve$.


Recall that we identify $(a, g_2)$ in $k^\times \times \GL(2, k)$ with 
$\blockdiag(a, g_2, (\det g_2) /a)$ in $M$.
\begin{lemma}\label{lemma:trG}
Let $\mu$ be a character of $k^\times$.  Let $\chi$ be a character of $\K^\times$.  
Let $A$ in $\Hom_G(\mu\rtimes\pi(\chi), \mu\rtimes\ve\pi(\chi))$
be the intertwining operator defined in {\rm Section \ref{sec:intop}}.
Then,
\begin{multline}\label{eq:trG}
\la \mu\rtimes\pi(\chi), f\ra_A =
\frac{1}{4} \int_{t = (a, t_2) \in \bar{T}^\reg}
\left[\mu(a) + \mu\lp\frac{\N_{\K/k} \gamma}{a}\rp\right] 
\left[\chi(\gamma) + \chi(\bar{\gamma})\right]\\
\cdot \lambda(\K/k, \psi)\ve(\gamma - \bar{\gamma})D^{1/2}_{T \bs G}(t) O_G^\ve(f, t) \;dt.
\end{multline}
Here, $t_2 = \phi(\gamma) \in T_\K$, $\gamma \in \K^\times$.
\end{lemma}
\begin{proof}
We keep the same notation used in the previous section, i.e. $\tau_2 = \pi(\chi)$, 
$(\pi, V) = \mu \rtimes \pi(\chi)$,
$\tau = \mu\otimes\tau_2,\dots$, etc.
Let $K$ be the hyperspecial maximal compact subgroup $\GSp(2, \mathcal{O})$ of $G$.
We have the Iwasawa decomposition $G = PK = NMK$; hence, $\bar{G} = N \bar{M} K$, where $\bar{M} = M / Z$.
\indexi{Iwasawa decomposition}%
We also have the corresponding measure decomposition $dg = \delta_M^{-1}(m)dn dm dk$, where
$dm$ is a measure on $\bar{M}$.
Recall that $\ve|_N = 1$.
For $h \in G$, $\fee \in V$,
\begin{multline*}
(\pi(f)A\fee)(h)
\\\begin{split}
&= \int_{\bar{G}} f(g)(\pi(g) (A\fee))(h) \;dg =\int_{\bar{G}} f(g) (A\fee)(hg) dg\\
& = \int_G f(h^{-1} g) (A\fee)(g) \;dg = \int_G f(h^{-1} g)\ve(g) A_2(\fee(g)) \;dg\\
&=\int_N \int_{\bar{M}} \int_K f(h^{-1}n_1 m k)\ve(n_1mk)A_2(\fee(n_1 m k)) \delta_M^{-1}(m)\;dk \;dm\; dn_1\\
&=\int_N \int_{\bar{M}} \int_K f(h^{-1}n_1 m k)\ve(m k) \delta_M^{-1/2}(m)A_2(\tau(m)\fee(k))\;dk\; dm\; dn_1.
\end{split}
\end{multline*}
Writing a general $m \in M$ in the form $(a, g_2)$, where $a \in k^\times$ and $g_2 \in \GL(2, k)$,
the above expression is equal to
\begin{multline*}
\int_N \int_{\bar{M}}\int_K f(h^{-1}n_1 m k)
\ve(g_2)\ve(k)\delta_M^{-1/2}(m)\mu(a)A_2\tau_2(g_2)\fee(k) \;dk \;dm \;dn_1
\\
\shoveleft{
=\!\! \int_N \int_{\bar{M}} \int_K\!
f(h^{-1}n_1 m k)
\ve(g_2)\ve(k)\delta_M^{-1/2}(m)\mu(a) \ve(g_2) \tau_2(g_2) A_2 \fee(k) dk dm dn_1}
\\= \int_N \int_{\bar{M}} \int_K f(h^{-1}n_1 m k)\ve(k)\delta_M^{-1/2}(m)\mu(a)\tau_2(g_2) A_2 \fee(k) \;dk \;dm \;dn_1.
\end{multline*}

For each $m \in \bar{M}$, we change variables $n_1 \mapsto n$, where $n$ is defined by $n^{-1}mnm^{-1} = n_1$.
The resulting Jacobian is $\abs{\det (1 - {\rm Ad}\; m)|_{\mathfrak{n}}}$.
Let
$D_{M\bs G}(m) = \abs{\det (1 - \text{Ad}\;m)|_{\mathfrak{g/m}}}$,
\indexs{D@$D_{M\bs G}$}%
where lower case gothic type denotes the Lie algebra of a group.
Since  \[D^{1/2}_{M\bs G}(m) = \delta_M^{-1/2}(m) \abs{\det (1 - \text{Ad}\;m)|_{\mathfrak{n}}},\]
the trace $\lp\pi(f)A\fee\rp(h)$ is equal to
\[
\int_N \int_{\bar{M}} \int_K\\ f(h^{-1}n^{-1} m n k)\ve(k)
D^{1/2}_{M \bs G}(m)\mu(a) \pi(g_2) A_2 \fee(k)\; dk\; dm\; dn.
\]
Since $G = PK$ and functions in $V$ transform under $P$ via $\tau$, we may
identify $V$ with a space of functions on $K$.  The operator $\fee \mapsto \pi(f)A\fee$ on $V$
is therefore an integral operator with kernel
\begin{multline*}
(h, k) \in K \times K \mapsto\\
\int_N \int_{\bar{M}}
f(h^{-1}n^{-1} m n k)\ve(k)
D^{1/2}_{M \bs G}(m)\mu(a) \pi(g_2) A_2\; dm\; dn.
\end{multline*}
This kernel is an operator on the vector space $V_2$ of $\tau_2$, for each $(h, k) \in K\times K$.

We compute $\la \pi, f\ra_A = \tr \pi(f)A$ by integrating the character of this kernel
(as an operator on $V_2$) over the diagonal $\{(k, k) : k \in K\}$.  We have:
\[\begin{split}
\la \pi, f\ra_A & = \int_K\int_N \int_{\bar{M}}
f(k^{-1}n^{-1}mnk)\ve(k)D^{1/2}_{M \bs G}(m)\mu(a)\chi_{\tau_2}^\ve(g_2)\;dm\; dn\;dk\\
&= \int_{\bar{M}} D^{1/2}_{M \bs G}(m) \mu(a)\chi_{\tau_2}^\ve(g_2)
\int_K \int_N f(k^{-1} n^{-1} m n k) \ve(k)\;dn\; dk\; dm,
\end{split}\]

Let $N_G(M)$ be the normalizer of $M$ in $G$.
Let $W_{M\bs G} = N_G(M) / M$.  Then, $\abs{W_{M\bs G}} = 2$, and the nontrivial element of $W_{M\bs G}$ sends $(a, g_2)$ \indexs{W@$W_{M\bs G}$}%
to $({(\det g_2)}/{a}, g_2)$.  We have:
\begin{multline*}
\la \pi, f\ra_A
= \frac{1}{\abs{W_{M\bs G}}}\int_{\bar{M}} \lp\mu(a)+\mu\lp\tfrac{\det g_2}{a}\rp\rp D^{1/2}_{M \bs G}(m)\chi_{\tau_2}^\ve(g_2)\\
 \cdot \Biggl[\int_K\int_N  f(k^{-1}n^{-1}m n k) \ve(k)\;dn\; dk\Biggr] dm.
\end{multline*}

Let $t_2$  be a regular element in $\GL(2, k)$.
For any $g_2 \in \GL(2, k)$, put $t_2^{g_2} := g_2^{-1} t_2 g_2$.
Recall that $\chi_{\tau_2}^\ve(t_2)$ is zero unless $t_2$ is conjugate to an element in $T_\K$,
and $\chi_{\tau_2}^\ve(t_2^{g_2}) = \ve(g_2)\chi_{\tau_2}^\ve(t_2)$ for all $g_2 \in \GL(2, k)$.
By the Weyl integration formula,
\indexi{Weyl integration formula}%
if a function $f$ on a group $H$ is zero outside of the conjugacy class of a maximal torus $T$, then
\begin{equation}\label{eq:Weylformula}
\int_H f(h)\; dh 
= \frac{1}{\abs{W_{T\bs H}}}\int_{T^\reg} D_{T\bs H}(t)\int_{T\bs H} f({\bar{h}}^{-1} t \bar{h})\;d\bar{h}\;dt.
\end{equation}
Here, $D_{T\bs H}(t) := \abs{\det(1 - \text{Ad }t)|_{\mf{h}/\mf{t}}}$.
\indexs{D@$D_{T\bs H}$}%
Viewing 
$\la \pi, f\ra_A$ as an integral over $\bar{M}$, the integrand vanishes outside of the union of the tori 
conjugate to $T$ in $M$.
Applying \eqref{eq:Weylformula}, we obtain: $\la \pi, f\ra_A =$
\begin{multline*}
\frac{1}{\abs{W_{T\bs M}}} \frac{1}{\abs{W_{M\bs G}}}
{\int_{t=(a, t_2) \in \bar{T}^\reg} D_{T\bs M}(t)}
\\\cdot \int_{\bar{m} = (b, g_2) \in T \bs M} D^{1/2}_{M \bs G}(\bar{m}^{-1}t\bar{m})
\left[\mu(a^b) + \mu\lp\tfrac{\det t_2^{g_2}}{a^b}\rp\right]\chi_{\tau_2}^\ve(t_2^{g_2})
\\\cdot \int_{K}\int_N f(k^{-1}n^{-1}\bar{m}^{-1}t\bar{m}n k)\ve(k)\;dn\; dk\; d\bar{m}\; dt
\end{multline*}
\begin{multline*}
=\frac{1}{\abs{W_{T\bs G}}}\int_{t=(a, t_2) \in \bar{T}^\reg}
D_{T\bs M}(t)D^{1/2}_{M \bs G}(t)\\
\int_{\bar{m} = (b, g_2) \in T \bs M}
\left[(\mu(a^b) + \mu\lp\tfrac{\det t_2^{g_2}}{a^b}\rp\right]
\chi_{\tau_2}^\ve(t_2^{g_2})\\
\int_{K}\int_N f(k^{-1}n^{-1}\bar{m}^{-1}t\bar{m}n k )\ve(k)\;dn\; dk\; d\bar{m}\; dt.
\end{multline*}
Here, $W_{T\bs G}$ denotes the Weyl group of $T$ in $G$.

For $m = (b, g_2) \in M$, we have
$\chi_{\tau_2}^\ve(t_2^{g_2}) =\ve(g_2)\chi_{\tau_2}^\ve(t_2) = \ve(m)\chi_{\tau_2}^\ve(t_2).$
Inserting $\ve(n) = 1$ ($n \in N$) at an appropriate place,
we express $\la \pi, f\ra_A$ as follows:
\begin{multline*}
\frac{1}{\abs{W_{T\bs G}}}\int_{t = (a, t_2) \in \bar{T}^\reg}D_{T \bs M}(t)D^{1/2}_{M \bs G}(t) 
\left[\mu(a) + \mu\lp\tfrac{\det{t_2}}{a}\rp\right]\chi_{\tau_2}^\ve(t_2)\\
\cdot \int_{T \bs M}\ve(\bar{m})
\int_K\int_N f(k^{-1}n^{-1}\bar{m}^{-1}t\bar{m} n k)\ve(nk)\;dn\; dk\;d\bar{m}\; dt
\end{multline*}
\begin{multline*}
= \frac{1}{\abs{W_{T\bs G}}}\int_{\bar{T}^\reg}D_{T \bs M}(t)D^{1/2}_{M \bs G}(t) 
\left[\mu(a) + \mu\lp\tfrac{\det{t_2}}{a}\rp\right]\chi_{\tau_2}^\ve(t_2)\\
\cdot \int_{T \bs M}
\int_K\int_N \ve(\bar{m}nk)f(k^{-1}n^{-1}\bar{m}^{-1}t\bar{m} n k)\;dn\; d\bar{m}\; dt
\end{multline*}
\begin{multline*}
= \frac{1}{\abs{W_{T\bs G}}}\int_{\bar{T}^\reg}
\left[\mu(a) + \mu\lp\tfrac{\det{t_2}}{a}\rp\right] D^{1/2}_{T\bs M}(t)\chi_{\tau_2}^\ve(t_2)\\
\cdot D^{1/2}_{T\bs M}(t)D^{1/2}_{M \bs G}(t) \int_{T \bs G} 
\ve(\bar{g})f(\bar{g}^{-1}t \bar{g})\;d\bar{g}\; dt.
\end{multline*}
Since $M \cong k^\times \times \GL(2, k)$,
$D_{T \bs M}((a, t_2))$ is equal to $D_{T_\K}(t_2) := D_{T_\K \bs \GL(2, k)}(t_2)$.
Moreover, $D_{T\bs M}D_{M\bs G} = D_{T\bs G}$.  Hence, 
\begin{multline}\label{eq:trpiAorig}
\la \pi, f\ra_A =
\frac{1}{\abs{W_{T\bs G}}}\int_{t = (a, t_2) \in \bar{T}^\reg}
\left[\mu(a) + \mu\lp\tfrac{\det{t_2}}{a}\rp\right]D^{1/2}_{T_\K}(t_2)\chi_{\tau_2}^\ve(t_2)\\
\cdot D^{1/2}_{T \bs G}(t) O_G^\ve(f, t)\;dt.
\end{multline}
By \eqref{eq:twistedchi2}, the above expression is equal to
\begin{multline*}
\frac{1}{\abs{W_{T\bs G}}}\int_{t = (a,\; \phi(\gamma)) \in \bar{T}^\reg}
\left[\mu(a) + \mu\lp\tfrac{\N_{\K/k} \gamma}{a}\rp\right]\cdot
\Big[\chi(\gamma) + \chi(\bar{\gamma})\Big]\\
\cdot \lambda(\K/k, \psi)\ve(\gamma - \bar{\gamma})D^{1/2}_{T \bs G}(t) O_G^\ve(f, t)\;dt.
\end{multline*}
Since $\abs{W_{T\bs G}} = 4$, the lemma follows.
\end{proof}
Let $P_1$ be the upper triangular parabolic subgroup of $H_1$. 
\indexs{P@$P_1$}%
Its Levi component is the maximal diagonal torus 
\[T_1 = \{(\diag(a, b), x) : a, b \in k^\times, x \in \K^\times, ab = \N_{\K/k}x\}.\]
\indexs{T@$T_1$}%
For characters $\mu_1, \mu_2$ of $k^\times$ and a character
$\chi$ of $\K^\times$, let $\mu_1\otimes\mu_2\otimes_1\chi$ denote the representation
of $P_1$ defined as follows:
\[
\mu_1\otimes\mu_2\otimes_1\chi : \lp \lp\lmx a & *\\0&b\rmx\rp, x\rp \mapsto \mu_1(a)\mu_2(b)\chi(x).
\indexs{t@$\;\otimes\;\otimes_1\;$}%
\]
Let $I(\mu_1, \mu_2)\otimes_1 \chi$
denote the representation of $H_1$ normalizedly induced from $\mu_1\otimes\mu_2\otimes_1\chi$.

Due to the condition that each element $(\diag(a, b), x) \in T_1$ satisfies $ab = \N_{\K/k} x$,
the $H_1$-module $I(\mu_1, \mu_2)\otimes_1 \chi$ is equal to 
$\eta^{-1} I(\mu_1, \mu_2)\otimes_1\lp\eta\circ\N_{\K/k}\rp\cdot\chi$ for
any character $\eta$ of $k^\times$.  In particular,
we have: \[I(\mu_1, \mu_2)\otimes_1 \chi = I(\mu_1\mu_2^{-1}, 1)\otimes_1\lp\mu_2\circ\N_{\K/k}\rp\cdot\chi.\]
Hence, any $H_1$-module induced from $P_1$ has the form
$I(\mu, 1)\otimes_1 \chi$, where $\mu$, $\chi$ are characters of $k^\times$, $\K^\times$, respectively.

Let $\bar{T}_1 = T_1 / Z_0(H_1)$.  
For $a \in k^\times$, $\gamma \in \K^\times$, recall that $(a, \phi(\gamma))$ denotes the element
$\diag(a, \phi(\gamma), (\N_{\K/k}\gamma)/a) \in T \subset G$.
The map 
\[
\lambda_1 : (a, \phi(\gamma)) \mapsto \lp\dmtwo{a}{\frac{\N_{\K/k}\gamma}{a}}, \gamma\rp
\indexs{lambda@$\lambda_1$}%
\]
defines an isomorphism from $\bar{T}$ to $\bar{T}_1$.  Let $dt_1$ be the Haar measure on $\bar{T}_1$ 
which is compatible with the measure $dt$ on $\bar{T}$ via $\lambda_1$.  That is,
for each measurable subspace $\bar{S} \subset \bar{T}$, we have
$
\int_{\lambda_1(\bar{S})}\; dt_1 = \int_{\bar{S}}\;dt.
$

For $t_1 \in T_1$ and $f_1 \in C(H_1)$, let 
\[
O_{H_1}(f_1, t_1) = \int_{T_1 \bs H_1} f_1(\bar{h}_1^{-1} t_1 \bar{h}_1)\;d\bar{h}_1.
\indexi{orbital integral}%
\]

Let $\mu, \chi$ be the characters of $k^\times, \K^\times$, respectively.
The split central character of \dt $I(\mu, 1)\otimes_1 \chi$ is $\chi|_{k^\times}\cdot\mu$.
\begin{lemma}\label{lemma:trH1}
For $f_1 \in C(H_1, \chi|_{k^\times}\cdot\mu)$, the following holds:
\begin{multline}\label{eq:trH1}
\la I(\mu, 1)\otimes_1 \chi, f_1 \ra\\
= \frac{1}{4}
\int_{t_1 = (\diag(a, b),\gamma) \in \bar{T}_1^\reg}
\Big[\mu(a) + \mu(b) \Big] \Big[ \chi(\gamma) + \chi(\bar{\gamma})\Big]\\
\cdot D^{1/2}_{T_1\bs H_1}(t_1)O_{H_1}(f_1, t_1)\;dt_1.
\end{multline}
\end{lemma}
\begin{proof}
We skip the proof, which 
is similar to that of Lemma \ref{lemma:trG} and is simpler because there is no twisting
by $\ve$.
\end{proof}

Recall the embedding $\phi : \K^\times \hookrightarrow \GL(2, k)$ which sends 
$a + b\sqrt{B}$ to $\lp\lsm a&bB\\b&a\rsm\rp$.  
Note that $\det \phi(\gamma) = \N_{\K/k}\gamma$ for all $\gamma \in \K^\times$.
Given $a, b \in k^\times$, $\gamma \in \K^\times - k^\times$, such that $ab = \N_{\K/k}\gamma$
and $a \neq b$,
the element $t = \blockdiag(a, \phi(\gamma), b) \in M$
is regular. 
Every norm of $t$ in $H_1$ is conjugate to either
$t_1 = (\diag(a, b), \gamma)$ or  ${\bar{t}_1} = \lp\diag(a, b),\bar{\gamma}\rp$.
In particular, $t_1$ and $\bar{t}_1$ are $\mb{G}$-regular.
We make the following choice of transfer factors for the pairs $(t_1, t)$, $(\bar{t}_1, t)$:
\indexi{transfer factor}%
\[
\Delta(t_1, t) = \Delta(\bar{t}_1, t) = 
\lambda(\K/k, \psi)\ve(a)\ve(\gamma - \bar{\gamma})
\left[\frac{D_{T \bs G}(t)}{D_{T_1 \bs H_1}(t_1)}\right]^{1/2}.
\]
\begin{define}
{\rm
We say that $f \in C(G)$ and $f \in C(H_1)$ are {\bf weakly matching} if
\indexi{matching functions!weakly}%
\[
O_{H_1}(f_1, t_1) = \Delta(t_1, t)O_G^{\ve}(f, t)
\]
for every $(t_1, t) \in T_1 \times T^\reg$ such that $t_1$ is a norm of $t$.
}
\end{define}
The assumption is that there exists a choice of transfer factors for $H_1 \times G$, compatible with 
$\Delta(t_1, t)$, such that matching implies weakly matching.
\begin{prop}\label{prop:trGH1}
Let $\mu$ be a character of $k^\times$.  Let $\chi$ be a character of $\K^\times$.
For weakly matching functions, the following character identity holds:
\[
\la \mu\rtimes\pi(\chi), f\ra_A = 
\la I(\ve\mu, 1)\otimes_1\chi, f_1\ra.
\]
\end{prop}
\begin{proof}

The elements in $H_1^\reg$ which are not $\mb{G}$-regular form a set of measure zero.
Hence, by Lemma \ref{lemma:trH1} and our choice of $dt_1$, $\la I(\ve\mu, 1)\otimes_1\chi, f_1\ra$ is equal to
\begin{multline*}
\frac{1}{4}
\int_{t = (a, \phi(\gamma)) \in \bar{T}^\reg}
\Big[\ve\mu(a) + \ve\mu\lp\tfrac{\N_{\K/k}\gamma}{a}\rp \Big] \Big[ \chi(\gamma) + \chi(\bar{\gamma})\Big]\\
\cdot D^{1/2}_{T_1\bs H_1}(\lambda_1(t))O_{H_1}(f_1, \lambda_1(t))\; dt
\end{multline*}
\begin{multline}\label{eq:trH1'}
= \frac{1}{4}
\int_{t = (a, \phi(\gamma)) \in \bar{T}^\reg}
\ve(a)\Big[\mu(a) + \mu\lp\tfrac{\N_{\K/k}\gamma}{a}\rp \Big] \Big[ \chi(\gamma) + \chi(\bar{\gamma})\Big]\\
\cdot D^{1/2}_{T_1\bs H_1}(\lambda_1(t))O_{H_1}(f_1, \lambda_1(t))\; dt.
\end{multline}

Since the functions $f_1 \in C(H_1)$ and $f \in C(G)$ are weakly matching, we have:
\[
\begin{split}
O_{H_1}(f_1, \lambda_1(t)) &= 
\lambda(\K/k, \psi)\ve(a)\ve(\gamma - \bar{\gamma})\left[\frac{D_{T \bs G}(t)}{D_{T_1 \bs H_1}(\lambda_1(t))}\right]^{1/2}
O_{G}^\ve(f, t)
\end{split}
\]
for all regular $t  = (a, \phi(\gamma)) \in T$.
The expression \eqref{eq:trH1'} is therefore equal to
\begin{multline}
\frac{1}{4}
\int_{t = (a, \phi(\gamma)) \in \bar{T}^\reg}
\Big[ \mu(a) + \mu\lp\tfrac{\N_{\K/k}\gamma}{a}\rp \Big] \Big[ \chi(\gamma) + \chi(\bar{\gamma})\Big]\\
\cdot\lambda(\K/k, \psi)\ve(\gamma - \bar{\gamma})D^{1/2}_{T\bs G}(t)O_{G}^\ve(f, t)\; dt.
\end{multline}
By Lemma \ref{lemma:trG}, the proposition follows.
\end{proof}
\subsection{A Trace Identity between $G$ and $H_2$}
Let $T_2 = \{(\diag(\alpha, \beta), c)_* : \alpha,\beta \in \K^\times, c\in k^\times\}$
\indexs{T@$T_2$}%
be the maximal diagonal torus of $H_2$.
Here, lower $*$ denotes the image of $(\diag(\alpha, \beta), c)$ in $H_2$.  Let $\bar{T_2} = T_2 / Z_0(H_2)$.
\indexs{T@$\bar{T_2}$}%

Recall that
$
\displaystyle
T = \left\{ (a, \phi(\gamma)) :a\in k^\times, \gamma \in \K^\times\right\} \subset G,
$
where $(a, \phi(\gamma))$ denotes the element $\blockdiag(a, \phi(\gamma), (\det\phi(\gamma))/a)$.
Any element in $T$ has the form \dt \[t = \lp c\N_{\K/k}\alpha, c\phi\!\lp\alpha\bar{\beta}\rp\rp\]
for some $c \in k^\times$, $\alpha, \beta \in \K^\times$.

The norm correspondence between $G$ and $H_2$ gives the following isomorphism from $\bar{T}$ to $\bar{T}_2$:
\begin{equation*}
\lambda_2 : \lp c\N_{\K/k}\alpha, c\phi\!\lp\alpha\bar{\beta}\rp\rp \mapsto
\indexs{lambda@$\lambda_2$}%
(\diag(\alpha, \beta), c)_*,\quad \forall\; \alpha, \beta\in \K^\times,\; c \in k^\times.
\end{equation*}
Let $dt_2$ be the Haar measure on $\bar{T}_2$ which is compatible with the measure $dt$ on $\bar{T}$ via
$\lambda_2$.

For $t_2 \in T_2$ and $f_2 \in C(H_2)$, let
\[
O_{H_2}(f_2, t_2) = \int_{T_2 \bs H_2}f_2(\bar{h}_2^{-1}t_2 \bar{h}_2)\;d\bar{h}_2.
\]

Let $c$ be an element in $k^\times$.  Let $\alpha, \beta$ be elements in $\K^\times$.
Suppose $\alpha\bar{\beta}\notin k^\times$,  then the element 
$t = \lp c\N_{\K, k}\alpha, c\phi\!\lp\alpha\bar{\beta}\rp\rp \in T$ is regular.
A norm of $t$ in $H_2$ is conjugate to either
$t_2 =  (\diag(\alpha, \beta), c)_*$ or
${\bar{t}_2} = (\diag(\bar{\alpha}, {\bar{\beta}}), c)_*$.
Extend the character $\ve$ of $k^\times$ to $\K^\times$.
We make the following choice of transfer factors for the pairs $(t_2, t)$, $(\bar{t}_2, t)$:
\indexi{transfer factor}%
\[
\begin{split}
\Delta(t_2, t) = \Delta(\bar{t}_2, t)
& =
\lambda(\K/k, \psi)\ve(c\N_{\K,k}\alpha)\ve(c\alpha\bar{\beta} - c\bar{\alpha}\beta)
\left[\frac{D_{T \bs G}(t)}{D_{T_2 \bs H_2}(t_2)}\right]^{1/2}\\
& =
\lambda(\K/k, \psi)\ve(\alpha\bar{\beta} - \bar{\alpha}\beta)
\left[\frac{D_{T \bs G}(t)}{D_{T_2 \bs H_2}(t_2)}\right]^{1/2}.
\end{split}
\]
\begin{define}{\rm 
We say that $f \in C(G)$ and $f_2 \in C(H_2)$ are {\bf weakly matching} if
\indexi{matching functions!weakly}%
\[
O_{H_2}(f_2, t_2) = \Delta(t_2, t)O_G^{\ve}(f, t)
\]
for every $(t_2, t) \in T_2 \times T^\reg$ such that $t_2$ is a norm of $t$.
}\end{define}
\begin{prop}\label{prop:trGH2}
Let $\mu, \chi$ be  characters of $k^\times, \K^\times$, respectively.
For weakly matching functions $f$ in $\C(G)$ and $f_2$ in $\C(H_2)$, the following character identity
holds:
\[
\la \frac{\ve\mu}{\chi|_{k^\times}}\rtimes \pi(\chi), f \ra_{\!\!A}
=
\la I\!\lp\frac{\mu\circ\N}{\chi},\; \chi\rp \otimes_2 \mu, f_2 \ra.
\]
Here, $A$ is the intertwining operator defined in {\rm Section \ref{sec:intop}}.
\end{prop}
\begin{proof}
For simplicity, let $\N$ denote $\N_{\K/k}$.
By Lemma \ref{lemma:trG}, $\la \frac{\ve\mu}{\chi|_{k^\times}}\rtimes \pi(\chi), f \ra_A$
is equal to
\begin{multline}\label{eq:trG1}
\frac{1}{4}
\int_{t = (c\N\alpha,\; c\phi(\alpha\bar{\beta})) \in \bar{T}^\reg}
\ve\mu(c)\\
{\cdot\bigg[
\mu(\N\alpha)\lp\chi\lp\tfrac{\beta}{\alpha}\rp + \chi\lp\tfrac{\bar{\beta}}{\bar{\alpha}}\rp\rp +
\mu(\N\beta) \lp\chi\lp\tfrac{\alpha}{\beta}\rp + \chi\lp\tfrac{\bar{\alpha}}{\bar{\beta}}\rp\rp
\bigg]}\\
\cdot\lambda(\K/k, \psi)\ve(c\alpha\bar{\beta} - c\bar{\alpha}\beta)D^{1/2}_{T \bs G}(t) O_G^\ve(f, t)\;dt.
\end{multline}

Employing the same technique used in the proof of Lemma \ref{lemma:trG}, we obtain:
\begin{multline}\label{eq:trH2}
\la I\!\lp\tfrac{\mu\circ\N}{\chi}, \chi\rp \otimes_2\mu, f_2 \ra =\\
\shoveleft{\frac{1}{\abs{\Gal(\K/k)}}\cdot\frac{1}{\abs{W_{T_2\bs H_2}}}}
\\\cdot\int_{t_2 = (\diag(\alpha, \beta), c)_* \in \bar{T}_2^\reg}
 \bigg[
\mu(\N\alpha)\lp\chi\lp\tfrac{\beta}{\alpha}\rp + \chi\lp\tfrac{\bar{\beta}}{\bar{\alpha}}\rp\rp +
\mu(\N\beta) \lp\chi\lp\tfrac{\alpha}{\beta}\rp + \chi\lp\tfrac{\bar{\alpha}}{\bar{\beta}}\rp\rp
\bigg]\\
\cdot \mu(c)
D^{1/2}_{T_2\bs H_2}(t_2)O_{H_2}(f_2, t_2)\;dt_2.
\end{multline}

By our choice of $dt_2$, and the fact that $\abs{\Gal(\K, k)} = \abs{W_{T_2 \bs H_2}} = 2$, the above expression is equal to
\begin{multline*}
\frac{1}{4}\int_{t = (c\N\alpha,\; c\phi(\alpha\bar{\beta})) \in \bar{T}^\reg}
\bigg[
\mu(\N\alpha)\lp\chi\lp\tfrac{\beta}{\alpha}\rp + \chi\lp\tfrac{\bar{\beta}}{\bar{\alpha}}\rp\rp +
\mu(\N\beta) \lp\chi\lp\tfrac{\alpha}{\beta}\rp + \chi\lp\tfrac{\bar{\alpha}}{\bar{\beta}}\rp\rp
\bigg]
\\\cdot \mu(c) D^{1/2}_{T_2\bs H_2}(\lambda_2(t))O_{H_2}(f_2, \lambda_2(t))\;dt.
\end{multline*}
The proposition now follows from the weakly matching condition on $f, f_2$.
\end{proof}
{\sc Remark.}
By the Weyl integration formula and the Iwasawa decomposition for real reductive groups, 
it is clear that the archimedean analogues of Propositions \ref{prop:trGH1}, \dts \ref{prop:trGH2} also hold.  
That is: If $k = \mbb{R}$, $\K = \mbb{C}$, and
$\mu$, $\chi$ are characters of $\mbb{R}^\times$, $\mbb{C}^\times$, respectively, then
\[
\begin{split}
\la \mu \rtimes \pi(\chi), f\ra_A &= \la I(\ve\mu, 1)\otimes_1 \chi, f_1\ra,\\
\la \tfrac{\ve\mu}{\chi|_{k^\times}}\rtimes \pi(\chi), f \ra_{\!\!A} &=
\la I\!\lp\tfrac{\mu\circ\N}{\chi},\; \chi\rp \otimes_2 \mu, f_2 \ra
\end{split}
\] 
for matching functions on $\mb{G}(\mbb{R})$, $\mb{H_1}(\mbb{R})$, $\mb{H_2}(\mbb{R})$.  
Here, $\ve(x)$ is the sign of $x$ for all $x \in \mbb{R}$.

\subsection{Examples}
For an admissible  representation $\tau$ of $\GL(2, k)$, let $B_{\K/k}\tau$ denote the representation of
$\GL(2, \K)$ which is obtained from $\tau$ via base change (\cite{L1}, \cite{F4}).  
\indexi{base change}%
In particular, for characters $\mu_1, \mu_2$ of $k^\times$, $B_{\K/k}I(\mu_1, \mu_2)$ is the induced representation
$I(\mu_1\circ\N_{\K/k},\; \mu_2\circ\N_{\K/k})$.  


For $x \in k$, let $\ord x$ denote the $p$-adic valuation of $x$, with $\ord 0 := \infty$.
Let $q$ be the cardinality of the residue field of $k$.
Let $\nu : x \mapsto q^{-\ord x}$ be the normalized absolute value function on $k$.
\begin{corollary}\label{corollary:trGH1'}
Let $\theta$ be a character of $\K^\times$.
The following holds for weakly matching functions $f \in C(G)$ and $f_1 \in C(H_1)\!:$
\[
\la \ve\nu\rtimes\nu^{-1/2}\pi(\theta), f\ra_A =
\la I(\nu^{1/2}, \nu^{-1/2})\otimes_1 \theta, f_1 \ra.
\]
\end{corollary}
\begin{proof}
Noting that
$
I(\nu^{1/2}, \nu^{-1/2})\otimes_1 \theta =
I(\nu, 1)\otimes_1 \lp\nu^{-1/2}\circ\N_{\K/k}\rp\cdot\theta,
$
the corollary follows from Proposition \ref{prop:trGH1}.
\end{proof}
\begin{corollary}\label{corollary:trGH1H2}
Let $\eta$ be a character of $k^\times$.
For weakly matching functions $f, f_1, f_2$ on $G, H_1, H_2$, respectively,
the following quantities are equal to one another:
\begin{itemize}
\item
$\la \ve\nu\times\ve\rtimes\nu^{-1/2}\eta, f\ra_A$,
\item
$\la \eta I(\nu^{1/2},\nu^{-1/2})\otimes_1 1, f_1 \ra$,
\item
$\la B_{\K/k}I(\eta\nu^{1/2}, \eta\nu^{-1/2})\otimes_2 \eta^2, f_2 \ra$.
\end{itemize}
\end{corollary}
\begin{proof}
The representation $\ve\nu\times\ve\rtimes\nu^{-1/2}\eta$ is equal to 
$\ve\nu\rtimes \pi\lp(\nu^{-1/2}\eta)\circ\N_{\K/k}\rp$.
The representation $\eta I(\nu^{1/2},\nu^{-1/2})\otimes_1 1$ of $H_1$ is equal to
$I(\nu, 1)\otimes_1 (\nu^{-1/2}\eta)\circ\N_{\K/k}$.  
The corollary follows from Corollary \ref{corollary:trGH1'} and Proposition \ref{prop:trGH2}.
\end{proof}

Let $\nu_\K$ be the normalized absolute value function on $\K$. 
\begin{corollary}\label{corollary:trGH2'}
Suppose the extension $\K/k$ is unramified.
Let $\theta$ be a character of $\K^\times$ which satisfies 
${{}^\sigma\theta}/{ \theta} = \ve'\circ\N_{\K/k}$ for some quadratic character $\ve'$ of
$k^\times$.  The following holds for weakly matching functions $f \in C(G)$ and $f_2 \in C(H_2)\!:$
\[
\la \ve\ve'\nu\rtimes\nu^{-1/2}\pi(\theta), f\ra_A =
\la \theta I(\nu_{\K}^{1/2}, \nu_{\K}^{-1/2})\otimes_2 \theta|_{k^\times}\cdot\ve',\; f_2 \ra.
\]
\end{corollary}
\begin{proof}
Since $\K/k$ is unramified, $\nu_\K$ is equal to $\nu\circ\N_{\K/k}$.
The corollary then follows from Proposition \ref{prop:trGH2}.
\end{proof}
\section{Induced Representations---Split Case}
\indexi{representation!induced}%
Let $F$ be a number field.  Let $\idc{F}$ denote the id\`ele class group of $F$.%
Let $\E$ be a nontrivial quadratic character of $\idc{F}$.  Let
$E$ be the quadratic extension of $F$ which corresponds to $\E$ via global class field theory.

Regard $\mb{G} = \GSp(2)$ as an $F$-group. Let $\mb{H_1}$ and $\mb{H_2}$ be the 
elliptic $\E$-endoscopic groups of $\mb{G}$ over $F$.
\indexi{endoscopic group}\indexi{H1@$\mb{H_1}$}\indexi{H2@$\mb{H_2}$}%

Let $v$ be a (possibly archimedean) place of $F$ which splits in $E$.  Then, $E_v := E\otimes_F F_v = F_v \oplus F_v$.
We have:
\begin{multline*}
\mb{H_1}(F_v)
= \lp \GL(2, F_v)\times F_v^\times \times F_v^\times\rp'\\
:= \left\{ (g, a, b) \in \GL(2, F_v) \times F_v^\times \times F_v^\times
: \det g = ab\right\},
\end{multline*}
\[\mb{H_2}(F_v) =
\lp \GL(2, F_v)\times\GL(2, F_v)\times F_v^\times\rp/
\left\{
\lp a I_2, b I_2, (ab)^{-1}\rp : a, b\in F_v^\times
\right\}.\]

Let $G_v = \mb{G}(F_v)$, $H_{i,v} = \mb{H}_i(F_v)$ ($i = 1, 2$).

Let $M_v$ denote the Levi component of the Heisenberg parabolic subgroup of $G_v$.
The norm correspondence between $\mb{G}$ and $\mb{H_1}$ gives the following isomorphism
from $M_v$ to $H_{1,v}$:
\[
\lambda_1 : (a, g_2) \mapsto (g_2, a, (\det g_2)/a), \quad \forall a \in k^\times,\; g_2 \in \GL(2, k).
\indexs{lambda@$\lambda_1$}%
\]
Here, $(a, g_2)$ denotes the element $\blockdiag(a, g_2, (\det g_2)/a)$ in $M_v$.
In particular, $\lambda_1$ defines a one to one correspondence between the conjugacy classes of the maximal tori in $M_v$
and those of the maximal tori in $H_{1,v}$.

Each maximal torus in $M_v$ has the form
\[
T_v = \{(a, t'): a \in F_v^\times, t' \in T_v'\},
\]
where $T_v'$ is a maximal torus in $\GL(2, F_v)$.  


A norm in $H_{1,v}$ of a regular element $t = (a, t') \in M_v$ is conjugate to either $t_1 = (t', a, (\det t') /a)$
or $\bar{t}_1 = (t', (\det t')/a, a)$.  Let $T_v, T_{1,v}$ be the maximal tori which contain $t, t_1$, respectively.
We define the transfer factors for the pairs $(t_1, t)$, $(\bar{t}_1, t)$ as follows:
\indexi{transfer factor}%
\[
\Delta(t_1, t) = \Delta(\bar{t}_1, t) = \left[\frac{D_{T_v\bs G_v}(t)}{D_{T_{1,v}\bs H_{1,v}}(t_1)}\right]^{1/2}.
\]

Let $Z_v$ be the center of $G_v$.
Let $C(G_v)$ denote the space of functions on $G_v$ which are smooth and compactly supported modulo $Z_v$.
Let $Z'(H_{1,v})$ denote the subgroup\dt
\[\{(\diag(z, z), z, z) : z \in F_v^\times\} \subset H_{1,v}.\]  
It is the group of $F_v$-points of the
maximal $F$-split component of the center of $\mb{H}_1$.
Let $C(H_{1,v})$ denote the space of functions on $H_{1,v}$ which are smooth and compactly supported modulo $Z'(H_{1,v})$.

For a maximal torus $T_{1,v}$ in $H_{1,v}$, put $\bar{T}_{1,v} := T_{1,v} / Z'(H_{1,v})$.
Put $\bar{G}_v := G_v / Z_v$, $\bar{H}_{1,v} := H_{1,v} / Z'(H_{1,v})$.
We fix a Haar measure $dt$ on each maximal torus $\bar{T}_v$ in $\bar{G}_v$.  It induces via $\lambda_1$ a measure $dt_1$
on the maximal torus $\bar{T}_{1,v} = \lambda_1(\bar{T}_v)$ in $\bar{H}_{1,v}$.

We say that $f_v \in C(G_v)$ and $f_{1,v} \in C(H_{1,v})$ are {\bf weakly matching} if
\indexi{matching functions!weakly}%
\[
O_{H_{1,v}}(f_{1,v}, t_1) = \Delta(t_1, t) O_{G_v}(f_v, t)
\]
for all $(t_1, t) \in H_{1,v}\times M_v^\reg$ such that  $t_1$ is a norm of $t$.

\begin{lemma}\label{lemma:indsplitH1}
Let $\tau$ be an admissible representation of $\GL(2, F_v)$.  Let $\mu_1, \mu_2$ be
characters of $F_v^\times$.  Then,
\[
\la \frac{\mu_1}{\mu_2}\rtimes \mu_2\tau, f_v \ra = \la \tau\otimes_1 \lp\mu_1\otimes\mu_2\rp, f_{1,v} \ra
\]
for weakly matching functions.
\end{lemma}
\begin{proof}
The proof is similar to that of Proposition \ref{prop:trGH1}.

By the Weyl integration formula,
$\la \tfrac{\mu_1}{\mu_2} \rtimes \mu_2\tau, f \ra$ is equal to
\begin{multline*}
\sum_{T_v} \frac{1}{\abs{W_{T_v\bs M_v}}}
{\int_{t = (a, t') \in \bar{T}_v^\reg} D_{T_v\bs M_v}(t)
\tfrac{\mu_1}{\mu_2}(a)\chi_{\mu_2\tau}(t')
D^{1/2}_{M_v\bs G_v}(t)O_{G_v}(f_v, t)\; dt} \\
\shoveleft{=
\sum_{T_v} \frac{1}{\abs{W_{T_v\bs M_v}}}}\\\cdot \int_{t = (a, t') \in \bar{T}_v^\reg}
D_{T_v\bs M_v}^{1/2}(t) \mu_1(a)\mu_2\!\lp \tfrac{\det t'}{a}\rp \chi_{\tau}(t')
D^{1/2}_{T_v\bs G_v}(t)O_{G_v}(f_v, t)\; dt.
\end{multline*}
Here, $T_v$ ranges over a set of representatives of the conjugacy classes of the maximal tori in $M_v$, and
$W_{T_v\bs M_v}$ is the Weyl group of $T_v$ in $M_v$.
The last equality holds because $D_{T_v\bs M_v} D_{M_v \bs G_v} = D_{T_v \bs G_v}$.

The trace $\la\tau\otimes_1 (\mu_1\otimes\mu_2),f_{1,v}\ra$ is equal to
\begin{multline*}\label{eq:traceh1split}
\sum_{{T}_{1,v}} \frac{1}{\abs{W_{T_{1,v}\bs H_{1,v}}}}
\int_{t_1 = \lp t',\; a,\; (\det t') /a\rp \in \bar{T}_{1,v}^\reg} 
D^{1/2}_{T_{1,v}\bs H_{1,v}}(t_1)\mu_1(a)\mu_2\lp \tfrac{\det t'}{a}\rp\chi_{\tau}(t')
\\ \shoveright{\cdot D^{1/2}_{T_{1,v}\bs H_{1,v}}(t_1) O_{H_{1,v}}(f_{1,v}, t_1)\; dt_1}
\\ 
\shoveleft{
= \sum_{\lambda_1({T}_{v})} \frac{1}{\abs{W_{\lambda_1(T_{v})\bs H_{1,v}}}}
\int_{t = (a, t') \in \bar{T}_{v}^\reg} 
D^{1/2}_{T_{1,v}\bs H_{1,v}}(\lambda_1(t))\mu_1(a)\mu_2\lp \tfrac{\det t'}{a}\rp\chi_{\tau}(t')
}
\\ \cdot D^{1/2}_{T_{1,v}\bs H_{1,v}}(\lambda_1(t)) O_{H_{1,v}}(f_{1,v}, \lambda_1(t))\; dt.
\end{multline*}
Here, $T_{1,v}, T_v$ range over sets of representatives of the conjugacy classes of the maximal tori in $H_{1,v}$,
$G_v$, respectively.

By assumption, $f_v,f_{1,v}$ are weakly matching functions; hence,
\[
D_{T_v\bs G_v}^{1/2}(t) O_{G_v}(f_v, t)
= D_{T_{1,v}\bs H_{1, v}}^{1/2}(\lambda_1(t)) O_{H_{1,v}}(f_{1,v}, \lambda_1(t))
\]
for all regular $t$ in $M_v$.  Moreover, $\abs{W_{T_v\bs M_v}} = \abs{W_{\lambda_1(T_v)\bs H_{1,v}}} = 2$ for all maximal
tori $T_v$ in $M_v$.
The lemma follows.
\end{proof}
Let $T_{2,v}$ be the maximal diagonal torus in $H_{2,v}$.
Let $T_v$ be the maximal diagonal torus in $G_v$.
Any element in $T_v$ has the form $t = \diag(ac, ad, bc, bd)$, $a, b, c, d \in F_v^\times$.
The norm correspondence between $\mb{G}$ and $\mb{H_2}$ gives the following isomorphism from
$T_v$ to $T_{2,v}$:
\[
\lambda_2 : \diag(ac, ad, bc, bd) \mapsto \lp\diag(a, b), \diag(c, d), 1\rp_*.
\indexs{lambda@$\lambda_2$}%
\]
Here, lower $*$ denotes image in $H_{2,v}$.

A norm in $H_{2, v}$ of a regular element $t = \diag(ac, ad, bc, bd) \in T_v$ is conjugate to either 
$t_2 = (\diag(a, b), \diag(c, d), 1)_*$ or $\bar{t}_2 = (\diag(c, d), \diag(a, b), 1)_*$.
We define the transfer factors for the pairs $(t_2, t)$, $(\bar{t}_2, t)$ as follows:
\indexi{transfer factor}%
\[
\Delta(t_2, t) = \Delta(\bar{t}_2, t) = \left[\frac{D_{T_v\bs G_v}(t)}{D_{T_{2,v}\bs H_{2,v}}(t_{2,v})}\right]^{1/2}.
\]

Let $Z'(H_{2,v}) = \{(\diag(1, 1), \diag(1, 1), z)_* : z \in F_v^\times\} \subset H_{2,v}$.
It is the group of $F_v$-points of the maximal $F$-split component of the center of $\mb{H_2}$.
Let $C(H_{2,v})$ denote the space of smooth, compactly supported modulo $Z'(H_{2,v})$ functions on $H_{2,v}$.

Put $\bar{T}_{2,v} := T_{2,v} / Z'(H_{2,v})$.
As usual, the Haar measure $dt_2$ on $\bar{T}_{2,v}$ is chosen such that it is compatible with the measure
$dt$ on $\bar{T}_v$ via $\lambda_2$.
We say that the functions $f_v \in C(G_v)$ and $f_{2,v} \in C(H_{2,v})$ are {\bf weakly matching} if
\indexi{matching functions!weakly}%
\[
O_{H_{2,v}}(f_{2,v}, t_2) = \Delta(t_2, t) O_{G_v}(f_v, t)
\]
for every $(t_2, t) \in T_{2,v} \times T_v^\reg$ such that $t_2$ is a norm of $t$.
\begin{lemma}\label{lemma:indsplitH2}
Let $\alpha_1, \alpha_2, \beta_1, \beta_2$ be characters of $F_v^\times$ such that
$\alpha_1\beta_1 = \alpha_2\beta_2$.  Then,
\[
\la \beta_1^{-1}I(\alpha_2, \beta_2)\rtimes \beta_1, f_v\ra
= \la I(\alpha_1, \beta_1) \times I(\alpha_2, \beta_2), f_{2,v}\ra
\]
for weakly matching functions $f_v \in \C(G_v)$, $f_{2,v} \in \C(H_{2,v})$.
\end{lemma}
\begin{proof}
The induced representations are supported only on the conjugacy classes of the maximal diagonal tori.
As in the proof of Lemma \ref{lemma:indsplitH1}, the character identity follows from the Weyl integration formula
and the weakly matching condition on the test functions.
\end{proof}
\section{$\ve$-Invariant Packets}
\indexi{packet!epsilon-invariant@$\ve$-invariant}\indexi{packet!local}%
Recall that $k$ is our fixed $p$-adic field, $\ve$ is a nontrivial quadratic character of $k^\times$, and
$\K$ is the quadratic extension of $k$ corresponding to $\ve$.
Let $G = \mb{G}(k)$, $H_i = \mb{H}_i(k)$ ($i = 1, 2$).

Let $\bar{G} = \PGSp(2, k)$.
\indexs{P@$\PGSp(2)$}%
Let $C_0 = \PGL(2, k)\times\PGL(2, k)$.
\indexs{C@${C_0}$}%
For representations $\tau_1, \tau_2$ of $\PGL(2, k)$, let $\tau_1\times\tau_2$ denote the $C_0$-module
where $(g, h) \in C_0$ acts by $\tau_1(g)\otimes\tau_2(h)$ on the tensor product of the spaces of $\tau_1, \tau_2$.


Let $\tau$ be a cuspidal or one dimensional representation of $\GL(2, k)$ with trivial central
character. 
Write $\dagger$ for $+$ if $\tau$ is cuspidal, or $\times$ if $\tau$ is one dimensional.
Following the same notation as in \cite{F1}, we let $\{\pi^\dagger, \pi^-\}$ denote 
the local (quasi-)packet of $\bar{G}$ which lifts to the representation
$I_{(2, 2)}(\tau, \ve\tau)$ of $\PGL(4, k)$.
In the terminology of \cite{F1}, which we also adopt,
$\{\pi^+, \pi^-\}$ is a packet, while $\{\pi^\times, \pi^-\}$ is a quasi-packet.
\indexi{packet}\indexi{packet!quasi-}\indexs{pi@$\pi^+$, $\pi^-$}\indexs{pi@$\pi^\times$, $\pi^-$}%

\begin{lemma}\label{lemma:veinvariantlocal}
The representations $\pi^\dagger$ and $\pi^-$ are $\ve$-invariant.
\end{lemma}
\begin{proof}
Let $f, f_4, f_{C_0}$ be matching functions on $\bar{G}$, $\PGL(4, k)$, $C_0$, respectively.
Let $\zeta = 1$ or $\ve$.  By \cite[Prop. V. 5]{F1} and \cite[Prop. V. 8.6]{F1}, we have:
\begin{equation*}
\begin{split}
\la \zeta\pi^\dagger, f\ra &= 
\frac{1}{2} \Big( \la \zeta I_{(2, 2)}(\tau, \ve\tau), f_4\ra + \la \zeta\tau\times\zeta\ve\tau, f_{C_0}\ra\Big)\\
\la \zeta\pi^-, f\ra &= 
\frac{1}{2} \Big( \la \zeta I_{(2, 2)}(\tau, \ve\tau), f_4\ra - \la \zeta\tau\times\zeta\ve\tau, f_{C_0}\ra\Big).
\end{split}
\end{equation*}
The representation  $I_{(2, 2)}(\tau, \ve\tau)$ is $\ve$-invariant.
By the norm correspondence between $\bar{G}$ and $C_0$, and the matching condition on the test functions,
the convolution operator $(\tau\times\ve\tau)(f_{C_0})$ is equal to $(\ve\tau\times\tau)(f_{C_0})$.
Hence, the right-hand sides of the above equations are independent of the choice of $\zeta$.
Consequently,
\[
\begin{split}
\la \pi^\dagger, f\ra &= \la \ve\pi^\dagger, f\ra,\\
\la \pi^-, f\ra &= \la \ve\pi^-, f\ra.
\end{split}
\]
By the linear independence of characters, we conclude that $\pi^\dagger$ and $\pi^-$ are $\ve$-invariant.
\end{proof}
\begin{corollary}\label{corollary:veinvariantlocal}
Let $F$ be a totally imaginary number field.  Let $\E$ be a nontrivial quadratic character of the \idele class group of $F$.
Let $\mct$ be a cuspidal or one dimensional automorphic representation of $\GL(2, \Af)$
with trivial central character.  Let $\{\Pi\}$ be the unstable $($quasi$)$-packet $[\mct, \E\mct]$ of $\mb{G}(\Af)$
{\rm (see Section \ref{sec:arthurunstable})}.
Then, every member of $\{\Pi\}$ is $\E$-invariant.
\end{corollary}
\begin{proof}
Let $\Pi$ be a member of $\{\Pi\}$. By Lemma \ref{lemma:veinvariantlocal}, $\Pi_v$ is $\E_v$-invariant for each finite place $v$.

Let $E$ be the quadratic extension of $F$ corresponding to $\E$.  Let $v$ be an archimedean place of $F$.
Since $F$ is totally imaginary, $F_v$ is algebraically closed.  Thus, $v$ splits in $E$, which implies $\E_v = 1$.
\end{proof}
\section{Character Identities for Unstable Packets}\label{sec:bc}
\indexi{packet!unstable|(}\indexi{packet!local|(}%

{\bf Notation}
\begin{itemize}
\item
For a number field $F$, let $\idc{F}$ denote the \idele class group of $F$.
\item
Suppose $S$ is a set of places of a number field $F$.
For any ad\`elic object (representation, test function, trace, \dots, etc) over $F$, 
let subscript $S$ denote the tensor product of local components over the places in $S$.
For example, if $\Pi$ denotes an irreducible representation of $\mb{G}(\Af)$, we put $\Pi_S := \otimes_{v \in S}\Pi_v$.
\item
For an algebraic group $\mb{H}$ over a number field $F$, put $H_v := \mb{H}(F_v)$ for each place
$v$ of $F$.
\end{itemize}

\subsection{One Dimensional Representations}\label{sec:bc1dim}
\indexi{representation!one dimensional}%
Let $1_{2,l}$  ($l = k, \K$) denote the trivial representation of $\GL(2, l)$.  It is nontempered.
\indexs{1@$1_2$}%
Let ${\rm St}_{2,l}$ denote the Steinberg representation of $\GL(2, l)$.  It is square integrable.
\indexs{S@${\rm St}_2$}%

Let $\xi$ be a character of $k^\times$ with $\xi^2 = 1$.
Let $\{\pi^\times, \pi^-\}$ be the local quasi-packet of $G$ which 
lifts to the representation $I_{(2, 2)}(\xi 1_{2,k},\; \ve \xi 1_{2,k})$ of $\PGL(4, k)$
according to {\rm \cite[Prop. V. 8.5]{F1}}.
In other words, $\pi^\times$ is the  nontempered subquotient $L(\nu\ve, \ve\rtimes \nu^{-1/2}\xi)$ of
\indexs{L@$L({\rm char., repn.})$}%
the induced representation $\nu\ve\times\ve\rtimes\nu^{-1/2}\xi$ (see \cite{ST}), and
$\pi^- = \delta^-(\ve\nu^{1/2}{\rm St}_{2, k}, \nu^{-1/2}\xi)$, the cuspidal member
\indexs{delta@$\delta^-({\rm repn., char.})$}%
of the local packet which contains the unique square integrable subrepresentation 
$\delta(\ve\nu^{1/2}{\rm St}_{2, k}, \nu^{-1/2}\xi)$ of $\nu\ve\times\ve\rtimes\nu^{-1/2}\xi$ (see \cite[Sect. V. 8]{F1}).
By Lemma \ref{lemma:veinvariantlocal}, $\pi^\times$ and $\pi^-$ are $\ve$-invariant representations.
\begin{prop}\label{prop:bc1dim}
Suppose $\ve$ is unramified.
There exist nontrivial intertwining operators $A^\times \in \Hom_G(\pi^\times, \ve\pi^\times)$,
$A^- \in \Hom_G(\pi^-, \ve\pi^-)$
such that the following identities hold for matching functions:
\[
\begin{split}
\la \pi^\times, f\ra_{A^\times} + \la \pi^-, f\ra_{A^-} &= \la \xi 1_{2, k}\otimes_1 1, f_1 \ra,\\
\la \pi^\times, f\ra_{A^\times} - \la \pi^-, f\ra_{A^-} &= \la (\xi\circ\N_{\K/k})1_{2, \K}\otimes_2 1, f_2\ra.
\end{split}
\]
\begin{remark}
It is likely that the proposition still holds if we remove the condition that $\ve$ is unramified; however, the
author is unable to prove it.
\end{remark}
\end{prop}
The proof will be given after Claim \ref{claim:beforebc1dimproof}.  

Using Krasner's lemma and the strong approximation theorem (see \cite{PR}), we
construct a totally imaginary number field $F$ and a quadratic extension
$E$ of $F$ such that:
\begin{itemize}
\item
There exists a finite place $w$ of $F$ such that $F_w = k$;
\item
the place $w$ is unramified, prime in $E$;
\item
$E_w = \K$.
\end{itemize}
Here, $w$ denotes both the place of $F$ and the unique place of $E$ which lies above it.  
Let $V$ be the set of places of $F$.  Let $\Vr$ be the set of places of $F$ which are unramified in $E$.

Let $\E$ be the quadratic character of  $\idc{F}$ which corresponds to the extension $E/F$ via global class field theory.
In particular, $\E_w = \ve$.

Consider $\mb{G}$ as an $F$-group.  Let $\mb{H_1}, \mb{H_2}$ be the elliptic $\E$-endoscopic groups of $\mb{G}$ over $F$.
For $i =1, 2$, we have $\mb{H}_i(F_w) = H_i$.
We fix Tamagawa measures (\cite{PR}) $dg = \otimes_{v}dg_v$ on $\mb{G}(\Af)$, and $dh_i = \otimes_vdh_{i,v}$ on $\mb{H}_i(\Af)$,
such that $dg_w$ and $dh_{i,w}$ coincide with the local Haar measures we have chosen for $G$ and $H_i$, respectively.

Let $1_2 $ be the trivial representation of $\GL(2, \Af)$.  Let $1_{2, E}$ be the trivial representation of $\GL(2, \Ae)$.
Let $\Z$ be a character of $\idc{F}$ such that $\Z^2 = 1$ and $\Z_w = \xi$.
Let $\pi_1$ be the one dimensional representation $\Z 1_2\otimes_1 1$ of $\mb{H_1}(\Af)$.
Let $\pi_2$ be the one dimensional representation 
$(\Z \circ \N_{E/F})1_{2, E}\otimes_2 1$
of $\mb{H_2}(\Af)$.
Let \[\{\Pi\} = \bigotimes_{v \in V} \{\Pi_v^\times, \Pi_v^-\}\]
be the unstable quasi-packet $[\Z 1_2, \E \Z 1_2]$ of $\mb{G}(\Af)$ (see Section \ref{sec:arthurunstable}).
It lifts to the induced representation
$I_{(2, 2)}(\Z 1_{2}, \E \Z 1_{2})$ of $\GL(4, \Af)$.
From \cite[Sect. V. 10.6]{F1}, $\Pi_v^\times = L(\ve_v\nu_v, \ve_v\rtimes\nu_v^{-1/2}\xi_v)$ for each place $v$, and
\[
\Pi^-_v = \begin{cases}
\delta^-(\ve_v\nu_v^{1/2}{\rm St}_{2,v}, \nu^{-1/2}\xi_v) & \text{ if } \ve_v \neq 1,\\
L(\nu_v^{1/2}{\rm St}_{2,v}, \nu_v^{-1/2}\xi_v) & \text{ if } \ve_v = 1.
\end{cases}
\]
Here, $\nu_v$ is the normalized absolute value function on $F_v$, 
and ${\rm St}_{2,v}$ denotes the Steinberg representation of $\GL(2, F_v)$.
In particular, $\Pi_w^\times = \pi^\times$ and $\Pi_w^- = \pi^-$.


Fix a place $u$ of $F$ which splits in $E$.
Then, $\E_u$ is trivial, and
\[
H_{1,u} = \mb{H}_1(F_u) = \{ (g, a, b) \in \GL(2, F_u) \times F_u^\times \times F_u^\times : 
\det g = ab\}.
\]
We have $\pi_{1 ,u} = \Z_u 1_{2,u} \otimes_1 (1 \otimes 1)$.  
Here, $1_{2, u} := 1_{\GL(2, F_u)}$, the trivial representation of $\GL(2, F_u)$.
The representations 
\[
\Pi_u^\times = L(\nu_v, 1\rtimes\nu_v^{-1/2}\xi_v), \quad \Pi_u^- = L(\nu_v^{1/2}{\rm St}_{2,v}, \nu_v^{-1/2}\xi_v)
\] 
are the two inequivalent nontempered subquotients of $1 \rtimes \Z_u 1_{2,u}$ (see \cite[Lemma 3.8]{ST}).  
\begin{lemma}\label{lemma:trGH1u}
For matching functions $f_u$ on $G_u$ and $f_{1,u}$ on $H_{1,u}$, the following holds:
\[
\la \Pi_u^\times, f_u\ra + \la\Pi_u^-, f_u\ra=
\la 1 \rtimes  \Z_u 1_{2,u},\; f_u\ra
= \la \pi_{1,u}, f_{1,u}\ra.
\]
\end{lemma}
\begin{proof}
The first equality is clear.  The second follows from Lemma \ref{lemma:indsplitH1}.
\end{proof}
The group of $F_u$-points of $\mb{H_2}$ is
\[
H_{2,u} = \lp\GL(2, F_u) \times \GL(2, F_u) \times F_u^\times\rp/
\{(a I_2 ,b I_2, (ab)^{-1}): a,b\in F_u^\times\}.
\]

For a character $\omega_u$ of the center $Z_u$ of $G_u$,
let $C(G_u, \omega_u)$ denote the space of smooth functions on $G_u$ which are compactly supported modulo
$Z_u$ and transform under $Z_u$ via $\omega_u^{-1}$.
If $\omega_u = 1$, then each $f_u \in C(G_u, 1)$ defines a test function $\bar{f}_u$ on $\PGSp(2, F_u)$.

The map $\diag(z, z, z, z) \mapsto (\diag(1, 1), \diag(1, 1), z)_*$, $z \in F_u^\times$, defines an isomorphism from $Z_u$ to
\[
Z'(H_{2,u}) := \{(\diag(1, 1), \diag(1, 1), z)_* : z\in F_u^\times\}\subset H_{2,u}.
\] 
Here, lower $*$ denotes image in $H_{2,u}$.
For a character $\omega_{2,u}$ of the group $Z'(H_{2,u})$, let
$C(H_{2,u}, \omega_{2,u})$ denote the space of smooth, compactly supported mod $Z'(H_{2,u})$ functions on $H_{2,u}$
which transform under $Z'(H_{2,u})$ via $\omega_{2,u}^{-1}$.
If a function $f_{2, u}$ in $C(H_{2,u}, \omega_{2,u})$ matches some $f_u \in C(G_u, 1)$, then $\omega_{2,u}$ is 
necessarily trivial.

The map $(g, h, c)_* \mapsto (g, ch)$ defines an isomorphism from $H_{2,u}$ to the group
\[
(\GL(2, F_u)\times\GL(2, F_u))/\{(z I_2, z^{-1}I_2):z\in F_u^\times\}.
\]
Hence, the representations of $H_{2,u}$ with trivial central characters may be identified
with the representations of $\PGL(2, F_u) \times \PGL(2, F_u)$, which is the 
group of $F_u$-points of the $F$-group $\mb{C_0} = \PGL(2)\times\PGL(2)$.
Moreover, if a function $f_{2,u}$ in $C(H_{2,u}, 1)$ matches a function $f_u$ in $C(G_u, 1)$, then
$f_{2,u}$ defines a function $\bar{f}_{2, u}$ on $C_{0, u}$  which matches $\bar{f}_u$ in the context of twisted endoscopy 
for $\PGL(4)$ (see \cite{F1}).

\begin{lemma}\label{lemma:trGH2u}
The following holds for matching functions $f_u \in C(G_u, 1)$ and $f_{2,u} \in C(H_{2,u}, 1)\!:$
\[
\la\Pi_u^\times,f_u\ra - \la\Pi_u^-,f_u\ra =
\la \Z_u 1_{2, u} \times \Z_u 1_{2, u},\; \bar{f}_{2, u}\ra = \la \pi_{2,u}, f_{2,u}\ra.
\]
\end{lemma}
\begin{proof}
The first equality follows from \cite[Prop. V. 8.6]{F1}.  
Let $u_1, u_2$ be the distinct places of $E$ which lie above $u$.
For $i = 1, 2$, we have $E_{u_i} = F_u$, $1_{2, E_{u_i}} = 1_{2, u}$, and $\lp\Z\circ\N_{E/F}\rp_{u_i} = \Z_u$.
Hence, $\pi_{2,u} = \lp \Z_{u} 1_{2, u}\times \Z_{u} 1_{2, u}\rp\otimes_2 1$, which
may be identified with the representation $\Z_u 1_{2,u} \times \Z_u 1_{2,u}$ of $C_{0, u}$.
The second equality follows.
\end{proof}
In summary, the following local identities hold for matching functions:
\begin{equation}\label{eq:leverage}
\begin{split}
\la \Pi_u^\times,f_u\ra + \la \Pi_u^-,f_u\ra &= \la \pi_{1,u},f_{1,u}\ra,\\
\la \Pi_u^\times,f_u\ra - \la \Pi_u^-,f_u\ra  &= \la\pi_{2,u},f_{2,u}\ra.
\end{split}
\end{equation}

Let $S'$ be the finite set of bad places for $\Z$, i.e. it is
the union of  $V - \Vr$ and the set of places $v$ where $\Z_v$ is ramified.
In particular, the fixed places $w$, $u$ lie outside of $S'$.
For each place $v \notin S'$, 
$\pi_{1,v}$ and $\pi_{2,v}$ are unramified representations.

Fix two places $w_1, w_2 \notin S' \cup \{w\}$ which are prime in $E$.
For a character $\omega$ of the center of $\mb{G}(\Af)$, recall the definition of
$\ELw{\mb{G}(\Af), \omega}$ (resp. $\ELw{\mb{H}_i(\Af), \omega}$):
It denotes the space of functions in $C(\mb{G}(\Af), \omega)$ (resp. $C(\mb{H}_i(\Af), \omega)$)
whose local components at $w_1, w_2$ are elliptic.

Let $S_0 = S'\cup\{w_1, w_2\}$.
Let $S$ be a finite set of places of $F$ containing $S_0$.

By Corollary \ref{corollary:veinvariantlocal}, each representation $\Pi'$ in $\{\Pi\}$ is $\E$-invariant;
hence, the twisted character
$\la \Pi', f\ra_\E$ is equal to a product $\prod_v \la \Pi'_v, f_v\ra_{\E_v}$ of twisted local characters.
Here, $\la \Pi'_v, f_v\ra_{\E_v}$ denotes the trace of $\Pi'_v(f_v)\rho(\E)_v$, 
where $\rho(\E)_v$ is the operator on the space of $\Pi'_v$ restricted from the operator 
$\rho(\E) : \phi \mapsto \E\phi$ on $L^2(\mb{G}(F)\bs\mb{G}(\Af))$.
Put $\la \Pi', f\ra_{\E, S}:= \prod_{v\in S} \la \Pi'_v, f_v\ra_{\E_v}$.

Let $f = \otimes_v f_v$ be a function in $\ELw{\mb{G}(\Af), 1}$ such that $f_v$ is spherical for all $v\notin S$.
For $i = 1, 2$, let $f_i = \otimes_v f_{i,v}$ be a function in $\ELw{\mb{H}_i(\Af), 1}$ which matches $f$
and has spherical local components at all $v\notin S$.
The following $\E$-trace identity follows from Proposition \ref{prop:globalbc}:
\begin{equation}\label{eq:pretrivialbc}
\sum_{\Pi' \in \{\Pi\}} m(\Pi')\la \Pi', f\ra_{\E, S} =
\frac{1}{2}\la \pi_1, f_1\ra_S + \frac{1}{2}\la \pi_2, f_2\ra_S.
\end{equation}
Here, $m(\Pi')$ is the multiplicity of $\Pi'$ in the discrete spectrum of $\mb{G}$.
From \cite[Sect. V. 10.6]{F1}, we have:
\[
m(\Pi') = \frac{1}{2}\lp 1 + (-1)^{n\lp\Pi'\rp}\rp,
\]
where $n\lp\Pi'\rp$ is the number of places  $v$ for which $\Pi'_v = \Pi_v^-$.  In other words,
$\Pi' \in \{\Pi\}$ appears (with multiplicity one) in the discrete spectrum if and only if 
$\Pi'_v = \Pi_v^-$ for an even number of places $v$.

Let $P^+(S)$ denote the set of representations $\Pi'$ in $\{\Pi\}$
such that $\Pi'_v$ is unramified for all $v \notin S$, and $\Pi'_v = \Pi_v^-$ for an {\it even} number of places $v$ in $S$.

Let $P^-(S)$ denote the set of $\Pi'$ in $\{\Pi\}$ such that $\Pi_v'$ is unramified for all $v\notin S$, and
$\Pi_v' = \Pi_v^-$ for an {\it odd} number of places $v$ in $S$.

By the multiplicity formula for $\{\Pi\}$, the equation \eqref{eq:pretrivialbc} is equivalent to
\begin{equation}\label{eq:trivialbc}
2\!\!\sum_{\Pi' \in P^+(S)}\la \Pi', f\ra_{\E, S} = \la \pi_1, f_1\ra_S + \la \pi_2, f_2\ra_S.
\end{equation}
\begin{claim}\label{claim:beforebc1dimproof}
The following identities hold for matching functions:
\[
\begin{split}
2\!\!\sum_{\Pi' \in P^+(S_0)}\la \Pi', f\ra_{\E, S_0} &= \la \pi_1, f_1\ra_{S_0} + \la \pi_2, f_2\ra_{S_0},\\
2\!\!\sum_{\Pi' \in P^-(S_0)}\la \Pi', f\ra_{\E, S_0} &= \la \pi_1, f_1\ra_{S_0} - \la \pi_2, f_2\ra_{S_0}.
\end{split}
\]
\end{claim}
\begin{proof}
Recall that we have fixed a place $u$ which splits in $E$.
The identity \eqref{eq:trivialbc}, when applied to $S = S_0 \cup \{u\}$, is equivalent to
\begin{multline}\label{eq:trivialbcS0u}
2\!\!\sum_{\Pi' \in P^+(S_0)}\la \Pi', f\ra_{\E, S_0}\la \Pi^\times_u, f_u\ra +
2\!\!\sum_{\Pi' \in P^-(S_0)}\la \Pi', f\ra_{\E, S_0}\la \Pi^-_u, f_u\ra\\
= \la \pi_1, f_1\ra_{S_0}\la \pi_{1,u}, f_{1,u}\ra
+ \la \pi_2, f_2\ra_{S_0}\la \pi_{2,u}, f_{2,u}\ra.
\end{multline}
By \eqref{eq:leverage}, the right-hand side of the above equation is equal to
\[
\left(\la \pi_1, f_1\ra_{S_0} + \la \pi_2, f_2\ra_{S_0}\right)\la \Pi_u^\times, f_u\ra +
\left(\la \pi_1, f_1\ra_{S_0} - \la \pi_2, f_2\ra_{S_0}\right)\la \Pi_u^-, f_u\ra.
\]
Since the representations $\Pi^\times_u$, $\Pi^-_u$ are inequivalent, the claim follows from
the linear independence of characters.
\end{proof}
\begin{proof}[Proof of Proposition \ref{prop:bc1dim}]
Recall that we have fixed a place $w$ of $F$ such that the local components of the
global objects (fields, groups, representations, \dots, etc.) at $w$ coincide with the desired local objects over $k$.
Applying \eqref{eq:trivialbc} to $S = S_0 \cup \{w\}$, we obtain the following identity for matching functions
whose local components at all $v \notin S$ are spherical:
\begin{multline}\label{eq:trivialbcS0w}
2\!\!\sum_{\Pi' \in P^+(S_0)}\la \Pi', f\ra_{\E, S_0}\la \Pi^\times_w, f_w\ra_{\E_w} +
2\!\!\sum_{\Pi' \in P^-(S_0)}\la \Pi', f\ra_{\E, S_0}\la \Pi^-_w, f_w\ra_{\E_w}\\
= \la \pi_1, f_1\ra_{S_0}\la \pi_{1,w}, f_{1,w}\ra
+ \la \pi_2, f_2\ra_{S_0}\la \pi_{2,w}, f_{2,w}\ra.
\end{multline}
By Claim \ref{claim:beforebc1dimproof}, the left-hand side of the above equation is equal to
\[
\la \pi_1, f_1\ra_{S_0} \la \Pi_w^\times + \Pi_w^-, f_w\ra_{\E_w} + 
\la \pi_2, f_2\ra_{S_0} \la \Pi_w^\times - \Pi_w^ -, f_w\ra_{\E_w}.
\]
Rearranging terms, \eqref{eq:trivialbcS0w} is equivalent to
\begin{multline}\label{eq:trivialbcseesaw}
\la \pi_1, f_1\ra_{S_0}\left[ \la \Pi_w^\times + \Pi_w^-, f_w\ra_{\E_w} - \la \pi_{1,w}, f_{1,w}\ra\right] \\=
- \la \pi_2, f_2\ra_{S_0}\left[ \la \Pi_w^\times - \Pi_w^-, f_w\ra_{\E_w} - \la \pi_{2,w}, f_{2,w}\ra\right].
\end{multline}

Pick a place $v \in S_0$ which does not split in $E$.  
Recall that an elliptic regular element in $G_v$ is of type $i$ ($i = 1, 2$) if its norms lie only in $\mb{H}_i(F_v)$.
\indexi{element!type}%
Choose $f_v$ to be a function on $G_v$
whose $\E_v$-twisted orbital integral is supported only on the set of elliptic regular elements of type $1$.
We may then set $f_2 = 0$ without violating the matching condition on the test functions.
For such choices of test functions, \eqref{eq:trivialbcseesaw} becomes
\[
\la \pi_1, f_1\ra_{S_0}\left[ \la \Pi_w^\times + \Pi_w^-, f_w\ra_{\E_w} - \la \pi_{1,w}, f_{1,w}\ra\right] = 0.
\]
Since there exists $f_1 \in \ELw{\mb{H_1}(\Af), 1}$ such that $\la \pi_1, f_1\ra_{S_0} \neq 0$,
we obtain:
\begin{equation}\label{eq:bc1dimpi1}
\displaystyle \la \Pi_w^\times, f_w\ra_{\E_w} + \la \Pi_w^-, f_w\ra_{\E_w} = \la \pi_{1,w}, f_{1,w}\ra.
\end{equation}

We may likewise pick a place $v$ in $S_0$, and choose the function $f$, such that
the twisted orbital integral of $f_v$ is supported only on the elliptic regular set of type 2.  
It then follows that
\begin{equation}\label{eq:bc1dimpi2}
\displaystyle \la \Pi_w^\times, f_w\ra_{\E_w} - \la \Pi_w^-, f_w\ra_{\E_w} = \la \pi_{2,w}, f_{2,w}\ra.
\end{equation}

Let $f$, $f_1$, $f_2$ now be matching test functions on $G$, $H_1$, $H_2$, respectively.
For $* = \times, -$, we have by construction $\Pi^*_w = \pi^*$.
Let $A^* \in \Hom_G(\pi^*, \ve\pi^*)$ be the intertwining operator
defined by $\rho(\E)_w$.  By \eqref{eq:bc1dimpi1} and \eqref{eq:bc1dimpi2}, we have:
\[
\begin{split}
\la \pi^\times, f\ra_{A^\times} + \la \pi^-, f\ra_{A^-} &= \la \xi 1_{2, k}\otimes_1 1, f_1 \ra,\\
\la \pi^\times, f\ra_{A^\times} - \la \pi^-, f\ra_{A^-} &= \la (\xi\circ\N_{\K/k})1_{2, \K}\otimes_2 1, f_2\ra.
\end{split}
\]

If one of $A^\times, A^-$ is zero, then
\[
\la \xi 1_{2, k}\otimes_1 1, f_1 \ra = \pm \la (\xi\circ\N_{\K/k})1_{2, \K}\otimes_2 1, f_2\ra.
\]
Let $f$ be an arbitrary test function on $G$ whose $\ve$-twisted orbital integral is 
supported only on the elliptic regular set of type $1$.  Let $f_2 = 0$.  Then,
$\la \xi 1_{2, k}\otimes_1 1, f_1 \ra$ is zero for each elliptic function $f_1$ on $H_1$ which matches $f$, 
which is a contradiction because $1_{2, k}\otimes_1 1$ is supported on elliptic elements.
Thus, neither $A^\times$ nor $A^-$ is trivial.  The proposition follows.
\end{proof}
{\sc Remark.}
Our proof of Proposition \ref{prop:bc1dim} relies on the multiplicity formula deduced in \cite{F1} for the global quasi-packet
\indexi{multiplicity formula}%
$[\Z 1_2, \E\Z 1_2]$.  Said multiplicity formula is proven only for the case of trivial central character, hence the assumption
that the character $\xi$ in the hypothesis of the proposition satisfies $\xi^2 = 1$.  If the same multiplicity formula
holds in general, then it is clear that we may remove the condition $\xi^2 = 1$.
\subsubsection{Steinberg Representations}\label{sec:steinberg}
Suppose the character $\ve$ of $k^\times$ is unramified, i.e. $\K/k$ is an unramified field extension.

Let $\nu, \nu_{\K}$ be the normalized absolute value functions on $k, \K$, respectively.
In particular, since $\K/k$ is unramified, $\nu_\K = \nu \circ\N_{\K/k}$.
Recall that ${\rm St}_{2,l}$ ($l = k, \K$) denotes the Steinberg square integrable representation of $\GL(2, l)$.

Let $\xi$ be a character of $k^\times$ with $\xi^2 = 1$.
By \cite[Lemma 3.6]{ST},
the induced representation $I = \ve\nu\times\ve\rtimes\nu^{-1/2}\xi$ of $G$ is reducible, with $4$ distinct subquotients:
\[\begin{split}
I =
\delta(\nu^{1/2}\ve{\rm St}_{2, k}, \nu^{-1/2}\xi)
& + L(\nu^{1/2}\ve{\rm St}_{2, k}, \nu^{-1/2}\xi)\\
& + L(\nu^{1/2}\ve{\rm St}_{2, k}, \nu^{-1/2}\ve\xi) + L(\nu\ve, \ve\rtimes\nu^{1/2}\xi).
\end{split}\]
Here, $\delta(\nu^{1/2}\ve{\rm St}_{2, k}, \nu^{-1/2}\xi)$
is the unique square integrable subrepresentation of $I$.
For $\zeta = \xi, \ve\xi$,
the representation $L(\nu^{1/2}\ve{\rm St}_{\GL(2, k)}, \nu^{-1/2}\zeta)$ is the unique nontempered
quotient of $\nu^{1/2}\ve{\rm St}_{2, k}\rtimes \nu^{-1/2}\zeta$.  The representation
$L(\nu\ve, \ve\rtimes\nu^{-1/2}\xi)$ is a nontempered quotient of $\nu^{1/2}\ve 1_{2,k}\rtimes \nu^{-1/2}\xi$.

Since the \rep $L(\nu^{1/2}\ve{\rm St}_{2, k}, \nu^{-1/2}\xi)$ is the unique nontempered quotient of
$\nu^{1/2}\ve{\rm St}_{2, k}\rtimes \nu^{-1/2}\xi$, and $L(\nu^{1/2}\ve{\rm St}_{2, k}, \nu^{-1/2}\ve\xi)$ 
that of 
\[
\nu^{1/2}\ve{\rm St}_{2, k}\rtimes \nu^{-1/2}\ve\xi = \ve\otimes(\nu^{1/2}\ve{\rm St}_{2, k}\rtimes \nu^{-1/2}\xi),
\]
we have $L(\nu^{1/2}\ve{\rm St}_{2, k}, \nu^{-1/2}\ve\xi) \cong \ve L(\nu^{1/2}\ve{\rm St}_{2, k}, \nu^{-1/2}\xi)$.
Since the subquotients of $I$ are distinct, neither $L(\nu^{1/2}\ve{\rm St}_{2, k}, \nu^{-1/2}\xi)$ 
nor $L(\nu^{1/2}\ve{\rm St}_{2, k}, \nu^{-1/2}\ve\xi)$ is $\ve$-invariant.

Recall that the representation $\delta(\nu^{1/2}\ve{\rm St}_{2, k}, \nu^{-1/2}\xi)$ 
belongs to a local packet of size two.   The other member of this local packet is a cuspidal
representation, denoted by $\delta^-(\nu^{1/2}\ve{\rm St}_{2, k}, \nu^{-1/2}\xi)$.

Let $\delta = \delta(\nu^{1/2}\ve{\rm St}_{2, k}, \nu^{-1/2}\xi)$,
$\pi^- = \delta^-(\nu^{1/2}\ve{\rm St}_{2, k}, \nu^{-1/2}\xi)$, and let $\pi^\times$ denote \dt $L(\nu\ve, \ve\rtimes\nu^{-1/2}\xi)$.
\begin{lemma}\label{lemma:trbcspecial}
There exist nontrivial intertwining operators $A^\delta \in \Hom_G(\delta, \ve\delta)$, $A^- \in \Hom_G(\pi^-,\ve\pi^-)$,
such that the following twisted character identities hold for matching functions:
\begin{equation}\label{eq:trbcspecial}
\begin{split}
\la \delta, f\ra_{A^\delta} - \la \pi^-, f\ra_{A^-} &= \la \xi\;{\rm St}_{2, k} \otimes_1 1, f_1 \ra,\\
\la \delta, f\ra_{A^\delta} + \la \pi^-, f\ra_{A^-} &= \la (\xi \circ\N_{\K/k}){\rm St}_{2, \K} \otimes_2  1, f_2 \ra.
\end{split}
\end{equation}
\end{lemma}
\begin{proof}
The representation $I = \ve\nu\times\ve\rtimes\nu^{-1/2}\xi$ is equal to $\ve\nu\rtimes I(\nu^{-1/2}\xi, \ve\nu^{-1/2}\xi)$,
and $I(\nu^{-1/2}\xi, \ve\nu^{-1/2}\xi)$ is equal to the monomial representation $\pi((\nu^{-1/2}\xi)\circ\N_{\K/k})$
associated with the character $(\nu^{-1/2}\xi)\circ\N_{\K/k}$ of $\K^\times$.  Thus,
\[
I = \ve\nu\rtimes \pi((\nu^{-1/2}\xi)\circ\N_{\K/k}).
\] 

Let $A \in \Hom_G(I, \ve I)$ be the operator defined in Section \ref{sec:intop}.
It intertwines each subquotient $L$ of $I$ to a representation which is equivalent to $\ve L$.
The subquotients \dt
\[
L(\nu^{1/2}\ve{\rm St}_{2, k}, \nu^{-1/2}\xi),\quad L(\nu^{1/2}\ve{\rm St}_{2, k}, \nu^{-1/2}\ve\xi)
\]
are obtained from each other by tensoring with $\ve$.  Since they are inequivalent,
$A$ swaps their spaces, which implies that
\[
\tr L(\nu^{1/2}\ve{\rm St}_{2, k}, \nu^{-1/2}\xi)(f)A =
\tr L(\nu^{1/2}\ve{\rm St}_{2, k}, \nu^{-1/2}\ve\xi)(f)A = 0.
\]

Since $\delta$ is the unique square integrable subrepresentation of 
$I$, the operator $A$ maps the space of $\delta$ to itself.
We have determined the image under $A$ of three of the four subquotients of $I$; hence,
$A$ must map the space of the unique irreducible quotient $\pi^\times$ of $I$ to itself.
Let $A^\delta$ denote the restriction of $A$ to $\delta$.  Let ${A}'$ denote the operator on 
$\pi^\times$ induced from $A$.  Then, $A^\delta$, ${A}'$ are nontrivial operators in $\Hom_G(\delta, \ve\delta)$, 
$\Hom_G(\pi^\times, \ve\pi^\times)$, respectively.  Moreover, we have:
\begin{equation}\label{eq:AtwistedI}
\la I, f \ra_A
= \la \delta, f\ra_{A'} + \la \pi^\times, f \ra_{A^\delta}.
\end{equation}

For the $\ve$-endoscopic groups, we have:
\begin{equation}\label{eq:steinbergH1H2}
\begin{split}
\xi I(\nu^{1/2}, \nu^{-1/2})\otimes_1 1 &= \xi 1_{2, k} \otimes_1 1 + \xi\;{\rm St}_{2, k} \otimes_1 1,\\
(\xi\circ\N_{\K/k})I(\nu_\K^{1/2}, \nu_\K^{-1/2})\otimes_2 1 &= (\xi\circ\N_{\K/k}) 1_{2, \K} \otimes_2 1 
+ (\xi\circ\N_{\K/k}){\rm St}_{2, \K} \otimes_2 1.
\end{split}
\end{equation}

By Corollary \ref{corollary:trGH1H2}, the following characters are equal to one another:
\begin{itemize}
\item
$\la \ve\nu\times\ve\rtimes\nu^{-1/2}\xi, f\ra_{A}$,
\item
$\la \xi I(\nu^{1/2},\nu^{-1/2})\otimes_1 1, f_1 \ra$,
\item
$\la B_{\K/k}\lp\xi I(\nu^{1/2}, \nu^{-1/2})\rp\otimes_2 1, f_2 \ra 
= \la (\xi \circ \N_{\K/k})I(\nu_\K^{1/2},\nu_\K^{-1/2})\otimes_2 1, f_2 \ra$.
\end{itemize}

By Proposition \ref{prop:bc1dim}, there exist nontrivial intertwining operators $A^\times, A^-$ such that
\begin{equation}\label{eq:localbc1dimrecap}
\begin{split}
\la \pi^\times, f\ra_{A^\times} + \la \pi^-, f\ra_{A^-} &= \la \xi 1_{2, k} \otimes_1 1, f_1 \ra,\\
\la \pi^\times, f \ra_{A^\times} - \la \pi^-, f\ra_{A^-} &= \la (\xi\circ\N_{\K/k})1_{2, \K}\otimes_2 1, f_2 \ra.
\end{split}
\end{equation}

Recall that $A^\times$ is restricted from a global operator on 
$L^2(\mb{G}(F)\bs\mb{G}(\Af))$ (see proof of Proposition \ref{prop:bc1dim}).  
The representation $\pi^\times$ is irreducible.  
Since both $A'$ and $A^\times$ have the property that their squares are the identity,  
by Schur's lemma we have $A^\times = d \cdot A'$, where $d = \pm 1$.
It follows from \eqref{eq:steinbergH1H2}, \eqref{eq:localbc1dimrecap}, and \eqref{eq:AtwistedI} that
\[
\la \xi\;{\rm St}_{2, k} \otimes_1 1, f_1 \ra
= -d \la \pi^\times, f\ra_{A'} - \la\pi^-, f\ra_{A^-} + \la \delta, f\ra_{A^\delta} + \la \pi^\times, f\ra_{A'}.
\]
The representation $\pi^\times$ is nontempered, 
while the rest of the representations in the above equation are square integrable.
The central exponents of square integrable representations all decay (\cite[Prop. 29]{Be}).
Invoking the linear independence of central exponents (see \cite[Sect. 21]{FK}), we conclude that
$-d \la \pi^\times, f\ra_{A'} + \la \pi^\times, f\ra_{A'} = 0$.  Hence, $d = 1$, and
\[
 \la \xi\;{\rm St}_{2, k} \otimes_1 1, f_1 \ra\\
= \la \delta, f\ra_{A^\delta} - \la\pi^-, f\ra_{A^-}.
\]

Applying the same argument to $\la (\xi\circ\N_{\K/k}){\rm St}_{2, \K}\otimes_2 1, f_2\ra$, we obtain:
\[
\la (\xi\circ\N_{\K/k}){\rm St}_{2, \K} \otimes_2  1, f_2 \ra = \la \delta, f\ra_{A^\delta} + \la \pi^-, f\ra_{A^-}.
\]
The lemma follows.
\end{proof}
\subsection{Cuspidal Representations}
\indexi{representation!cuspidal|(}%
In this section, the character $\ve$ is not assumed to be unramified.

In what follows, we prove local character identities by contructing global objects whose
local components at a chosen place coincide with specific local conditions.
For instance, we contruct automorphic representations whose components at a chosen finite place are
the $p$-adic representations we are interested in.   To do so,
we make use of the following generalization of \cite[Prop. III. 3]{F2}:
\begin{prop}\label{prop:construct}
Let $F$ be a number field.  Let $\mb{H}$ be a reductive connected $F$-group.
Fix a nonarchimedean place $w$ of $F$.  Let $\tau_w$ be a cuspidal $H_w$-module.
Let $\{v_i\}_{i \in I}$ be a finite set of nonarchimedean places.
For each $i \in I$, let $\tau_{v_i}$ be a square integrable $H_{v_i}$-module.
Let $S$ be a finite set of places which contains the union of
$w$, $\{v_i\}_{i \in I}$, and all the archimedean places. Then, there exists an automorphic
representation $\pi$ of $\mb{H}(\Af)$ such that: {\rm (i)} $\pi_w = \tau_w$; {\rm (ii)} $\pi_{v_i} = \tau_{v_i}$
for all $i \in I$; {\rm (iii)} $\pi_v$ is unramified for all $v \notin S$.
\end{prop}  

Let $\tau$ be an irreducible, cuspidal, non-$\K$-monomial representation of $\GL(2, k)$ with trivial central character $\omega_\tau$.  
Let $\{\pi^+, \pi^-\}$ be the local packet of $\PGSp(2, k)$
which lifts to the induced representation $I_{(2, 2)}(\tau, \ve\tau)$ of $\PGL(4, k)$, according to
\cite[Prop. V. 5]{F1}.
\begin{prop}\label{prop:bccuspidal}
There exist nontrivial operators $A^+$ in $\Hom_G(\pi^+, \ve\pi^+)$, $A^-$ in \dt $\Hom_G(\pi^-, \ve\pi^-)$,
such that the following twisted character identities hold for matching functions $f \in \C(G)$, $f_i \in \C(H_i)$ $(i = 1, 2)\!:$
\[
\begin{split}
\la \pi^+, f\ra_{A^+} &= \frac{1}{2} \la \tau\otimes_1 1, f_1 \ra + \frac{1}{2}\la B_{\K/k}\tau\otimes_2 \omega_\tau, f_2\ra,\\
\la \pi^-, f\ra_{A^-} &= \frac{1}{2} \la \tau\otimes_1 1, f_1 \ra - \frac{1}{2}\la B_{\K/k}\tau\otimes_2 \omega_\tau, f_2\ra.
\end{split}
\]
\begin{remark}
Note that by assumption $\omega_\tau$ is equal to $1$.  We nonetheless write $\omega_\tau$ in the above equations, for we
expect the same character identities to hold in the cases where $\omega_\tau \neq 1$.
\end{remark}
\end{prop}

\begin{proof}
Construct a totally imaginary number field $F$ and a quadratic extension
$E$ of $F$ with the following properties:
\begin{itemize}
\item
There exists a finite place $w$ of $F$, prime in $E$, such that $F_w = k$.
\item
$E_w = \K$.  
Here, $w$ denotes both the place of $F$ and the unique place of $E$ which lies above it.
\end{itemize}
Let $V$ be the set of places of $F$.  
Let $\Vr \subset V$ be the set of finite places which are unramified in $E$.
Let $w_1, w_2 \in V$ be two (finite) places, different from $w$, which are unramified, prime in $E$.

Using Proposition \ref{prop:construct}, we construct a cuspidal automorphic representation $\mct$ of \dt $\GL(2, \Af)$ 
with trivial central character such that:
\begin{itemize}
\item
$\mct_w = \tau$;
\item
$\mct_{w_i}$ ($i = 1, 2$) is the Steinberg representation ${\rm St}_{2, w_i}$ of $\GL(2, F_{w_i})$;
\item
$\mct_v$ is unramified for each finite place $v \notin \{w, w_1, w_2\}$.
\end{itemize}
In particular, since $\tau$ is non-$\K$-monomial, $\mct$ is non-$E$-monomial.

Let $\E$ be the quadratic character of $\idc{F}$ which corresponds to the extension $E/F$.
In particular, $\E_w = \ve$.
View $\mb{G} = \GSp(2)$ as an $F$-group, and $\mb{H_1}$, $\mb{H_2}$ as the elliptic $\E$-endoscopic groups of $\mb{G}$
over $F$.

Let $\{\Pi\} = \otimes_{v} \{\Pi_v\}$ be the unstable packet $[\mct, \E\mct]$ of $\mb{G}(\Af)$ (see Section \ref{sec:arthurunstable}).
It lifts to the induced representation $I_{(2, 2)}(\mct, \E\mct)$ of $\GL(4, \Af)$.
For each place $v$ where $\mct_v$ is square integrable,
the local packet $\{\Pi_v\}$ consists of two (square integrable) representations $\Pi^+_v$, $\Pi^-_v$.
In particular, $\Pi_w^+ = \pi^+$ and $\Pi_w^- = \pi^-$.

At a place $v$ where $\mct_v$ is fully induced (for instance, where $\mct_v$ is unramified), 
the packet $\{\Pi_v\}$ consists of a single representation, denoted by  $\Pi_v^+$.  
For convenience, we put $\Pi_v^- := 0$, and write $\{\Pi_v\} = \{\Pi_v^+, \Pi_v^-\}$.

At the place $w_i$ ($i = 1, 2$), where $\mct_{w_i} = {\rm St}_{2, w_i}$, the local packet
$\{\Pi_{w_i}\}$ consists of the square integrable representation 
$\Pi_{w_i}^+ = \delta(\ve_{w_i}\nu_{w_i}^{1/2}{\rm St}_{2, w_i}, \nu_{w_i}^{-1/2})$ and
the cuspidal $\Pi_{w_i}^- = \delta^-(\ve_{w_i}\nu_{w_i}^{1/2}{\rm St}_{2, w_i}, \nu_{w_i}^{-1/2})$.

By Corollary \ref{corollary:veinvariantlocal}, each representation $\Pi'$ in $\{\Pi\}$ is $\E$-invariant.
Hence, as in the case of one dimensional representations (Section \ref{sec:bc1dim}), the twisted character
$\la \Pi', f\ra_\E$ is equal to a product $\prod_v \la \Pi'_v, f_v\ra_{\E_v}$ of twisted local characters.

Let $\pi_1$ be the automorphic representation $\mct\otimes_1 1$ of $\mb{H_1}(\Af)$.
Let $\pi_2$ be the automorphic representation $B_{E/F}\mct\otimes_2 1$ of $\mb{H_2}(\Af)$.


Let $S_0 = V - \Vr$.  Note that $w$ may or may not lie in $S_0$. 
Let $S_1 = S_0 \cup \{w\}$.  Let $S = S_1\cup\{w_1, w_2\}$.
By Corollary \ref{corollary:globalbc}, the following $\E$-trace identity holds for matching functions
$f$ in $\ELw{\mb{G}(\Af), 1}$ and $f_i$ in $\ELw{\mb{H}_i(\Af), 1}$ ($i = 1, 2$) whose local components at all
$v \notin S$ are spherical:
\begin{equation}\label{eq:precuspidalbc}
\sum_{\Pi' \in \{\Pi\}} m(\Pi')\la \Pi', f\ra_{\E, S} =
\frac{1}{2}\la \pi_1, f_1\ra_S + \frac{1}{2}\la \pi_2, f_2\ra_S.
\end{equation}
Here, $m(\Pi')$ is the multiplicity of $\Pi'$ in the discrete spectrum of $\mb{G}$.
From \cite[Sect. V. 10.3]{F1}, we have:
\[ 
m(\Pi') = \frac{1}{2}\lp 1 + (-1)^{n\lp\Pi'\rp}\rp,
\]
where $n\lp\Pi'\rp$ is the number of places $v$ for which $\Pi'_v = \Pi_v^-$.
In other words,
$\Pi' \in \{\Pi\}$ appears (with multiplicity one) in the discrete spectrum if and only if 
$\Pi'_v = \Pi_v^-$ for an even number of places $v$.

Let $P^+(S)$ denote the set of representations $\Pi'$ in $\{\Pi\}$
such that $\Pi'_v$ is unramified for all $v \notin S$, and $\Pi'_v = \Pi_v^-$ for an {\it even} number of places $v$ in $S$.

Let $P^-(S)$ denote the set of $\Pi'$ in $\{\Pi\}$ such that $\Pi_v'$ is unramified for all $v\notin S$, and
$\Pi_v' = \Pi_v^-$ for an {\it odd} number of places $v$ in $S$.
By the multiplicity formula, the equation \eqref{eq:precuspidalbc} is equivalent to
\begin{equation}\label{eq:cuspidalbc}
2\!\!\sum_{\Pi' \in P^+(S)}\la \Pi', f\ra_{\E, S} = \la \pi_1, f_1\ra_S + \la \pi_2, f_2\ra_S.
\end{equation}

If a place $v$ of $F$ is archimedean, then $\GL(2, F_v)$ is equal to $\GL(2, \mbb{C})$ because $F$ is totally imaginary.
Hence, $\mct_v$ is a fully induced representation.  

If $v$ is a finite place not in $S$, then by construction $\mct_v$ is unramified.  Since $\mct$ is cuspidal, $\mct_v$ is fully induced.

Hence, for all $v \in S_1 - \{w\}$, $\mct_v$ is fully induced.
Suppose $\mct_v = I(\alpha_v, \beta_v)$ for some characters $\alpha_v, \beta_v$ of $F_v^\times$. 
Then,
\[
\begin{array}{|c|c|c|}
\hline &v \text{ prime in }E& v \text{ splits in } E\\
\hline
\pi_{1,v} =& I(\alpha_v, \beta_v)\otimes_1 1 
& I(\alpha_v, \beta_v)\otimes_1 (1\otimes 1)\\
\hline
\pi_{2,v} =& B_{E/F}I(\alpha_v, \beta_v) \otimes_2 (\alpha_v\beta_v)\circ\N_{E_v/F_v} & 
\lp I(\alpha_v, \beta_v)\times I(\alpha_v, \beta_v)\rp \otimes_2 \alpha_v\beta_v\\
\hline
\end{array}
\]
Let $\Pi_v''$ be the induced representation $\ve_v{\alpha_v}/{\beta_v}\rtimes I(\beta_v, \beta_v\ve_v)$ of $G_v$.
By Propositions \dts \ref{prop:trGH1}, \ref{prop:trGH2}, and Lemmas \ref{lemma:indsplitH1}, \ref{lemma:indsplitH2},
$\la \pi_{1,v}, f_{1,v}\ra$ and $\la \pi_{2,v}, f_{2,v}\ra$ are equal to $\la \Pi_v'', f_v\ra_{A_v}$
Here, $A_v$ is the intertwining
operator defined in Section \ref{sec:intop} if $v$ is prime in $E$; otherwise, $A_v = 1$.

For each $v \in S_1 - \{w\}$, regard the local components of the test functions at all places
$u \neq v$ as fixed. 
By \eqref{eq:cuspidalbc} and the linear independence of characters, we conclude that
$\{\Pi_v\}$ is the singleton consisting of $\Pi_v^{''}$.
Hence, the $(S_1 - \{w\})$-components of both sides of \eqref{eq:cuspidalbc} cancel one another.  

Let $S_3 = S - (S_1 - \{w\}) = \{w, w_1, w_2\}$.  After cancellation, \eqref{eq:cuspidalbc} becomes
\begin{equation}\label{eq:cuspidalbcS3}
2\!\!\sum_{\Pi' \in P^+(S_3)}\ep(\Pi')\la \Pi', f\ra_{\E, S_3}
= \prod_{v \in S_3}\la \pi_{1, v}, f_{1, v}\ra + \prod_{v \in S_3}\la \pi_{2, v}, f_{2, v}\ra.
\end{equation}
Here, $\ep(\Pi') = \pm 1$.  It appears in the equation because $\la \Pi_v'', f_v\ra_{A_v} = \pm \la \Pi_v'', f_v\ra_{\E_v}$
for $v \in S_1 - \{w\}$.
Absorbing the signs $\ep(\Pi')$ into the local intertwining operators $\rho(\E)_v$ at $v \in S_3$ if necessary, we assume
that $\ep(\Pi') = 1$ for all $\Pi'$.

By construction, $\pi_{1,w_i} = {\rm St}_{2, w_i}\otimes_1 1$ and $\pi_{2,w_i} = B_{E/F}{\rm St}_{2, w_i}\otimes_2 1$
for $i = 1, 2$.
Hence, by Lemma \ref{lemma:trbcspecial}, there exist constants $\ep^+_i, \ep^-_i = \pm 1$ such that
\[
\begin{split}
\ep^+_i \la \Pi^+_{w_i}, f_{w_i}\ra_{\E_{w_i}} - \ep^-_i\la \Pi^-_{w_i}, f_{w_i}\ra_{\E_{w_i}} &= \la \pi_{1, w_i}, f_{1, w_i}\ra,\\
\ep^+_i \la \Pi^+_{w_i}, f_{w_i}\ra_{\E_{w_i}} + \ep^-_i\la \Pi^-_{w_i}, f_{w_i}\ra_{\E_{w_i}} &= \la \pi_{2, w_i}, f_{2, w_i} \ra.
\end{split}
\]
The equation \eqref{eq:cuspidalbcS3} is therefore equivalent to
\begin{multline*}
\lp 2\!\!\sum_{\Pi' \in P^+(\{w, w_1\})}\!\!\la \Pi', f\ra_{\E, \{w, w_1\}}\rp 
\la \Pi^+_{w_2}, f_{w_2}\ra_{\E_{w_2}}\\
\shoveleft{
+ \lp 2\!\!\sum_{\Pi' \in P^-(\{w, w_1\})}\!\!\la \Pi', f\ra_{\E, \{w, w_1\}}\rp
\la \Pi^-_{w_2}, f_{w_2}\ra_{\E_{w_2}}
}\\
= \ep^+_2 \lp\prod_{v \in \{w, w_1\}}\!\!\la \pi_{1, v}, f_{1, v}\ra + \!\prod_{v \in \{w, w_1\}}\!\!\la \pi_{2, v}, f_{2, v}\ra\rp
\la \Pi^+_{w_2}, f_{w_2}\ra_{\E_{w_2}}\\
+ \ep^-_2 \lp \prod_{v \in \{w, w_1\}}\!\!\la \pi_{2, v}, f_{2, v}\ra - \prod_{v \in \{w, w_1\}}\!\!\la \pi_{1, v}, f_{1, v}\ra\rp
\la \Pi^-_{w_2}, f_{w_2}\ra_{\E_{w_2}}.
\end{multline*}
By the linear independence of characters, we have:
\begin{multline}\label{eq:cuspidalbc+}
2\!\!\sum_{\Pi' \in P^+(\{w, w_1\})}\!\!\la \Pi', f\ra_{\E, \{w, w_1\}}\\
= \ep^+_2 \lp\prod_{v \in \{w, w_1\}}\!\!\la \pi_{1, v}, f_{1, v}\ra + \!\prod_{v \in \{w, w_1\}}\!\!\la \pi_{2, v}, f_{2, v}\ra\rp,
\end{multline}
\begin{multline}\label{eq:cuspidalbc-}
2\!\!\sum_{\Pi' \in P^-(\{w, w_1\})}\!\!\la \Pi', f\ra_{\E, \{w, w_1\}}\\
= \ep^-_2 \lp \prod_{v \in \{w, w_1\}}\!\!\la \pi_{2, v}, f_{2, v}\ra - \prod_{v \in \{w, w_1\}}\!\!\la \pi_{1, v}, f_{1, v}\ra\rp.
\end{multline}

Equation \eqref{eq:cuspidalbc+} has the same form as \eqref{eq:cuspidalbcS3}, with $w_2$ removed from the set $S_3 = \{w, w_1, w_2\}$.
By applying the same argument inductively on the set $\{w_1, w_2\}$, we conclude that
\begin{equation*}
2\!\!\sum_{\Pi' \in P^+(\{w\})}\!\!\la \Pi', f\ra_{\E, \{w\}}
= \ep^+_1 \ep^+_2  \lp\la \pi_{1, w}, f_{1, w}\ra + \la \pi_{2, w}, f_{2, w}\ra\rp.
\end{equation*}
Hence, there exists a constant $\ep^+ = \pm 1$ such that
\begin{equation}\label{eq:cuspidalbcw}
2 \la \Pi^+_w, f_w\ra_{\E_w} = 
2\!\!\sum_{\Pi' \in P^+(\{w\})}\!\!\la \Pi', f\ra_{\E, \{w\}} =
\ep^+ \la \pi_{1, w}, f_{1, w}\ra +  \ep^+\la \pi_{2, w}, f_{2, w}\ra.
\end{equation}

Similarly, \eqref{eq:cuspidalbc-} implies that there exists $\ep_- = \pm 1$ such that
\begin{equation}\label{eq:cuspidalbcw-}
2 \la \Pi^-_w, f_w\ra_{\E_w} = 
2\!\!\sum_{\Pi' \in P^-(\{w\})}\!\!\la \Pi', f\ra_{\E, \{w\}} =
\ep^-\la \pi_{1, w}, f_{1, w}\ra -  \ep^-\la \pi_{2, w}, f_{2, w}\ra.
\end{equation}

For $* = +$, $-$,
let $\rho(\E)_w^*$ denote the operator in $\Hom_{G_w}(\Pi^*_w,\; \E_w \Pi^*_w)$ restricted from $\rho(\E)$.
Let $A^* = \ep^* \cdot \rho(\E)_w^*$.  The proposition then follows from \eqref{eq:cuspidalbcw} and
\eqref{eq:cuspidalbcw-}.
\end{proof}
{\sc Remark.}  As in Proposition \ref{prop:bc1dim}, we may extend Proposition \ref{prop:bccuspidal}
to the case of nontrivial central characters if the multiplicity formula for the global packet $\{\Pi\}$ holds in general.
\subsubsection{Monomial Representations}\label{sec:localmonomial}
\indexi{representation!monomial|(}%
Let $\chi$ be a character of $\K^\times$.
The representation $\pi = 1 \rtimes \pi(\chi)$ of $G$ is reducible with two inequivalent tempered constituents
$\pi^+$, $\pi^-$.  
\begin{lemma}\label{lemma:localmonomialinduced}
There exist intertwining operators $A^* \in \Hom_G(\pi^*, \ve\pi^*)$ $(*= +, -)$ such that the following holds
for matching functions $f$ in $C(G)$ and $f_1$ in $C(H_1)\!:$
\[
\la \pi^+, f\ra_{A^+} + \la \pi^-, f\ra_{A^-} = \la \pi(\chi)\otimes_1 1, f_1\ra.
\]
\end{lemma}
\begin{proof}
Construct a totally imaginary number field $F$ and a quadratic extension $E$ of $F$ such that:
(i) There is a place $w$ of $F$ such that $F_w = k$; (ii) $w$ is prime in $E$, $E_w = \K$.

Let $w_1$, $w_2$ be two (finite) places, different from $w$, which are prime in $E$.

Let $\E$ be the quadratic character of $\idc{E}$ corresponding to the extension $E/F$.
We identify the generator $\sigma$ of $\Gal(\K/k)$ with that of $\Gal(E/F)$.

Construct a character $\X$ of $\idc{E}$ such that: (i) $\X_{w_i}$ ($i = 1, 2$) is not fixed by the action of $\sigma$;
(ii) $\X_v$ is unramified for each finite place $v \notin \{w_1, w_2\}$.

Let $\pi(\X)$ be the cuspidal $E$-monomial representation of $\GL(2, \Af)$ associated with $\X$.
At each finite place $v$ where $\X_v$ is unramified, $\pi(\X)_v = \pi(\X_v)$ is fully induced.  
Since $F$ is totally imaginary,
$\pi(\X)_v$ is a fully induced representation of $\GL(2, \CC)$ for each archimedean place $v$.

Let $\pi_1$ be the automorphic representation $\pi(\X)\otimes_1 1$ of $\mb{H}(\Af)$.
In particular, $\pi_{1, v}$ is parabolically induced for each $v \notin \{w_1, w_2\}$.

Let $\mct$ be the following representation of the Heisenberg parabolic subgroup $P$ of $\mb{G}(\Af)$:
\[
\mct = 1 \otimes \pi(\X) : \lp\lsm a &*\;* &*\\& g_2 &*\\&&\tfrac{\det g_2}{a}\rsm\rp \mapsto \pi(\X)(g_2).
\]
Let $\Pi = 1 \rtimes \pi(\X)$, the representation of $\mb{G}(\Af)$ parabolically induced from $\mct$.
Since $\E\otimes (1 \rtimes \pi(\X)) = 1 \rtimes \E\pi(\X) \cong 1\rtimes\pi(\X)$, the \rep
$\Pi$ is $\E$-invariant.
At each place $v$ where $\pi(\X)_v$ is fully induced, $\Pi_v$ is irreducible.
If $\pi(\X)_v$ is cuspidal, then $\Pi_v$ has two inequivalent tempered constituents $\Pi_v^+$, $\Pi_v^-$.

The operator $I_{P, \mct}(\E)$, defined in Section \ref{sec:finechiexpoverview}, intertwines $\Pi$ with $\E\Pi$.
For each place $v$, let $I_{P, \mct}(\E)_v \in \Hom_{G_v}(\Pi_v, \E_v\Pi_v)$ be the local component of $I_{P, \mct}(\E)$
at $v$.
For an irreducible constituent $\Pi'$ of $\Pi$, and a test function $f$ on $\mb{G}(\Af)$, 
put $\la \Pi', f\ra_{\E} := \tr \Pi'(f)I_{P, \mct}(\E)$.
It is equal to the product $\prod_{v \in V}\la \Pi'_v, f_v\ra_{\E_v}$, where 
$\la \Pi'_v , f_v\ra_{\E_v} := \tr \Pi'_v(f_v)I_{P, \mct}(\E)_v$.

Let $S_0$ be a finite set of places of $F$ which contains the union of $\{w_1, w_2\}$, 
\dt the archimedean places, and the finite places which are ramified in $E$.
Let $f$, $f_1$ be matching functions on $\mb{G}(\Af)$, $\mb{H_1}(\Af)$, with elliptic components 
at $w_1$, $w_2$ and spherical components at all the places outside of $S_0$.
By \eqref{eq:globalmonomialinduced} in Section \ref{sec:somegloballifts}, we have:
\[
\prod_{v \in S_0} \lp\la \Pi_v^+, f_v\ra_{\E_v} - \la \Pi_v^-, f_v\ra_{\E_v}\rp = \prod_{v \in S_0} \la \pi_{1,v}, f_{1, v}\ra.
\]
Since $\pi_{1,v}$ is fully induced for all $v \neq w_1, w_2$, 
by Proposition \ref{prop:trGH1} and Lemma \dts \ref{lemma:indsplitH1} the $(S_0 - \{w_1, w_2\})$-components of the above equation cancel.  
Hence,
\begin{equation}\label{eq:localmonomialinducedw1w2}
\prod_{v \in \{w_1, w_2\}} \lp\la \Pi_v^+, f_v\ra_{\E_v} - \la \Pi_v^-, f_v\ra_{\E_v}\rp 
= \prod_{v \in \{w_1, w_2\}} \la \pi_{1,v}, f_{1, v}\ra.
\end{equation}

Now, construct a character $\X'$ of $\idc{E}$ such that: (i) $\X'_w = \chi$; (ii) $\X'_{w_i} = \X_{w_i}$ for $i = 1, 2$;
(iii) $\pi(\X')_v$ is fully induced for each $v \neq w_1, w_2$.

Let $\Pi = 1\rtimes \pi(\X')$.  Let $\pi_1 = \pi(\X')\otimes_1 1$.
Let $S = S_0 \cup \{w\}$.
By \eqref{eq:globalmonomialinduced} and the usual cancellation, we have:
\[
\prod_{v \in \{w_1, w_2, w\}}\lp \la \Pi_v^+, f_v\ra_{\E_v} - \la \Pi_v^-, f_v\ra_{\E_v}\rp 
= \prod_{v \in \{w_1, w_2, w\}} \la \pi_{1,v}, f_{1, v}\ra.
\]
By \eqref{eq:localmonomialinducedw1w2}, the $\{w_1, w_2\}$-components of the above equation cancel each other; hence,
\[
\la \Pi_w^+, f_w\ra_{\E_w} - \la\Pi_w^-, f_w\ra_{\E_w} = \la \pi_{1,w}, f_{1, w}\ra.
\]

Since $\Pi^+_w$, $\Pi^-_w$ are inequivalent, $I_{P, \mct}(\E)_w$ either swaps their spaces, or it
maps each space to itself.  Since the right-hand side of the above character identity is nonzero in general, we
conclude that $I_{P, \mct}(\E)_w$ restricts to $I_{P, \mct}(\E)_w^+$ in \dt $\Hom_{G_w}(\Pi_w^+,\, \E_w\Pi_w^+)$
and $I_{P, \mct}(\E)_w^-$ in $\Hom_{G_w}(\Pi_w^-,\, \E_w\Pi_w^-)$.

Let $A^+ = I_{P, \mct}(\E)_w^+$.  Let $A^- = -I_{P, \mct}(\E)_w^-$.  The lemma follows.
\end{proof}
Let $\theta, \chi$ be two distinct characters of $\K^\times$ such that 
none of $\theta, \chi, \theta\chi, \theta\,{}^\sigma\!\chi$ is fixed by the generator $\sigma$ of $\Gal(\K/k)$.



Suppose $\theta\chi|_{k^\times} = \theta\,{}^\sigma\!\chi|_{k^\times} = \ve$.  
Then, the central characters of $\pi(\theta\chi)$ and $\pi(\theta\,{}^\sigma\!\chi)$ are trivial.
Let $\{\pi^+, \pi^-\}$ be the local packet
of $G$ which lifts to the induced representation $I_{(2, 2)}(\pi(\theta\chi), \pi(\theta\;{}^\sigma\!\chi))$ of $\PGL(4, k)$
according to \cite[Prop. V. 5]{F1}.
\begin{prop}\label{prop:H1monomial}
There exist operators $A^* \in \Hom_G(\pi^*\!,\, \ve\pi^*)$ $(* = +, -)$, and a constant $\ep = \pm 1$,
such that the following character identities hold for matching functions:
\[
\begin{split}
\la \pi^+, f\ra_{A^+} + \la\pi^-, f \ra_{A^-} &= \la \pi(\chi)\otimes_1\theta, f_1\ra,\\
\la \pi^+, f\ra_{A^+} - \la\pi^-, f \ra_{A^-} &= \ep\la \pi(\theta)\otimes_1\chi, f_1\ra.
\end{split}
\]
\end{prop}
\begin{proof}
Construct a totally imaginary number field $F$ and a quadratic extension $E$ of $F$ such that:
(i) There is a place $w$ of $F$ such that $F_w = k$; (ii) $w$ is prime in $E$, $E_w = \K$.

Let $w_1$, $w_2$ be two (finite) places, different from $w$, which are prime in $E$.

Let $\E$ be the quadratic character of $\idc{E}$ corresponding to the extension $E/F$.
We identify the generator $\sigma$ of $\Gal(\K/k)$ with that of $\Gal(E/F)$.

Construct characters $\X$, $\mco$ of $\idc{E}$ such that: 
(i) $\mco_w = \theta$, $\X_w = \chi$; 
(ii) $\X_{w_i}$ ($i = 1, 2$) is not fixed by the action of $\sigma$;
(iii) $\mco_{w_i} = \U_{w_i}\circ\N_{E/F}$ ($i = 1, 2$) for some character $\U_{w_i}$ of $F_{w_i}^\times$;
(iv) $\X_v$, $\mco_v$ are unramified for each finite place $v \neq w_1, w_2, w$.

Let $\{\Pi\}$ be the global unstable packet $[\pi(\OO\X), \pi(\OO\,{}^\sigma\!\X)]$ of $\mb{G}(\Af)$, described as follows:
\begin{itemize}
\item
At a place $v$ where $\pi(\OO\X)_v$ and $\pi(\OO\,{}^\sigma\!\X)_v$ are fully induced, $\{\Pi_v\}$ is the singleton
consisting of an irreducible fully induced representation, denoted by $\Pi_v^+$.
For convenience, we put $\Pi_v^- := 0$ and write $\{\Pi_v\} = \{\Pi_v^+, \Pi_v^-\}$.

In particular, $\{\Pi_v\}$ is a singleton for each $v$ which splits in $E$.  Since $F$ is totally imaginary, the archimedean
places all split in $E$.
\item
At a place $v$, prime in $E$, 
such that $\OO_v = {}^\sigma \OO_v$ but $\X_v \neq {}^\sigma \X_v$, the local packet $\{\Pi_v\}$ consists of the
two inequivalent tempered constituents $\Pi_v^+$, $\Pi_v^-$ of $1 \rtimes \pi(\OO_v\X_v)$.
\item
At a place $v$, prime in $E$,  where none of $\OO_v, \X_v, \OO_v\X_v, \OO_v\,{}^\sigma\!\X_v$ is fixed by $\sigma$,
the local packet $\{\Pi_v\}$ consists of two inequivalent cuspidal representations $\Pi_v^+$, $\Pi_v^-$.
In particular, $\{\Pi_w\} = \{\pi^+, \pi^-\}$.
\end{itemize}

From \cite[Sect. 10.3]{F1},
a representation $\Pi' \in \{\Pi\}$ appears in the discrete spectrum of $\mb{G}$ with multiplicity
\[
m(\Pi') = \frac{1}{2}\lp 1 + (-1)^{n(\Pi')}\rp,
\]
where $n(\Pi')$ is the number of places $v$ for which $\Pi'_v = \Pi_v^-$.  In other words, 
$\Pi'$ appears with multiplicity one in the discrete spectrum if and only if $\Pi_v' = \Pi_v^-$
for an even number of places $v$.

Let $\pi_1$ be the representation $\pi(\X)\otimes_1 \OO$ of $\mb{H}(\Af)$.  Let $\cp_1 = \pi(\OO)\otimes_1 \X$.
By the construction of $\X$ and $\OO$, we have $\pi_{1,w} = \pi(\chi)\otimes_1 \theta$ and $\cp_{1,w} = \pi(\theta)\otimes_1 \chi$.
Moreover, the representations $\pi_{1, v}$, $\cp_{1, v}$ are parabolically induced for each place $v \neq w, w_1, w_2$.

For a finite set of places $S$,
let $P^+(S)$ denote the set of representations $\Pi'$ in $\{\Pi\}$
such that $\Pi'_v$ is unramified for each $v \notin S$, and $\Pi'_v = \Pi_v^-$ for an {\it even} number of places $v$ in $S$.

Let $P^-(S)$ denote the set of $\Pi'$ in $\{\Pi\}$ such that $\Pi_v'$ is unramified for each $v\notin S$, and
$\Pi_v' = \Pi_v^-$ for an {\it odd} number of places $v$ in $S$.

Let $S$ be the finite set of places  $\{w, w_1, w_2\}$.  
Let $f$, $f_1$ be matching functions in $\ELw{\mb{G}(\Af), 1}$, $\ELw{\mb{H_1}(\Af), 1}$,
respectively, whose local components at all the places outside of $S$ are spherical.
By \eqref{eq:monomial1} and the usual cancellation, we obtain:
\begin{equation}\label{eq:H1monomialS}
2\!\!\sum_{\Pi' \in P^+(S)}\!\!\la \Pi', f\ra_{\E, S} = \prod_{v \in S}\la \pi_{1, v}, f_{1, v}\ra + \prod_{v \in S}\la \cp_{1, v}, f_{1, v}\ra.
\end{equation}

For $v = w_1, w_2$, we have by construction $\OO_v = \U_v\circ\N_{E/F}$ for some character $\U_v$ of $F_v^\times$.
Hence,
\[
\begin{split}
\pi_{1,v} &= \pi(\X_v)\otimes_1 \U_v\circ\N_{E/F} = \U_v\pi(\X_v)\otimes_1 1,\\
\cp_{1, v} &= \U_v\, I(1, \E_v)\otimes_1 \X_v.
\end{split}
\]
The local packet $\{\Pi_v\}$ consists of the inequivalent tempered constituents $\Pi_v^+$, $\Pi_v^-$ of
$1 \rtimes \U_v\pi(\X_v)$.
By Lemma \ref{lemma:localmonomialinduced},
we have:
\begin{equation}\label{eq:localmonomindlev}
\begin{split}
\la \pi_{1,v}, f_{1, v} \ra &= \ep_v^+\la \Pi_v^+, f_v\ra_{\E_v} + \ep_v^- \la \Pi_v^-, f_v\ra_{\E_v}.
\end{split}
\end{equation}
Here, $\ep_v^+,\ep_v^- = \pm 1$.  They appear because each twisted character in Lemma \ref{lemma:localmonomialinduced}
may differ from the corresponding character here by a sign.
Since $f_{1,v}$ is elliptic and $\cp_{1,v}$ is parabolically induced,
the term $\prod_{v \in S}\la \cp_{1, v}, f_{1, v}\ra$ in \eqref{eq:H1monomialS} vanishes.

Let $v$ be still $w_1$ or $w_2$.  Let $\Pi_v = 1 \rtimes \pi(\OO_v\X_v) = 1\rtimes \U_v\pi(\X_v)$.
Let $A$ be the intertwining operator in $\Hom_{G_v}(\Pi_v,\,\E_v\Pi_v)$ defined in Section \ref{sec:intop}.
By Lemma \ref{lemma:trG}, the twisted character $\la \Pi_v, f_v\ra_A$ is zero for all elliptic test function $f_v$
on $G_v$.
Since $\Pi_v = \Pi_v^+ + \Pi_v^-$, 
the twisted character $\la \Pi_v^+, f_v\ra_A$ must be equal to $-\la \Pi_v^-, f_v\ra_A$ for all elliptic $f_v$.
Hence, there exists a constant $\xi_v = \pm 1$ such that $\la \Pi_v^+, f_v\ra_{\E_v} = \xi_v \la \Pi_v^-, f_v\ra_{\E_v}$,
and
\[
\la \pi_{1,v}, f_{1, v} \ra = \ep_v^+\la \Pi_v^+, f_v\ra_{\E_v} + \ep_v^- \xi_v\la \Pi_v^+, f_v\ra_{\E_v}.
\]

Since $\la \pi_{1,v}, f_{1,v}\rp \neq 0$ for a general elliptic test function $f_{1,v}$, 
the twisted character $\la \Pi_v^+, f_v\ra_{\E_v}$ is nonzero.
Moreover, the sign $\ep_v^+$ must be equal to $\ep_v^-\xi_v$.  There is therefore a constant $\ep = \pm 1$ such that
\[
\la \pi_{1,v}, f_{1, v} \ra = 2\;\ep_v\la \Pi_v^+, f_v\ra_{\E_v}.
\]

We conclude that the equation \eqref{eq:H1monomialS} is equivalent to
\begin{multline}\label{eq:H1monomialww1w2}
2\la \Pi_{w_2}^+, f_{w_2}\ra_{\E_{w_2}}\\
\cdot \lp \sum_{\Pi' \in P^+(\{w, w_1\})}\prod_{\;\;v \in \{w, w_1\}}\!\!\la \Pi'_v, f_v\ra_{\E_v}
+\; \xi_{w_2}\!\! \sum_{\Pi' \in P^-(\{w, w_1\})}\prod_{\;\;v \in \{w, w_1\}}\!\!\la \Pi'_v, f_v\ra_{\E_v} \rp
\\ = 2\;\ep_{w_2}\la \Pi_{w_2}^+, f_{ w_2}\ra_{\E_{w_2}}\cdot\!\!\prod_{v \in \{w, w_1\}}\!\!\la \pi_{1, v}, f_{1, v}\ra .
\end{multline}
Hence,
\begin{multline*}
\sum_{\Pi' \in P^+(\{w, w_1\})}\prod_{\;\;v \in \{w, w_1\}}\!\!\la \Pi'_v, f_v\ra_{\E_v}
+\; \xi_{w_2}\!\! \sum_{\Pi' \in P^-(\{w, w_1\})}\prod_{\;\;v \in \{w, w_1\}}\!\!\la \Pi'_v, f_v\ra_{\E_v}\\
= \ep_{w_2}\!\!\prod_{v \in \{w, w_1\}}\la \pi_{1, v}, f_{1, v}\ra.
\end{multline*}

By repeating the same argument for the place $w_1$, we conclude that there are constants
$\xi^+\!,\, \xi^- = \pm 1$ such that
\begin{equation*}
\xi^+ \la \Pi_w^+, f_w\ra_{\E_w} +\; \xi^- \la \Pi_w^-, f_w\ra_{\E_w} = \la \pi_{1, w}, f_{1, w}\ra.
\end{equation*}
By the symmetry between $\pi_{1,w} = \pi(\chi)\otimes_1 \theta$ and $\cp_{1,w} = \pi(\theta)\otimes_1 \chi$,
we conclude that there are constants $\tilde{\xi}^+, \tilde{\xi}^- = \pm 1$ such that
\begin{equation*}
\tilde{\xi}^+ \la \Pi_w^+, f_w\ra_{\E_w} +\; \tilde{\xi}^- \la \Pi_w^-, f_w\ra_{\E_w} = \la \cp_{1, w}, f_{1, w}\ra.
\end{equation*}

Since $\pi_{1,w} \ncong \cp_{1,w}$, the ordered couple $(\xi^+, \xi^-)$ is equal
to neither $(\tilde{\xi}^+, \tilde{\xi}^-)$ nor \dt $(-\tilde{\xi}^+, -\tilde{\xi}^-)$, 
or else we would have $\la \pi_{1, w}, f_{1, w}\ra = \pm \la \cp_{1, w}, f_{1, w}\ra$.

For $* = +/-$, let $\rho(\E)_w^*$ be the operator in $\Hom(\Pi_w^*,\,\E_w\Pi_w^*)$
restricted from $\rho(\E)$.  The proposition now follows on letting $A^* = \xi^*\,\rho(\E)_w^*$ and $\ep = \xi^+/\tilde\xi^+$.
\end{proof}


\indexi{packet!unstable|)}\indexi{packet!local|)}\indexi{representation!cuspidal|)}\indexi{representation!monomial|)}%
\section{Character Identities for Stable Packets}
\indexi{packet!local|(}\indexi{packet!stable|(}\indexi{lifting!local|(}%
\subsection{Lifting from $H_1$}\label{sec:localH1}
Let $\K/k$ be the quadratic field extension corresponding to $\ve$ via local class field theory.

Let $\chi$ be a character of $\K^\times$ such that $\chi \neq \sig\chi$.
From \cite[Prop. 4.8]{ST}, the induced (local) representation $\nu\ve\rtimes\nu^{-1/2}\pi(\chi)$ of $G$
is reducible of length 2.  Its constituents consist of 
a nontempered quotient, denoted by  $L\!\lp\nu\ve, \nu^{-1/2}\pi(\chi)\rp$, and
a square integrable subrepresentation, which we denote here by $\delta\!\lp\nu\ve, \nu^{-1/2}\pi(\chi)\rp$.
\begin{lemma}\label{lemma:nonbcH11dimlocal}
Let $L$ denote the local representation $L\!\lp\nu\ve, \nu^{-1/2}\pi(\chi)\rp$ of $G$.
There exists an intertwining operator $A \in \Hom_G(L, \ve L)$, and a nonzero constant $c$,
such that the following identity holds for matching functions $f \in \C(G)$, $f_1 \in \C(H_1)\!:$
\[
\la L, f\ra_A = c \la 1_{2, k}\otimes_1 \chi, f_1 \ra.
\]
\end{lemma}
\begin{proof}
Construct a totally imaginary number field $F$ and a quadratic extension
$E$ of $F$ such that: 
(i) There exists a finite place $w$ of $F$ such that $F_w = k$, (ii)
$w$ is prime in $E$, (iii) $E_w = \K$.
We identify the generator $\sigma$ of $\Gal(\K/k)$ with that of $\Gal(E/F)$.

Let $V$ be the set of places of $F$.
Let $\Vr \subset V$ be the set of finite places which are unramified in $E$.

Let $\E$ be the quadratic character of $\idc{F}$ corresponding to the extension $E/F$.
In particular, $\E_w = \ve$.
We consider $\mb{G}$ as an $F$-group, and $\mb{H_1}$, $\mb{H_2}$ as the elliptic $\E$-endoscopic groups of $\mb{G}$
over $F$.

Construct a character $\X$ of $\idc{E}$ with $\X_w = \chi$.

Let $\{\Pi\} = \bigotimes_v \{\Pi_v\}$ be the stable quasi-packet \[\{L\!\lp\nu\E, \nu^{-1/2}\pi(\X)\rp\}\] of $\mb{G}(\Af)$
\indexs{L@$L({\rm char., repn.})$}%
(see Section \ref{sec:arthurunstable}, \cite[Prop. V. 10.10]{F1}).
It lifts to the Langlands quotient \rep
$J\!\lp\nu^{1/2}\pi\lp\X\rp, \nu^{-1/2}\pi\lp\X\rp\rp$ of $\GL(4, \Af)$.
\indexi{Langlands quotient}%
Here, by abuse of notation we let $\nu$ denote also the normalized absolute value function on $\Af$.
The local components of $\{\Pi\}$ are as follows:
\begin{itemize}
\item
If $\E_v \neq 1$ and $\X_v \neq \sig \X_v$, then
the local quasi-packet $\{\Pi_v\}$ consists of \dt
\[
L(\nu_v\E_v, \nu_v^{-1/2}\pi(\X_v)),
\]
the unique nontempered quotient of 
$\nu_v\E_v\rtimes\nu_v^{-1/2}\pi(\X_v)$.
In particular, $\{\Pi_w\}$ is the singleton consisting of $L$.
Since $L(\nu_v\E_v, \nu_v^{-1/2}\pi(\X_v))$ is the unique nontempered quotient of the $\E_v$-invariant
$\nu_v\E_v\rtimes\nu_v^{-1/2}\pi(\X_v)$, it is $\E_v$-invariant.
\item
If $\E_v \neq 1$, and $\X_v = \U_v\circ\N_{E/F}$ for some character $\U_v$ of $\idc{F}$, then
$\{\Pi_v\}$ is the quasi-packet
\[
\{L(\nu_v\E_v, \E_v\rtimes \U_v\nu_v^{-1/2}),\; \delta^-(\E_v\nu_v^{1/2}{\rm St}_{2, v}, \nu_v^{-1/2}\U_v)\}.
\]

By Lemma \ref{lemma:veinvariantlocal}, each member of $\{\Pi_v\}$ is $\E_v$-invariant.
\item
If $v$ splits into two places $v_1$, $v_2$ of $E$, then $\{\Pi_v\}$ is the singleton consisting of
the irreducible representation $\X_1/\X_2 \rtimes \X_2\, 1_{2, F_v}$.  Here, $\X_i := \X_{v_i}$ ($i = 1, 2$), and
$1_{2, F_v}$ is the trivial representation of $\GL(2, F_v)$.
\end{itemize}

It follows from the above remarks that each member of $\{\Pi\}$ is $\E$-invariant.
Hence, for each $\Pi' \in \{\Pi\}$,
$\la \Pi', f\ra_\E$ is the product $\prod_{v \in V}\la \Pi'_v, f_v\ra_{\E_v}$ of twisted local characters.
Here, $\la \Pi'_v, f_v\ra_{\E_v} := \tr \Pi'_v(f_v)\rho(\E)_v$, where $\rho(\E)_v$ is the local intertwining
operator in $\Hom_{G_v}(\Pi'_v,\,\E_v\Pi'_v)$ restricted from $\rho(\E)$.

Let $w_1$, $w_2$ be two (finite) places, different from $w$, which are prime in $E$.

Let $S_0$ be the (finite) set of bad places for $\X$.  That is, it is the union of $V - \Vr$ and the 
places $v \in \Vr$ for which $\X_v$ is ramified.  In particular, since $\X_w \neq \sig\X_w$, the place $w$
lies in $S_0$.  Let $S = S_0 \cup \{w_1, w_2\}$.

The split central character of $\pi_1$ is $\X|_{\idc{F}}$.
Let $f \in \ELw{\mb{G}(\Af), \X|_{\idc{F}}}$, $f_1 \in \ELw{\mb{H_1}(\Af), \X|_{\idc{F}}}$ be matching functions whose
local components at the places outside of $S$ are spherical.

Since $\{\Pi\}$ is stable, each representation $\Pi' \in \{\Pi\}$ occurs with multiplicity one
in the discrete spectrum of $\mb{G}$.  Here, we are assuming that the multiplicity one property holds for $\GSp(2)$.

By Claim \ref{claim:globalsoudryH1} and the stability of $\{\Pi\}$, 
the following twisted trace identity holds:
\begin{equation}\label{eq:globalsoudryH1local}
\sum_{\Pi' \in \{\Pi\}} 
\la \Pi', f\ra_{\E, S} = \frac{1}{2}\la 1_2\otimes_1 \X, f_1\ra_{S}
+ \frac{1}{2}\la 1_2\otimes_1 \sig\X, f_1\ra_{S}.
\end{equation}

By the matching condition on the functions $f$ and $f_1$, the terms $\la 1_2\otimes_1 \X, f_1\ra_{S}$ and
$\la 1_2\otimes_1 \sig\X, f_1\ra_{S}$ have equal contribution.  
For each $\Pi' \in \{\Pi\}$, $\Pi'_w$ is equal to \dt $L = L(\nu\ve, \nu^{-1/2}\pi(\chi))$.
Hence, 
\eqref{eq:globalsoudryH1local} is equivalent to
\[
\lp 
\sum_{\Pi' \in \{\Pi\}}\la \Pi', f\ra_{\E, S - \{w\}}\rp\cdot \ep(\Pi')\cdot\la L, f_w\ra_{\E_w}
= \prod_{v \in S}\la 1_{2, F_v}\otimes_1 \X_v,\, f_{1,v}\ra.
\]
Here, $\ep(\Pi') = \pm 1$.  It appears in the equation because, for two distinct members $\Pi'$, $\Pi''$ of
$\{\Pi\}$, the local intertwining operators $\rho(\E)_w$ for $L$ as a component of $\Pi'$ and for
$L$ as a component of $\Pi''$ may differ by a sign.

Fixing the local components of the test functions at each place in $S - \{w\}$, the lemma follows.
\end{proof}
Let $\chi$ be a character of ${\K}^\times$  such that $\chi \neq \sig\chi$.
Let $L = L(\nu\ve, \nu^{-1/2}\pi(\chi))$.  Let $\delta = \delta(\nu\ve, \nu^{-1/2}\pi(\chi))$.
\begin{prop}\label{prop:nonbcH1steinberg}
There exist intertwining operators $A^L$ in $\Hom_G(L,\, \ve L)$, and $A^\delta$ in $\Hom_G(\delta,\,\ve\delta)$, 
such that the following identities hold for matching functions:
\[\begin{split}
\la L, f\ra_{A^L} &= \la 1_{2, k}\otimes_1 \chi, f_1\ra,\\
\la \delta, f\ra_{A^\delta} &= \la {\rm St}_{2, k}\otimes_1 \chi, f_1\ra.
\end{split}\]
\end{prop}
\begin{proof}
Let $I = \nu\ve\rtimes\nu^{-1/2}\pi(\chi)$.  Thus, $I = \delta + L$.

Let $A$ be the operator in $\Hom_G(I, \ve I)$ defined in Section \ref{sec:intop}.
By Proposition \dts \ref{prop:trGH1}, we have:
\begin{equation*}\label{eq:nonbcH1st1}
\la I, f\ra_A = \la I(\nu^{1/2}, \nu^{-1/2})\otimes_1 \chi, f_1\ra.
\end{equation*}

Since $\delta$ is the unique subrepresentation of $I$, and $L$ is the unique quotient,
the operator $A$ induces operators $A|_\delta$ in $\Hom_G(\delta, \ve\delta)$ and
$A|_L$ in $\Hom_G(L, \ve L)$.  We have therefore:
\begin{equation}\label{eq:nonbcH1stI}
\la I, f\ra_A = \la \delta, f\ra_{A|_\delta} + \la L, f\ra_{A|_L} = \la I(\nu^{1/2}, \nu^{-1/2})\otimes_1 \chi, f_1\ra.
\end{equation}

By Lemma \ref{lemma:nonbcH11dimlocal}, there exist a constant $c$ and an
operator $A^L \in \Hom_G(L, \ve L)$ such that 
\begin{equation}\label{eq:nonbcH1st2}
\la L, f\ra_{A^L} = c \la 1_{2, k}\otimes_1 \chi, f_1\ra.
\end{equation}
The operator $A^L$ differs from $A|_L$ by at most a sign.
Absorbing any negative sign into the constant $c$ if necessary, we may assume $A^L = A|_L$.

Since
\[
I(\nu^{1/2}, \nu^{-1/2})\otimes_1 \chi = 1_{2, k}\otimes_1\chi + {\rm St}_{2, k}\otimes_1 \chi,
\]
by \eqref{eq:nonbcH1stI} and \eqref{eq:nonbcH1st2} we have:
\[
\la \delta, f\ra_{A|_{\delta}} - \la {\rm St}_{2, k}\otimes_1 \chi, f_1\ra = (1 - c)\la 1_{2, k}\otimes_1 \chi, f_1\ra.
\]
The right-hand side of the above equation consists of a nontempered representation, \dt while the representations which appear
on the left are square integrable.  
The central exponents of a square integrable representation all decay.
By the linear independence of central exponents (see \cite[Sect. 21]{FK}), we conclude that $c = 1$.

Hence,
\begin{align*}
\la L, f\ra_{A^L} &= \la 1_{2, k}\otimes_1 \chi, f_1\ra, &
\la \delta, f\ra_{A|_\delta} &= \la {\rm St}_{2, k}\otimes_1 \chi, f_1\ra.
\end{align*}
Let $A^\delta = A|_\delta$.  The proposition follows.
\end{proof}
\begin{prop}\label{prop:nonbcH1local}
Let $\tau$ be an irreducible, cuspidal, non-$\K$-monomial, admissible representation of $\GL(2, k)$,
$\chi$ be a character of $\K^\times$, such that: There does not exist a character $\ve'$ of $k^\times$ such that
${}^\sigma\chi/\chi = \ve'\circ\N_{\K/k}$ and $\tau \cong \ve'\tau$.
There exist a collection $[\pi]$ of distinct $\ve$-invariant, cuspidal, admissible representations of $G$, 
\indexs{b@$[{\rm repn.}]$}%
and an operator $A(\pi')$ in $\Hom_G(\pi', \ve\pi')$ for each $\pi' \in [\pi]$, 
such that the following holds for matching functions:
\[
\sum_{\pi' \in [\pi]}\la \pi', f \ra_{A(\pi')} = \la \tau\otimes_1 \chi,\, f_1\ra.
\]
\end{prop}
\begin{proof}
Construct a totally imaginary number field $F$ and a quadratic extension $E$ of $F$ such that:
(i) There exists a place $w$ of $F$ such that $F_w = k$;
(ii) $w$ is prime in $E$, with $E_w = \K$.
We identify the generator $\sigma$ of $\Gal(\K/k)$ with that of $\Gal(E/F)$.

Let $\E$ be the quadratic character of $\idc{F}$ corresponding to $E/F$.  Thus, $\E_w = \ve$.
We view $\mb{G}$ as an $F$-group, and $\mb{H_1}$, $\mb{H_2}$ as the elliptic $\E$-endoscopic groups of $\mb{G}$ over $F$.

Let $w_1$, $w_2$ be two (finite) places, different from $w$, which are prime in $E$.

Construct a character $\X$ of $\idc{E}$ such that:
\begin{itemize}
\item
$\X_w = \chi$;
\item
$\X_{w_i} \neq \sig \X_{w_i}$ ($i = 1, 2$);
\item
$\X_v$ is unramified for each finite place $v \neq w, w_1, w_2$.
\end{itemize}

Construct a cuspidal automorphic representation $\mct$ of $\GL(2, \Af)$
such that:
\begin{itemize}
\item
$\mct_w = \tau$;
\item
$\mct_{w_i} = {\rm St}_{2, w_i}$ ($i = 1, 2$);
\item
$\mct_v$ is unramified for each finite place $v \neq w, w_1, w_2$.
\end{itemize}
In particular, since $\tau$ is not $\K$-monomial, $\mct$ is not $E$-monomial.

Let $\pi_1$ be the automorphic representation $\mct\otimes_1 \X$ of $\mb{H_1}(\Af)$.
Let $\sig \pi_1$ be $\mct\otimes_1 \sig \X$.

Let $V$ be the set of places of $F$.

Let $\{\Pi\} = \bigotimes_{v \in V}\{\Pi_v\}$ be the global packet of $\mb{G}(\Af)$ whose global datum is the lift of
the global datum parametrizing $\pi_1$ (see Section \ref{sec:globaldata}).  
By Lemma \ref{lemma:nbcH1stable}, $\{\Pi\}$ is a stable packet.  Hence, each representation $\Pi' \in \{\Pi\}$ occurs
with multiplicity one in the discrete spectrum of $\mb{G}$.  Here, we are assuming that $\GSp(2)$ has the multiplicity one
property.

Let $\{\Pi\}^\E$ denote the set of $\E$-invariant representations in $\{\Pi\}$.
For each $\Pi' \in \{\Pi\}^\E$, the character $\la \Pi', f\ra_\E$ is a product $\prod_{v \in V}\la \Pi'_v, f_v\ra_{\E_v}$
of twisted local characters.  

For a finite set of places $S$,
let $\{\Pi\}^\E_S$ denote the set of representations in $\{\Pi\}^\E$ which have the same 
unramified local component $\Pi_v^0$ at each finite place $v \notin S$.

Let $\Vr \subset V$ be the set of finite places which are unramified in $E$.
Let $S_0 = V - \Vr$.  Let $S = S_0 \cup \{w, w_1, w_2\}$.  
Then, $\pi_{1,v}$ is unramified for each place $v \notin S$.

Let $\omega$ be the split central character of $\pi_1$.
Let $f$ and $f_1$  be matching functions in \dt $\ELw{\mb{G}(\Af), \omega}$ and $\ELw{\mb{H_1}(\Af), \omega}$, respectively, such that
their local components at the finite places outside of $S$ are spherical.

By Proposition \ref{prop:nonbcH1} and the stability of $\{\Pi\}$, 
the following twisted trace identity holds:
\begin{equation}\label{eq:prenonbcH1S}
\sum_{\Pi' \in \{\Pi\}^\E_S} 
\la\Pi', f\ra_{\E, S} = \frac{1}{2}\la \pi_1,f_1\ra_S 
+ \frac{1}{2}\la \sig\pi_1, f_1\ra_S.
\end{equation}

The matching condition on $f$ and $f_1$ implies that $\la \pi_1, f_1\ra_S$ and $\la \sig \pi_1, f_1\ra_S$ 
have equal contribution to the trace identity.
Consequently, \eqref{eq:prenonbcH1S} is equivalent to
\begin{equation}\label{eq:nonbcH1S}
\sum_{\Pi' \in \{\Pi\}^\E_S}\la \Pi', f\ra_{\E, S} = \la \pi_1, f_1\ra_S = \prod_{v \in S}\la \pi_{1,v}, f_{1,v}\ra.
\end{equation}

Since $F$ is totally imaginary, and $\mct_v$ is unramified for each finite place $v \neq w, w_1, w_2$, 
the representation $\pi_{1,v}$ is parabolically induced for each place $v$ in $S - \{w, w_1, w_2\}$.
For $i = 1, 2$, the representation $\pi_{1, w_i}$ is equal to ${\rm St}_{2, w_i}\otimes_1 \X_{w_i}$.
By Propositions \ref{prop:trGH1}, \ref{prop:nonbcH1steinberg}, Lemma \ref{lemma:indsplitH1}, and 
the linear independence of characters, 
the local components at each place $v \in S - \{w\}$ of the representations in $\{\Pi\}^\E_S$ are all equal to the
same representation $\Pi_{0,v}$.

Hence, the $(S - \{w\})$-components of \eqref{eq:nonbcH1S} cancel up to a sign.
That is, there exists a set $\left\{\ep(\Pi') = \pm 1 : \Pi' \in \{\Pi\}^\E_S\right\}$ such that
\[
\sum_{\Pi' \in \{\Pi\}^\E_S}\ep(\Pi')\la \Pi'_w, f_w\ra_{\E_w} = \la \pi_{1, w}, f_{1, w}\ra.
\]
The signs $\ep(\Pi')$ appear because, for a place $v \in S - \{w, w_1, w_2\}$ (resp. $v = w_1, w_2$),
the twisted character in Proposition \ref{prop:trGH1} (resp. Proposition \ref{prop:nonbcH1steinberg})
may differ from $\la \Pi'_v , f_v\ra_{\E_v}$ by a sign.

For $\Pi' \in \{\Pi\}^\E$, let $\rho(\E)|_{\Pi'\!,\,w}$ denote the operator in $\Hom_{G_w}(\Pi'_w,\,\E_w \Pi'_w)$ 
which is defined by the restriction of $\rho(\E)$ to $\Pi'$.
Let
\[
[\pi] = \left\{\Pi'_w : \Pi' \in \{\Pi\}^\E_S\right\}.
\]
Since the local components at each place $v \neq w$ of the representations in $\{\Pi\}^\E_S$ are all equal to
the same representation, the elements in $[\pi]$ are distinct.

Since $\pi_{1,w}$ is cuspidal, its central exponents are equal to zero.
Hence, by the linear independence of central exponents (\cite[Sect. 21]{FK}),
every representation in $[\pi]$ is cuspidal.
For each $\pi' = \Pi'_w$ in $[\pi]$, let $A(\pi') = \ep(\Pi')\cdot \rho(\E)|_{\Pi'\!,\,w}$.  The proposition follows.
\end{proof}

\subsubsection{Some Additional Twisted Character Identities for $H_1$}
Using the global character identities listed in Section \ref{sec:somegloballifts}, we list
the following local character identities.  We skip their proofs, which are no different in nature from
what we have encountered so far.  
For an $\ve$-invariant representation $\pi$ of $G$, we let $A$ be the (unique up to a sign) nontrivial intertwining
operator in $\Hom_G(\pi, \ve\pi)$.
\begin{enumerate}
\item
Let $\theta, \chi$ be characters of $\K^\times$ such that  $\theta\chi = \mu\circ\N_{\K/k}$ 
  for some character $\mu$ of $k^\times$.
    
Suppose ${\sig{\chi}}/{\chi}$ does not factor through $\N_{\K/k}$.
    Then, $\pi\lp{\sig{\chi}}/{\chi}\rp$ is cuspidal, 
    and the following holds for matching functions:
    \[
    \la \pi\!\lp \frac{\sig{\chi}}{\chi}\rp \rtimes \mu,\, f\ra_{\!\!A} = \la \pi(\theta)\otimes_1 \chi,\, f_1\ra.
    \]
    
If ${\sig{\chi}}/{\chi} = \ve'\circ\N_{\K/k}$ for some nontrivial quadratic character
    $\ve'$ of $k^\times$, then the following holds for matching functions:
    \[
    \la \ve'\times\ve'\ve\rtimes\mu,\, f\ra_A = \la \pi(\theta)\otimes_1 \chi,\, f_1\ra.
    \]
 \item
  Let $\K' \neq \K$ be a quadratic extension of $k$.  Let $\ve'$ be the quadratic character
  of $k^\times$ corresponding to $\K'/k$ via local class field theory.  

  Let $\chi$ be a character of $\K^\times$ such that 
  ${\sig \chi}/{\chi} = \ve'\circ\N_{\K/k}$.  In particular, $\chi \neq \sig \chi$.
  
  Let $\theta$ be a character of ${\K'}^\times$ which is not fixed by the action of $\Gal(\K'/k)$.  
  Let $\mu$ be a character of ${\K'}^\times$ which satisfies $\mu\circ\N_{\K\K'/\K'} = \chi \circ \N_{\K\K'/\K}$.
  Here, $\K\K'$ is the compositum of $\K$ and $\K'$.  
  
  Let $\pi_{\K'}(\theta)$, $\pi_{\K'}(\mu\theta)$ be the $\K'$-monomial representations associated with
  $\theta$, $\mu\theta$, respectively.
  
  The following holds for matching functions:
  \[
  \la \ve\ve'\rtimes\pi_{\K'}(\mu\theta),\, f\ra_A = \la \pi_{\K'}(\theta)\otimes_1 \chi,\, f_1\ra.
  \]
 \end{enumerate}
\subsection{Lifting from $H_2$}\label{sec:localH2}
Let $\nu, \nu_\K$ be the normalized absolute value functions on the $p$-adic fields $k, \K$, respectively.
Let $1_{2, l}$ ($l = k, \K$) denote the trivial representation of $\GL(2, l)$.
Let ${\rm St}_{2, l}$ denote the Steinberg representation of $\GL(2, l)$.

For a nontrivial quadratic character $\zeta$ of $k^\times$, and a cuspidal $\zeta$-invariant representation $\tau$ of
$\GL(2, k)$, let $L = L(\zeta\nu, \nu^{-1/2}\tau)$ denote the unique
\indexs{L@$L({\rm char., repn.})$}%
nontempered quotient of the induced representation $\zeta\nu\rtimes\nu^{-1/2}\tau$ of $G$.
Let $\delta = \delta(\zeta\nu, \nu^{-1/2}\tau)$ denote the unique square integrable subrepresentation of
\indexs{delta@$\delta({\rm char., repn.})$}%
$\zeta\nu\rtimes\nu^{-1/2}\tau$.

Let $\chi$ be a character of $\K^\times$ such that $\chi \neq \sig \chi$ and
${\sig \chi}/{\chi} = \ve'\circ\N_{\K/k}$ for some nontrivial quadratic character $\ve'$ of $k^\times$.
Thus, the $\ve$-invariant representation $\pi(\chi)$ is also $\ve'$-invariant (\cite{LL}).
\begin{lemma}\label{lemma:nonbcH21dimlocal}
Let $\zeta = \ve'$ or $\ve'\ve$.  
There exists an operator $A$ in \dts $\Hom_G(L, \ve L)$,
and a nonzero constant $c$, such that the following holds for matching functions $f \in \C(G)$, $f_2 \in \C(H_2)\!:$
\[
\begin{split}
\la L, f\ra_A &= c \la \chi 1_{2, \K}\otimes_2 \chi|_{k^\times}\cdot\zeta,\, f_2\ra.
\end{split}
\]
\end{lemma}
\begin{proof}
This is due to  Claim \ref{claim:globalsoudryH2}.  The proof is otherwise identical with 
that of Lemma \ref{lemma:nonbcH11dimlocal}, and we skip it.
\end{proof}
\begin{prop}\label{prop:nonbcH2steinberg}
Let $\zeta = \ve'$ or $\ve'\ve$.  
If the extension $\K/k$ is unramified, then
there exist operators $A^L$ in $\Hom_G(L, \ve L)$, and $A^\delta$ in $\Hom_G(\delta, \ve\delta)$, such that
the following identities hold for matching functions:
\[
\begin{split}
\la L, f\ra_{A^L} &= \la \chi 1_{2, \K}\otimes_2 \chi|_{k^\times}\cdot\zeta,\, f_2\ra,\\
\la \delta, f\ra_{A^\delta} &= \la \chi {\rm St}_{2,\K}\otimes_2 \chi|_{k^\times}\cdot\zeta,\, f_2\ra.
\end{split}
\]
\end{prop}
\begin{proof}
Since $\K/k$ is unramified, $\nu_\K$ is equal to $\nu\circ\N_{\K/k}$.
By Proposition \ref{prop:trGH2}, we have: 
\begin{multline*}
\la \ve\zeta\nu\rtimes \nu^{-1/2}\pi(\chi), f\ra_A
\\=  \la L, f\ra_A + \la \delta, f\ra_A =
\la \chi I(\nu_\K^{1/2}, \nu_\K^{-1/2})\otimes_2 \chi|_{k^\times}\cdot\zeta,\, f_2\ra.
\end{multline*}
Here, $A$ is the intertwining operator defined in Section \ref{sec:intop}.
Using Lemma \ref{lemma:nonbcH21dimlocal}, the proposition follows from the same argument used in the proof of 
Proposition \ref{prop:nonbcH1steinberg}.
\end{proof}

For a representation $\tau$ of $\GL(2, \K)$, let $\sig \tau$ denote the following $\GL(2, \K)$-module
on the space of $\tau$:
\[
\sig \tau : (g_{ij}) \mapsto \tau((\sigma g_{ij})),\quad \forall (g_{ij}) \in \GL(2, \K).
\]

\begin{prop}\label{prop:nonbcH2local}
Let $\tau$ be an irreducible, cuspidal, admissible representation of \dts $\GL(2, \K)$
such that: {\rm (i)} $\tau \ncong \sig \tau$, and {\rm (ii)} 
$\omega_\tau = \mu \circ\N_{\K/k}$ for some character $\mu$ of $k^\times$.
There exist a collection $[\pi]$ of distinct $\ve$-invariant, cuspidal, admissible $G$-modules, and an operator $A(\pi')$ 
\indexs{b@$[{\rm repn.}]$}%
in $\Hom_G(\pi', \ve\pi')$ for each $\pi' \in [\pi]$, 
such that the following holds for matching functions:
\[
\sum_{\pi' \in [\pi]}\la \pi', f\ra_{A(\pi')} = \la \tau\otimes_2 \mu,\, f_2 \ra.
\]
\end{prop}
\begin{proof}
Construct a totally imaginary number field $F$ and a quadratic extension $E$ of $F$ such that:
(i) There is a finite place $w$ of $F$ such that $F_w = k$; (ii) $w$ is prime in $E$, with $E_w = \K$.
We identify the generator $\sigma$ of $\Gal(\K/k)$ with that of $\Gal(E/F)$.

Let $\E$ be the quadratic character of $\idc{F}$ corresponding to $E/F$.
View $\mb{G}$ as an $F$-group, and $\mb{H_1}$, $\mb{H_2}$ as the elliptic $\E$-endoscopic groups of $\mb{G}$ over $F$.

Let $V$ be the set of places of $F$.  
Let $\Vr \subset V$ be the set of finite places which are unramified in $E$.

Let $w_1$, $w_2$ be two finite places, different from $w$, which are unramified, prime in $E$.

Construct a character $\U$ of $\idc{F}$ such that: (i) $\U_w = \mu$; (ii) $\U_{w_1} = \U_{w_2} = 1$;
(iii) $\U_v$ is unramified for each finite place $v \neq w, w_1, w_2$.

Construct a cuspidal automorphic representation $\mct$ of $\GL(2, \Ae)$ with the following properties:
\begin{itemize}
\item
$\omega_{\mct} = \U\circ\N_{E/F}$;
\item
$\mct_w = \tau$;
\item
For $i = 1, 2$, the representation $\mct_{w_i}$ is equal to $\chi_{w_i}{\rm St}_{2, w_i}$,
where $\chi_{w_i}$ is a character of $E_{w_i}^\times$ such that: 
(i) $\chi_{w_i} \neq \sig \chi_{w_i}$, and (ii) $\chi_{w_i}^2 = 1$;
\item
$\mct_v$ is unramified for each finite place $v \neq w, w_1, w_2$.
\end{itemize}
Let $\pi_2$ be the cuspidal automorphic representation $\mct\otimes_2 \U$ of $\mb{H_2}(\Af)$.
Its local component at $w$ is  $\pi_{2,w} = \tau\otimes_2 \mu$.  Let $\sig \pi_2 = {}^\sigma\mct \otimes_2 \U$. 

Since the number field $F$ is totally imaginary, $\pi_{2, v}$ and $\sig \pi_{2,v}$ are fully induced for each archimedean place $v$.

Let $\{\Pi\} = \bigotimes_{v \in V}\{\Pi_v\}$ be the global packet of  $\mb{G}(\Af)$ whose global datum is the lift of 
the global datum parametrizing $\pi_2$.
Since all the cases of lifting to the unstable (quasi-)packets of $\mb{G}(\Af)$ have been accounted for, 
the packet $\{\Pi\}$ is stable.  Hence, each representation $\Pi' \in \{\Pi\}$ occurs with 
multiplicity one in the discrete spectrum of $\mb{G}$ (that is, assuming that the multiplicity one property is valid for $\GSp(2)$).

Let $\{\Pi\}^\E$ denote the set of representations in $\{\Pi\}$ which are $\E$-invariant.
For each $\Pi'$ in $\{\Pi\}^\E$, the twisted character $\la \Pi', f\ra_\E$ is a product $\prod_{v \in V}\la\Pi'_v, f_v\ra_{\E_v}$ 
of twisted local characters.

For a finite set of places $S$, let $\{\Pi\}^\E_S$ denote the set of  representations in
$\{\Pi\}^\E$ whose local components at each place $v \notin S$ are the same unramified representation $\Pi^0_v$.

Let $S_0 = V - \Vr$.
Let $S = S_0 \cup \{w, w_1, w_2\}$.  Let $\omega$ be the split central character of $\pi_2$.
Let $f$ and $f_2$ be matching functions in $\ELw{\mb{G}(\Af), \omega}$ and $\ELw{\mb{H_2}(\Af), \omega}$, respectively, such that
$f_v, f_{2,v}$ are spherical for each place $v \notin S$.

By Proposition \ref{prop:nonbcH2} and the stability of $\{\Pi\}$, the following trace identity holds:
\begin{equation}\label{eq:prenonbcH2S}
\sum_{\Pi' \in \{\Pi\}^\E_S}\la\Pi', f\ra_{\E, S} = \frac{1}{2}\la \pi_2,f_2\ra_S 
+ \frac{1}{2}\la \sig\pi_2, f_2\ra_S.
\end{equation}

The matching condition on $f, f_2$ implies that $\la \pi_2,f_2\ra_S = \la \sig\pi_2, f_2\ra_S$.
Hence, \eqref{eq:prenonbcH2S} is equivalent to
\begin{equation}\label{eq:nonbcH2S}
\sum_{\Pi' \in \{\Pi\}^\E_S} \la\Pi', f\ra_{\E, S} = \la \pi_2,f_2\ra_S =\prod_{v \in S}\la \pi_{2,v}, f_{2,v}\ra.
\end{equation}

For each place $v \in S - \{w, w_1, w_2\}$, the representation $\pi_{2,v}$ is fully induced.  
For $v = w_1, w_2$, it is equal to $\chi_v {\rm St}_{2, v}\otimes_2 1$.
By Propositions \ref{prop:trGH2}, \ref{prop:nonbcH2steinberg}, and
Lemma \ref{lemma:indsplitH2}, and by the linear independence of characters, we conclude that the local components 
at each place $v \in S - \{w\}$ of the representations in $\{\Pi\}^\E_S$ are all equal to the same representation
$\Pi_{0, v}$.  

Hence, as in the proof of Proposition \ref{prop:nonbcH1local},
the  $(S - \{w\})$-components of \eqref{eq:nonbcH2S} cancel up to a sign.
That is, there exists a set $\left\{\ep(\Pi') = \pm 1 : \Pi' \in \{\Pi\}^\E_S\right\}$ such that
\[
\sum_{\Pi' \in \{\Pi\}^\E_S} \ep(\Pi')\la \Pi'_w, f_w\ra_{\E_w} = \la \pi_{2, w}, f_{2, w}\ra.
\]

For each $\Pi'$ in $\{\Pi\}^\E_S$, let $\rho(\E)|_{\Pi'\!,\,w}$ denote the operator in 
$\Hom_{G_w}(\Pi'_w,\,\E_w \Pi'_w)$ which is defined by the restriction of $\rho(\E)$ to $\Pi'$.
Let
\[
[\pi] = \left\{\Pi'_w : \Pi' \in \{\Pi\}^\E_S\right\}.
\]
Since the local components at each place $v \neq w$ of the representations in $\{\Pi\}^\E_S$ are all equal to
the same representation, the elements in $[\pi]$ are distinct.

Since $\pi_{2,w}$ is cuspidal, its central exponents are equal to zero.
Hence, by the linear independence of central exponents (\cite[Sect. 21]{FK}),
every representation in $[\pi]$ is cuspidal.
For each $\pi' = \Pi'_w$ in $[\pi]$, let $A(\pi') = \ep(\Pi')\cdot \rho(\E)|_{\Pi'\!,\,w}$.  The proposition follows.
\end{proof}

\subsubsection{Some Additional Twisted Character Identities for $H_2$}
Let $\mu$ be a character of $k^\times$.  
The induced representation $\ve\rtimes\mu 1_{2, k}$ of $G$ is irreducible (\cite[Prop. 4.8]{ST}).
It is $\ve$-invariant because $\ve\otimes(\ve\rtimes\mu 1_{2,k}) = \ve \rtimes \ve \mu\, 1_{2, k}$ is equivalent to 
$\ve^{-1}\rtimes \ve\ve\mu\, 1_{2, k} = \ve \rtimes \mu 1_{2, k}$.
\begin{claim}
Suppose the characters $\mu$ and $\ve$ are unramified.
There exists an intertwining operator $A$ in $\Hom_G(\ve\rtimes \mu 1_{2,k},\; \ve \rtimes \ve\mu\, 1_{2,k})$
such that the following holds for matching functions $f \in \C(G), f_2 \in \C(H_2)\!:$
\[
\la \ve \rtimes \mu  1_{2, k},\, f\ra_A = \la \lp\mu\circ\N_{\K/k}\rp\! 1_{2, \K}\otimes_2 \mu^2\ve,\, f_2 \ra.
\]
\end{claim}
\begin{proof}
Construct a totally imaginary number field $F$ and a quadratic extension $E$ of $F$ such that:
(i) There exists a finite place $w$ of $F$ such that $F_w = k$; (ii) $w$ is prime in $E$, with $E_w = \K$.
We identify the generator $\sigma$ of $\Gal(\K/k)$ with that of $\Gal(E/F)$.

Let $V$ be the set of places of $F$.  
Let $\Vr \subset V$ be the set of finite places which are unramified in $E$.

Let $\E$ be the quadratic character of $\idc{F}$ corresponding to $E/F$.
We view $\mb{G}$ as an $F$-group, and $\mb{H_1}$, $\mb{H_2}$ as the elliptic $\E$-endoscopic groups of $\mb{G}$
over $F$.

Let $w_1$, $w_2$ be two (finite) places, different from $w$, which are prime in $E$.

Construct a character $\U$ of $\idc{F}$ such that: (i) $\U_w = \mu$, and (ii) $\U_v$ is unramified for each
finite place $v \neq w, w_1, w_2$.

Let $1_2$ be the trivial representation of $\GL(2, \Af)$.
Let $\mct$ be the following representation of the Heisenberg parabolic subgroup $P$ of $\mb{G}(\Af)$:
\[
\mct = \E \otimes \U\, 1_2 : \lp\lsm a &*\;* &*\\& g_2 &*\\&&\tfrac{\det g_2}{a}\rsm\rp \mapsto \E(a)\, \U(\det g_2).
\]
Let $\Pi = \E\rtimes \U 1_2$, the representation of $\mb{G}(\Af)$ parabolically induced from $\mct$.
\begin{itemize}
\item
At a place $v$ where $\E_v$ is trivial, $\Pi_v$ is reducible, with two nontempered constituents 
$\Pi_v^+ = L(\nu_v, 1\rtimes \nu_v^{-1/2}\U_v)$ and 
$\Pi_v^- = L(\nu_v^{1/2}{\rm St}_{2, F_v}, \nu_v^{-1/2}\U_v)$ (\cite[Lemma 3.8]{ST}).  
\item
If $\E_v \neq 1$, then $\Pi_v$ is irreducible.  We put $\Pi_v^+ := \Pi_v$ and $\Pi_v^- := 0$.
\end{itemize}
In particular, $\Pi_w$ is by construction irreducible and unramified.

Let $\{\Pi\} = \bigotimes_{v \in V}\{\Pi_v^+, \Pi_v^-\}$ be the set of irreducible {\it automorphic} constituents of $\Pi$.
That is, an irreducible constituent $\Pi'$ of $\Pi$ is a member of $\{\Pi\}$ if and only if  $\Pi'_v$ is unramified for 
almost all $v$.

The representation $\Pi$ is $\E$-invariant.  The operator $I_{P, \mct}(\E)$, defined in 
Section \ref{sec:finechiexpoverview}, intertwines $\Pi$ with $\E\Pi$.
For an irreducible constituent $\Pi'$ of $\Pi$, and a test function $f$ on $\mb{G}(\Af)$, 
put $\la \Pi', f\ra_{\E} := \tr \Pi'(f)I_{P, \mct}(\E)$.  It is a product $\prod_{v \in V}\la \Pi'_v, f_v\ra_{\E_v}$
of twisted local characters.

Let $1_{2, E}$ be the trivial representation of $\GL(2, \Ae)$.
Let $\pi_2$ be the representation \dt $\lp\U\circ\N_{E/F}\rp\!1_{2, E}\otimes_2 \U^2\E$ of $\mb{H_2}(\Af)$.
In particular, $\pi_{2,w}$ is equivalent to $(\mu\circ\N_{\K/k})1_{2,\K}\otimes_2 \mu^2\ve$.

Let $S_0 = V - \Vr$. 
Let $S_1 = S_0\cup\{w_1, w_2\}$.  Then, $\Pi_v$ and $\pi_{2,v}$ are unramified for each place $v \notin S_1$.
The split central character of $\pi_2$ is equal to $\U^2\E$.
Let $f \in \ELw{\mb{G}(\Af), \U^2\E}$ and $f_2 \in \ELw{\mb{H_2}(\Af), \U^2\E}$ be matching functions such that
$f_v, f_{2,v}$ are spherical for each $v \notin S_1$.
By Claim \ref{claim:bcH2veeq}, the following holds:
\begin{equation}\label{eq:bcH2veS}
\prod_{v \in S_1}\lp \la \Pi_v^+, f_v\ra_{\E_v} - \la \Pi_v^-, f_v\ra_{\E_v} \rp
\\= \prod_{v \in S_1} \la \pi_{2, v}, f_{2,v}\ra.
\end{equation}

For each $\Pi'$ in $\{\Pi\}$, $\Pi_w'$ is equal to $\Pi_w = \ve\rtimes \mu 1_{2,k}$.  
Let $S = S_1 \cup \{w\}$. 
Let $f$, $f_2$ be matching functions with elliptic components at $w_1$, $w_2$, and
spherical components at all the places outside of $S$.  Applying \eqref{eq:bcH2veS} to $S$, we obtain:
\begin{equation}\label{eq:bcH2veSw}
\la \Pi_w, f_v\ra_{\E_w} \prod_{v \in S_1}\lp \la \Pi_v^+, f_v\ra_{\E_v} - \la \Pi_v^-, f_v\ra_{\E_v} \rp
= \la \pi_{2,w}, f_{2,w}\ra \prod_{v \in S_1}\la \pi_{2,v}, f_{2,v}\ra_{\E_v}.
\end{equation}
By \eqref{eq:bcH2veS}, the $S_1$-components of the above equation cancel; hence,
\[
\la \pi_w, f_w\ra_{\E_w} = \la \pi_{2,w}, f_{2,w}\ra.
\]
Let $A$ be the operator in $\Hom_{G_w}(\Pi_w,\, \E_w\Pi_w)$ which is induced from $I_{P, \mct}(\E)$.
The claim follows.
\end{proof}

\begin{corollary}\label{corollary:bcH2vesteinberg}
Suppose the extension $\K/k$ is unramified.  Let $\mu$ be an unramified character of $k^\times$.
Let $\pi$ be the induced representation $\ve \rtimes \mu\, {\rm St}_{2, k}$ of $G$.
There exists an operator $A$ in $\Hom_G(\pi, \ve\pi)$ such that the following holds for matching functions:
\[
\la \pi, f\ra_A 
= \la \lp\mu\circ\N_{\K/k}\rp {\rm St}_{2, \K}\otimes_2 \mu^2\ve,\, f_2\ra.
\]
\end{corollary}
\begin{proof}
We have:
\[
\lp\mu\circ\N_{\K/k}\rp I(\nu_\K^{1/2}, \nu_\K^{-1/2})
= \lp\mu\circ\N_{\K/k}\rp 1_{2, \K} + \lp\mu\circ\N_{\K/k}\rp {\rm St}_{2, \K}.
\]
Since the extension $\K/k$ is unramified, $\ve$ is an unramified character of $k^\times$, and
$\nu_\K$ is equal to $\nu\circ\N_{\K/k}$.  The corollary then follows from Proposition \ref{prop:trGH2} and the previous claim.
\end{proof}
\begin{prop}
Let $\tau$ be an irreducible, cuspidal, non-$\K$-monomial representation of $\GL(2, k)$.
Let $\pi$ be the induced representation $\ve\rtimes \tau$ of $G$.
There exists an operator $A$ in $\Hom_G(\pi, \ve\pi)$ such that the following holds for matching functions:
\[
\la \pi, f\ra_A = \la B_{\K/k}\tau\otimes_2 \omega_\tau \ve,\, f_2\ra.
\indexi{base change}%
\]
\end{prop}
\begin{proof}
Construct a totally imaginary number field $F$ and a quadratic extension $E$ of $F$ such that:
(i) There exists a finite place $w$ of $F$ such that $F_w = k$; (ii) $w$ is prime in $E$, with $E_w = \K$.
We identify the generator $\sigma$ of $\Gal(\K/k)$ with that of $\Gal(E/F)$.

Let $\E$ be a quadratic character of $\idc{F}$ corresponding to $E/F$.
We consider $\mb{G}$ as an $F$-group, and $\mb{H_1}$, $\mb{H_2}$ as the elliptic $\E$-endoscopic groups of $\mb{G}$
over $F$.

Let $V$ be the set of places of $F$.  
Let $\Vr \subset V$ be the set of finite places which are unramified in $E$.

Let $w_1$, $w_2$ be two (finite) places, different from $w$, which are unramified, prime in $E$.

Construct a cuspidal automorphic representation $\mct$ of $\GL(2, \Af)$ such that:
\begin{itemize}
\item
$\mct_w = \tau$;
\item
$\mct_{w_i} = {\rm St}_{2, F_{w_i}}$ ($i = 1, 2$);
\item
$\mct_v$ is unramified for each finite place $v$.
\end{itemize}

Let $\pi_2 = B_{E/F}\mct\otimes_2 \omega_\mct\E$.  It is a cuspidal automorphic representation of $\mb{H_2}(\Af)$,
with $\pi_{2,w} = B_{\K/k}\tau\otimes_2 \omega_\tau\ve$.

Let $\Pi$ be the $\mb{G}(\Af)$-module $\E \rtimes \mct$.  
It is parabolically induced from the representation $\Omega = \E\otimes\mct$
of the Heisenberg parabolic subgroup $P$ of $\mb{G}(\Af)$.  
In particular, $\Pi_w = \ve \rtimes \tau$.

Since $F$ is totally imaginary, $\mct_v$ is a fully induced representation of $\GL(2, \CC)$ for each \dt  archimedean place $v$.
For each finite place $v \neq w, w_1, w_2$, the representation $\mct_v$ is by construction unramified, hence fully induced.
Therefore, by \cite[Thm 4.2]{ST}, \cite[Lemma 3.4]{ST}, and \cite[Prop. 4.8]{ST},
the representation $\Pi_v$ is irreducible for each place $v$.  Consequently, $\Pi$ is irreducible.

The operator $I_{P, \Omega}(\E)$, defined in Section \ref{sec:finechiexpoverview}, intertwines $\Pi$ with $\E\Pi$.
For a test function $f$ on $\mb{G}(\Af)$, put $\la \Pi, f\ra_{\E_v} := \tr \Pi(f)I_{P, \Omega}(\E)$.
It is a product $\prod_{v \in V}\la \Pi_v, f_v\ra_{\E_v}$ of twisted local characters.


Let $S_0 = V - \Vr$.  Let $S = S_0 \cup \{w, w_1, w_2\}$.
The central character of $\Pi$ and the split central character of $\pi_2$ are both equal to $\omega_\mct\E$.
Let $f$ in $\ELw{\mb{G}(\Af), \omega_{\mct}\E}$ and $f_2$ in $\ELw{\mb{H_2}(\Af), \omega_\mct\E}$ 
be matching functions such that $f_v, f_{2,v}$ are spherical for each place $v \notin S$.

By Claim \ref{claim:bcH2veeq}, we have:
\begin{equation}\label{eq:bcH2vecuspidalS}
\prod_{v \in S}\la \Pi_v, f_v\ra_{\E_v} = \prod_{v \in S} \la \pi_{2, v}, f_{2, v}\ra.
\end{equation}
Note that since $\Pi$ is irreducible, the term $\Pi_v^-$ in Claim \ref{claim:bcH2veeq} is zero for
each place $v$.

By Proposition \ref{prop:trGH2}, Lemma \ref{lemma:indsplitH2}, and Corollary \ref{corollary:bcH2vesteinberg},
the \dts $(S_0 \cup \{w_1, w_2\})$-components of \eqref{eq:bcH2vecuspidalS} cancel up to a sign $\ep = \pm 1$.
Hence,
\begin{equation}
\la \Pi_w, f_w \ra_{\E_w} = \ep \la \pi_{2,w}, f_{2,w}\ra.
\end{equation}
Let $I_{P, \Omega}(\E)_w \in \Hom_{G_w}(\Pi_w,\,\E_w\Pi_w)$ be the operator induced from $I_{P, \Omega}(\E)$.
Let $A = \ep\cdot I_{P, \Omega}(\E)_w$.  The proposition follows.
\end{proof}
\indexi{packet!local|)}\indexi{packet!stable|)}\indexi{lifting!local|)}%

\begin{appendix}
\chapter{Summary of Global Lifting}\label{chap:globaltables}
\indexi{packet!global|(}\indexi{packet!quasi-|(}\indexi{packet!unstable|(}\indexi{packet!stable|(}%
\indexi{lifting!global|(}\indexi{base change|(}\indexi{restriction of scalars|(}%
\subsubsection{Notation:}
\begin{itemize}
\item 
$E/F$ a quadratic extension of number fields, $\sigma$ the generator of $\Gal(E/F)$.
\item
For a number field $Q$, $\idc{Q}$ denotes the \idele class group of $Q$,
${\idcc{Q}}$ the set of (unitary) characters of $\idc{Q}$, and
$1_{2, Q}$ the trivial representation of $\GL(2, \mbb{A}_Q)$.
\item
$\ve$ a character of $\idc{F}$, with $\ve \neq 1 = \ve\circ\N_{E/F}$.
\item
For an automorphic representation $\tau$ of $\GL(2, \Ae)$, 
$\pi(\tau)$ denotes the automorphic $\GL(4, \Af)$-module
associated with $\tau$ via the twisted endoscopic lifting from \dt $\rR_{E/F}\GL(2)$ to $\GL(4)$.
$\rR_{E/F}$ denotes restriction of scalars from $E$ to $F$.
\item
$\mc{A}(\GL(2))$ denotes the set of equivalence classes of automorphic representations of $\GL(2,\Af)$.
For $\tau \in \mc{A}(\GL(2))$, $B_{E/F}\tau$ denotes the automorphic $\GL(2, \Ae)$-module obtained from
$\tau$ via base change.
\end{itemize}
In the following tables, each entry in the second row is a (pseudo-)(quasi-)packet of the group above it in the first row.
Two (quasi-)packets belong to the same row if the global datum parametrizing one is the lift of the global
datum parametrizing the other via the $L$-group embedding associated with their corresponding groups
(see Section \ref{sec:globaldata}).\\
\quad\\
\section{Unstable (quasi-)packets of $G$}\label{sec:tableunstable}
\begin{enumerate}
\item
  $\tau \in \mc{A}(\GL(2))$ cuspidal non-$E$-monomial, or one dimensional:
  \begin{center}
    \begin{tabular*}{.8\textwidth}{@{\extracolsep{\fill}}|c|c|c|c|}
      \hline
      $\GL(4)$ & $G$ & $H_1$ & $H_2$\\
      \hline
      $I(\tau, \ve\tau)$&
      $[\tau, \ve\tau]$&
      $\tau\otimes_1 1$ & 
      $B_{E/F}\tau\otimes_2\omega_{\tau}$\\
      \hline
    \end{tabular*}
  \end{center}
  $[\tau, \ve\tau]$ is the unstable (quasi-)packet which lifts to $I(\tau, \ve\tau)$ of $\GL(4, \Af)$
  (see Section \ref{sec:arthurunstable}).
  
  {\sc Remark:} The above table also covers representations of the form $\tau\otimes_1\chi$, with $\chi = {}^\sigma\chi$, 
  because by Hilbert 90 $\chi = \mu\circ\N_{E/F}$ for some $\mu \in {\idcc{F}}$, and
  $\tau\otimes_1\mu\circ\N_{E/F} = \mu\tau\otimes_1 1$.
\item
  $\theta$, $\chi \in {\idcc{E}}$\;; none of $\theta$, $\chi$, $\theta\chi$,
  $\theta\,\sig\chi$ is fixed by $\sigma$.
  
  $\{\pi(\theta)\otimes_1 \chi\}$ denotes the set of the packets
  $\pi(\theta)\otimes_1 \chi$, $\pi(\theta)\otimes_1\sig\chi$, $\pi(\chi)\otimes_1 \theta$,
  $\pi(\chi)\otimes_1\sig\theta$.
  \begin{center}
    \begin{tabular*}{.8\textwidth}{@{\extracolsep{\fill}}|c|c|c|}
      \hline
      $\GL(4)$ & $G$ & $H_1$\\
      \hline
      $I(\pi(\theta\chi), \pi(\theta\,\sig\chi))$&
      $[\pi(\theta\chi), \pi(\theta\,\sig\chi)]$&
      $\{\pi(\theta)\otimes_1 \chi\}$\\
      \hline
    \end{tabular*}
  \end{center}  
\end{enumerate}
\section{Stable (quasi-)packets}\label{sec:tablestable}
Every $\ve$-invariant stable (quasi-)packet of $\mb{G}(\Af)$ is lifted {\bf exclusively} from either
$\mb{H_1}$ or $\mb{H_2}$.
\subsection{Stable (quasi-)packets lifted from  $H_1$ via $\xi_1 : {}^LH_1\hookrightarrow {}^L G$}
\begin{enumerate}
\item
  $\tau \in \mc{A}(\GL(2))$ cuspidal;  $\chi \in \idcc{E}$ with $\chi \neq {}^\sigma \chi$.  
  \begin{center}
    \begin{tabular*}{.8\textwidth}{@{\extracolsep{\fill}}|c|c|c|}
      \hline
      $\GL(4)$& $G$ & $H_1$\\
      \hline
      $\pi\!\lp\chi B_{E/F}\tau\rp$&
      $\xi_1^*(\tau\otimes_1\chi)$&
      $\tau\otimes_1 \chi$,\; $\tau\otimes_1\sig\chi$\\
      \hline
    \end{tabular*}
  \end{center}
  $\xi_1^*(\tau\otimes_1\chi)$ denotes the packet of $\mb{G}(\Af)$ which is the lift of $\tau\otimes_1 \chi$.
  All of its members are cuspidal.
\item
  $\chi \in \idcc{E}$ with $\chi \neq {}^\sigma \chi$.
  \begin{center}
    \begin{tabular*}{.8\textwidth}{@{\extracolsep{\fill}}|c|c|c|}
      \hline
      $\GL(4)$& $G$ & $H_1$\\
      \hline
      $J(\nu^{1/2}\pi(\chi), \nu^{-1/2}\pi(\chi))$&
      $L(\ve\nu,\nu^{-1/2}\pi(\chi))$ &
      \parbox{.7in}{\flushright$1_{2, F}\otimes_1 \chi$,\\ $1_{2, F}\otimes_1 \sig\chi$}\\
      \hline
    \end{tabular*}
  \end{center}
      {\sc Remark:} Since $\mu 1_{2, F}\otimes_1 \chi = 1_{2, F}\otimes_1 (\mu\circ\N_{E/F})\chi$,
      $\forall \mu \in {\idcc{F}}$, the table above covers all cases of $\tau\otimes_1 \chi$ where
      $\tau$ is a one dimensional representation of $\GL(2, \Af)$.
\end{enumerate}
\subsection{Stable (quasi-)packets lifted from $H_2$ via $\xi_2 : {}^LH_2\hookrightarrow {}^L G$}
\begin{enumerate}
\item
  $\tau \in \mc{A}(\rR_{E/F}\GL(2))$ cuspidal; $\mu \in \idcc{F}$ with $\omega_\tau = \mu \circ \N_{E/F}$.
  \begin{center}
    \begin{tabular*}{.8\textwidth}{@{\extracolsep{\fill}}|c|c|c|}
      \hline
      $\GL(4)$&$G$ & $H_2$\\
      \hline
      $\pi(\tau)$&
      $\xi_2^*(\tau\otimes_2\mu)$& 
      \parbox{.7in}{\flushright$\tau\otimes_2\mu$,\\ ${}^\sigma\tau\otimes_2\mu$}\\
      \hline
    \end{tabular*}
  \end{center}
  $\xi_2^*(\tau\otimes_2\mu)$ denotes the packet of $\mb{G}(\Af)$ which is the lift of $\tau\otimes_2\mu$.
  All of its members are cuspidal.
\item       
  $\chi \in \idcc{E}$, $\ve' \neq \ve \in \idcc{F}$, with
  (i) $\ve' \neq (\ve')^2 = 1$, (ii) ${\sig\chi}/{\chi} = \ve'\circ\N_{E/F}$;
  $\zeta \in \{\ve', \ve\ve'\}$.
  \begin{center}
    \begin{tabular*}{.8\textwidth}{@{\extracolsep{\fill}}|c|c|c|}
      \hline
      $\GL(4)$& $G$& $H_2$\\
      \hline
      $J(\nu^{1/2}\pi(\chi), \nu^{-1/2}\pi(\chi))$&
      $L(\ve\zeta\nu, \nu^{-1/2}\pi(\chi))$&
      $\chi 1_{2, E}\otimes_2 \chi|_{\idc{F}}\cdot\zeta$,\\
      && $\sig\chi 1_{2, E}\otimes_2 \chi|_{\idc{F}}\cdot\zeta$\\
      \hline
    \end{tabular*}
  \end{center}
\end{enumerate}
Each $\ve$-invariant, discrete spectrum, automorphic representation of $\mb{G}(\Af)$ {\it with two elliptic local components}
appears in exactly one of the tables in Sections \ref{sec:tableunstable} and \ref{sec:tablestable}.
\section{Induced representations}\label{sec:tableinduced}
\indexi{representation!induced|(}%
\begin{enumerate}
\item
  $\mu \in \idcc{F}$.
  \begin{center}
    \begin{tabular*}{.8\textwidth}{@{\extracolsep{\fill}}|c@{}|c@{}|c@{}|c@{}|}
      \hline
      $\GL(4)$& $G$ & $H_1$ & $H_2$\\
      \hline
      $I(\mu,\mu,\ve\mu,\ve\mu)$&
      $\ve\times\ve\rtimes\mu$&
      $\mu I(1, 1)\otimes_1 1$ & 
      $\lp \mu\circ\N_{E/F}\rp I(1, 1)\otimes_2\mu^2$\\
      \hline
    \end{tabular*}
  \end{center}
\item
  $\mu \in \idcc{F}$;\; $\theta$, $\chi \in {\idcc{E}}$\;; with
  (i) $\theta \neq {}^\sigma \theta$, (ii) $\chi \neq {}^\sigma \chi$,
  (iii) $\theta \chi = \mu \circ \N_{E/F}$.
  
  (a) $\chi^2 \neq \sig\chi^2$.
  \begin{center}
    \begin{tabular*}{.8\textwidth}{@{\extracolsep{\fill}}|c|c|c|}
      \hline
      $\GL(4)$&$G$ & $H_1$\\
      \hline
      $I\lp\mu\,\pi\!\lp{\sig\chi}/{\chi}\rp, \mu, \mu\ve\rp$&
      $\pi\!\lp{\sig\chi}/{\chi}\rp\rtimes \mu$&
      $\pi(\theta)\otimes_1 \chi$,\; $\pi(\theta)\otimes_1 \sig\chi$\\
      \hline
    \end{tabular*}
  \end{center}
  (b) $\exists\;\ve' \in \idcc{F}$ with (i) $\ve' \neq (\ve')^2 = 1$, (ii) ${\sig\chi}/{\chi} = \ve'\circ\N_{E/F}$.
  \begin{center}
    \begin{tabular*}{.8\textwidth}{@{\extracolsep{\fill}}|c|c|c|}
      \hline
      $\GL(4)$&$G$ & $H_1$\\
      \hline
      $I(\mu, \mu\ve'\ve, \mu\ve, \mu\ve')$&
      $\ve'\times\ve'\ve\rtimes \mu$&
      $\pi(\theta)\otimes_1 \chi$\\
      \hline
    \end{tabular*}
  \end{center}
\item
  $\theta \in {\idcc{E}}$ with $\theta \neq {}^\sigma \theta$.
  \begin{center}
    \begin{tabular*}{.8\textwidth}{@{\extracolsep{\fill}}|c|c|c|}
      \hline
      $\GL(4)$&$G$ & $H_1$\\
      \hline
      $I(\pi(\theta), \pi(\theta))$&
      $1\rtimes \pi(\theta)$&
      $\pi(\theta)\otimes_1 1$\\
      \hline
      $I(\pi(\theta), \pi(\theta))$&
      $\ve\rtimes \pi(\theta)$&
      $I(1, 1)\otimes_1 \theta$,\; $I(1, 1)\otimes_1 {}^\sigma \theta$\\
      \hline
    \end{tabular*}
  \end{center}
\item
  $[E', F] = 2$, $E' \neq E$, with $\Gal(E'/F) = \la \sigma'\ra$.
  
  $\ve' \in \idcc{F}$ with $\ve' \neq \ve' \circ {\N_{E'/F}} = 1$.
  
  $\chi \in \idcc{E}$;\; $\mu, \theta \in \idcc{E'}$\,;
  (i) ${\sig\chi}/{\chi} = \ve'\circ\N_{E/F}$,
  (ii) $\mu\circ\N_{EE'/E'} = \chi \circ \N_{EE'/E}$,
  (iii) $\mu\theta \neq {}^{\sigma'}\!(\mu\theta)$.
  \begin{center}
    \begin{tabular*}{.8\textwidth}{@{\extracolsep{\fill}}|c|c|c|}
      \hline
      $\GL(4)$&$G$ & $H_1$\\
      \hline
      $I\!\lp\pi_{E'}(\mu\theta), \ve\pi_{E'}(\mu\theta)\rp$&
      $\ve\ve'\rtimes\pi_{E'}(\mu\theta)$&
      $\pi_{E'}(\theta)\otimes_1 \chi$\\
      \hline
    \end{tabular*}
  \end{center}
\item
  $\tau \in \mc{A}(\GL(2))$ cuspidal non-$E$-monomial.
  \begin{center}
    \begin{tabular*}{.8\textwidth}{@{\extracolsep{\fill}}|c|c|c|}
      \hline
      $\GL(4)$&$G$ & $H_2$\\
      \hline
      $I(\tau, \ve\tau)$&
      $\ve\rtimes\tau$&
      $B_{E/F}\tau\otimes_2\omega_{\tau}\ve$\\
      \hline
    \end{tabular*}
  \end{center}
\item
  $\mu \in \idcc{F}$.
  \begin{center}
    \begin{tabular*}{.8\textwidth}{@{\extracolsep{\fill}}|c|c|c|}
      \hline
      $\GL(4)$&$G$ & $H_2$\\
      \hline
      $I(\mu,\mu,\ve\mu,\ve\mu)$&
      $1 \times \ve \rtimes \mu$&
      $\lp\mu\circ\N_{E/F}\rp I(1, 1)\otimes_2\mu^2\ve$\\
      \hline
    \end{tabular*}
  \end{center}
\end{enumerate}
Each parabolically induced $\mb{G}(\Af)$-module which contributes discretely to the fine $\chi$-expansion of
the $\ve$-twisted trace formula of $\mb{G}(\Af)$ appears in exactly one of the tables in Section \ref{sec:tableinduced}.
\indexi{packet!global|)}\indexi{packet!quasi-|)}\indexi{packet!unstable|)}\indexi{packet!stable|)}%
\indexi{lifting!global|)}\indexi{representation!induced|)}\indexi{base change|)}\indexi{restriction of scalars|)}%
\chapter{Fundamental Lemma}\label{chap:fundlemma}
Let $F$ be a $p$-adic field.  
Let $\ve$ be a character of $F^\times$ with $\ve^2 = 1$.
To prove the Fundamental Lemma in the context of $\mb{G}$ and $\mb{H_1}$, $\mb{H_2}$,
\indexi{Fundamental Lemma}%
by \cite{H} it suffices to prove it for 
the case where $F$ has odd residual characteristic, and then only for the unit elements
in the Hecke algebras of the groups, when $\ve \neq 1$ and also when $\ve = 1$.
\indexi{Hecke algebra}%
Moreover, we may assume that $\ve$ is an unramified character.
The case of $\ve = 1$ is straightforward, given the results of \cite{F}.
Hence, in this chapter we assume that $F$ has odd residual characteristic, and that
$\ve$ is nontrivial.  

We fix once and for all an element $A \in F^\times - {F^\times}^2$ such that $F_A := F(\sqrt{A})$
is the unramified quadratic extension of $F$ which corresponds to $\ve$ via local class field theory.
Consider $\mb{G} = \GSp(2)$ as an $F$-group.  Let $G = \mb{G}(F)$.
Recall the definitions in Chapter \ref{chap:endoscopy} of the elliptic $\ve$-endoscopic groups of $\mb{G}$:
\begin{enumerate}
\item
$\displaystyle
H_1 = \mb{H_1}(F) = \{(g, \xi) \in \GL(2, F)\times \EA^\times : \det g = \N_{\EA/F}x\}.
$
\item
$H_2 = \mb{H_2}(F)$
\[= \{(g, c) \in \GL(2, \EA)\times F^\times\}/\{(\diag(z, z), \N_{\EA/F}z^{-1}): z \in \EA^\times\}.\]
\end{enumerate}

Let $R$ be the ring of integers in $F$.  Fix a uniformizer $\vp$ of $F$.
Let $K$ be the maximal compact, hyperspecial subgroup $\mb{G}(R)$ of $G$.
Let $1_K$ in $C_c^\infty(G)$ be the characteristic function of $K$.
Let $dg$ be the Haar measure on $G$ for which the volume of $K$ is one.
For a function $f \in C_c^\infty(G)$, define the \emph{$\ve$-twisted orbital integral} of $f$ at an
\indexi{orbital integral}%
element $\delta \in G$ as follows:
\[
O_\delta(f) = \int_{Z_G(\delta)\bs G}f(g^{-1}\delta g)\ve(g)\,dg.
\]
Here, $Z_G(\delta)$ is the centralizer of $\delta$ in $G$.

For $i = 1, 2$, the group $\mb{H}_i$ is defined over $R$.
Let $K_i = \mb{H}_i(R)$.  Let $dh_i$ be the Haar measure on $H_i$ for which $K_i$ has volume one.
Let $1_{K_i} \in C_c^\infty(H_i)$ be the characteristic function of $K_i$.  
For $H = H_1$ or $H_2$, $f_H \in C_c^\infty(H)$, and $\delta_H \in H$, let
\[
\begin{split}
O_{\delta_H}(f_H) &= \int_{Z_H(\delta_H)\bs H}f_H(h^{-1} \delta_H h)\,dh,\\
SO_{\delta_H}(f_H) &= \sum_{\delta_H'} O_{\delta_H}(f_H).
\end{split}
\]
\indexi{orbital integral!stable}%
Here, the sum in the second expression is over representatives $\delta_H'$ of the $F$-conjugacy classes 
within the stable conjugacy class of $\delta_H$.
\begin{thm}\label{thm:fundlemma}
Let $i = 1, 2$.  For a regular element $\delta \in G$ and a $\mb{G}$-regular $\delta_i$ in $H_i$
which is a norm of $\delta$,
there exists a choice of transfer factor $\Delta(\delta_i, \delta) \in \mbb{C}$ such that the following holds:
\indexi{transfer factor}%
\[
SO_{\delta_i}(1_{K_i}) = \sum_{\delta'} \Delta(\delta_i, \delta')O_{\delta'}(1_K).
\]
The sum on the right is over representatives $\delta'$ of the $F$-conjugacy classes
within the the stable conjugacy class of $\delta$.
\begin{remark}
Since $O_{\delta'}(1_{K})$ is equal to zero if $\delta$ is not stably conjugate to an element in $K$, it suffices to prove the theorem for
the elements $\delta$ which lie in $K$.
\end{remark}
\end{thm}
We shall prove the theorem only for the \emph{elliptic} regular elements.  For elements which lie
in parabolic subgroups, the theorem reduces via Harish-Chandra descent (see, for example, \cite{H1})
to the Fundamental Lemma for lower rank groups.

An element $g \in G$ is called {\bf absolutely semisimple} if $g^m = 1$ for some integer $m$ prime to the residual characteristic of $F$.
\indexi{absolutely semisimple}%
An element $g$ is called {\bf topologically unipotent} if  $g^{(q^N)}$ converges to $1$ as the integer $N$ approaches infinity.
\indexi{topologically unipotent}%
Here, $q$ is the cardinality of the residue field of $F$.
By an extension of a theorem of Kazhdan's (\cite[Prop. I. H. 1]{F}), an element $\delta \in K$ decomposes uniquely as $\delta = su = us$, 
\indexi{Kazhdan decomposition}%
where $s$ is absolutely semisimple, $u$ is {topologically unipotent}, and $s, u \in K$.

With regard to proving Theorem \ref{thm:fundlemma},
we restrict further our attention to the elliptic regular elements which are topologically unipotent.  For elements whose absolutely
semisimple parts are nontrivial, the method of semisimple reduction (\cite[Part. III]{F}) reduces the theorem to the 
Fundamental Lemma for lower rank groups.

\subsubsection{Notation}\label{sec:flnotations}
\begin{itemize}
\item 
For a field $k$ and an element $B \in k^\times$, let
$k_B$ denote the possibly trivial field extension $k(\sqrt{B})$ of $k$.
\item
Let $\GL(2, k)^{k_B}$ denote the subgroup $\{g \in \GL(2, k) : \det g \in \N_{k_B/k} k_B^\times\}$
of $\GL(2, k)$.
\item
For $\xi = a + b\sqrt{B} \in k_B^\times$, $a, b \in k$, put
\[
\phi^{k : B}(\xi) := \phi^B(\xi) := \lp\begin{matrix} a & bB\\b&a \end{matrix}\rp.
\indexs{phi@$\phi^{\;:\;}$, $\phi^{\;:\;}_{\rho}$, $\phi^{\;:\;:}$}%
\]
Observe that if $[k_B: k] = 2$, then $\det \phi^{B}(\xi) = \N_{k_B/k}\xi$ for all $\xi \in k_B^\times$.
\item
For $\rho \in k^\times$, $\xi \in k_B$, put
\[
\phi^{k:B}_\rho(\xi) := \phi^B_\rho(\xi) 
:= \lp\begin{matrix} a & bB\rho\\b/\rho&a \end{matrix}\rp.
\]
\item
Let $T_B^\rho$ denote the maximal torus $\{\phi^B_\rho(\xi):\xi \in k_B^\times\}$ in $\GL(2, k)$.  It is elliptic if $[k_B:k] = 2$.
\item 
Let $L = k_B(\sqrt{\xi})$, where $\xi \in k_B^\times$.
For $l = \alpha + \beta \sqrt{\xi} \in L^\times$, $\alpha, \beta \in k_B$, put
\[
\phi^{k:B:\xi}(l) :=
\lp\begin{matrix}
\phi^B(\alpha) & \phi^B(\beta)\phi^B(\xi)\\
\phi^B(\beta) & \phi^B(\alpha)
\end{matrix}\rp.
\]
\item
For $\lp\lsm a &b\\c&d\rsm\rp, \lp\lsm x&y\\z&t\rsm\rp \in \GL(2)$, put
\[
[\lp\lsm a &b\\c&d\rsm\rp, \lp\lsm x&y\\z&t\rsm\rp]
:= 
\lp\lsm a &&& b\\&x&y&\\&z&t&\\c&&&d\rsm\rp \in \GL(4).
\]
\end{itemize}
\section{Norm Correspondence---Elliptic Elements}\label{sec:normmappingell}
\indexi{norm correspondence|(}\indexi{element!elliptic|(}\indexi{maximal torus!elliptic|(}%
We now describe the norm correspondence between the elliptic regular elements of $G$ and 
the elliptic $\mb{G}$-regular elements of  $H_1$, $H_2$.
But first, we need to classify the stable conjugacy classes of the elliptic regular elements of the groups.
This is equivalent to classifying the stable conjugacy classes of the elliptic maximal tori.
We list below the representatives of these stable conjugacy classes.  Moreover, we list the conjugacy
classes within each stable conjugacy class.
\subsection{Elliptic Tori of $G$}
The $\ve$-twisted orbital integral of a function on $G$ at an element $\delta \in G$ is zero 
if $Z_\delta(G)$ is not contained in the kernel of $\ve$.  In fact, as we shall see, elliptic $\mb{G}$-regular
elements of the $\ve$-endoscopic groups can be the norms of only those elliptic regular elements in 
$G$ whose centralizers are contained in
$\ker \ve$.  Thus, for our purpose, it suffices to list only the stable 
conjugacy classes of the elliptic tori which lie in $\ker \ve$.

Recall that $F_A$ is the quadratic extension of $F$ which corresponds to $\ve$ via local class field theory.
The stable conjugacy classes of the elliptic tori of $G$ are computed
in \cite[Sect. I. C]{F}.  The representatives of those classes which lie in the kernel of $\ve$ are as follows:
\begin{itemize}
\item
$\displaystyle T_{\text{I}, A} =
\boxed{
\left\{
\left[\phi^A(\gamma), \phi^A(\zeta)\right]: \gamma,\zeta \in \EA^\times;\; \N_{{\EA}/F} \gamma = \N_{{\EA}/F} \zeta
\right\}
}\,;$
\item
$T_{\text{III}, A, AD} =$
\[
\boxed{
\left\{
\left[\phi^{A}(\gamma), \phi^{AD}(\zeta)\right] :\\ \gamma \in \EA^\times, \zeta \in \EAD^\times;\;
\N_{{\EA}/F}\gamma = \N_{{\EAD}/F} \zeta\right
\}
}\,,
\]
$D, AD \in F^\times - A{F^\times}^2$;
\item
$T_{\text{III}, D, AD} =$
\[
\boxed{
\left\{
\left[\phi^{D}(\gamma), \phi^{AD}(\zeta)\right] : \gamma \in \ED^\times, \zeta \in \EAD^\times;\;
\N_{{\ED}/F}\gamma = \N_{{\EAD}/F} \zeta
\right\}
}\,,
\]
$D, AD \in F^\times - A{F^\times}^2$;
\item
$\displaystyle T_{\text{IV}, A, D} =
\boxed{
\left\{ \phi^{F:A:D}(l) : l \in L^\times;\;\N_{L/{\EA}} l \in F^\times \right\}}\,,
$\\
$L = \EA(\sqrt{D})$, $D \in \EA^\times - {\EA^\times}^2$.
\end{itemize}
Note that $T_{\text{III}, A, AD}$ and $T_{\text{III}, D, AD}$ are distinct cases, for
$A$ depends on $\ve$ and  is not arbitrary.

With the exception of $T_{{\rm I}, A}$, the conjugacy class of each elliptic maximal torus listed above
are stable.  In the case of $T_{{\rm I}, A}$, there are two conjugacy classes within the stable conjugacy
class, described as follows:

For $\xi = a + b \sqrt{A} \in \EA^\times$ ($a, b \in F$), 
put $\displaystyle \phi^A_\vp(\xi) := \lp\lsm a & b \vp A\\ b/\vp&a \rsm\rp.$
Put 
\[
T_{{\rm I}, A}^\vp := \dmfour{1}{1}{\vp}{1}^{-1} T_{{\rm I}, A}\dmfour{1}{1}{\vp}{1}
= \left\{[\phi^A(\gamma), \phi_\vp^A(\zeta)] : \N_{\EA/F}\gamma = \N_{\EA/F}\zeta\right\}.
\]  
From \cite[Sect. I. C]{F}, the tori
$T_{{\rm I}, A}$ and $T_{{\rm I}, A}^\vp$ are representatives of the conjugacy classes within the stable conjugacy class
of $T_{{\rm I}, A}$.
\subsection{Elliptic Tori of $H_1$}
There are two stable conjugacy classes of elliptic tori in $H_1$.  They are represented by the following tori:
\begin{itemize}
\item
\boxed{
T_{1,A} = \{(\phi^A(\gamma),\zeta):\gamma,\zeta\in \EA^\times;\; \N_{\EA/F}\gamma = \N_{\EA/F}\zeta\}
}\,;
\item
\boxed{
T_{1, AD} = \{(\phi^{AD}(\gamma),\zeta) : \gamma \in \EAD^\times,\, \zeta \in \EA^\times;\;
\N_{\EAD/F} \gamma = \N_{\EA/F} \zeta\}
}\,,\\ 
$D, AD \in F^\times - A{F^\times}^2$.
\end{itemize}

An element of ${H_1}$ has the form $(g, \zeta) \in \GL(2, F)\times \EA^\times$, with
$\det g = \N_{\EA/F}\zeta$.
Suppose two elements $(g, \zeta), (g',\zeta')$ in $H_1$ are stably conjugate to each other.
Then, $\zeta = \zeta'$.  Since 
$\mb{H_1}(\bar{F}) = \GL(2, \bar{F})\times  {\bar{F}}^\times\times \bar{F}^\times$, 
and there is no distinction between
stable conjugacy and $F$-conjugacy for $\GL(2)$, the elements $g, g'$ are conjugate via some element of $\GL(2, F)$.

Let $\GL(2, F)^{\EA}$ denote the subgroup of elements in $\GL(2, F)$ whose determinants lie in $\N_{\EA/F}\EA^\times$.
We have an exact sequence:
\[
1 \rightarrow \GL(2, F)^{\EA} \rightarrow \GL(2, F) \xrightarrow{\det}
F^\times / \N_{\EA/F} \EA^\times \rightarrow 1.
\]
Since $[F^\times : \N_{\EA/F}\EA^\times] = 2$,
$\GL(2, F)^{\EA}$ is a subgroup of index 2 in $\GL(2, F)$.  Thus,
there are at most two conjugacy classes within the stable conjugacy class
of each element $g \in \GL(2, F)^{\EA}$.
\begin{claim}
A regular element $\delta_1 = (\phi^A(\gamma), \zeta)\in T_{1, A}$ is stably conjugate but not conjugate to 
$\delta_1^\vp = (\phi_\vp^A(\gamma),\zeta)$.
\end{claim}
\begin{proof}
Since $\phi^A_\vp(\gamma) = \lp\lsm 1&\\&\vp\rsm\rp^{-1}\phi^A(\gamma)\lp\lsm 1&\\&\vp\rsm\rp$, it is clear that
$\delta_1$ is stably conjugate to $\delta_1^\vp$.
Suppose $g$ is an element in $\GL(2, F)$ such that
$g^{-1} \phi^A(\gamma) g = \phi^A_\vp(\gamma)$.  Then, $g\lp\lsm 1&\\&\vp\rsm\rp$ commutes with $\phi^A(\gamma)$;
hence, $g\lp\lsm 1&\\&\vp\rsm\rp = \lp\begin{smallmatrix}x & Az\\z&x\end{smallmatrix}\rp$ for some $x, z \in F$.
The determinant of $g$ is therefore equal to
$\vp^{-1}(x^2 - z^2 A)$, 
which does not lie in $\N_{{\EA}/F} \EA^\times$, for $x^2 - z^2 A$ lies in $\N_{\EA/F}\EA^\times$, while $\vp$ does not
because $\EA/F$ is unramified.
Hence, $g \notin \GL(2, F)^{\EA}$, which implies that $\delta_1^\vp$ is not conjugate to ${\delta_1}$ in $H_1$.
\end{proof}
\begin{corollary}
The conjugacy class of the tori in $H_1$ which are stably conjugate but not conjugate to $T_{1, A}$
is represented by
\[
T_{1, A}^\vp := \left\{\lp\phi^A_\vp(\gamma), \zeta\rp: \gamma,\zeta \in \EA^\times;\;
\N_{\EA/F}\gamma = \N_{\EA/F}\zeta\right\}.
\]
\end{corollary}
\begin{claim}
The conjugacy class of $T_{1, AD}$ is stable.
\end{claim}
\begin{proof}
Suppose an element $\delta' = (t', \zeta)$ in $H_1$ is stably conjugate to the element
$\delta = (\phi^{AD}(\xi), \zeta) \in T_{1, AD}$, where
$\N_{\EAD/F}\xi = \N_{\EA/F}\zeta$.  There exists $g$ in $\GL(2, F)$ such that $t' = g^{-1}\phi^{AD}(\xi)g$.
Since $\GL(2, F)$ is the disjoint union of $\GL(2, F)^{\EA}$ and \dt $\GL(2, F)^{\EA}\diag(1, \vp)$,
we may assume that $g = \diag(1, \vp)$.
Since by assumption the residual characteristic of $F$ is odd, and $\EA/F$ is unramified, the extension $\EAD/F$ is ramified.
Hence, there exists an element $\gamma \in \EAD^\times$ such that $\N_{\EAD/F}\gamma = \vp^{-1}$.  Since $\phi^{AD}(\xi)$ commutes with
$\phi^{AD}(\gamma)$, $t'$ is equal to 
\[
\lp\phi^{AD}(\gamma) g\rp^{-1}\cdot\phi^{AD}(\xi)\cdot\lp\phi^{AD}(\gamma) g\rp.
\]  
The element $\phi^{AD}(\gamma) g$,
having determinant $1$, lies in $\GL(2, F)^{\EA}$; hence, $\delta'$ is conjugate to $\delta$ in $H_1$.
\end{proof}
\subsection{Elliptic Tori of $H_2$}
For an element $D \in \EA^\times$, and $l = \alpha + \beta\sqrt{D} \in \EA(\sqrt{D})^\times$ ($\alpha, \beta \in \EA$), 
put
\[
\phi^{\EA : D}(l) := \lp\lsm \alpha & \beta D\\\beta&\alpha\rsm\rp \in \GL(2, \EA).
\]

A maximal elliptic torus in $H_2$ has the following form:
\[
\boxed{
T_{2, D} = \left\{\lp\phi^{{\EA}:D}\lp l \rp, c \rp : l \in E^\times\!,\; c \in F^\times\right\}}\;,
\]
where $D \in \EA^\times - {\EA^\times}^2$ and $E = \EA(\sqrt{D})$.
There are 3 possibilities for $D$:
\begin{itemize}
\item
$D \in F^\times -A{F^\times}^2$; hence, $E/F$ is Galois, biquadratic.  
\item
$D \notin F^\times$, $E/F$ cyclic of degree $4$.
\item
$D \notin F^\times$,  $E/F$ non-Galois, with Galois closure $\tilde{E}$ of $E$ over $F$.  
Then, the Galois group $\Gal(\tilde{E}/F)$ is equivalent to $D_4$, the dihedral group of
order 8.
\end{itemize}

Since $H_2$ is a quotient group of $\GL(2, \EA)\times F^\times$, there is no distinction here between
conjugacy and stable conjugacy.
\subsection{Norm Correspondence---Explicit Description}
We now use the results in Section \ref{sec:normmapping} to describe explicitly the norm correspondence
between elliptic regular elements.  

The groups of $\bar{F}$-points of $\mb{H_1}$, $\mb{H_2}$ are as follows:
\[\begin{split}
\mb{H_1}(\bar{F}) &= \{(g, \alpha, \beta) \in \GL(2, \bar{F}) \times \bar{F}^\times \times \bar{F}^\times:
\det g = \alpha\beta\},\\
\mb{H_2}(\bar{F}) &= \lp\GL(2, \bar{F}) \times \GL(2, \bar{F}) \times \bar{F}^\times\rp /
\{(\alpha I_2, \beta I_2, (\alpha\beta)^{-1}) : \alpha, \beta \in \bar{F}^\times\}.
\end{split}\]

Recall the definitions of $\N_1, \N_2$:
\[
\begin{split}
\N_1 &: \mb{H_1}(\bar{F}) \ni (\diag(a, \lambda/a);b,\lambda/b) \mapsto
\diag(b,a, \lambda/a, \lambda/b) \in \mb{G}(\bar{F}),\\
\N_2 &: \mb{H_2}(\bar{F}) \ni (\diag(a, b), \diag(c, d), 1) \mapsto 
\diag(ac, ad, bc, bd) \in \mb{G}(\bar{F}).
\end{split}
\]
By definition, an element $\delta_i \in \mb{H}_i(F)$ ($i = 1, 2$) is a norm of an element $\delta \in \mb{G}(F)$
if $\delta_i$ is conjugate to a diagonal element $t_i$ in $\mb{H}_i(\bar{F})$, and $\delta$ to a 
diagonal element $t$ in $\mb{G}(\bar{F})$, such that $t = \N_i(t_i)$.

For an algebraic group $\mb{H}$, and elements $g, h \in \mb{H}(\bar{F})$, we write
$g \sim h$ if $g, h$ are conjugate under $\mb{H}(\bar{F})$.
\subsubsection{Matching $T_{1, A}$}\label{sec:matT1A}
Let $\sigma$ be the generator of $\Gal(\EA/F)$.
Put $\bar{\gamma} := \sigma \gamma$ for $\gamma \in \EA$.
\begin{claim}
Let $\gamma, \zeta$ be elements in $\EA^\times - F^\times$ such that $\N_{\EA/F}\gamma = \N_{\EA/F}\zeta$.
The element $\delta_1 = \lp \phi^A(\gamma), \zeta\rp$ in $T_{1, A}$ is a norm of
$[\phi^A(\zeta), \phi^A(\gamma)] \in T_{{\rm I}, A}$.
\end{claim}
\begin{proof}
The element $(\phi^A(\gamma),\zeta)$ is conjugate to
$(\diag(\gamma, \bar{\gamma});\zeta, \bar{\zeta})$ in $\mb{H_1}(\bar{F})$.
We have:
\[
\N_1 : (\diag(\gamma, \bar{\gamma});\zeta, \bar{\zeta}) \mapsto \diag(\zeta, \gamma, \bar{\gamma}, \bar{\zeta}).
\]
The element $\diag(\zeta, \gamma, \bar{\gamma}, \bar{\zeta})$ is conjugate to $\delta$ in $\mb{G}(\bar{F})$.
\end{proof}
It is clear from the proof of the claim that
\[
\lp\phi^A(\gamma),\zeta\rp,\;\lp\phi^A(\gamma),\bar{\zeta}\rp,\;\lp\phi^A(\zeta),\gamma\rp,\;
\lp\phi^A(\zeta),\bar{\gamma}\rp \in T_{1, A}
\]
are all norms of $\delta$.  If $\gamma \neq \zeta, \bar{\zeta}$, the elements listed above are
not stably conjugate to one another.
This is an example of a norm correspondence which is not one to one.
\subsection{Matching $T_{1, AD}$}\label{sec:matT1AD}
Let $D$ be an element of $F^\times - A {F^\times}^2$.  
Recall our notation: ${\EAD} = F(\sqrt{AD})$,
${\ED} = F(\sqrt{D})$.
Let $E$ be the biquadratic extension $F(\sqrt{D}, \sqrt{A})$ of $F$.

Suppose $\Gal(E/F) = \langle \sigma, \tau \rangle$, such that
$\EA = E^\tau$, $\ED  = E^\sigma$, and ${\EAD} = E^{\sigma\tau}$.
Here $E^\mu$ denotes the subfield of $E$ fixed by $\mu \in \Gal(E/F)$.
Thus, $\Gal({\ED} / F) = \la \tau \ra = \la \sigma\tau\ra $, $\Gal({\EA}/F) = \la \sigma \ra = \la \tau\sigma\ra$, and
$\Gal({\EAD}/F) = \la \sigma \ra = \la \tau \ra$.
\begin{claim}
Let $\zeta$ be an element of $\EAD^\times - F^\times$, $\gamma$ an element of $\EA^\times - F^\times$, such that
$\N_{\EAD/F}\zeta = \N_{\EA/F}\gamma$.
The element $\delta_1 =  (\phi^{AD}(\zeta),\gamma)$ in $T_{1,AD}$
is a norm of \dt $\delta = [\phi^A(\gamma),\phi^{AD}(\zeta)] \in T_{{\rm III}, A, AD}$.
\end{claim}
\begin{proof}  
Consider $\delta_1$, $\delta$ as elements in $\mb{H_1}(\bar{F})$, $\mb{G}(\bar{F})$, respectively.  We have:
\[
\delta_1 \sim
(\diag(\zeta, \tau \zeta);\gamma, \sigma\gamma) \stackrel{\N_1}{\longmapsto}
\diag(\gamma,\zeta,\tau\zeta,\sigma\gamma) \sim \delta.
\]
The element $\delta$ lies in $T_{\text{III}, A, AD}$ because
\[\det \phi^{AD}(\zeta) = \N_{\EAD/F}\zeta = \N_{\EA/F}\gamma = \det \phi^{A}(\gamma).\]
\end{proof}
Observe that the elements $(\phi^{AD}(\zeta),\gamma)$ and $(\phi^{AD}(\zeta),\sigma\gamma)$ in $T_{1, AD}$, though not
stably conjugate to each other, are both norms of $[\phi^A(\gamma),\phi^{AD}(\zeta)] \in T_{{\rm III}, A, AD}$.

\subsection{Matching $T_{2,D}$, Biquadratic Case}
Let $D$ be an element in $F^\times - A{F^\times}^2$.
Then, $E = \EA(\sqrt{D})$ is a Galois, biquadratic extension of $F$.
Suppose $\Gal(E/F) = \langle \sigma, \tau \rangle$, such that ${\ED} = E^\sigma$, 
$E_{A} = E^\tau$, and ${\EAD} = E^{\sigma\tau}$.
Then,
$\Gal({\ED} / F) = \la \tau \ra = \la \sigma\tau\ra$, $\Gal({\EA}/F) = \la \sigma \ra = \la \tau\sigma\ra$, and
$\Gal({\EAD}/F) = \la \sigma \ra = \la \tau \ra$.

\begin{claim}
For $c \in F^\times$ and $l \in E^\times$,
the element $\delta_2 = (\phi^{\EA:D}(l),c)$ in $T_{2, D}$ is a norm
of \mbox{$\delta = [\phi^{F:D}(c\zeta'),\phi^{F:AD}(c\zeta)] \in T_{{\rm III},D,AD}$,}
where $\zeta' = \N_{E/\ED}l = l\sigma l$ and $\zeta = \N_{E/\EAD}l = l\sigma\tau l$.
\end{claim}
\begin{proof}
Recall that $\mb{H_2}(\bar{F})$ a quotient of $\GL(2, \bar{F})\times\GL(2, \bar{F})\times \bar{F}^\times$.
The image of $\delta_2$ in $\mb{H_2}(\bar{F})$ is represented by
$(\phi^{\EA:D}(l),c\,\sigma\phi^{\EA:D}(l), 1)$.

Suppose $l = \alpha + \beta \sqrt{D}$ for some $\alpha, \beta \in \EA$.
Since $D \in F^\times$ and $\sigma\sqrt{D} = \sqrt{D}$, the element $(\phi^{\EA:D}(l),c\,\sigma\phi^{\EA:D}(l), 1)$
is equal to
\[
\lp \lp\lsm  \alpha & \beta D \\ \beta & \alpha\rsm \rp,
c\lp\lsm \sigma\alpha & \sigma\beta D \\ \sigma\beta&\sigma\alpha\rsm\rp, 1\rp \sim
(\diag(l,\tau l),\diag(c\sigma l,c\tau\sigma l), 1) \in \mb{H_2}(\bar{F}),
\]
which maps via $\N_2$ to
\begin{multline*}
\diag(cl\sigma l, cl\tau\sigma l,c\tau l\cdot \sigma l, \tau l\cdot \tau\sigma l)
= \diag(c\zeta', c\zeta, c\tau\zeta, c\tau \zeta') \\\sim \delta 
= [\phi^{F:D}(c \zeta'),\phi^{F:AD}(c \zeta)].
\end{multline*}
Since $\det \phi^{F:D}(c\zeta') = c^2\N_{E/F} l = \det \phi^{F:AD}(c\zeta)$, the element $\delta$ lies in
$T_{\text{III}, D, AD}$.
\end{proof}
If $l \notin \ED^\times$, the regular elements $(\phi^{\EA:D}(l),c), (\sigma\phi^{\EA:D}(l),c) \in T_{2, D}$,
which are not stably conjugate to each other,
are both norms of the (possibly singular) element \dt $[\phi^{F:D}(c \zeta'),\phi^{F:AD}(c \zeta)] \in T_{{\rm III}, D, AD}$.
\subsection{Matching $T_{2,D}$, Non-Biquadratic Case}
Let $D$ be an element in $\EA^\times - {\EA^\times}^2$.
Let $\tilde{E}$ be the Galois closure of $E = F(\sqrt{D})$ over $F$.  Let $\sigma$ be the element 
of order $4$ in $\Gal(\tilde{E}/F)$ defined by
$\sigma\A = -\A$, $\sigma \D = \sqrt{\sigma D}$, $\sigma^2\D = -\D$, $\sigma^3\D = -\sqrt{\sigma D}$.
In particular, $\Gal(\EA/F)$ is generated by the restriction of $\sigma$ to $\EA$.
\begin{claim}
For $c \in F^\times$ and $l \in E^\times$, the element $\delta_2 = (\phi^{\EA:D}(l), c)$ in $T_{2, D}$ is a norm of
$\delta = \phi^{F:A:D}(cl\sigma l) \in T_{{\rm IV}, A, D}$.
\end{claim}
\begin{proof}
The image of $\delta_2$ in $\mb{H_2}(\bar{F})$ is represented by
$(\phi^{\EA:D}(l),c\,\sigma\phi^{\EA:D}(l), 1)$.
Suppose $l = \alpha + \beta \sqrt{D}$ for some $\alpha, \beta \in \EA$.
Then, $\delta_2$ is equal to
\[
 \lp \lp\lsm  \alpha & \beta D \\ \beta & \alpha\rsm \rp,
c\lp\lsm \sigma\alpha & \sigma\beta \sigma D \\ \sigma\beta&\sigma\alpha\rsm\rp,1\rp
\sim (\diag(l,\sigma^2 l),\diag(c\sigma l,c\sigma^3 l),1) \in \mb{H_2}(\bar{F}).
\]
We have:
\begin{multline*}
\N_2 : (\diag(l,\sigma^2 l),\diag(c\sigma l,c\sigma^3 l),1)
\mapsto \diag(cl\sigma l,cl\sigma^3 l,c\sigma^2 l\cdot\sigma l, c\sigma^2 l\cdot \sigma^3 l)\\
\sim \delta = \phi^{F:A:D}(cl\sigma l).
\end{multline*}
Since $\N_{E/\EA}l\sigma l = l\sigma l\cdot\sigma^2l\cdot\sigma^3l \in \EA^\sigma = F$, 
$\delta$ lies in $T_{{\rm IV}, A, D}$.
\end{proof}
If $l \notin \ED^\times$, the elements $(\phi^{\EA:D}(l),c)$ and $(\sigma\phi^{\EA:D}(l),c) \in T_{2, D}$ are not stably conjugate to each
other, but they are both norms of the element $\phi^{F:A:D}(cl\sigma l)$ in $T_{{\rm IV}, A, D}$.

\subsubsection{Summary}
In summary, the norm correspondence between the stable conjugacy classes of the elliptic regular elements is as follows
(bar denotes action of $\sigma$):
\begin{enumerate}
\item Case of $H_1$
\[
\begin{split}
T_{1, A} \ni
\left\{
\lsm
(\phi^A(\gamma),\zeta),\\
(\phi^A(\gamma),\bar{\zeta}),\\
(\phi^A(\zeta),\gamma),\\
(\phi^A(\zeta),\bar{\gamma})
\rsm
\right\}
&\longleftrightarrow
[\phi^A(\gamma),\phi^A(\zeta)] \in T_{\text{I}, A}.\\
T_{1, AD} \ni \left\{ \lsm (\phi^{AD}(\zeta),\gamma), \\ (\phi^{AD}(\zeta),\bar{\gamma})\rsm\right\}
& \longleftrightarrow [\phi^{A}(\gamma),\phi^{AD}(\zeta)] \in T_{\text{III}, A, AD}.
\end{split}
\]
\item{Case of $H_2$}

$D \in F^\times - A{F^\times}^2$:
\[
T_{2, D} \ni 
\left\{\lsm (\phi^{\EA:D}(l),\,c),\\(\sigma\phi^{\EA:D}(l),\,c)\rsm\right\}
\longleftrightarrow [\phi^{D}(cl\sigma l),\phi^{AD}(cl\sigma\tau l)] \in T_{\text{III}, D, AD}.
\]
$D \in \EA^\times - {\EA^\times}^2$:
\[
T_{2, D} \ni 
\left\{\lsm (\phi^{\EA:D}(l),\,c),\\ (\sigma{\phi^{\EA:D}(l)},\,c)\rsm\right\}
\longleftrightarrow \phi^{F:A:D}(cl\sigma l) \in T_{\text{IV}, A, D}.
\]
\end{enumerate}
\indexi{norm correspondence|)}\indexi{element!elliptic|)}\indexi{maximal torus!elliptic|)}%
\section{Comparison of Orbital Integrals}
Recall that $\ve$ is the nontrivial, unramified, quadratic character of $F^\times$, and $\EA$
is the quadratic, unramified extension of $F$.  It corresponds to $\ve$ via local class
field theory.  We fix a uniformizer $\vp$ in the ringer of integers $R$ of $F$.
Since $\EA/F$ is unramified, $\vp$ is also a uniformizer for $\EA$.

Let $K = \mb{G}(R)$,  $K_{i} = \mb{H}_i(R)$ ($i = 1, 2$).

For $\mb{H} = \mb{G}$, $\mb{H_1}$, or $\mb{H_2}$, fix a pair $(\mb{B}, \mb{T})$, 
where $\mb{T}$ is the maximal diagonal torus of $\mb{H}$, and 
$\mb{B}$ is the upper triangular Borel subgroup of $\mb{H}$ containing $\mb{T}$.
Let $\Delta$ be the set of roots of $\mb{H}$ with respect to  $(\mb{B}, \mb{T})$.
We define the {\bf discriminant} $D_H(\delta)$ of an element $\delta \in \mb{H}(F)$ as follows:
\indexi{discriminant}%
\[
D_H(\delta) := \prod_{\alpha \in \Delta}\abs{\alpha(t_\delta) - 1}_{\bar{F}}.
\]
Here, $\abs{\;\cdot\;}_{\bar{F}}$ is the unique  absolute value on $\bar{F}$ extending the normalized $p$-adic absolute value on $F$, and 
$t_\delta$ is an  element in $\mb{T}(\bar{F})$ which is conjugate to $\delta$ in $\mb{H}(\bar{F})$.  The discriminant
$D_H(\delta)$ is independent of the choice of $t_\delta$.  Note that $D_H(\delta) = 0$ if $\delta$ is not regular.

For a $p$-adic field $k$ with ring of integers $R_k$, each element $x \in k^\times$ may be written
as $x = u\,\vp_k^{{\rm ord}_k x}$, where $u$ is a unit in $R_k$, $\vp_k$ is a uniformizer of $k$, and ${\rm ord}_k x$
is an integer.  The function ${\rm ord}_k : k^\times \rightarrow \mbb{Z}$ is independent
of the choice of $\vp_k$.
Let $q$ be the cardinality of the residue field $R / (\vp)$ of $F$.  Put ${\rm ord} := {\rm ord}_F$.
Let $R_{\EA}$ be the ring of integers of $\EA$.  

For elements $g, h$ in a group, put 
\[\displaystyle
\Int(g) h := g h g^{-1}.  
\indexs{I@$\Int$}%
\]
\subsection{Comparison for $T_{1, A}$, $T_{{\rm I}, A}$}\label{sec:compT1ATIA}
\subsubsection{Stable Orbital Integral---$T_{1, A}$}\label{sec:oiT1A}
\begin{claim}\label{claim:oiT1A}
Let $\zeta, \xi$ be elements in $R_{\EA}^\times - R^\times$ such that $\N_{\EA/F}\zeta = \N_{\EA/F}\xi$.
Suppose $\zeta = a + b \sqrt{A}$, where $a, b \in F$.
Let $\delta_1 = (\phi^A(\zeta), \xi) \in T_{1,A}$.
The following equality holds:
\[
SO_{\delta_1}(1_{K_1}) =
\frac{q+1}{q-1}q^{\ord b} - \frac{2}{q-1}.
\]
\end{claim}
\begin{proof}
Let $\GL(2, F)^{\EA}$ be the subgroup  of elements in $\GL(2, F)$ whose determinants  lie in $\N_{\EA/F}\EA^\times$.
Let $K' = \GL(2, R)$.  Since $\EA/F$ is unramified, $K'$ is a maximal compact subgroup in $\GL(2, F)^{\EA}$.
The group $K_1 = \mb{H_1}(R)$ is then equal to
$
\lp K' \times R_{\EA}^\times\rp' :=
\{(g, \gamma) \in K'\times R_{\EA}^\times : \det g = \N_{\EA/F}\gamma\}.
$

For $\rho \in F^\times$, put $\phi^A_\rho(\alpha) := \Int(\diag(1, \rho^{-1}))\phi^A(\alpha)$.
There are two conjugacy classes within the stable conjugacy class of $\delta_1$.
One is represented by $\delta_1^1 := \delta_1$, and the other by  $\delta_1^\vp := \lp\phi^A_\vp(\zeta), \xi\rp$.

For $\rho \in \{1, \vp\}$, put
\[
T_1^\rho := Z_{H_1}(\delta_1^\rho) = 
\{(\phi^A_\rho(\alpha), \beta) : \N_{\EA/F}\alpha = \N_{\EA/F}\beta\}.
\]
Put $T_1 := T_1^1 = T_{1, A}$.  Let $dh$ denote the fixed Haar measure $dh_1$ on $H_1$.

By definition, 
\[
SO_{\delta_1}(1_{K_1}) = 
\sum_{\rho \in \{1, \vp\}}\int_{T_1^\rho\bs H_1}1_{K_1}(h^{-1}\delta_1^\rho h)\,dh.
\]
The quotient $T_1^\rho \bs H_1$ is equal to
$T_A^\rho \bs \GL(2, F)^{\EA}$, where $T_A^\rho := \{\phi_\rho^A(\gamma):\gamma\in \EA^\times\}$.

From \cite[p. 30]{F}, we have the decomposition
\[
\GL(2, F) =
 \dot{\bigcup}_{j \geq 0} 
T_A^\rho \lp\begin{smallmatrix}1&0\\0&\vp^{j - \ord\rho}\end{smallmatrix}\rp K',
\]
where the disjoint union is over all nonnegative integers $j$.
Since $K'$ is a subgroup of \dt $\GL(2, F)^{\EA}$, and $\vp \notin \N_{\EA/F}\EA^\times$, the above decomposition implies that
\begin{equation*}
\GL(2, F)^{\EA} =
\dot{\bigcup}_{j - \ord{\rho} \text{ even } \geq 0} 
T_A^\rho \lp\begin{smallmatrix}1&0\\0&\vp^{j - \ord{\rho}}\end{smallmatrix}\rp K'.
\end{equation*}

Let $\delta_1'$ be the $\GL(2, F)$-component $\phi^A(\zeta)$ of $\delta_1$.
Put ${\delta_1'}^\vp := \phi^A_\vp(\zeta)$.
For any $h = (h', \gamma) \in H_1$,
the element $\Int\lp h^{-1}\rp\delta_1$ lies in  $K_1$ if and only if 
$\Int\lp {h'}^{-1}\rp\delta_1'$ lies in $K'$.  
The stable orbital integral $SO_{\delta_1}(1_{K_1})$ is therefore equal to
\begin{multline*}
\int_{T_A \backslash \GL(2,F)^{\EA}} 1_{K'}(\Int(h'^{-1})\delta_1')\, dh' + 
\int_{T_A^\vp \backslash \GL(2,F)^{\EA}} 1_{K'}(\Int(h'^{-1}){\delta_1'}^\vp)\, dh'\\
= \sum_{j \text{ even } \geq 0} 
\int_{ K' \cap \Int(\diag(1, \vp^{-j}))T_A \bs K'}
1_{K'}\lp\Int\lp k^{-1}\lp\lsm 1 & \\ & \vp^{-j}\rsm\rp\rp{\delta_1'}\rp dk\\
+ \sum_{j \text{ odd } \geq 0} 
\int_{ K' \cap \Int(\diag(1, \vp^{1-j})) T_A^\vp \bs K'}
1_{K'}\lp\Int\lp k^{-1}\lp\lsm 1 & \\ & \vp^{-j}\rsm\rp\rp{\delta_1'}\rp dk.
\end{multline*}
Here, $dk$ is the Haar measure on $K'$ such that $K'$ has volume one.

For a nonnegative integer $j$, let 
\[
R_{\EA}(j) = \{a + b \sqrt{A} : a, b \in F;\,\abs{a} \leq 1,\, \abs{b} \leq \abs{\vp^j} = q^{-j}\}.
\]
Then, the group $K' \cap \Int\lp\lp\lsm 1 & \\ & \vp^{ - j}\rsm\rp \rp T_A$ is isomorphic to
$R_{\EA}(j)^\times$.

Hence, for $\rho \in \{1, \vp\}$ we have:
\[
\left[K' \cap T_A : K' \cap \Int\lp\lp\lsm 1 & \\ & \vp^{\ord\rho - j}\rsm\rp \rp T_A^\rho\right] 
= [R_{\EA}^\times : R_{\EA}(j)^\times].
\]
Consequently,
\begin{equation}\label{eq:oiT1A}
\begin{split}
SO_{\delta_1}(1_{K_1})
=& \sum_{j \text{ even } \geq 0}  [R_{\EA}^\times : R_{\EA}(j)^\times] 
1_{K'}\lp\lp\begin{smallmatrix}a & b A \vp^j\\b\vp^{-j}&a\end{smallmatrix}\rp\rp\\
&+ \sum_{j \text{ odd } \geq 0} [R_{\EA}^\times : R_{\EA}(j)^\times] 
1_{K'}\lp\lp\begin{smallmatrix}a & b A \vp^j\\b\vp^{-j}&a\end{smallmatrix}\rp\rp\\
=& \sum_{j \geq 0} [R_{\EA}^\times : R_{\EA}(j)^\times]
1_{K'}\lp\lp\begin{smallmatrix}a & b A \vp^j\\b\vp^{-j}&a\end{smallmatrix}\rp\rp.
\end{split}
\end{equation}

Since $\EA/F$ is unramified, the cardinality $q_{\EA}$ of the residue field $R_{\EA} / (\vp)$ is $q^2$.
The index $[R_{\EA}^\times : R_{\EA}(j)^\times]$ is therefore equal to
\begin{multline*}
[R_{\EA}^\times : 1 +\vp^j R_{\EA}]\;/\;[R_{\EA}(j)^\times : 1 + \vp^j R_{\EA}] \\
= \frac{(q_{\EA} - 1)q_{\EA}^{j - 1}}{(q - 1)q^{j - 1}} =
\begin{cases}
1 & \text{if } j = 0,\\
(q+1)q^{j-1} & \text{if } j > 0.
\end{cases}
\end{multline*}

Since $1_{K'}\lp\lp\begin{smallmatrix}a & b A \vp^j\\b\vp^{-j}&a\end{smallmatrix}\rp\rp$ is zero unless $j \leq \ord{b}$,
we conclude that
\[
SO_{\delta_1}(1_K) = 1 + \frac{q + 1}{q}\sum_{j = 1}^{\ord b}q^j
= \frac{q+1}{q-1}q^{\ord b} - \frac{2}{q-1}.
\]
\end{proof}
\subsubsection{Twisted Orbital Integral---$T_{\text{I}, A}$}\label{sec:oiTIA}
\quad

{\bf Notation:}
\begin{itemize}
\item
$\displaystyle
\left[\lp\lsm a&b\\c&d\rsm\rp, \lp\lsm x&y\\z&t\rsm\rp\right]:=
\lp\lsm a&&&b\\&x&y&\\&z&t&\\c&&&d \rsm\rp.
$
\item
$C_0 := \{[g_1, g_2] : g_1, g_2 \in \GL(2, F);\, \det g_1 = \det g_2\} \subset G$.
\indexs{C@$C_0$, $C_0'$}%
\item
$C_0' = (\GL(2, F) \times \GL(2, F))'$
\[:= \left\{ (g_1, g_2) \in \GL(2, F)\times\GL(2, F) : \det g_1 = \det g_2\right\}.\]
\item
$T_A := \{\phi^A(\xi) : \xi \in \EA^\times\}$.  
\item
For $\rho \in F^\times$, put $T_A^\rho := \Int(\diag(1, \rho^{-1}))T_A,$ and
\[
T_\rho' := (T_A \times T_A^\rho)' := \{(t_1, t_2) \in T_A \times T_A^\rho : \det t_1 = \det t_2\}
\subset C_0'.\]
\item
For $\rho \in F^\times$, put $\phi^A_\rho(\gamma) := \Int(\diag(1, \rho^{-1}))\phi^A(\gamma)$.
Put
\[
T_\rho := T_{\text{I}, A}^\rho :=
\left\{\left[\phi^A(\alpha), \phi^A_\rho(\beta)\right]: \N_{\EA/F} \alpha = \N_{\EA/F} \beta \right\} \subset G.\]
\item
Put $K_0 := (\GL(2, R) \times \GL(2, R))'$, where the prime denotes equal determinants.
\item
For an integer $m$, put
$K_m := \{(g, h) \in C_0' : g - ehe \in \vp^m M_2(R)\}$.  Here, $M_2(R)$ denotes
\indexs{K@$K_m$}%
the set of all $2\times 2$ matrices with coefficients in $R$, and $e := \lp\lsm 1 &\\&-1\rsm\rp \in \GL(2, F)$.
\end{itemize}

%
Let $\xi, \zeta$ be elements in $R_{\EA}^\times - R^\times$ such that $\xi \neq \zeta$ or $\sigma\zeta$, and 
$\N_{\EA/F}\xi = \N_{\EA/F}\zeta$.
Let $\delta$ be the regular element $[\phi^A(\xi), \phi^A(\zeta)]$ in $T_{{\rm I}, A}$.
Let $\delta_\vp$ be the element $[\phi^A(\xi), \phi^A_\vp(\zeta)]$.  It is stably conjugate but not conjugate to $\delta$.
We assume that $\delta$ is topologically unipotent.

Suppose $\xi = a_1 + b_1 \sqrt{A}$, $\zeta = a_2 + b_2\sqrt{A}$, where 
$a_i, b_i$ ($i = 1, 2$) are elements in $F$, with $b_i \neq 0$.  Let  $N_i = \ord b_i$.  Let $X = \ord  (a_1 - a_2)$.
\begin{claim}\label{claim:oiTIA}
The sum $O_{\delta}(1_K) + O_{\delta_\vp}(1_K)$ of $\ve$-twisted orbital integrals is equal to
\[
(-q)^{N_1}(-q)^X\left[\frac{(q+1)q^{N_2} - 2}{q-1}\right].
\]
\end{claim}
\begin{proof}


We make use of the following  decomposition of $G$ (\cite[Lemma I. J. 1]{F}, see also \cite{WE}):
\[
\GSp(2, F) = \dot{\bigcup}_{m\geq 0} C_0 z(m) K.
\] 
Here, the disjoint union is over nonnegative integers $m$, and
\[
z(m) := \lp\lmx 1 & 0 & \vp^{-m} &0 \\ 0 &1&0&\vp^{-m}\\0&0&1&0\\0&0&0&1\rmx\rp.
\]

Put $K_m^{C_0} := C_0 \cap z(m) K z(m)^{-1}$.  For $\rho \in \{1, \vp\}$, put
$\delta_\rho := [\phi^A(\xi), \phi^A_\rho(\zeta)]$.
we have:
\begin{multline*}
O_{\delta_\rho}(1_K) =
\int_{T_\rho \bs G} 1_K(g^{-1} \delta_\rho g) \ve(g)\,dg\\
\begin{split}
& = \sum_{m \geq 0}\abs{K} \int_{T_\rho \bs C_0 / C_0 \cap z(m)Kz(m)^{-1}
}1_{K}(z(m)^{-1}g^{-1}\delta_\rho gz(m))\ve(gz(m))\,dg\\
& = 
\sum_{m \geq 0} \int_{T_\rho \bs C_0 / C_0 \cap z(m)Kz(m)^{-1}
}1_{K}(z(m)^{-1}g^{-1}\delta_\rho gz(m))\ve(g)\,dg\\
& =
\sum_{m \geq 0} \int_{T_\rho \bs C_0 / K_m^{C_0}
}1_{K_m^{C_0}}(g^{-1}\delta_\rho g)\ve(g)\,dg.
\end{split}
\end{multline*}
Note that $\ve$ is trivial on $K$ because it is unramified.  Also, $\abs{K} = 1$ by our choice of measure.

Let 
\[
C_0' =  \left\{ (g_1, g_2) \in \GL(2, F)\times\GL(2, F) :
\det g_1 = \det g_2\right\}.
\]
Let $e =\lp\lsm 1&\\&-1\rsm\rp$, $w = \lp\lsm &1\\1&\rsm\rp \in \GL(2, F)$.  For an integer $m$,
define an isomorphism $\phi_m : C_0' \rightarrow C_0$ as follows:
\[
\phi_m((g_1, g_2)) = \dmfour{1}{1}{\vp^m}{\vp^m}[g_1, ew\, g_2\, we]\dmfour{1}{1}{\vp^m}{\vp^m}^{-1}.
\]

For $\rho \in \{1, \vp\}$, let  $\delta_\rho' = (\phi^A(\xi), \phi^A_\rho(\zeta)) \in C_0'$.  
Let ${T_{\rho}'} = Z_{C_0'}(\delta_{\rho}')$.
For an integer $m \geq 0$, let $\delta_{\rho, m} = \phi_m^{-1}(\delta_\rho) \in C_0'$.
Let ${T_{\rho, m}} = Z_{C_0'}(\delta_{\rho,m})$.


Let
\[
K_m = \{(g_1, g_2) \in C_0':g_1 - e\,g_2\,e \in \vp^m M_2(R)\},
\]
where $M_2(R)$ is the set of $2 \times 2$ matrices with coefficients in $R$.
Then, $\phi_m$ maps $K_m$ isomorphically onto $K_m^{C_0}$ (\cite[p. 38]{F}).

Let $K_0 = (\GL(2, R) \times \GL(2, R))'$
\[:= \{(g_1, g_2) \in \GL(2, R) \times \GL(2, R) : \det g_1 = \det g_2\}.\]
Then, $[K  : K_m^{C_0}] = [K_0 : K_m]$.  We have:
\begin{multline}\label{eq:oiTIA0}
O_{\delta_\rho}(1_K)\\
=\sum_{m\geq 0}[K_0 : K_m] \int_{{T_{\rho, m}}\bs C_0'} 
1_{\phi_m(K_m)}\lp\phi_m(g')^{-1}\delta_{\rho}\phi_m(g')\rp\ve(\phi_m(g'))\,dg'.
\end{multline}

By abuse of notation, for $g'=(g_1, g_2) \in C_0'$,
let $\ve(g')$ denote the value of $\ve(\phi_m(g'))$.  It is equal
to $\ve(\det g_1) = \ve(\det g_2)$ and is independent of $m$.

We change variables $g'\mapsto \lp\vpmm, \vpmm we\rp  g'$.  
Since
\[
\ve\lp\lp\vpmm, \vpmm we\rp\rp = \ve(\det \vpmm) = \ve(\vp)^m = (-1)^m, 
\]
the expression \eqref{eq:oiTIA0} is equal to
\begin{multline}\label{eq:oiTIA1}
\sum_{m\geq 0}[K_0:K_m]\int_{{T_\rho}' \bs C_0'} 1_{K_m}({g'}^{-1} \delta_\rho'\, g')
\ve\Bigl(\lp\vpmm,\vpmm we\rp g'\Bigr) dg'\\
= \sum_{m\geq 0}(-1)^m [K_0:K_m]\int_{{T_\rho}' \bs C_0'} 1_{K_m}({g'}^{-1} \delta_\rho'\, g')\ve(g')\, dg'.
\end{multline}

For $\rho \in \{1, \vp\}$, we have the following decomposition of $C_0'$ (\cite[p. 45]{F}):
\[
C_0' = \dot{\bigcup}_{r \in \RR_\rho} T_\rho' \cdot r \cdot K_0\,,
\]
where
\[
\RR_\rho :=
\left\{
r_{j_1, j_2, \epsilon} = \Big(\vp^{-\floor{j_1/2}} \lp\begin{smallmatrix}1&0\\0&\vp^{j_1}\end{smallmatrix}\rp
,\,
\phi^{A}_\rho(\epsilon') \vp^{-\floor{(j_2-\ord{\rho})/2}}
\lp\begin{smallmatrix}1&0\\0&\epsilon\vp^{j_2-\ord{\rho}}\end{smallmatrix}\rp\Big)
\right\}.
\]
Here, (i) $j_1, j_2$ are nonnegative integers, with $j_1 - j_2 + \ord{\rho}$ even; (ii)
$\epsilon$ is a representative in $R^\times$ of a class in $R^\times / {R^\times}^2$;
(iii) $\epsilon'$ is an element of $R_{\EA}^\times$ such that $\N_{\EA/F}\epsilon' = \epsilon$.

Let $R_{T_\rho'} = K_0 \cap T_{\rho}'$.  The sum \eqref{eq:oiTIA1} is then equal to
\begin{multline*}
\sum_{m \geq 0} (-1)^m [K_0 : K_m]
\sum_{r \in \RR_\rho} [R_{T_\rho'} : T_\rho' \cap r K_0 r^{-1}]\ve(r)
\int_{K_0} 1_{K_m}(k^{-1} r^{-1} \delta_\rho' rk)dk.
\end{multline*}

For $r = r_{j_1, j_2, \epsilon} \in \RR_\rho$, we have
\[\begin{split}
\ve(r) &= 
\ve\lp\Big(\vp^{-\floor{j_1/2}} \lp\begin{smallmatrix}1&0\\0&\vp^{j_1}\end{smallmatrix}\rp
\;,\;
\phi^{A}_\rho(\epsilon') \vp^{-\floor{(j_2-\ord{\rho})/2}}
\lp\begin{smallmatrix}1&0\\0&\epsilon\vp^{j_2-\ord{\rho}}\end{smallmatrix}\rp\Big)\rp\\
& = \ve\lp\det\lp\vp^{-\floor{j_1/2}} \lp\begin{smallmatrix}1&0\\0&\vp^{j_1}\end{smallmatrix}\rp\rp\rp
=(-1)^{j_1}.
\end{split}\]
Thus, $O_{\delta_\rho}(1_K)$ is equal to
\begin{multline}\label{eq:oiTIAR}
\sum_{m \geq 0} (-1)^m [K_0 : K_m]
\sum_{\substack{j_1, j_2 \geq 0;\\ j_1 - j_2 + \ord{\rho} \text{ even}}} 
(-1)^{j_1} [R_{T_\rho'} : T_\rho' \cap r K_0 r^{-1}]\\\cdot \int_{K_0} 1_{K_m}(k^{-1} r^{-1} \delta_\rho' rk) dk.
\end{multline}

Hence, the sum describing the $\ve$-twisted orbital integral $O_{\delta_\rho}(1_K)$ is quite similar to 
that describing its nontwisted counterpart in \cite{F}.  
They are identical but for the factors $(-1)^m$, $(-1)^{j_1}$ which appear above in the $\ve$-twisted case.

The quantities $[R_{T_\rho'} : T_\rho' \cap r K_0 r^{-1}]$ and $[K_0 : K_m]
\int_{K_0} 1_{K_m}(k^{-1} r^{-1} \delta_\rho' rk)dk$ are computed in \cite[Sect. II. C]{F}.
Recall that $\xi = a_1 + b_1 \sqrt{A}$, $\zeta = a_2 + b_2\sqrt{A}$, and
$N_i = \ord b_i$ ($i = 1, 2$), $X = \ord(a_1 - a_2)$.
We follow the notation in \cite{F} and let $N = \min(N_1, N_2)$, $\nu_i = N_i - j_i$.
Let $\delta(m = 0)$ be $0$ if $m \neq 0$, and $1$ if $m \neq 0$.  Let
$\delta(m \geq 1)$ be $0$ if $m < 1$, and $1$ if $m \geq 1$.
It follows from \cite[Sect. II. C]{F} that 
$O_{\delta}(1_K) +  O_{\delta_\vp}(1_K)$ is the sum of the following three quantities:
\begin{itemize}
\item[I.]
\begin{multline}\label{eq:tIA1}
\sum_{m = 0}^N (-1)^m\left(\delta(m = 0) + \delta(m \geq 0)(1 - q^{-2})q^{3m}\rp\\
\begin{split}
& \left[ \sum_{\nu_1 = m}^{N_1 - 1} \sum_{\nu_2 = m}^{N_2 - 1} \left( \frac{q + 1}{q} \right)^2 \left( - q \right)^{N_1 - \nu_1} q^{N_2 - \nu_2} \right.\\
&\left. + \sum_{\nu_1 = m}^{N_1 - 1} \frac{q + 1}{q}
\left( - q \right)^{N_1 - \nu_1} + \sum_{\nu_2 = m}^{N_2 - 1} \frac{q + 1}{q}
q^{N_2 - \nu_2}  + 1\right],
\end{split}
\end{multline}
\item[II.]
\begin{equation}\label{eq:tIA2}
\sum_{0 = \nu}^{N - 1} \sum_{m = \nu + 1}^{2 N - \nu} (-1)^m\frac{1}{2} \left(
\frac{q + 1}{q} \right)^2 (-q)^{N_1 - \nu} q^{N_2 - \nu} 2 q^{m + 2 \nu},
\end{equation}
\item[III.]
\begin{multline}\label{eq:tIA3}
\sum_{m = N + 1}^{X - N} (-1)^m\frac{q + 1}{q} q^{2 N + m} + \\
\sum_{\nu = 0}^{N - 1} \sum_{m = 2 N - \nu + 1}^{X - \nu} (-1)^m \frac{1}{2} \left( \frac{q + 1}{q}
\right)^2 (-q)^{N - \nu}q^{N - \nu} 2 q^{m + 2 \nu}.
\end{multline}
\end{itemize}
The sum I + II + III has the same form as the sum on \cite[p. 60]{F} 
(even though they describe orbital integrals for different groups).  
There, it is computed to be
\[
(-q)^{N_1} (-q)^X\left[\frac{((q + 1)q^{N_2} -2}{q-1}\right].
\]
The proof of Claim \ref{claim:oiTIA} is now complete.
\end{proof}
\subsubsection{Transfer Factor---$T_{1,A}, T_{{\rm I}, A}$}
Let $\xi, \zeta$ be distinct elements in $\EA^\times - F^\times$ such that $\xi \neq \zeta, \sigma{\zeta}$.
Let $\delta$ be the element $[\phi^A(\xi), \phi^A(\zeta)]$ in $T_{{\rm I},A} \subset G$.
Let $\delta_1 = (\phi^A(\xi),\zeta) \in T_{1, A} \subset H_1$.
Then, $\delta$ is regular and $\delta_1$ is a norm of $\delta$.

There are two conjugacy classes in the stable conjugacy class of $\delta$.
One is represented by $\delta$, and the other by $\delta_\vp := [\phi^A(\xi), \phi_\vp^A(\zeta)]$. 
Since we are proving the matching of orbital integrals for the unit elements in the Hecke algebras,
we may assume without loss of generality that $\delta_1 \in K_1$ and $\delta \in K$, or in other words
$\xi, \zeta \in R_{\EA}^\times - R^\times$.

Extend $\ve$ to $\EA^\times$ by letting 
$\ve(x) = (-1)^{\text{ord}_{\EA}x}$ for all $x \in \EA^\times$.
We define the transfer factors for the pairs $(\delta_1, \delta)$ and $(\delta_1, \delta_\vp)$ as follows:
\[
\Delta(\delta_1, \delta) = \Delta(\delta_1, \delta_\vp) = 
\ve
\lp\!
\Big((\xi + \sigma \xi) - (\zeta + \sigma \zeta)\Big)
\lp \frac{\xi - \sigma\xi}{2\sqrt{A}}\rp
\!\rp
\left(\frac{D_G(\delta)}{D_{H_1}(\delta_1)}\right)^{1/2}.
\indexi{transfer factor}%
\]
We do not explain here why the above expression constitutes a transfer factor.  For those interested
in further details, \cite{LS} and \cite{H} are both useful sources on transfer factors.

Suppose $\xi = a_1 + b_1 \sqrt{A}$, $\zeta = a_2 + b_2 \sqrt{A}$, where $a_i, b_i \in F$ ($i = 1, 2$).
Let $N_i = \ord b_i$.  Let $X = \ord (a_1 - a_2)$.
\begin{claim}\label{claim:deltaT1ATIA}
The following holds:
\[
\Delta(\delta_1, \delta) = \Delta(\delta_1, \delta_\vp) = (-q)^{-N_1 - X}.
\]
\end{claim}
\begin{proof}
The proof is a straightforward computation, and we skip it.
\end{proof}
\begin{prop}\label{prop:oiT1ATIA}
The following identity holds:
\[
SO_{\delta_1}(1_{K_1}) =  \Delta(\delta_1, \delta)O_{\delta}(1_{K}) +  \Delta(\delta_1, \delta_\vp)O_{\delta_\vp}(1_{K}).
\]
\end{prop}
\begin{proof}
This follows from Claims \ref{claim:oiT1A}, \ref{claim:oiTIA}, and \ref{claim:deltaT1ATIA}.
\end{proof}

Let $\tilde{\delta}_1 = (\phi^A(\xi), \zeta) \in T_{1, A}$, i.e. we swap $\xi$ and $\zeta$ in $\delta_1$.
Then, $\tilde{\delta}_1$ is not stably conjugate to $\delta_1$, but it is also a norm of $\delta = [\phi^A(\xi), \phi^A(\zeta)]$.
We define the transfer factors for the pairs $(\tilde{\delta}_1, \delta)$ and $(\tilde{\delta}_1, \delta_\vp)$ as follows:
\begin{align*}
\Delta(\tilde{\delta}_1, \delta) &= \Delta(\delta_1,\delta), &
\Delta(\tilde{\delta}_1, \delta_\vp) &= - \Delta(\delta_1, \delta).
\end{align*}
\begin{prop}\label{prop:oiT1ATIAus}
The following holds:
\[
SO_{\tilde{\delta}_1}(1_{K_1}) = \sum_{\rho \in \{1, \vp\}}\Delta(\tilde{\delta}_1,\delta_\rho)O_{\delta_\rho}(1_K).
\]
\end{prop}
\begin{proof}

Put $O^{{\rm us}}_\delta(1_K) := O_\delta(1_K) - O_{\delta_\vp}(1_K)$.
\indexs{O@$O^{\rm us}(\;)$}%
Thus,
\[
\sum_{\rho \in \{1, \vp\}}\Delta(\tilde{\delta}_1, \delta_\rho)O_{\delta_\rho}(1_K) = 
\Delta(\tilde{\delta}_1, \delta) O^{\rm us}_{\delta}(1_K).
\]

It follows from \eqref{eq:oiTIAR} that $O^{\rm us}_{\delta}(1_K)$ is equal to
\begin{multline*}
= \sum_{\substack{r_{j_1, j_2, \epsilon};\\j_1 + j_2 \text{ even}}} \ve(r) [R_{T'} : T' \cap r K_0 r^{-1}]  \sum_{m \geq 0} (-1)^m [K_0 : K_m]\\
\shoveright{\cdot\int_{K_0} 1_{K_m}(k^{-1} r^{-1} \delta' rk)\,dk}\\
\shoveleft{- \sum_{\substack{r_{j_1, j_2,\ve};\\j_1 + j_2 \text{ odd}}}\ve(r)
[R_{T_\vp'} : T_\vp' \cap r K_0 r^{-1}] \sum_{m \geq 0} (-1)^m [K_0 : K_m]}\\
\shoveright{\cdot \int_{K_0} 1_{K_m}(k^{-1} r^{-1} \delta_\vp' rk)\,dk}\\
= \sum_{\rho \in \{1, \vp\}}\!\!
\sum_{\substack{r_{j_1, j_2, \epsilon};\\j_1 + j_2 +\ord{\rho}\text{ even}}} (-1)^{j_1 + j_2}\ve(r)
[R_{T_\rho'} : T_\rho' \cap r K_0 r^{-1}]  
\sum_{m \geq 0} (-1)^m [K_0 : K_m]\\
\cdot \int_{K_0} 1_{K_m}(k^{-1} r^{-1} \delta' rk)\,dk.
\end{multline*}

Hence, the quantity $O^{\rm us}_\delta(1_K)$ is obtained from $\sum_{\rho\in\{1, \vp\}}O_{\delta_\rho}(1_K)$ on
multiplying each summand in the expressions  \eqref{eq:tIA1} through \eqref{eq:tIA3} by 
\[(-1)^{j_1 + j_2} = (-1)^{N_1 - \nu_1 + N_2 - \nu_2}.\]
In other words, $O^{\rm us}_\delta(1_K)$ is equal to the sum of the following three quantities:
\begin{itemize}
\item[I.]
\begin{multline}
\sum_{m = 0}^N (-1)^m\left(\delta(m = 0) + \delta(m \geq 0)(1 - q^{-2})q^{3m}\rp\\
\begin{split}
& \left[\sum_{\nu_1 = m}^{N_1 - 1} \sum_{\nu_2 = m}^{N_2 - 1} \left( \frac{q + 1}{q} \right)^2 
q^{N_1 - \nu_1} (-q)^{N_2 - \nu_2} \right.\\
& \left.+ \sum_{\nu_1 = m}^{N_1 - 1} \frac{q + 1}{q}
q^{N_1 - \nu_1} + \sum_{\nu_2 = m}^{N_2 - 1} \frac{q + 1}{q}
(-q)^{N_2 - \nu_2}  + 1\right],
\end{split}
\end{multline}
\item[II.]
\begin{equation}
\sum_{0 = \nu}^{N - 1} \sum_{m = \nu + 1}^{2 N - \nu} (-1)^m\frac{1}{2} \left(
\frac{q + 1}{q} \right)^2 q^{N_1 - \nu} (-q)^{N_2 - \nu} 2 q^{m + 2 \nu},
\end{equation}
\item[III.]
\begin{multline}
\sum_{m = N + 1}^{X - N} (-1)^m\frac{q + 1}{q} q^{2 N + m} + \\
\sum_{\nu = 0}^{N - 1} \sum_{m = 2 N - \nu + 1}^{X - \nu} (-1)^m \frac{1}{2} \left( \frac{q + 1}{q}
\right)^2 q^{N - \nu}(-q)^{N - \nu} 2 q^{m + 2 \nu}.
\end{multline}
\end{itemize}
The sum of the above quantities is the same as $\sum_{\rho \in \{1, \vp\}}O_{\tilde{\delta}_\rho}(1_K)$, where
\[
\tilde{\delta} = [\phi^A(\zeta),\phi^A(\xi)],
\] 
i.e., $\delta$ with the two $\GL(2, F)$-components swapped.

By Claims \ref{claim:oiTIA} and \ref{claim:deltaT1ATIA}, we have
\[
O^{\rm us}_{\delta}(1_K) =
\sum_{\rho \in \{1, \vp\}} O_{\tilde{\delta}_\rho}(1_K) = 
(-q)^{N_2 + X}\left[\frac{((q + 1)q^{N_1} -2}{q-1}\right],
\]
\[
\Delta(\tilde{\delta}_1, \delta) = (-q)^{-N_2 - X}.
\]
By Claim \ref{claim:oiT1A}, the stable orbital integral $SO_{\tilde{\delta}_1}(1_{K_1})$ is equal to
\[
\frac{q+1}{q-1}q^{N_1} - \frac{2}{q-1}.
\]
The proposition follows.
\end{proof}
{\sc Remark:}  Since $\delta = [\phi^A(\xi), \phi^A(\zeta)]$ is conjugate to $\tilde{\delta} = [\phi^A(\zeta), \phi^A(\xi)]$, 
it seems contradictory that $\sum_{\rho \in \{1, \vp\}}O_{\delta_\rho}(1_K)$ should differ from 
$\sum_{\rho \in \{1, \vp\}}O_{\tilde{\delta}_\rho}(1_K)$.  However, note that we have chosen $[\phi^A(\xi), \phi^A_\vp(\zeta)]$
to be the representative of the conjugacy class, different from that of $\delta$, in the stable conjugacy class of $\delta$.
In other words, the choice of the representative involves a bias towards the second component of $[\phi^A(\xi), \phi^A(\zeta)]$.  
More concretely, there exists an element $g \in G$ such that $\tilde{\delta} = g^{-1}\delta g$, 
but $[\phi^A(\zeta), \phi^A_\vp(\xi)]$ is not conjugate to $[\phi^A(\xi), \phi^A_\vp(\zeta)]$ via the same $g$.
\subsection{Comparison for $T_{1, AD}$, $T_{{\rm III}, A, AD}$}
\label{sec:compT1ADTIIIAAD}
Let $D$ be an element in $F^\times - A{F^\times}^2\cup {F^\times}^2$.  Note that since by assumption 
$\EA/F$ is unramified, $A$ lies in $R^\times - {R^\times}^2$.
Let $E = \EA(\sqrt{D}) = F(\sqrt{A}, \sqrt{D})$.  
Suppose $\Gal(E/F) =  \langle \sigma, \tau \rangle$, such that
$\EA = E^\tau$, $\EAD = E^\sigma$, and $\Gal(\EA/F) = \langle \sigma \rangle$, $\Gal(\EAD/F) = \langle \tau \rangle$.  
Since by assumption the residual characteristic of $F$ is odd, the quadratic extensions $\ED/F$, $\EAD/F$ are ramified.

For a $p$-adic field $k$, let $R_k$ denote the ring of integers in $k$.
Let $\vp$ be a uniformizer of $F$.

Since $\EAD/F$ and $\ED/F$ are ramified, without loss of generality we assume that $AD = {-\vp}$, $D = {-A^{-1}\vp}$. 
We let $\vp_{AD} = \sqrt{-\vp}$ be the uniformizer of $\EAD$, thus $\N_{\EAD/F}\,\vp_{AD} = \vp$.
(We may also let $\sqrt{-A^{-1}\vp}$ be the uniformizer of $\ED$, but we do not use it in this section).

Let $\xi$ be an element in $R_{\EA}^\times$, $\zeta$ an element in $R_{\EAD}^\times$, such that $\xi, \zeta \notin R^\times$,
and $\N_{\EA/F}\xi = \N_{\EAD/F}\zeta$.
Let
\[\begin{split}
\delta_1 &= (\phi^{AD}(\zeta),\xi) \in T_{1, AD} \cap K_1,\\
\delta &= [\phi^{A}(\xi), \phi^{AD}(\zeta)] \in T_{{\rm III}, A, AD} \cap K.
\end{split}\]
Then, $\delta$ is regular, and $\delta_1$ is a norm of $\delta$.
Recall from Section \ref{sec:normmappingell} that the respective conjugacy classes of
$\delta_1, \delta$ are stable.  We assume that $\delta$ is topologically unipotent.

Suppose $\xi = a_1 + b_1 \sqrt{A}$, $\zeta = a_2 + b_2 \sqrt{AD}$, where
$a_i, b_i \in F$ ($i = 1, 2$).  Let $N_i = \ord b_i$.  Let $X = \ord (a_1 - a_2)$.
\subsubsection{Stable Orbital Integral---$T_{1, AD}$}
\begin{claim}\label{claim:oiT1AD}
The stable orbital integral $O_{\delta_1}(1_{K_1})$ is equal to
\[
\frac{q^{N_2 + 1} - 1}{q-1}.
\]
\end{claim}
\begin{proof}
Recall that $\GL(2, F)^{\EA}$ denotes $\{g \in \GL(2, F):\det g \in \N_{\EA/F}\EA^\times\}$.
Since the field extension $\EA/F$ is unramified,
$K' = \GL(2, R)$ is a maximal compact subgroup in $\GL(2, F)^{\EA}$.

Let
\[
T_{AD}^{\EA} = \{ \phi^{AD}(\gamma): \gamma \in  \EAD^\times;\, \N_{\EAD/F}\gamma \in \N_{\EA/F}\EA^\times\} 
\subset \GL(2, F)^{\EA}.
\]
For a nonnegative  integer $j$, put $r_j = \phi^{AD}(\vp_{AD}^{-j})\lp\lsm 1&\\&\vp^j\rsm\rp$.
Lemma I. I. 1 of \cite{F} implies that $\GL(2, F)^{\EA}$ decomposes as follows:
\begin{equation}\label{eq:decompGL2TAD}
\GL(2, F)^{\EA} = 
\dot{\bigcup}_{j\geq 0} T_{AD}^{\EA} \cdot  r_j \cdot K'.
\end{equation}

Let
\[
R_{\EAD}(j) = \{a + b \sqrt{AD} : a, b \in F;\,\abs{a} \leq 1,\, \abs{b} \leq \abs{\vp^j} = q^{-j}\}.
\]
We have $K' \cap r_j^{-1} T_{AD}^{\EA} r_j  = K' \cap r_j^{-1} T_{AD} r_j \cong R_{\EAD}(j)^\times$, and
$[R_{\EAD}^\times : R_{\EAD}(j)^\times]$ is equal to $q^j$ for all nonnegative integers $j$ (\cite[Sect. I. I.]{F}).
Hence,
\[\begin{split}
O_{\delta_1}(1_{K_1}) &= \int_{T_{1, AD} \bs H_1} 1_{K_1}(h^{-1} \delta_1 h)\, dh\\
&= \int_{T_{AD}^{\EA} \bs \GL(2, F)^{\EA}} 1_{K'}(h^{-1} \phi^{AD}(\xi) h)\,dh\\
&= \sum_{j \geq 0}
\int_{ K' \cap\, \Int(r_j^{-1})T_{AD}^{\EA} \bs K'} 1_{K'}\lp\Int\!\lp k^{-1} r_j^{-1}\rp{\delta_1'}\rp dk\\
&= \sum_{j \geq 0}[R_{\EAD}^\times : R_{\EAD}(j)^\times] 
1_{K'}\!\lp \lp\lsm 1&\\&\vp^j\rsm\rp^{\!-1}\phi^{AD}(\xi) \lp\lsm1&\\&\vp^j\rsm\rp\rp\\
&= \sum_{j \geq  0} q^j \cdot 1_{K'}\!\lp\lp\lsm a&b AD \vp^j\\ b\vp^{-j}&a \rsm\rp\rp
= \sum_{j = 0}^{N_2} q^j = \frac{q^{N_2 + 1} - 1}{q-1}.
\end{split}\]
\end{proof}
\subsubsection{Twisted Orbital Integral---$T_{\text{III}, A, AD}$}\label{claim:oiTIIIAAD}
Recall: $\xi = a_1 + b_1 \sqrt{A} \in R_{\EA}^\times$, $\zeta = a_2 + b_2 \sqrt{AD} \in R_{\EAD}^\times$,
$N_i = \ord b_i$ ($i = 1, 2$), $X = \ord (a_1 - a_2)$.  Recall that $\delta = [\phi^A(\xi), \phi^{AD}(\zeta)] \in T_{{\rm III}, A, AD}$.
\begin{claim}
The $\ve$-twisted orbital integral $O_{\delta}(1_K)$ is equal to
\[
(-q)^{N_1 + X} \lp\frac{q^{N_2 + 1} - 1}{q-1}\rp.
\]
\end{claim}
\begin{proof}
As before, $C_0'$ denotes the group $(\GL(2, F)\times\GL(2, F))'$, and 
$K_0$ denotes \dt $(\GL(2, R)\times\GL(2, R))'$, where the prime indicates equal determinants.
Recall that $K_m := \{(g, h) \in C_0' : g - ehe \in \vp^m M_2(R)\}$, $e := \diag(1, -1)$.

Let $\delta' = (\phi^A(\xi),\phi^{AD}(\zeta))$.
Let 
\[ 
T' = Z_{C_0'}(\delta') = 
\{(\phi^A(\alpha), \phi^{AD}(\beta)) : \alpha \in \EA^\times, \beta \in \EAD^{\times};\, \N_{\EA/F}\alpha = \N_{\EAD/F}\beta\}.
\]
Repeating the procedure used to compute the orbital integral in the case of $T_{{\rm I}, A}$, we obtain:
\begin{equation}\label{eq:oiTIIIAAD1}
O_{\delta}(1_{K}) =
\sum_{m\geq 0}(-1)^m[K_0, K_m] \int_{T' \bs C_0'}1_{K_m}(g^{-1} \delta' g)\,dg.
\end{equation}

Recall that  $\vp_{AD}$ is a uniformizer of $\EAD$ such that  $\N_{\EAD/F}\, \vp_{AD} = \vp$.
From \cite[Sect. II. G]{F}, we have the decomposition
\[
C_0' = \dot{\bigcup}_{r \in \mc{R}} T' \cdot r \cdot K_0.
\]
Here, $\RR$ is a set indexed by the couplets of nonnegative integers $j_1, j_2$, and the classes in $R^\times/{R^\times}^2$.
For a representative $\ep \in R^\times$ of a class in $R^\times/{R^\times}^2$,
let $\ep'$ be an element in $\EA^\times$ such that $\N_{\EA/F}\,\ep' = \ep^{-1}$.
For an integer $j$, let $s_{j}$ be $0$ if $j$ is even, $1$ if $j$ is odd.
The set $\RR$ consists of the elements
\[
r_{j_1, j_2, \ep} = \Big(\phi^A(\ep')\vp^{-\floor{j_1/2}} \lp\begin{smallmatrix}1&0\\0&\ep \vp^{j_1}\end{smallmatrix}\rp
,\, \phi^{AD}(\vp_{AD}^{-j_2 + s_{j_1}}) \lp\begin{smallmatrix}1&0\\0&\vp^{j_2}\end{smallmatrix}\rp\!\Big)
\in C_0',
\]

Let $R_{T'} = T' \cap K_0$.  Then, \eqref{eq:oiTIIIAAD1} is equal to
\begin{equation}\label{eq:oiTIIIAAD2}
\sum_{r \in \mc{R}} [R_{T'} : T' \cap r K_0 r^{-1}]\ve(r) \sum_{m \geq 0} (-1)^m [K_0 : K_m]
\int_{K_0} 1_{K_m}(k^{-1} r^{-1} \delta' rk)\,dk.
\end{equation}

The quantities $[R_{T'} : T' \cap r K_0 r^{-1}]$ and $[K_0 : K_m] \int_{K_0} 1_{K_m}(k^{-1} r^{-1} \delta' rk)\,dk$ are computed in 
\cite[p. 64]{F}.  Moreover,
\begin{multline*}
\ve\lp\phi^A(\ep')\vp^{-\floor{j_1/2}} \lp\begin{smallmatrix}1&0\\0&\ep\vp^{j_1}\end{smallmatrix}\rp,\,
\phi^{AD}(\vp_{AD}^{-j_2 + s_{j_1}}) \lp\begin{smallmatrix}1&0\\0&\vp^{j_2}\end{smallmatrix}\rp\rp\\
= \ve\lp\det \lp\vp^{-\floor{j_1/2}} \lp\begin{smallmatrix}1&0\\0&\vp^{j_1}\end{smallmatrix}\rp\rp\rp
= (-1)^{j_1}.
\end{multline*}
It follows that \eqref{eq:oiTIIIAAD2} has the same form as $1/2$ times the sum on page 78 of \cite{F}, 
which is computed there to be 
\[
2\, (-q)^{N_1 + X}\lp \frac{q^{N_2 + 1} - 1}{q-1}\rp.
\]
\end{proof}
\subsubsection{Transfer Factor---$T_{1, AD}$, $T_{{\rm III}, A, AD}$}
Recall: $\xi = a_1 + b_1 \sqrt{A} \in R_{\EA}^\times$, $\zeta = a_2 + b_2 \sqrt{AD} \in R_{\EAD}^\times$,
$N_i = \ord b_i$ ($i = 1, 2$), $X = \ord (a_1 - a_2)$.  Recall that $\delta_1 = (\phi^{AD}(\zeta), \xi) \in T_{1, AD}$
and $\delta = [\phi^A(\xi), \phi^{AD}(\zeta)] \in T_{{\rm III}, A, AD}$.

Recall that $E = F(\sqrt{A}, \sqrt{D})$.
We extend $\ve$ to $E^\times$ by letting
\[
\ve(x) = \lp\sqrt{-1}\rp^{{\rm ord}_E\,x}, \quad \forall x \in E^\times.
\]
We define the transfer factor $\Delta(\delta_1, \delta)$ as follows:
\[
\ve\!\lp\!
\Big((\xi + \sigma\xi) - (\zeta + \tau\zeta)\Big)\cdot
\lp \frac{\xi - \sigma\xi}{2\sqrt{A}}\rp
\!\rp
\lp\frac{D_G(\delta)}{D_{H_1}(\delta_1)}\rp^{1/2}.
\indexi{transfer factor}%
\]
\begin{prop}\label{prop:oiT1ADTIIIAAD}
The following equality holds:
\[
O_{\delta_1}(1_{K_1}) = \Delta(\delta_1, \delta) O_{\delta}(1_K).
\]
\end{prop}
\begin{proof}
This follows from Claims \ref{claim:oiT1AD}, \ref{claim:oiTIIIAAD}, and the fact that
\[
\Delta(\delta_1, \delta) = 
(-q)^{-N_1 - X}.
\]
\end{proof}
\subsection{Comparison for $T_{2, D}$, $T_{{\rm III}, D, AD}$}
Let $D$ be an element in $F^\times - A{F^\times}^2\cup {F^\times}^2$.  The quadratic extensions $\ED/F$, $\EAD/F$ are
ramified, for the residual characteristic of $F$ is odd and $\EA/F$ is unramified.

Let $E$ be the biquadratic extension $F(\sqrt{A}, \sqrt{D})$ of $F$.
The quadratic extensions $E/\ED$, $E/\EAD$ are unramified.
Suppose $\Gal(E/F) = \langle \sigma, \tau\rangle$, such that
$\EA= E^\tau$, $\ED = E^{\sigma}$, and $\EAD = E^{\sigma\tau}$.  
Thus, $\Gal(\EA/F) = \langle\sigma\rangle = \la \sigma\tau\ra$, 
$\Gal(\ED/F) = \langle\tau \rangle = \la \tau\sigma\ra$, and
$\Gal(\EAD/F) = \langle \sigma \rangle = \langle \tau \rangle$.

For a $p$-adic field $k$, let $R_k$ denote the ring of integers in $k$.  Let $R = R_F$.

Let $l$ be an element in $R_E^\times - R^\times$ such that $\N_{E/\ED}l$, $\N_{E/\EAD} l \notin F^\times$.
Let $\xi = \N_{E/\ED} l$, $\zeta = \N_{E/\EAD} l$.
Let $c$ be an element in $R^\times$.

Suppose $l = \alpha + \beta \sqrt{D}$ for some $\alpha, \beta \in \EA$.
Put
\[
\phi^{\EA:D}(l) := \lp\lsm \alpha&\beta D\\ \beta&\alpha\rsm\rp \in \GL(2, \EA).
\]

Recall that $K_2 = \mb{H_2}(R)$.  Let
\[
\delta_2 = \lp \phi^{\EA:D}(l), c\rp_* \in T_{2, D} \cap K_2,
\]
where $(\phi^{\EA:D}(l), c)_*$ denotes the image of $(\phi^{\EA:D}(l), c)$ in $H_2$.
It is a norm of the regular element
\[
\delta = c\left[\phi^D(\xi), \phi^{AD}(\zeta) \right] \in T_{\text{III}, D, AD}\cap K.
\]
The conjugacy classes of $\delta_2$ and $\delta$ are stable.  We assume that $\delta$ is topologically unipotent.

We depart from the previous notation and let $q$ be the cardinality of the residue field of $\EA$.
Let $q_0$ be the cardinality of the residue field of $F$.
Since $\EA/F$ is unramified, $q = q_0^2$.

Let $n = {\rm ord}_{\EA} \beta$.
Suppose $\xi = \N_{E/\ED}l = a_1 + b_1 \sqrt{D}$, and
$\zeta = \N_{E/\EAD}l = a_2 + b_2 \sqrt{AD}$, where $a_i, b_i \in F$ ($i = 1, 2$).
Let $N_i = {\rm ord}_F\, b_i$.  Let $B_i$ be the unit in $R^\times$ such that $b_i = B_i \vp^{N_i}$.
Let $N = \min(N_1, N_2)$.  Let $X = \ord(a_1 - a_2)$.

\subsubsection{Stable Orbital Integral---$T_{2, D}$}
\begin{claim}\label{claim:oiT2D}
The stable orbital integral $O_{\delta_2}(1_{K_2})$ is equal to
\[
\frac{q^{N + 1} - 1}{q-1}.
\]
\end{claim}
\begin{proof}
By \cite[Lemma I. I. 2]{F}, we have:
\[
O_{\delta_2}(1_{K_2}) = \frac{q^{n + 1} - 1}{q-1}.
\]
By \cite[Lemma II. L. 1]{F}, $n$ is equal to $N$.
\end{proof}
\subsubsection{Twisted Orbital Integral---$T_{\text{III}, D, AD}$}

Let $\kappa$ be the (unique) nontrivial quadratic character
of the group $R^\times / {R^\times}^2$.
\begin{claim}\label{claim:oiTIIIDAD}
The $\ve$-twisted orbital integral $O_\delta(1_K)$ is equal to
\[
2 \kappa(B_1 B_2) (-1)^{N_1}q_0^{1 + N_1 + N_2}
\lp \frac{q^{N + 1} - 1}{q - 1} \rp.
\]
\end{claim}
\begin{proof}
Fix an element $\epsilon_0 \in R^\times - {R^\times}^2$.
Without loss of generality, we may assume that $A = \ep_0$, $D = -\ep_0\vp$, and $AD = -\vp$.
We fix uniformizers
$\vp_{D} = \sqrt{-\epsilon_0\vp}$, $\vp_{AD} = \sqrt{-\vp}$ for $\ED, \EAD$, respectively.
Thus, $\N_{\ED/F}\vp_{D} = \epsilon_0\vp$, and $\N_{\EAD/F}\vp_{AD} = \vp$.

Let 
\[\begin{split}
C_0' &= \{ (g_1, g_2) \in \GL(2, F)^2 : \det g_1 = \det g_2\},\\
K_0 &= \{ (k_1, k_2) \in \GL(2, R)^2 : \det k_1 = \det k_2\}.
\end{split}\]
Let $K' = \GL(2, R)$.  For $D' = D, AD$, let
\[
T_{D'} = \{\phi^{D'}(\gamma) : \gamma \in F_{D'}^\times\} \subset \GL(2, F).
\]

We have the decompositions
\[
\GL(2, F) = \dot{\bigcup}_{j_1\geq 0}T_D\lp \lsm 1 & \\ & \vp^{j_1}\rsm\rp K'
= \dot{\bigcup}_{j_2\geq 0}T_{AD}\lp \lsm 1 & \\ & \vp^{j_2}\rsm\rp K'.
\]
Hence, each element $g \in C_0'$ may be written as
\[
g = (t_1, t_2)\cdot \Big(\phi^D(\vp_{D}^{-j_1})\lp\lsm1 & \\ & (\epsilon_0\vp)^{j_1}\rsm\rp,\,
\phi^{AD}(\vp_{AD}^{-j_2})\lp \lsm 1 & \\ & \vp^{j_2}\rsm\rp\Big)
\cdot (k_1, k_2),
\]
where $(t_1, t_2) \in T_D \times T_{AD}$, and $(k_1, k_2) \in \GL(2, R)^2$, such that $\det (t_1k_2)$ is equal to 
$\det (t_2k_2)$.
If $\det k_1 \neq \det k_2$, we may multiply $k_1$ by $\diag(x,x)$ from $T_{D}$, for some $x \in R^\times$,
such that  $\det k_1$ and $\det k_2$ differ by $\epsilon_0 \in R^\times - {R^\times}^2$.

Let $s(\ep)$ be $0$ if $\ep = 1$, $1$ if $\ep = \ep_0$.
Let $\R$ be the set of elements in $C_0'$ of the form:
\[
r_{j_1, j_2, \epsilon} = 
\Big(\phi^D\!\!\lp\vp_{D}^{-j_1 - s(\epsilon)}\rp\lp\lsm1 & \\ & \epsilon\cdot(\epsilon_0\vp)^{j_1}\rsm\rp,\,
\phi^{AD}\!\!\lp\vp_{AD}^{-j_2 - s(\epsilon)}\rp\lp \lsm 1 & \\ & \vp^{j_2}\rsm\rp\Big),
\]
where  $j_1$, $j_2$ range over nonnegative integers, and $\epsilon$ over the set $\{1, \epsilon_0\}$.
Then, $C_0'$ decomposes as follows:
\begin{equation}\label{eq:decompTIIIDAD}
C_0' = \dot{\bigcup}_{\epsilon \in \{1, \epsilon_0\}}
\dot{\bigcup}_{j_1, j_2 \geq 0}T\cdot r_{j_1, j_2, \epsilon} \cdot K_0.
\end{equation}


Let $T = T_{\text{III}, D, AD}$.  Recall that we put
\[
C_0 := \{[g_1, g_2]:g_1,g_2 \in \GL(2, F);\, \det g_1 = \det g_2\} \subset \GSp(2, F).
\]
It is isomorphic to $C_0'$ via $[g_1, g_2] \rightarrow (g_1, g_2)$.
By abuse of notation, we let $T$ and $\delta$ denote also their images in $C_0'$.

For $m \in \mbb{Z}$, recall that
\[
K_m := \{(g, h) \in C_0' : g - ehe \in \vp^m M_2(R)\},\quad e = \diag(1, -1).
\]
Let $R_T = T \cap K$.  By the same technique used in the case of $T_{{\rm I}, A}$, we have:
\begin{multline}\label{eq:tIIIDAD}
O_\delta(1_K) = \sum_{m\geq 0} (-1)^m
\sum_{j_1, j_2 \geq 0} \sum_{\,\,\epsilon=1, \epsilon_0} \ve(r)\\
\cdot [R_T : T \cap r K_0 r^{-1}] \int_{K_m\bs K_0}1_{K_m}(k^{-1}r^{-1}\delta r k)\,dk.
\end{multline}
Here, $r$ is the abbreviation of $r_{j_1, j_2, \ep}$.

Recall that $\xi = a_1 + b_1 \sqrt{D}$, $\zeta = a_2 + b_2\sqrt{AD}$,
where $b_i = B_i \vp^{N_i}$, $B_i \in R^\times$, $N_i \in \mathbb{Z}$.
For nonnegative integers $j_1$ and $j_2$, let $\nu_1 = N_1 - j_1$ and $\nu_2 = N_2 - j_2$.
Let $D_1 = D$ and $D_2 = AD$.  For $r = r_{j_1, j_2, \epsilon}$, we have:
\[
r^{-1}\delta r
=\lp
\lp\lsm
 a_1 & b_1' D_1' \\
b_1' & a_1
\rsm\rp, 
\lp\lsm
a_2 & b_2' D_2'\\
b_2' & a_2
\rsm\rp
\rp,
\]
where
\begin{align*} 
b_1' &= \lp{B_1\vp^{\nu_1}}\rp/\lp{\epsilon\epsilon_0^{j_1}}\rp, & D_1' &= D_1\ep^2\cdot(\ep_0\vp)^{2 j_1}.\\
b_2' &= B_2\vp^{\nu_2}, & D_2' &= D_2\vp^{2 j_2}.
\end{align*}

Put $R_m := R/\lp\vp^m R\rp$.  For $a \in R$, let $\bar{a}$ denote its image in $R_m$.
Following the same argument in \cite[p. 64]{F}, 
$\int_{K_m \bs K_0}1_{K_m}(k^{-1}r^{-1}\delta r k)\,dk$ is equal to the cardinality of
\[
L_{m, r} := 
\left\{x = \lp\lsm x_1&x_2\\x_3&x_4\rsm\rp \in SL(2, R_m):
\lp\lsm\ol{a_1} & \ol{b_1'}\ol{D_1'}\\ \ol{b_1'} & \ol{a_1'}\rsm\rp
\lp\lsm x_1&x_2\\x_3&x_4\rsm\rp
= \lp\lsm x_1&x_2\\x_3&x_4\rsm\rp 
\lp\lsm\ol{a_2} & \ol{b_2'}\ol{D_2'}\\ \ol{b_2'} & \ol{a_2'}\rsm\rp
\right\}.
\]

Recall that $X = \ord(a_1 - a_2)$.
Transcribing Lemma II. G. 3 of \cite{F} 
to the current situation, the following holds:
\begin{claim}\label{claim:T2DTIIIDADconditions}
The set $L_{m,r}$ $(r = r_{j_1, j_2, \epsilon})$ is nonempty precisely when the following conditions are satisfied:
\begin{enumerate}
\item
$0 \leq m \leq X.$
\item
$\nu_1, \nu_2 \geq 0$.
\item
$m \leq \nu_1$ if and only if $m \leq \nu_2$.
\item
If $\nu_1 < m$ or $\nu_2 < m$, then $\nu_1 = \nu_2$; we denote the common value by $\nu$.
\item
If $\nu < m$, then $(B_1/\ep\epsilon_0^{j_1})/B_2 \in {R^\times}^2$.
\item
If $\nu < m$, then $m \leq 2 N_i - \nu + {\rm ord}\, D_i$ $(i = 1, 2)$.
\end{enumerate}
\begin{remark}
The only differences between the above conditions and those in Lemma II. G. 3
of \cite{F} are item 5 and the definition of $D_i'$ ($i = 1, 2$).
\end{remark}
\end{claim}

Recall that $q_0$ is the cardinality of the residue field of $F$.
\begin{claim}\label{claim:T2DTIIIDADcardLmr}
If $L_{m, r}$ is nonempty, then
\[
\#L_{m, r} = 
\begin{cases}
1 & \text{ if } m = 0,\\
(q_0^2 - 1)q_0^{3m - 1} & \text{ if } 1 \leq m \leq \nu_1\, (\text{equivalently: } 1 \leq m \leq \nu_2),\\
2q_0^{m + 2\nu} & \text{ if } \nu < m.
\end{cases}
\]
\end{claim}
\begin{proof}
This follows from Lemma II. G. 3 of \cite{F}.
\end{proof}
\begin{claim}
For $r = r_{j_1, j_2, \epsilon}$, 
the following holds:
\[
[R_T : T \cap r K_0 r^{-1}] = q_0^{j_1 + j_2}.
\]
\end{claim}
\begin{proof}
For $k = \ED$ or $\EAD$, and a nonnegative integer $j$, put $R_k(j) := R + \vp^j R_k$.
The group $R_T = T\cap K$ is isomorphic to $\lp R_{\ED}^\times \times R_{\EAD}^\times\rp'$, 
where the prime indicates equal determinants. 
The computation of $[R_T : T \cap r K_0 r^{-1}]$ is equivalent to computing the cardinality of the 
kernel in the short exact sequence
\begin{equation*}\begin{split}
1 & \rightarrow (R_{\ED}^\times \times R_{\EAD}^\times)'/(R_{\ED}(j_1)^\times \times R_{\EAD}(j_2)^\times)'\\
& \rightarrow \lp R_{\ED}^\times \times R_{\EAD}^\times\rp/\lp R_{\ED}(j_1)^\times \times R_{\EAD}(j_2)^\times\rp
\\
& \rightarrow \lp R_{\ED}^\times \times R_{\EAD}^\times\rp/
\Big[\lp R_{\ED}^\times \times R_{\EAD}^\times\rp' \lp R_{\ED}(j_1)^\times \times R_{\EAD}(j_2)^\times\rp\Big]
\rightarrow 1.
\end{split}\end{equation*}
As shown in \cite[Lemma I. I. 2]{F}, for $k = \ED$, $\EAD$,
\[
[R_k^\times : R_k(j)^\times] = 
[R_k^\times : 1 + \vp^j R_k] / [R^\times : 1 + \vp^j R] = q_0^j.
\]
Hence, the cardinality of the middle term of the sequence is $q_0^{j_1 + j_2}$.
The image in the short exact sequence is isomorphic via $\N_{\ED/F} \times \N_{\EAD/F}$ to 
\begin{multline}\label{eq:TIIIDADshortexact}
\lp \N_{\ED/F}R_{\ED}^\times \times \N_{\EAD/F}R_{\EAD}^\times\rp/\\
\left\{ x\cdot \N_{\ED/F}R_{\ED}(j_1)^\times\times x\cdot \N_{\EAD/F}R_{\EAD}(j_2)^\times :\right.\\
\left. x \in \N_{\ED/F}R_{\ED}^\times \cap \N_{\EAD/F}R_{\EAD}^\times\right\}.
\end{multline}
Since the field extensions $\ED/F$ and $\EAD/F$ are ramified, we have
\[
\N_{\ED/F}R_{\ED}^\times =  \N_{\EAD/F}R_{\EAD}^\times = {R^\times}^2.
\]
Consequently, the set \eqref{eq:TIIIDADshortexact} has cardinality $1$.  
The claim follows.
\end{proof}
Rearranging the terms in \eqref{eq:tIIIDAD}, and noting that $\ve(r_{j_1, j_2, \epsilon})$ is equal to
\[
\ve\lp \det \lp\phi^D\!(\vp_D^{-j_1})\,\phi^D\!\!\lp\vp_D^{-s(\epsilon)}\rp\rp \cdot (\epsilon_0\vp)^{j_1}\rp
= \ve\lp(\epsilon_0\vp)^{s(\epsilon)}\rp = (-1)^{s(\epsilon)}, 
\]
the expression \eqref{eq:tIIIDAD} is the sum of the following two quantities:
\begin{gather}
\label{eq:tIIIDAD1}
\sum_{\epsilon = 1, \epsilon_0} (-1)^{s(\epsilon)}
\sum_{\substack{0\leq\nu_1\leq N_1\\ 0\leq \nu_2 \leq N_2}}q_0^{N_1 - \nu_1 + N_2 - \nu_2}
\lp 1 + (1 - q_0^{-2})\sum_{m = 1}^{\min (\nu_1,\nu_2)} q_0^{3m}\rp,\\
\label{eq:tIIIDAD2}
\sum_{\epsilon = 1, \epsilon_0} (-1)^{s(\epsilon)}  q_0^{N_1 + N_2 - 2\nu} \sum_{0 \leq \nu \leq N} \sum_{\nu < m \leq X - \nu} (-1)^m
L_{m, r_{j_1, j_2, \epsilon}}.
\end{gather}
Here,  $\nu$ is the common value of $\nu_1$ and $\nu_2$ when $\nu_1 = \nu_2 < m$. 

The second sum in \eqref{eq:tIIIDAD1} is independent of $\ep$; hence, \eqref{eq:tIIIDAD1} is
equal to zero.

In \eqref{eq:tIIIDAD2}, suppose $L_{m, r} \neq 0$.
Since $\nu < m$, by Claim \ref{claim:T2DTIIIDADconditions} we have \dts $B_1/(B_2 \ep \epsilon_0^{j_1}) \in {R^\times}^2$.
If $B_1/ B_2$ lies in ${R^\times}^2$, then since $\epsilon_0 \in R^\times - {R^\times}^2$,
$j_1$ is even or odd depending on (in that order) $\epsilon$ being $1$ or $\epsilon_0$.
If $B_1/ B_2 \notin {R^\times}^2$, then $j_1$ is odd or even depending on (in that order) $\epsilon$ being $1$ or $\epsilon_0$.
Recall that $\kappa$ is the character of $R^\times$ which is $1$ at squares and $-1$ otherwise.
Thus, $(-1)^{s(\ep)}$ is equal to $\kappa(B_1 B_2)(-1)^{j_1}$ if $L_{m, r_{j_1, j_2, \ep}} \neq 0$.

By the comments above and Claim \ref{claim:T2DTIIIDADcardLmr}, we conclude that
\begin{equation*}\begin{split}
O_\delta(1_K)
&=\kappa(B_1 B_2) \sum_{0 \leq \nu \leq N} (-1)^{N_1 - \nu} q_0^{N_1 + N_2 - 2\nu}\sum_{\nu < m \leq X - \nu} (-1)^m
2q_0^{m+ 2\nu}\\
&= 2 \kappa(B_1 B_2) (-1)^{N_1}q_0^{N_1 + N_2} \sum_{0 \leq \nu \leq N} (-1)^\nu \sum_{\nu < m \leq X - \nu}
(-q_0)^m\\
&= 2  \kappa(B_1 B_2) (-1)^{N_1}q_0^{1 + N_1 + N_2}\frac{1}{q_0 + 1}\sum_{0\leq \nu \leq N}
\left[(-q_0)^X q_0^{-\nu} - q_0^\nu \right]\\
&= 2 \kappa(B_1 B_2) (-1)^{N_1}q_0^{1 + N_1 + N_2} 
\frac{1}{q_0^2 - 1}((-1)^X q_0^{X - N} - 1)(q_0^{N+1} - 1).
\end{split}\end{equation*}

By Lemma II. L. 1 of \cite{F}, $X = 2N + 1$.
Hence,
\[
\begin{split}
O_\delta(1_K)
&= 
2 \kappa(B_1 B_2) (-1)^{N_1}q_0^{1 + N_1 + N_2}
\left(
\frac{(q_0^2)^{N + 1} - 1}{q_0^2 - 1}
\right)\\
&= 
2 \kappa(B_1 B_2) (-1)^{N_1}q_0^{1 + N_1 + N_2}
\left(
\frac{q^{N + 1} - 1}{q - 1}
\right).
\end{split}
\]
\end{proof}
\subsubsection{Transfer Factor---$T_{2, D}, T_{{\rm III}, D, AD}$}
We define the transfer factor $\Delta(\delta_2, \delta)$ to be
\[
\frac{1}{2}\,
\kappa\!\lp\frac{\xi - \tau\xi}{2\sqrt{D}}\cdot \frac{\zeta - \sigma\tau\zeta}{2\sqrt{AD}}\rp
\ve\! \lp \frac{\xi - \tau\xi}{2\sqrt{A}}\rp
\lp \frac{D_G(\delta)}{D_{H_2}(\delta_2)}\rp^{1/2}.
\indexi{transfer factor}%
\]
\begin{claim}\label{claim:deltaT2DTIIIDAD}
The following holds:
\[\Delta(\delta_2, \delta) = \frac{1}{2}\, \kappa(B_1 B_2) (-1)^{N_1}q_0^{-1 - N_1 - N_2}.\]
\end{claim}
\begin{proof}
Note that by hypothesis $\xi = \N_{E/\ED}l = l\sigma l$ and $\zeta = \N_{E/\EAD}l = l\sigma\tau l$, 
where $l \in R_E^\times$.
The claim then follows from a straightforward computation.
\end{proof}
\begin{prop}\label{prop:oiT2DTIIIDAD}
The following identity holds:
\[
O_{\delta_2}(1_{K_2}) = \Delta(\delta_2, \delta) O_{\delta}(1_K).
\]
\end{prop}
\begin{proof}
This follows from Claims \ref{claim:oiT2D}, \ref{claim:oiTIIIDAD}, and \ref{claim:deltaT2DTIIIDAD}.
\end{proof}
\subsection{Comparison for $T_{2, D}$, $T_{{\rm IV}, A, D}$}
Let $D$ be an element in $\EA^\times - {\EA^\times}^2 \cup F^\times$.  Let $E = \EA(\sqrt{D}) = F(\sqrt{D})$.
Let $\tilde{E}$ be the Galois closure of $E$ over $F$.  Then, $\Gal(\tilde{E}/F) = \Gal(E/F)$
is cyclic of order $4$ if $E/F$ is Galois, and $\Gal(\tilde{E}/F) = D_4$ if $E/F$ is non-Galois.

Let $\sigma$ be the element in $\Gal(\tilde{E}/F)$ defined by
$\sigma\A = -\A$, $\sigma \D = \sqrt{\sigma D}$, $\sigma^2\D = -\D$, $\sigma^3\D = -\sqrt{\sigma D}$.
Note that $\sigma$ is of order $4$, and it generates $\Gal(E/F)$ if $E/F$ is Galois.
Moreover, $\Gal(\EA/F)$ is generated by the restriction of $\sigma$ to $\EA$.



Let $q$ be the cardinality of the residue field of $\EA$. 
Let $q_0$ be the cardinality of the residue field of $F$.  Since by assumption $\EA/F$ is unramified, $q = q_0^2$.
We fix a uniformizer $\vp$ for both  $F$ and $\EA$.  

Let $l$ be an element  in $R_E^\times - R_{\EA}^\times$, and $c$ an element in $R^\times$.
Suppose $l = \alpha + \beta \sqrt{D}$, where $\alpha, \beta \in \EA$.
Let
\[
\delta_2 = (\phi^{\EA:D}(l), c)_* := \lp\lp\lsm \alpha &\beta D\\\beta &\alpha\rsm\rp, c\rp_* \in T_{2, D} \cap K_2.
\]
Here, lower $\ast$ indicates image in $H_2$.

For an element $\gamma = a + b/\sqrt{A} \in \EA^\times$, with $a, b \in F$, put
$\displaystyle \bg{\gamma} := \lp\lsm a & b/A \\ b & a\rsm\rp.$
Suppose $l\sigma l = \xi_1 + \xi_2\sqrt{D}$, with $\xi_1, \xi_2\in \EA$.
Let
\[
\delta = c\lp\lsm \bg{\xi_1} & \bg{\xi_2}\bg{D}\\\bg{\xi_2} & \bg{\xi_1}\rsm\rp \in T_{{\rm IV}, A, D} \cap K.
\]

The element $\delta$ is regular, and $\delta_2$ is a norm of $\delta$.
The conjugacy classes of $\delta_2$ and $\delta$ are stable.
We assume that $\delta$ is topologically unipotent.

Let $n = {\rm ord}_{\EA} \beta$, and $N = {\rm ord}_{\EA} \xi_2$.  
Let $B$ be the element in $R_{\EA}^\times$ such that $\xi_2 = B\vp^N$.

\subsubsection{Stable Orbital Integral---$T_{2, D}$}
\begin{claim}\label{claim:oiT2D'}
The stable orbital integral $O_{\delta_2}(1_{K_2})$ is equal to
\begin{itemize}
\item
$\displaystyle \frac{(q+1)q^{n} - 2}{q - 1}$\; if $E/\EA$ is unramified,
\item
$\displaystyle \frac{q^{n + 1} - 1}{q - 1}$\; if $E/\EA$ is ramified.
\end{itemize}
\end{claim}
\begin{proof}
This follows from Lemma I. I. 2 of \cite{F}.
\end{proof}
\subsubsection{Twisted Orbital Integral---$T_{\text{IV}, A, D}$}
Let $\kappa$ be the nontrivial quadratic character of the group $R_{\EA}^\times / {R_{\EA}^\times}^{\!\!\!2}$.
\begin{claim}\label{claim:oiTIVAD}
The $\ve$-twisted orbital integral $O_{\delta}(1_K)$ is equal to
\begin{itemize}
\item
$\displaystyle (-q)^N \lp \frac{(q+1)q^N - 2}{q - 1} \rp$ if $E/\EA$ is unramified,
\item
$\displaystyle -2\kappa(B)\,(-q)^Nq_0\lp\frac{q^{N + 1}  - 1}{q - 1}\rp$ if $E/\EA$ is ramified.
\end{itemize}
\end{claim}
\begin{proof}



Let
\[
\begin{split}
C_A &= \GL(2, \EA)' := \{g \in \GL(2, \EA) : \det g \in F^\times\},\\
\mb{C}_A &= \left\{\lp\begin{smallmatrix} \bg{\alpha} & \bg{\beta} \\ \bg{\gamma} & \bg{\eta}\end{smallmatrix}\rp
: \lp\lsm \alpha & \beta \\ \gamma&\eta\rsm\rp \in C_A \right\}.
\indexs{C@$C_A$, $\mb{C}_A$}%
\end{split}
\indexs{C@$C_A$, $\mb{C}_A$}%
\]

For a nonnegative integer $m$, let
$u_m = \lp\begin{smallmatrix}1&0&0&\vp^{-m}/A\\0&1&0&0\\0&0&1&0\\0&0&0&1\end{smallmatrix}\rp$.
Let
 \[\begin{split} 
K_m &= \GL(2, R_{\EA}(m))',\\
K_m^A &= \mb{C}_A \cap u_m K u_m^{-1}.
\indexs{K@$K_m^A$}%
 \end{split}\]  
Here, $R_{\EA}(m) := R + (\vp^m \sqrt{A}) R$, and the prime indicates determinant in $F^\times$.

For $a = a_1 + a_2/\sqrt{A}, b = b_1 + b_2/\sqrt{A}, c= ... \in \EA^\times$, put
\[
\phi\begin{pmatrix}a&b\\c&d\end{pmatrix} := 
\begin{pmatrix}
a_1 & a_2/A & b_2 & b_1\\ 
a_2 & a_1 & b_1 A &b_2\\
c_2/A & c_1 /A&d_1&d_2/A\\
c_1&c_2/A&d_2&d_1
\end{pmatrix}.
\indexs{phi@$\phi(\;)$}%
\]
We define an isomorphism {$\phi_m : C_A \rightarrow \mb{C}_A$} as follows:
\[
\phi_m : \lp\begin{smallmatrix} \alpha&\beta\\\gamma&\eta\end{smallmatrix}\rp \mapsto
\phi
\lp
\lp\begin{smallmatrix}1&0\\0&A\vp^m\end{smallmatrix}\rp
\lp\begin{smallmatrix}1&0\\-1&1\end{smallmatrix}\rp
\lp\begin{smallmatrix}1&0\\0&\sqrt{A}\end{smallmatrix}\rp
\lp\begin{smallmatrix}\alpha&\beta\\ \gamma&\eta\end{smallmatrix}\rp
\lp\begin{smallmatrix}1&0\\0&1/\sqrt{A}\end{smallmatrix}\rp
\lp\begin{smallmatrix}1&0\\1&1\end{smallmatrix}\rp
\lp\begin{smallmatrix}1&0\\0&1/\lp A\vp^m\rp\end{smallmatrix}\rp
\rp.
\]
In particular, $\phi_m$ maps $K_m$ isomorphically onto $K_m^A$ (\cite[Lemma I. J. 3]{F}).

We have the disjoint decomposition (\cite[Lemma I. J. 1]{F}):
\begin{equation}\label{eq:decompTIVAD}
G = \dot{\bigcup}_{m\geq 0} \mb{C}_A u_m K.
\end{equation}
Let $T = T_{\text{IV}, A, D}$.  
Noting that $\ve(u_m) = 1$,
we have:
\begin{equation}\label{eq:oiTIVADsum}
\begin{split}
O_\delta(1_K) &= \int_{T\bs G} 1_K(g^{-1}\delta g)\ve(g)\,dg\\
&= \sum_{m\geq 0} \abs{K}_G \int_{T \bs \mb{C}_A / \mb{C}_A \cap u_m K u_m^{-1}} 1_K(u_m^{-1}h^{-1} \delta h u_m) \ve(h u_m)\,dh\\
&= \sum_{m\geq 0} \abs{K}_G\int_{T \bs \mb{C}_A / K_m^A} 1_{K_m^A}(h^{-1} \delta h) \ve(h)\,dh.
\end{split}
\end{equation}

Let $\delta'$ be the element in $C_A$ such that $\phi(\delta') = \delta$.  Let $T' = Z_{C_A}(\delta')$.
For $m \geq 0$, let $\delta_m = \phi_m^{-1}(\delta) \in C_A$.  Let $T_m = Z_{C_A}(\delta_m) = \phi_m^{-1}(T)$.
We rewrite \eqref{eq:oiTIVADsum} as follows:
\begin{equation}\label{eq:oiTIVAD1}
O_\delta(1_K) 
=\sum_{m\geq 0} \abs{K}_G\int_{T_m \bs C_A / K_m} 1_{\phi_m(K_m)}(\phi_m(h')^{-1}\delta \phi_m(h'))\ve(\phi_m(h'))dh'.
\end{equation}
Here, $h'$ denotes the element in $C_A$ such that $\phi_m(h') = h \in \mb{C}_A$.

The similitude factor of $\phi(h')$ is $(\det h')^{-1}$; hence, 
\[
\ve(\phi_m(h')) = \ve((\det h')^{-1}) = \ve(\det h').
\]
We change variables
$h' \mapsto \lp\begin{smallmatrix}1&0\\1&1\end{smallmatrix}\rp
\lp\begin{smallmatrix}1&0\\0&1/\lp A\vp^m\rp\end{smallmatrix}\rp h'$.  The expression \eqref{eq:oiTIVAD1} becomes
\[
\sum_{m\geq 0}\abs{K}_G\abs{K_m}_{C_A}^{-1} \int_{T' \bs C_A} 1_{K_m}({h'}^{-1}\delta'h')\ve(A\vp^m)\ve(\det h')\,dh'.
\]

Since $\EA/F$ is unramified, $\ve(A)$ is equal to one.  Consequently, $\ve(A\vp^m) = (-1)^m$ for each nonnegative integer $m$.

Let $\tilde{T}'$ be the centralizer of $\delta'$ in $\GL(2, \EA)$.  It contains $T'$.
Let 
\[\begin{split}
K' &= \GL(2, R_{\EA})' = \left\{ k \in \GL(2, R_{\EA}) : \det k \in F^\times\right\},\\
T_0' &= T' \cap K'.
\end{split}\]

Suppose $E/\EA$ is unramified.
Since 
\[\GL(2, \EA) = \dot{\bigcup}_{j \geq 0}\tilde{T}'\cdot \dmtwo{1}{\vp^j} \cdot\GL(2, R_{\EA}),\]
the group $C_A$ decomposes as follows (\cite[Lemma II. M. 3]{F}): 
\[
C_A = \dot{\bigcup}_{j, \ep} T'\cdot r_{j, \ep}\cdot  K'.
\]
The disjoint union is over nonnegative integers $j$ and $\ep \in R_{\EA}^\times/{R_{\EA}^\times}^{\!\!2}$ if $j \geq 1$.
Here, $r_{j,\epsilon} = t_\epsilon\lp\begin{smallmatrix}1&0\\0&\epsilon\vp^j\end{smallmatrix}\rp$, where
$t_\ep$ is an element in $\tilde{T}'$ such that $\det t_\ep = \ep^{-1}$.

By the above decomposition of $C_A$, the twisted orbital integral $O_{\delta}(1_K)$ is equal to
\begin{equation}\label{eq:oiTIVAD2}
\sum_{m\geq 0}(-1)^m\sum_{r = r_{j,\ep}} \ve(\det r)
[T_0' : T' \cap r K' r^{-1}]\cdot [K_0 : K_m] \int_{K_0} 1_{K_m}(k^{-1}r^{-1}\delta'rk)dk.
\end{equation}
Since $\EA/F$ is unramified, $\ve(r_{j,\ep})$ is equal to $(-1)^j$.  
Thus, $O_{\delta}(1_K)$ is equal to
\begin{equation}\label{eq:oiTIVAD3}
\sum_{m\geq 0}(-1)^m\sum_{r = r_{j, \ep}} (-1)^j
[T_0' : T' \cap r K' r^{-1}]\cdot[K_0 : K_m] \int_{K_0} 1_{K_m}(k^{-1}r^{-1}\delta'rk)\,dk.
\end{equation}

Recall: $\delta = c \lp\lsm \bg{\xi_1}&\bg{\xi_2}\bg{D}\\ \bg{\xi_2}&\bg{\xi_1}\rsm\rp \in T$,
$N = {\rm ord}_{\EA} \xi_2$.  Let $X = {\rm ord}_{\EA}(\xi_1 - \sigma\xi_1)$.

By the computation of $[T_0' : T' \cap r K' r^{-1}]$, $[K_0 : K_m] \int_{K_0} 1_{K_m}(k^{-1}r^{-1}\delta'rk)\,dk$
in \cite[Sect. II. M]{F}, and \eqref{eq:oiTIVAD3}, the twisted orbital integral $O_\delta(1_K)$ 
is equal to the sum of the following two quantities if $E/\EA$ is unramified:
\begin{itemize}
\item[I.]
$\displaystyle
\lp\sum_{m = 0} 1 + \frac{q+1}{q}\sum_{1\leq m \leq N} (-q_0)^{3m}\rp
\lp \sum_{\nu = m} 1 + \frac{q+1}{q}\sum_{m \leq \nu < N}(-q)^{N - \nu}\rp,
$
\item[II.]
$\displaystyle
\sum_{0 \leq \nu < N} \frac{q + 1}{q}(-q)^{N - \nu}\sum_{\nu < m \leq X - \nu} (-q_0)^m q^\nu.
$
\end{itemize}
The sum I is equal to 
\begin{equation}\label{eq:TIVADunramsum1}
\frac{(-q)^N}{1 - q_0}\lp1 - \lp\frac{q+1}{q_0}\rp q_0^N + q_0^{-1}\rp,
\end{equation}
and the sum II is equal to
\begin{equation}\label{eq:TIVDunramsum2}
\frac{(q+1)(-q)^N}{(q_0 + 1)(1 - q_0)q_0}\lp q_0^{X - N + 1} - q_0^{X + 1} + q_0^N - 1\rp.
\end{equation}

From \cite[Sect. II. M]{F}, since $\delta$ is topologically unipotent, $X$ is equal to $2N$.
Hence, $O_{\delta}(1_K) =$ I + II is equal to
\[\begin{split}
&= \frac{(-q)^N}{1 - q_0}\lp\substack{ 1 - \lp\frac{q+1}{q_0}\rp q_0^N + q_0^{-1} +
\frac{q+1}{(q_0 + 1) q_0} q_0^{X -N + 1}\\ - \frac{q+1}{(q_0 + 1)q_0}q_0^{X + 1} +
\frac{q+1}{(q_0+1)q_0}q_0^N - \frac{q+1}{(q_0 + 1) q_0}}\rp\\
&= \frac{(-q)^N}{1 - q_0}\frac{1}{(q_0 + 1)q_0}((q+1)(-q_0^{2N + 1} - 1) + q_0^2 + 2q_0 + 1)\\
&= \frac{(-q)^N}{-(q - 1)}\frac{1}{q_0}((q+1)(-q^N q_0) - q - 1 + q + 2q_0 + 1)\\
&= (-q)^N \left[ \frac{(q+1)q^N - 2}{q - 1} \right].
\end{split}\]
This completes the proof of the claim for the case where $E/\EA$ is unramified.

Now, suppose $E/\EA$ is ramified.  
Recall that we have fixed a uniformizer $\vp$ for both $F$ and $\EA$.
Since $E/\EA$ is ramified, and $D \in \EA^\times - F^\times$, 
there exists an element $\ep_0 \in  R_{\EA}^\times - {R_{\EA}^\times}^{\!\!\!2}\cup F^\times$ such that $D = -\ep_0\vp$.
Then, $\N_{E/\EA}\sqrt{D} = \ep_0\vp$, and the element $t_0 = \phi^{\EA:D}(\sqrt{D}) \in  \tilde{T}'$ 
has determinant $\ep_0\vp$.

The group $C_A$ has a disjoint decomposition as follows:
\[
C_A = \dot{\bigcup}_{j,\ep} T'\cdot r_{j, \ep}\cdot  K',
\]
where $j$ ranges over the nonnegative integers, and $\ep$ over the set $\{1, \ep_0\}$.
Here, $r_{j,\ep}$ is an element in $\tilde{T}'\lp\lsm 1 &\\ &(\ep_0\vp)^j/\ep\rsm\rp$ such that 
$\det r_{j, \ep} = \begin{cases} 1 & \text{if } \ep = 1,\\ \vp & \text{if }\ep = \ep_0.\end{cases}$

The twisted orbital integral $O_\delta(1_K)$ may then be expressed as the sum
\begin{equation}
\sum_{m\geq 0}(-1)^m\sum_{r = r_{j,\ep}} (-1)^{s(\ep)}
[T_0' : T' \cap r K' r^{-1}]\cdot[K_0 : K_m] \int_{K_0} 1_{K_m}(k^{-1}r^{-1}\delta'rk)dk,
\end{equation}
where $s(\ep)$ is equal to $0$ if $\ep = 1$, $1$ if $\ep = \ep_0$.

Recall that $\delta = c\lp\lsm \bg{\xi_1} & \bg{\xi_2}\bg{D}\\\bg{\xi_2} & \bg{\xi_1}\rsm\rp \in \mb{C}_A$.
Its preimage under $\phi$ in $C_A$ is the element $\delta' = c\, \phi^{\EA:D}(\xi_1 + \xi_2\sqrt{D})$.
Recall that $B$ is the element in $R_{\EA}^\times$ such that $\xi_2 = B\vp^N$ ($N = {\rm ord}_{\EA}\xi_2$).

For a nonnegative integer $j$, let $\nu = N - j$.  Then,
\[
r_{j, \ep}^{-1}\delta' r_{j, \ep} =
\lp\lsm
\xi_1 & \xi_2' D' \\
\xi_2' & \xi_1
\rsm\rp,
\]
where
$ \xi_2' = (B\ep/\ep_0^j)\vp^\nu$, $D' = (D/\ep^2) (\ep_0\vp)^{2j}$.

Let $X = {\rm ord}_{\EA}(\xi_1 - \sigma\xi_1)$.
By the same argument used in the proof of Lemma II. M. 6 in \cite{F}, the integral
$\int_{K_0} 1_{K_m}(k^{-1}r_{j,\ep}^{-1}\delta'r_{j,\ep}k)\,dk$ is nonzero precisely when
$0 \leq \nu \leq N$ and $0 \leq m \leq X$.  Moreover, if $m > \nu$, then $\nu + m \leq X$ and
$\ep_0^j$ lies in $B\ep{R_{\EA}^\times}^{\!\!\!\!2}$.

For $\nu < m \leq X - \nu$, the condition $\ep_0^j \in B\ep {R_{\EA}^\times}^{\!\!\!2}$ implies that $j$ is 
even if $B\ep \in {R_{\EA}^\times}^{\!\!\!2}$, and odd otherwise.  Thus, $\ep \in \{1, \ep_0\}$ depends on $j$.
Moreover, $(-1)^{s(\ep)}$ is equal to $\kappa(B)(-1)^{j}$, where $\kappa$ is the nontrivial quadratic character of the group 
$R_{\EA}^\times / {R_{\EA}^\times}^{\!\!\!2}$.

The quantities $[T_0' : T' \cap r K' r^{-1}]$ and 
$[K_0 : K_m] \int_{K_0} 1_{K_m}(k^{-1}r^{-1}\delta'rk)\,dk$
are computed in \cite[Sect. II. M]{F}.  
Recall that $q_0$ is the cardinality of the residue field of $F$, and $q = q_0^2$ that of $\EA$.
The twisted orbital integral $O_\delta(1_K)$ is equal to the sum of the following:
\begin{enumerate}
\item
$\displaystyle \sum_{\ep = 1,\, \ep_0}(-1)^{s(\ep)}\sum_{0 \leq \nu \leq N} q^{N - \nu}$,
\item
$\displaystyle \sum_{\ep = 1,\,\ep_0}(-1)^{s(\ep)}\sum_{0\leq \nu \leq N}\frac{q^{N - \nu}(q + 1)}{q}\sum_{1 \leq m \leq \nu} (-q_0)^{3m}$,
\item
$\displaystyle \sum_{0 \leq \nu \leq N} \kappa(B) (-q)^{N - \nu} \sum_{\nu < m \leq X - \nu} 2(-q_0)^m q^\nu$

$= 2\kappa(B)(-q)^N \sum_{0 \leq \nu \leq N} (-1)^\nu \sum_{\nu  < m \leq X - \nu}(-q_0)^m$.
\end{enumerate}
The first two sums are zero, for their inner sums are independent of $\ep$.
The twisted orbital integral $O_\delta(1_K)$ is therefore equal to the third sum, which is
\[
2\kappa(B)(-q)^N(-q_0)\lp\frac{(-1)^X q_0^{X - N} - (-1)^X q_0^{X +1} + q_0^{N+1} - 1}{q - 1}\rp.
\]
From \cite[Sect. II. M]{F}, $X = 2N + 1$; hence,
\[
O_\delta(1_K) 
= - 2\kappa(B)(-q)^N q_0\lp\frac{q_0^{2 N + 2}  - 1}{q - 1}\rp = -2\kappa(B)(-q)^Nq_0\lp\frac{q^{N + 1}  - 1}{q - 1}\rp.
\]
\end{proof}

\subsubsection{Transfer Factor---$T_{2, D}$, $T_{{\rm IV}, A, D}$}
Recall: $l = \alpha + \beta \sqrt{D} \in R_E^\times$, $c \in R^\times$, $l \sigma l = \xi_1 + \xi_2 \sqrt{D}$.
We let $\delta_2$ be the element $\lp \phi^{\EA:D}(l), c\rp_*$ in $T_{2, D}$.  It is a norm of
$\delta = c \lp \lsm \bg{\xi_1} & \bg{\xi_2}\mb{D}\\ \bg{\xi_2} & \bg{\xi_1}\rsm\rp \in T_{{\rm IV},  A, D}$.
Recall that $n = {\rm ord}_{\EA}\beta$, $N = {\rm ord}_{\EA}\xi_2$, and $\xi_2 = B\vp^N$, with $B \in R_{\EA}^\times$.

Suppose $E/\EA$ is unramified.
We extend $\ve$ to $E^\times$ by letting $\ve(x) = (-1)^{\text{ord}_E x}$ for all $x \in E^\times$.
We define the transfer factor for $(\delta_2, \delta)$ as follows:
\[
\Delta(\delta_2, \delta) =
\ve\!\lp\frac{\xi - \sigma^2 \xi}{2\sqrt{D}}\rp \lp\frac{D_G(\delta)}{D_{H_2}(\delta_2)}\rp^{\!\!1/2}.
\indexi{transfer factor}%
\]
\begin{prop}\label{prop:oiT2DTIVADunram}
The following identity holds:
\[
O_{\delta_2}(1_{K_2}) = \Delta(\delta_2, \delta)O_{\delta}(1_K).
\]
\end{prop}
\begin{proof}
By Lemma M. II. 1 of \cite{F}, 
we have $n = N$.  The proposition then follows from Claims \ref{claim:oiT2D'}
and \ref{claim:oiTIVAD}.
\end{proof}
Suppose $E/\EA$ is ramified.  
Let $\ve_{E/\EA}$ be the quadratic character associated with the extension $E/\EA$.
We define the transfer factor for $(\delta_2, \delta)$ as follows:
\[
\Delta(\delta_2, \delta) = 
-\frac{1}{2}\,\ve_{E/\EA}\!\lp \frac{l \sigma l - \sigma^2l \sigma^3 l}{2 \sqrt{D}}\rp
\lp\frac{D_G(\delta)}{D_{H_2}(\delta_2)}\rp^{\!\!1/2}.
\indexi{transfer factor}%
\]
\begin{prop}
The following identity holds:
\[
O_{\delta_2}(1_{K_2}) = \Delta(\delta_2, \delta)O_{\delta}(1_K).
\]
\end{prop}
\begin{proof}
Note that $({l \sigma l - \sigma^2l \sigma^3 l})/(2 \sqrt{D})$ is equal to $\xi_2 = B\vp^N$.
Recall that $D = -\ep_0\vp$, where $\ep_0 \in R_{\EA}^\times - {R_{\EA}^\times}^{\!\!\!2}$.
The element $\xi_2 \in \EA^\times$ may be written as $\xi_2 = B\ep_0^{-N} (\ep_0\vp)^N$.
The element $\ep_0\vp = \N_{E/\EA} \sqrt{D}$ is a norm from $E^\times$ but $\ep_0$ is not.  Moreover, since
$R_{\EA}^\times \cap \N_{E/\EA}E^\times = {R_{\EA}^\times}^{\!\!\!2}$, the restriction of $\ve_{E/\EA}$ to $R_{\EA}^\times$
coincides with $\kappa$.  Hence, $\ve_{E/\EA}(\xi_2) = (-1)^N\kappa(B)$.

The Jacobian factor $\lp{D_G(\delta)}/{D_{H_2}(\delta_2)}\rp^{1/2}$ is equal to 
\[
\vert\xi_2 \sqrt{D}\vert_{\bar{F}} = \abs{\xi_2}_{\EA}\abs{D}^{1/2}_{\EA} = q^{-N} q^{-1/2} = q^{-N}q_0^{-1}.
\]

By Lemma II. M. 1 of \cite{F}, $n = N$.
The proposition then follows from Claims \ref{claim:oiT2D'} and \ref{claim:oiTIVAD}.
\end{proof}

\end{appendix}
\backmatter
\bibliographystyle{amsalpha}

\printindex{symbols}{List of Symbols}
\printindex{index}{Index}
\end{document}